\documentclass[usenames,dvipsnames,12pt]{article}
\usepackage{hyperref}
\usepackage{amssymb, amsthm,amsmath}
\usepackage{stmaryrd,ulem}
\usepackage{float}
\usepackage{amsmath}
\usepackage{graphicx,psfrag,epsf}
\usepackage{enumerate}
\usepackage{natbib}
\usepackage{url} 
\usepackage{tikz}
\usepackage{pgfplots}
\usepackage{subcaption}
\usepackage{a4wide}
\usepackage[linecolor=blue!60!,backgroundcolor=blue!10!,textwidth=2.5cm,textsize=scriptsize]{todonotes}
\usepackage{comment}

\usepackage{booktabs}
 \usepackage{tabularx}
 \usepackage{array}
 \usepackage{makecell}
\usepackage[table]{xcolor}
\input{macro}

\usepackage[utf8]{inputenc}
\usetikzlibrary{patterns}
\usetikzlibrary{shapes}
\usepackage[linesnumbered,ruled]{algorithm2e}

\usepackage{bm}
\usepackage{dsfont}
\usepackage{amssymb}
\usepackage{amsfonts}
\usepackage{amsthm}
\usepackage{mathtools}
\usepackage{xcolor}
\usepackage{multirow}

\newtheorem{theorem}{Theorem}[section]
\newtheorem{lemma}[theorem]{Lemma}
\newtheorem{proposition}[theorem]{Proposition}
\newtheorem{corollary}[theorem]{Corollary}
\newtheorem{remark}{Remark}[section]

\newtheorem{assumption}{Assumption}

\newtheorem{definition}{Definition}

\newcommand{\mP}{\ensuremath{\mathbb P}}

\newcommand{\OSCI}{{\mbox{{\tiny SCI}}}}
\newcommand{\bench}{0}

\newcommand{\reaction}{\mathrm{react}}
\newcommand{\lfdr}{\mathrm{lfdr}}

\newcommand{\et}[1]{\textcolor{black}{#1}}

\newcommand{\pie}[1]{\textcolor{black}{#1}}
\newcommand{\ind}[1]{{\mathbf{1}\{#1\}}}
\newcommand{\range}[1]{[#1]}
\newcommand{\R}{\mathbb{R}}
\newcommand{\E}{\ensuremath{\mathbb E}}
\newcommand{\V}{\ensuremath{\mathbb V}}
\renewcommand{\P}{\ensuremath{\mathbb P}}

\newcommand{\wh}{\widehat}

\newcommand{\IER}{\mathrm{IER}}

\newcommand{\FDP}{\mathrm{FDP}}
\newcommand{\FP}{\mathrm{FCP}}

\newcommand{\Err}{\mathrm{err}}

\newcommand{\FCP}{\mathrm{FCP}}

\newcommand{\dd}{\operatorname{d}\!}

\newcommand{\Ec}[2]{\E\brac{#1\middle | #2}}
\newcommand{\Pc}[2]{\probp{#1 \mid #2}}

\newcommand{\scal}[3][a]{\inner[#1]{#2;#3}}

\newcommand{\method}{\texttt{OnlineSCI}}

\SetKwComment{Comment}{/* }{ */}

\title{
Online selective conformal inference: adaptive scores, convergence rates and optimality
}
\date{}

\author{Humbert, P., Gazin, U., Heller, R., Roquain, E.}

\begin{document}

\maketitle

\begin{abstract}
In a supervised online setting, quantifying uncertainty has been proposed in the seminal work of Gibbs and Cand\`es \cite{gibbs2021adaptive}. For any
 given point-prediction algorithm, their method (ACI) produces a conformal prediction set with an average miscoverage  getting close to a prespecified level $\alpha$ 
 for a long time horizon. We introduce an extended version of this algorithm, called \method, allowing the user to additionally select times where such an inference should be made.  \method\ encompasses several prominent online selective tasks, such as building prediction intervals for extreme outcomes,  classification with abstention, and online testing. 
 \method\  controls the false coverage proportion among selected times via a pathwise bound for arbitrary sequences, as well as the instantaneous error rate conditional on selection, up to a non-asymptotic remainder term, under stochastic assumptions.  
Importantly, our theory covers the case where \method\  updates the point-prediction algorithm at each time step, a property which we refer to as {\it adaptive} capability. We show that the adaptive versions of \method\   can converge to an optimal solution and provide an explicit convergence rate in several application settings, under model- and estimator-specific regularity conditions.  The favorable behavior of \method\  in practice is illustrated by numerical experiments. 
\end{abstract}

\section{Introduction} \label{sec:intro}

\subsection{Background}

The problem of online uncertainty quantification has received much attention in recent years.
The sequential data stream is  $\set{(X_t, Y_t), t \geq 1}$, where $X_t \in \cX$ is the covariate and $Y_t \in \cY$ is the outcome. The goal is typically to produce a prediction set at each time $t$,  $\cC_t(X_t)$ for $Y_t$, using the covariate $X_t$ and  all previous information $\set{(X_s, Y_s), s <t}$. We consider two data settings:
\begin{enumerate}
    \item {\textit{Adversarial}:} $(X_t, Y_t)$, $t\geq 1$, is an arbitrary sequence of random elements in $\cX \times \cY$, including also the deterministic case.
    \item {\textit{Identically distributed}:} the variables $(X_t, Y_t)$, $t\geq 1$, are identically distributed. The prototypical example is the case where these variables are independent, in which case they form an independent and identically distributed (iid) sequence. We will also consider the case of a dependent, autoregressive model.
\end{enumerate}
The adversarial setting imposes no distributional assumptions: the sequence may be deterministic or stochastic and may exhibit arbitrary distribution shift over time.  By contrast, the identically distributed setting (and in particular the iid setting) is a favorable situation for which we can expect optimality results to hold.

For the adversarial setting, \cite{gibbs2021adaptive} suggested the {\it adaptive conformal inference} (ACI) method: a gradient descent method that adjusts the width of the prediction set.  The sensitivity of ACI to its step size parameter has led to modifications that tune the step size adaptively \citep{zaffran2022adaptive,gibbs2024conformal}. 
  The aforementioned methods guarantee  that the average miscoverage, $\sum^t_{k=1}  \ind{Y_k \notin \cC_k}/t$, is in the long run  at most a prespecified target $\alpha$. 
 \et{However, the  instantaneous miscoverage probability, {$\mP(Y_t \notin \cC_t\mid \mathcal{F}_{{t-1}})$, where $ \mathcal{F}_{{t-1}}$ denotes the past information at time $t$}, may deviate substantially from $\alpha$,  even if $\set{(X_t, Y_t), t \geq 1}$ are iid. Addressing this problem, \cite{angelopoulos2024online} suggested requiring 
 that the algorithm achieves asymptotic non-coverage probability $\alpha$ when the online examples are iid. Thus, asymptotically, the coverage guarantee is the same as for a single example in the classic split conformal literature \citep{vovk2005algorithmic}.  
 More formally, \cite{angelopoulos2024online} suggested modifications of ACI that ensure: (1) $\lim_{t\rightarrow \infty} \sum^t_{k=1}  \ind{Y_k \notin \cC_k}/t=\alpha$ for any adversarial data setting, and (2) in the iid  data setting, $\lim_{t\rightarrow \infty} \mP(Y_t \notin \cC_t) = \alpha$. Their decaying-step ACI methods provide both guarantees when inference is produced at every time point. The question addressed here is whether both guarantees can be retained when inference is reported only at selected times. 
}

\subsection{Selection effect}

The guarantees above concern inference at every time point. In many applications, however, the analyst is interested only in a selected subset of the inferences, for example, when the outcome is expected to be extreme  (say extreme weather forecast, where arguably accurate forecast when the weather is not extreme is not of great importance), or when the covariates indicate immediate attention is needed (say, for high-risk patients). A typical setting where selection occurs is when the stream of examples is dense (i.e., new observations arrive frequently) and reporting prediction sets for all examples is not meaningful, so  selection  greatly improves the relevance of the reported results. 
With selection, a natural cumulative error criterion  is  the false coverage proportion (FCP) among selected times, rather than on all examples.  Control of the expected FCP under selection  has been addressed thus far in the batch setting, see \cite{bao2024selective, gazin2024selecting}. 
\et{Under stochastic assumptions,  a natural instantaneous error criterion is  the non-coverage probability {\it conditionally on being selected and on the past/history}\et{, called here the {\it instantaneous error rate } (IER)}. This criterion prevents satisfactory cumulative error control from masking severe miscalibration at particular times. Selection-conditional coverage has been studied in the offline batch setting  \cite{jin2024confidence} and, more recently, in online settings \cite{bao2024cas, sale2025online, zheng2026online}; \cite{sale2025online, zheng2026online} additionally condition on past decisions. To the best of our knowledge, however, coverage conditional on selection and the full observed history  has not been previously investigated. 
}

\begin{figure}[t]
		\centering
		\includegraphics[width=1.\linewidth]{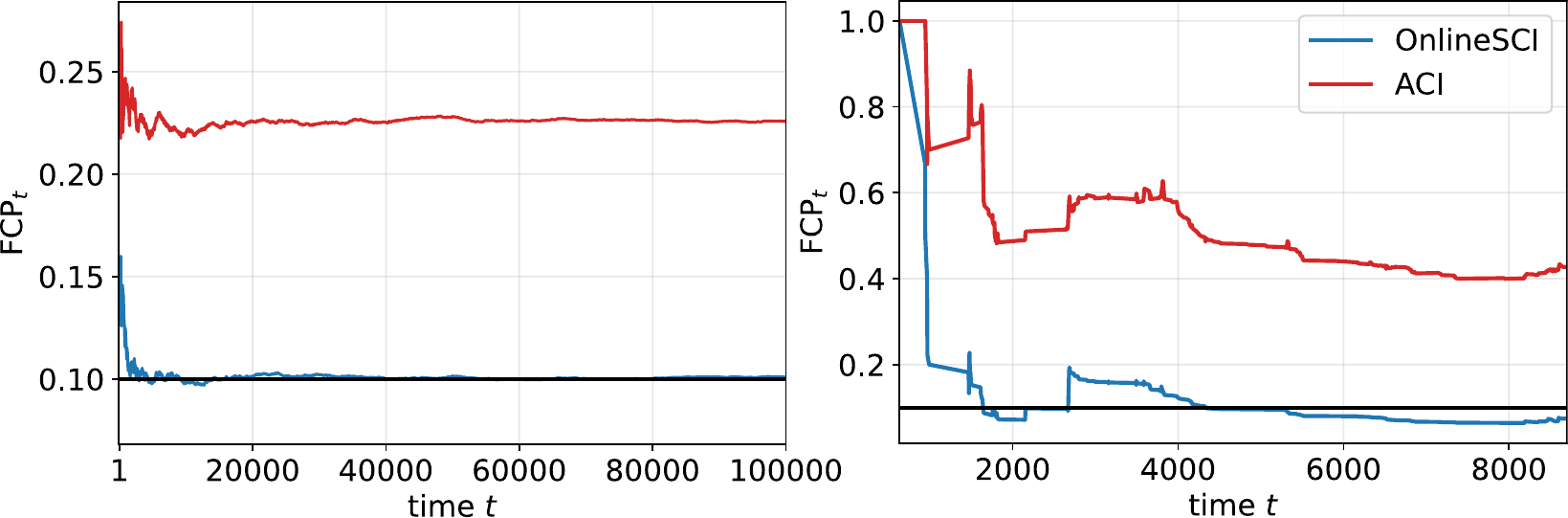}
        \caption{Cumulative FCP among selected times  versus  time $t$ for  ACI with decaying step sizes \citep{angelopoulos2024online} (red) and  the  proposed method \method\  (blue). \et{The horizontal black line corresponds to the target level $\alpha=0.1$.}
        Left panel:  online selective classification for $K=2$ classes (see \S~\ref{subsec-simClassification} for more details). Right panel:  online selective prediction intervals for a time series of  electricity prices (selection takes place if the  outcome at the previous time point is below a threshold, see \S~\ref{sec:elec}   for more details).}
        \label{fig:intro_Selective_interval}
\end{figure}

\et{
Selection invalidates methods that build on ACI by potentially inflating the FCP among the selected. 
Consider, for example,  selective classification, where the outcome is a class ($\cY=[K]:=\{1,\ldots,K \}$) and the analyst is interested only in prediction sets of size one (hence corresponding to an `actual' classification when selected).  Figure~\ref{fig:intro_Selective_interval} (left panel) shows  that for $K=2$ classes classic conformal prediction sets \citep{vovk2005algorithmic}, which guarantee marginal miscoverage probability at level $\alpha$, are invalidated by selection:  while the FCP {(misclassification proportion in this case)} is acceptable without selection, it is unacceptably high among the selected examples, where the prediction set is necessarily of size one.  The inflation follows since we expect a fraction $\alpha$ of examples to have errors  across all time points, and  these errors are necessarily among the selected. Consequently,  $\mP(Y_t \notin \cC_t \mid S_t = 1)$ (where $S_t=\ind{|\cC_t|=1}$ denotes the indicator of selecting at time $t$) can substantially exceed $\alpha$ with ACI, as it does in Figure~\ref{fig:intro_Selective_interval}.  }
Instead,  \method\ tailored for online selective classification in   Algorithm~\ref{alg:BinaryClassification} updates its calibration threshold on the selected subsequence, and the resulting  FCP remains close to the target  level. Algorithm~\ref{alg:BinaryClassification} adds to the available methods for  classification with abstention (or with a rejection option), introduced by \cite{chow70}. 
In the machine learning literature,  abstention is often accounted for by adding a term to the risk that penalizes any rejection \cite{herbei06,bartlett08, wegkamp2011support, cortes2016learning, ni2019calibration, cao2022generalizing}. \et{Our suggestion, on the other hand, aims to control the misclassification error rate on the selected, as done in \cite{zhao2023controlling, gazin2024selecting} for the batch setting, and we also investigate the IER which corresponds to the instantaneous classification error conditionally on being selected and on the history. }

\et{
The inflation in FCP for ACI-type methods in the setting of constructing selective prediction intervals is illustrated for a time series of electricity prices in the right panel of Figure~\ref{fig:intro_Selective_interval} (where selection occurs if the outcome observed at the previous time point is below a prespecified value). The FCP is far smaller, and  close to the target level $\alpha$, for \method\ tailored for online selective prediction intervals in   Algorithm~\ref{alg:PredictXoriented}.
}

\begin{figure}[t]
		\centering
		\includegraphics[trim={0 0 0 55}, clip, width=1.\linewidth]{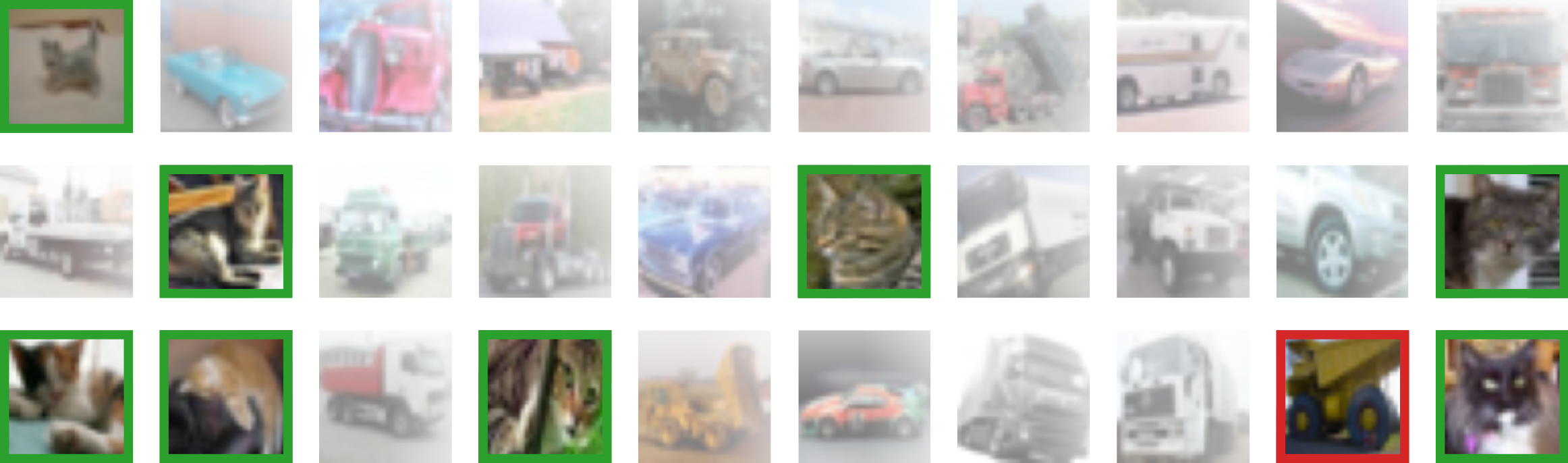}
        \caption{Illustration of the online conformal testing application for \texttt{CIFAR-10} \et{with $\alpha=0.1$}. Classes `\textit{\et{car}}' and `\textit{truck}' correspond to $Y=0$ (null) and `\textit{cat}' to $Y=1$ (alternative). A red square indicates a false discovery and a green square a true discovery. Images without a colored border are not selected.
}
        \label{fig:intro_ND}
\end{figure}

\et{
A prominent related setting is online multiple testing, where selected times correspond to rejected hypotheses.  Online testing  focuses 
on inference over an infinite sequence of fixed, unknown parameters ordered in time, see \cite{FS2008,Javanmard18,ramdas2017online,gang2023structure}, among others. At time $t$, let  $Y_t\in \{0,1 \}$ denote the hypothesis state ($Y_t=0$ if the null hypothesis is true, $Y_t=1$ otherwise). Let $S_t\in\{0,1\}$ be the selection rule ($S_t=1$ if the null hypothesis is rejected, $S_t=0$ otherwise). In this setting, the error corresponds to a false rejection, and the FCP reduces to the  false discovery proportion (FDP) and IER to the `instantaneous' false discovery rate (probability that a rejected null is true conditional on rejection and the history). 
In classical online testing, $(Y_t)_{t\geq 1}$ is unobserved. 
In the feedback setting studied here, $Y_t$ is revealed after each decision, so subsequent decisions may use past labels.  This framework includes  online novelty detection (also referred to as out-of-distribution testing or outlier detection,  \cite{blanchard2010semi, liang2022integrative, pimentel2014review}), where the goal is to test whether a new observation comes from the same unknown distribution as a reference null sample.   
In the CIFAR-10 illustration in Figure~\ref{fig:intro_ND}, \method\ makes eight discoveries, one of which is false, whereas  LORD++ \citep{ramdas2017online}, SAFFRON \citep{ramdas2018saffron} and Feedback LORD \citep{lu2025feedback} make no discoveries on the 20 images.  
}

\et{
These examples motivate the central question of the paper: can a single online procedure retain a pathwise error guarantee on the selected subsequence while also calibrating the history-conditional error of each selected decision under regular stochastic regimes? 
}

\subsection{Related work and main contributions}

Selective variants of the ACI algorithm have recently been proposed \cite{bao2025cap,sale2025online}, with the threshold  updated  only at selected times. 
This modification is appealing because it directly leads to a valid pathwise FCP bound for the selected observations 
(see Theorem~\ref{thangelinfo} below for a formal statement in a more general setting). 
 \et{Existing methods, however, do not jointly provide, within one recursion, a pathwise selected-error guarantee for arbitrary data streams, explicit rates for history-conditional one-step calibration, adaptively updated scores, and selection rules coupled to the current threshold. 
 }
 
\et{\method\,, with its single recursion, jointly provides all of the above. It uses decaying step sizes indexed by the number of selections rather than by calendar time. Selection therefore replaces the usual deterministic stochastic-approximation clock by a random effective sample size $J(t)$ \eqref{def:Jt}. For informative selection, the decision to update is itself coupled to the evolving threshold. Adaptive scores enhance the good practical and theoretical properties. The same recursion applies to selective regression, selective classification, and testing, including settings in which an inference is reported only when its current output is sufficiently informative.  We also propose a practical data-driven rule for tuning the method's hyperparameters\footnote{See the interactive website \url{https://onlinesciapp-srjzthcznjslvzjiavpahm.streamlit.app/}}.
}

\et{
In the study of our novel method  \method, we show: 
\begin{itemize}
\item[(i)] \textit{Joint cumulative and instantaneous guarantees.} A single \method\ recursion provides a non-asymptotic pathwise bound on the realized FCP among selected times for arbitrary data streams and, under stochastic regularity conditions, explicit convergence rates for the history-conditional IER.
\item[(ii)] \textit{A general convergence analysis for selective and adaptive procedures.} Our bounds explicitly separate the effects of stochastic-approximation error, the randomness in the number of selections, and estimation error from adaptively updated predictors or test statistics. To our knowledge, the resulting high-probability rate is new even in the nonselective decaying-step setting. 
\item[(iii)] \textit{Informative selection.} Beyond covariate-oriented selection, our theory covers threshold-coupled selection rules that report an inference only when its current output is sufficiently informative. This includes length-restricted prediction intervals, prediction intervals lying entirely  above an indifference zone, selective classification, and online testing. The analysis identifies  an intrinsic  critical-level phenomenon and an unavoidable failure-probability for such rules.
\item[(iv)] \textit{Oracle convergence and optimality.} Under model- and estimator-specific regularity conditions,  adaptive  \method\ converges at explicit rates to the corresponding oracle rules. We establish optimality results, up to explicit remainder terms, in regression, classification, and testing.
\end{itemize}
}

 \et{
The principal contribution is therefore not only the selected-time recursion itself, but a unified theory regarding the inference.
After describing \method, we proceed with the following outline.  In \S~\ref{sec:FDPbound}, we establish a  non-asymptotic FCP 
bound, valid in the adversarial setting.}  In  \S~\ref{sec:rate} we provide the main results:  general convergence rate results for the IER of \method, both with high probability (Theorem~\ref{th:rateconvfull}) and in expectation (Theorem~\ref{th:genL2-first}). These results come from general non-asymptotic approximation bounds provided in  \S~\ref{sec:genapproxbounds}, which are valid under weak distributional assumptions. In \S~\ref{sec:incompresscritic} we address the special case of an informative selection, where the selection may depend on the procedure threshold itself. We show that, in that case, a failure probability $\delta\geq 0$, defined in  \eqref{deltaprobaJtfinite}, is unavoidable in our high-probability bound.

 \et{The main results are applied to online selective regression in \S~\ref{sec:baoselect} for various selection rules: time-constant (\S~\ref{sec:appliiidregression_timeconstant}), large prediction (\S~\ref{sec:adaptiveXorientedSelection}), length restriction (\S~\ref{sec:lengthrestriction}), informative lower bound (\S~\ref{rk:plb}). While these results are in the iid regression model, we also consider the stationary autoregressive model in \S~\ref{sec:AR}. In each setting, we provide a specific convergence rate to the oracle which explicitly depends on the estimation quality of the regression mean (and variance).}
 
\et{The case of online selective classification, where prediction sets are selected only if they are of size one, is explored in  \S~\ref{sec:selectiveclassifiid}, while the case of online conformal testing, where selection corresponds to  rejecting the null hypothesis, is investigated in \S~\ref{sec:selectiveNDiid}. In both cases,  we show that \method\: converges to an optimal solution. As a byproduct, we provide an online method for selection by prediction (\S~\ref{sec:selectpredict}), originally introduced by \cite{Jin2023selection} in the batch setting. }

 Numerical illustrations are provided in \S~\ref{sec:num}.  They support the theoretical findings:  the FCP and IER of \method\  approach the  nominal level $\alpha$, whereas naive post-selection use of  the ACI method of \cite{angelopoulos2024online} can fail. Also, the simulations show that the adaptation stage (i.e., the online updating of the  score or test statistic computations) in \method\  is important for convergence to the optimal procedure.

\section{Preliminaries}\label{subsec:prelim}

We introduce the setting, notation, and the error criteria we aim to control, present our method \method\ , and derive its error upper bound.

\subsection{Setting}\label{subsec:prelim:setting}

To model the information that is available as data arrive sequentially, we introduce the filtration
\begin{equation}\label{filtration}
\mathcal{F}_{t} := \sigma(\cD_0, \set{(X_k, Y_k), 1 \leq k \leq t}), \quad t \geq 1 ,
\end{equation}
where $\cD_0 := (X_{i}, Y_{i})_{-(n-1)\leq i\leq 0}$ is a holdout set of size $n$ available at time $1$. 
Hence, $\mathcal{F}_{t-1}$ represents information available at time $t$, before observing $X_t$.
\et{Then, as $X_t$ is observed, a decision (selection plus inference) is made, $Y_t$ is observed and this feedback is used to adjust the next decision.}

A {\it procedure} $\cR=((S_t)_{t\geq 1}, (\cC_t)_{t\geq 1})$  is any sequence such that, for each  $t\geq 1$,  it comprises a selection decision $S_t\in \{0,1\}$ and a prediction set $\cC_t \in \mathcal{P}(\cY)$, both measurable functions of $\mathcal{F}_{t-1}$ and $X_t$, so  that each event $\{Y_t\in \cC_t\}$ is measurable.
Here, $\mathcal{P}(\cY)$ denotes the power set of $\cY$.
Throughout the paper, we focus on the following class of procedures. 

\begin{definition}\label{def:thr_proc}
A threshold-sequence based procedure is a procedure of the form
$\cR=((S_t(X_t, \qt_t))_{t\geq 1}, $ $ (\cC_t(X_t, \qt_t))_{t\geq 1}),$ where \et{the selection rules $(\St_t(\cdot,\cdot))_{t\geq 1}$ and prediction set rules $(\cC_t(\cdot,\cdot))_{t\geq 1}$ are as follows:}
\begin{itemize}
    \item[(i)] for all $t\geq 1$:  
 each random $\qt_t\in \R$ is a $\mathcal{F}_{t-1}$-measurable threshold;  each random $\St_t: (x,q)\in \cX \times \mbr \mapsto \St_t(x,q)\in \set{0, 1}$ is a selection rule such that $\St_t(x,q)$ is $\mathcal{F}_{t-1}$-measurable; each random $\cC_t: (x,q)\in \cX \times \mbr \mapsto \cC_t(x,q) \in \mathcal{P}(\cY)$ is a prediction set rule such that $\cC_t(x,q)$ is $\mathcal{F}_{t-1}$-measurable \et{and each event $\{Y_t\in \cC_t\}$ is measurable.}
 \item[(ii)] For some $B>0$, at least one of the following conditions must apply for all $x\in \cX$, $t\geq 1$, and $q\geq B$:  the prediction set is the entire space, $\cC_t(x,q)=\cY$ ; or there is no selection, $\St_t(x,q)=0$.
\end{itemize}
\end{definition}

Typically, the prediction set rule will take the following classical form \citep{gibbs2021adaptive, angelopoulos2024online}: 
\begin{align}
    \cC_t(x, \qt) &= \set{y \in \cY\; : \; V_t(x, y) \leq \qt} \mbox{ , $x\in \cX$, $q\geq 0$,}\; \label{eq:set_pred_funct}
\end{align}
with $\cC_t(x, \qt)$ left arbitrary for $q<0$, where $V_t: (x, y) \in \cX \times \cY \mapsto V_t(x, y)\in [0,B]$ 
is a given $\mathcal{F}_{t-1}$-measurable {\it non-conformity score function} which is large when $y$ is far from its point-prediction at $x$. 

In (two-sided) prediction intervals for regression  $\cY=\R$, $\cX=[0,1]^d$, we can choose the score function to be the rescaled residual \citep{papadopoulos2008normalized, papadopoulos2011regression}: for $x\in \cX$, $y\in \cY$,
    \begin{equation}\label{regressionVadaptive}
        V_t(x,y)=2B\Phi \Big(|y-\hat{\mu}_t(x)|/\hat{\sigma}_t(x)\Big)-B \in [0,B),
    \end{equation}
    for some estimators $\hat{\mu}_t(x)$ (resp. $\hat{\sigma}_t^{2}(x)$) of the conditional expectation (resp. variance) of $Y_{t}$ given $X_{t}=x$ (each estimator being $\mathcal{F}_{t-1}$ measurable). 
   The rescaling is done here by the function $2B\Phi(\cdot)-B$, where $\Phi$ is the c.d.f. of the standard normal distribution, although this choice is arbitrary and not important for the sequel.
    In this case, \eqref{eq:set_pred_funct} reduces to the prediction interval
\begin{equation}
\mathcal{C}_t(x,q) = \bigg[ \hat{\mu}_t(x) \pm \hat{\sigma}_t(x) \Phi^{-1}\Big(\frac{q+B}{2B}\Big)\bigg]\label{equCtregression},
\end{equation}
if $q\in ( 0,B)$, with $\mathcal{C}_t(x,q) =\{\hat{\mu}_t(x)\}$ if $q\leq 0$ and $\mathcal{C}_t(x,q) =\R$ if $q\geq B$.

\et{The selection rule $\St_t(\cdot, \cdot)$ at time $t$ can be arbitrarily general. There are  two main types of selection rules. An {\it $X$-oriented selection rule}  does not depend on the threshold decision $q$, i.e., $S_t(x,q) = S_t(x),$ but can depend on $t$ and on $\mathcal F_{t-1}$. Such rules have been considered in \cite{angelopoulos2024online, bao2025cap, sale2025online} and include, for instance, the time-constant selection $S_t(x)=S(x)=\ind{x\in \mathcal{A}}$ for some user-fixed specific measurable set $\mathcal{A}\subset \cX$ for the covariate vector. In the regression case above, it also encompasses $S_t(x)=\ind{\hat{\mu}_t(x)\geq y_0}$ or $S_t(X_t)=\ind{Y_{t-1}\geq y_0}$ for some benchmark value $y_0\in \R$ and estimator $\hat{\mu}_t(\cdot)$, which corresponds to situations where a particularly extreme event is occurring. 
Note that since $S_t(x,q)$ does not depend on $q$, a typical choice is $\cC_t(x,q)=\cY$ for $q\geq B$ to satisfy  Definition~\ref{def:thr_proc} (ii).
These selections will be considered in \S~\ref{sec:appliiidregression_timeconstant},\S~\ref{sec:adaptiveXorientedSelection} and~\S~\ref{sec:AR}.
}

\et{By contrast, an {\it informative selection rule} depends on the threshold decision $q$ and is typically of the following form:
\begin{align}\label{eq:selectionstat}
\St_t(x, \qt) = \ind{W_t(x) >q}, \:\:\: q\in\R,
\end{align}
where $W_t: x \in \cX \mapsto W_t(x)\in [0,B]$ 
is a given $\mathcal{F}_{t-1}$-measurable `informativeness' score. 
It satisfies the constraint $\St_t(x, \qt)=0$ (no selection) for $q\geq B$ in  Definition~\ref{def:thr_proc} (ii) and $\St_t(x, \qt)=1$ (always selection) when $q<0$.
Such selection rules are inspired from \cite{gazin2024selecting} (developed therein in the batch setting), and are called 
 `informative' in the sense that they are able to only select the examples for which the {\it actually produced prediction set} $\mathcal C_t$ is `interesting'. For instance, in the regression case above, we can prescribe a minimum accuracy by imposing that $\mathcal C_t$ should be of Lebesgue measure below $2\lambda_0>0$, 
 see \S~\ref{sec:lengthrestriction}.
  Another useful choice is to consider that the prediction set 
 $\mathcal C_t $ (or more precisely its one sided version) is informative only if its lower bound is above some benchmark value $y_0\in \R$, see \S~\ref{rk:plb}. 
}

\et{For classification outcomes, $\mathcal{Y}=[K]$, and selective classification and conformal testing correspond to informative selection with specific choices of $V_t(x,y)$ and $W_t$ in \eqref{regressionVadaptive} and \eqref{eq:selectionstat}, see \S~\ref{sec:selectiveclassifiid}~and~\S~\ref{sec:selectiveNDiid}.}

\subsection{Error rates}\label{prelim-errorrates}

To evaluate the performance of a given procedure $\cR=((S_t)_{t\geq 1}, (\cC_t)_{t\geq 1})$, we measure whether the decision taken at time $t$ for the response $y$ is an error or not with an \textit{error function}: 
\begin{equation}\label{eq:err_funct}
    \Err_t: (y,C)\in \cY \times \mathcal{P}(\cY) \mapsto \Err_t(y,C)\in  [0, 1] \, ,
\end{equation}
with the conventions $\Err_t(y,\cY):=0$ and $\Err_t(y,\emptyset):=1$ for all $y\in \cY$. 
 While the most classical choice is the coverage error \citep{gibbs2021adaptive}
\begin{equation}
    \label{coverror}
    \Err_t(y, C) = \ind{y \notin C},
\end{equation}
other choices are covered by our theory, see Remark~\ref{rem:weightederror}. Given this error function, the average error rate over selected times of the procedure $\cR=((S_t)_{t\geq 1}, (\cC_t)_{t\geq 1})$, referred to as the FCP, is defined\footnote{Note that we implicitly assume here that each $\Err_{t}(Y_{t}, \cC_t)$ is measurable.}  as follows at time   $t$: 
\begin{align}
    \qquad \FP_t(\cR) &:= \dfrac{\sum_{k=1}^t \Err_{k}(Y_{k}, \cC_k)  \cdot S_k}{1\vee \sum_{k=1}^t S_k  } \, . \label{equFCP} 
\end{align}
In the setting with  $\cC_k\equiv \{1\}$ and $\cY=\{0,1\}$, $ \FP_t(\cR) = \big(\sum_{k=1}^t \ind{Y_{k}=0} \cdot S_k \big)/\big(1\vee \sum_{k=1}^t S_k \big)$ corresponds to the FDP.

Algorithms guaranteeing  that the average error over the selected  is at most $\alpha$ may nevertheless be unsatisfactory, since we can have average non-coverage at most $\alpha$ even if the instantaneous non-coverage probability at each time step is either much larger, or much smaller, than $\alpha$. If it  is smaller than $\alpha$, the  prediction set is too wide (and hence less accurate); if it is larger than $\alpha$ then the prediction set fails to provide the desired confidence, since it is too narrow.

This work aims to control the probability that the output prediction set makes an error, conditional on being selected and on the past. 
We refer to this quantity as the \textit{instantaneous error rate} at time $t \geq 1$ for a given procedure $\cR = ((S_t)_{t \geq 1}, (\cC_t)_{t \geq 1})$, and define it as:
\begin{align}
\label{equ-instaneousgen}
\IER_t(\cR):= \E[\Err_t(Y_t, \cC_t) \:|\:  S_t=1,\mathcal{F}_{t-1}] \; ,
\end{align}
\et{with the convention $\IER_t(\cR):=0$ if $\P(S_t=1\:|\:\mathcal{F}_{t-1})=0$. Note that $\IER_t(\cR)$ is random and $\mathcal{F}_{t-1}$-measurable. }

\begin{remark}\label{rem:weightederror}
Our work can accommodate non-symmetric errors. For instance, in classification, we may consider $\Err_t(y, C) = w_y \ind{y \notin C}$, 
for some weights $w_y\geq 0$ with $\sum_{y\in \range{K}} w_y=1$. In regression, we can consider  $\Err_t(y,C)\in [0,1]$ measuring the distance 
between $y$ and the set $C$. 
\end{remark}

\subsection{OnlineSCI}\label{sec:method}

We introduce our general method, called \method.

\begin{definition}\label{defmainproc}
 The \method~procedure with selection rules $(\St_t(\cdot,\cdot))_{t\geq 1}$ and prediction set rules $(\cC_t(\cdot,\cdot))_{t\geq 1}$ as in Definition~\ref{def:thr_proc} is the threshold-sequence based procedure $\cR^{\OSCI} = ((S_t(X_t, \qt^{\OSCI}_t))_{t\geq 1},$ $ (\cC_t(X_t, \qt^{\OSCI}_t))_{t \geq 1})$, where the threshold sequence $(\qt^{\OSCI}_t)_{t\geq 1}$ is given by the following recursion:
    \begin{align} \label{equACIquantile}
        \qt_{t+1}^{\OSCI} &=  \qt_{t}^{\OSCI} + \gamma_{J(t)} \paren{\reaction_{t}(Y_t, \cC_t(X_t, \qt_t^{\OSCI}),q_t^{\OSCI})-\alpha} \cdot S_t(X_t, \qt_t^{\OSCI}) \; ,\:\:\:t\geq 1\;,
    \end{align}
    where $\qt_1^{\OSCI} \in [0,B)$, $(\gamma_j)_{j \geq 1}$ is a positive step size sequence, $J(t)$ given by 
    \begin{equation}\label{def:Jt}
        J(t):=\sum_{i=1}^{t-1} \St_{i}(X_i, \qt_i^{\OSCI}) + 1
    \end{equation}
  is the number of selections before time $t$ plus one, and the reaction function  is given by:
    \begin{equation}\label{equreact}
        \reaction_{t}(y, C,q) := 
        \et{\Err_t(y, C)\vee \ind{\qt < 0} }
\; ,
    \end{equation}
    with $\Err_t(\cdot)$ defined by \eqref{eq:err_funct}. 
\end{definition}

The recursion \eqref{equACIquantile} works by updating $q_t^{\OSCI}$ only if there is a selection, that is, $\St_t(X_{t}, \qt_t^{\OSCI})=1$, which leads to a piecewise constant sequence $\qt_t^{\OSCI}$, whose value is updated only at each selection time. A decision $\cC_t(X_t, \qt_t^{\OSCI})$ is produced at the selected time and the value of $Y_t$ is revealed. Then, the update increases  $q_t^{\OSCI}$ by $\gamma_{J(t)}(1-\alpha)$ when an error occurs (or when $q_t^{\OSCI} < 0$) and decreases it by $\gamma_{J(t)}\alpha$ otherwise. Importantly, the step size $\gamma_{J(t)}$ depends on the  number selected $J(t)$, which means that it is  updated only at each selection. This avoids artificially decreasing the step size when not needed. 
Since $q_t^{\OSCI}$ is $\mathcal{F}_{t-1}$-measurable,  \method~is a threshold-sequence based procedure (see Definition \ref{def:thr_proc}). 

\method\ is an extension of ACI \citep{gibbs2021adaptive} that presents several innovations. First, it follows the approach of \cite{angelopoulos2024online}, by using a step size sequence $(\gamma_j)_{j \geq 1}$ that tends to zero. 
As proved therein in case of no selective inference, this entails an interesting convergence behavior of $\qt_t^{\OSCI}$ in the iid case, which limits the oscillations of the original ACI threshold sequence (that used a constant $(\gamma_j)_j$ sequence), in particular reducing drastically the occurrences of the threshold below $0$ and above $B$. Second, \method\ allows for incorporating a selection rule by updating the threshold only at the selected time. Doing so, it inherits the classical FCP control of \cite{angelopoulos2024online} in restriction to the selection, as we formally prove in \S~\ref{sec:FDPbound}. Third, the proposed formulation possesses sufficient flexibility to encompass a range of settings, see Algorithms~\ref{alg:PredictXoriented}, \ref{alg:BinaryClassification}, and~\ref{alg:onlineBH}. The theoretical results established for  \method\ provide a unified basis for deriving the specific results associated with each variant, with the main results in \S~\ref{sec:rate} yielding corollaries in \S~\ref{sec:baoselect},  \S~\ref{sec:selectiveclassifiid} and \S~\ref{sec:selectiveNDiid}. 

The procedure \method\ has a specific behavior when $q^\OSCI_t$ is outside $[0,B]$. \et{Assume that there is a selection at time $t$.  }
If $q^\OSCI_t<0$\footnote{This  case corresponds to a 'favorable' situation where many errors have not been made in a row. While it is rarely met in practice, a relevant inference should be proposed in that case too.}, the reaction function forces $q^\OSCI_{t+1}$ to be larger than $q^\OSCI_t$, regardless of the rule $\cC_t(X_t, q)$ for $q<0$. As a consequence, the threshold sequence value $q^\OSCI_t$ does not depend on the choice of $\mathcal{C}_t(X_t,q^\OSCI_t)$ when $q^\OSCI_t<0$, so the analyst is free to choose what to do in this case.  
This differs from the method of \cite{angelopoulos2024online}, which imposes $\mathcal{C}_t(x,q)=\emptyset$ for $q<0$ (hence incurring an error whenever $q^\OSCI_t<0$). 
Choosing $\mathcal{C}_t(x,q)\neq \emptyset$ can only decrease the FCP value so is in general desirable, if it makes sense in the specific setting considered.

If $q^\OSCI_t\geq B$, according to Definition~\ref{def:thr_proc}, we should have either $\mathcal{C}_t(X_t,q^\OSCI_t)=\mathcal{Y}$ or $S_t(X_t,q^\OSCI_t)=0$. In the case where $S_t(x,q)=S_t(x)=1$ does not depend on $q$, we are in the former case and $q^\OSCI_{t+1}$ decreases below $q^\OSCI_t$ until going below $B$, and the algorithm goes on until it converges to the appropriate quantity (as we show in \S~\ref{sec:baoselect}, in which it is proven that $q^\OSCI_t\geq B$ occurs with a probability vanishing when $t$ grows, under mild assumptions). In the case where the selection is of the form \eqref{eq:selectionstat}, we are in the latter case and $S_t(X_t,q^\OSCI_t)=0$, which entails $q^\OSCI_{T}=q^\OSCI_t$ for all $T>t$. Hence, the algorithm 'is frozen' above $B$ so cannot converge to the appropriate quantity. 
This behavior is not a problem specific to \method\:. It is unavoidable if one seeks both an adversarial FCP bound and a non-trivial prediction set, see Proposition~\ref{prop:unavoidable}.
Importantly, the probability that $q^\OSCI_t\geq B$ is empirically small, but non-zero, and theoretically quantified by the term $\delta$ defined in \eqref{deltaprobaJtfinite} in \S~\ref{sec:incompresscritic}.

\begin{remark}[Restart strategy]\label{rem:rerun}
  \et{  In practice, when $q^\OSCI_{t_0}\geq B$ for a selection of the form \eqref{eq:selectionstat}, we suggest to restart \method, incorporating all the information up to $t_0$   \et{to adjust the parameters of the method, typically by using the approach of \S~\ref{sec:select_process}.} 
All our iid results apply with the new holdout set $\mathcal{D}$ corresponding to $\mathcal{F}_{t_0}$ of size $n:=n+t_0$,  provided that $t_0$ is a stopping time (the future stream is independent of  $\mathcal{F}_{t_0}$).  
This restart strategy can be combined with other restart rules (e.g., change-point detection) when the data distribution shifts over time. 
}
\end{remark}

\subsection{FCP bound}\label{sec:FDPbound}

Following the argument in \cite{angelopoulos2024online} restricted on the selection, we show the procedure \method\ controls  $\FP_t$ in \eqref{equFCP} in the adversarial setting. Similar bounds on the FCP in the selective case have been proposed by \cite{bao2025cap,sale2025online} for other variants of ACI. 
\begin{theorem}\label{thangelinfo}
For any sequence $(X_{t},Y_{t})_{t\geq 1}$, denoting \et{$N(t):=\sum_{k=1}^{t} S_{k}(X_{k}, \qt_{k}^{\OSCI}) = J(t+1)-1$} for all $t\geq 1$, see \eqref{def:Jt}, the procedure \method\ is such that 
for all $t\geq 1$,
\begin{align}\label{FCPbound}
    \FP_t(\mathcal{R}^{\OSCI})&\leq \alpha + \dfrac{B + \max_{1 \leq j \leq J(t)}(\gamma_{j})}{N(t)} \left(\gamma^{-1}_{1} + \sum_{j=2}^{N(t)} \lvert \gamma^{-1}_{j}-\gamma^{-1}_{j-1} \rvert \right) \; 
\end{align}
\et{if $N(t)\geq 1$, and $\FP_t(\mathcal{R}^{\OSCI})=0$ otherwise.} In particular,  further assuming that the sequence of $(\gamma_j)_{j\geq 1}$ is nonincreasing, we have for all $t\geq 1$,
\begin{align}\label{FCPbound_decre}
     \FP_t(\mathcal{R}^{\OSCI})&\leq \alpha + \frac{B+\gamma_{1}}{N(t)\gamma_{N(t)}} \; 
\end{align}
\et{if $N(t)\geq 1$, and $\FP_t(\mathcal{R}^{\OSCI})=0$ otherwise.}
\end{theorem}
    See \S~\ref{proof:thangelinfo} for the proof. 
The bound in \eqref{FCPbound} provides an upper bound on  $\FP_t$ and shares similarities with the one proposed in  \cite{angelopoulos2024online}. However, there are two main differences. First, our algorithm includes a selection part, which makes the problem more challenging and slows down the convergence of the upper bound if there are only few selections. 
\et{Second, a lower bound on  $\FP_t$ is also provided for the algorithm of \cite{angelopoulos2024online}, which is not possible here. This comes from the difference in OSCI algorithm: our formulation  
does not impose  an error whenever $\qt_t^{\OSCI}<0$, see \eqref{equreact}, which means that the case $\FP_t = 0$ cannot be formally excluded in \eqref{FCPbound_decre}. Since  the inequality should be true for any sequence $(X_{t},Y_{t})_{t\geq 1}$, the largest lower bound possible is $0$.   }

\begin{remark}
    In conformal testing, Theorem~\ref{thangelinfo} says that  $\FP_t =\FDP_t$ is bounded close to $\alpha$ if enough discoveries have been made before time $t$ in an adversarial setting, which is a new contribution. 
\end{remark}

\section{Convergence rates}\label{sec:advers}\label{sec:thresholdapprox}\label{sec:iid}\label{sec:rate}

This section presents approximation results for the \method\  procedure. For brevity, the explicit finite-sample approximation bounds are deferred to \S~\ref{sec:genapprox}, and we focus here on the resulting convergence rates.

\subsection{Benchmark IER function}

The instantaneous error rate \eqref{equ-instaneousgen} of \method\  corresponds, at each time step, to the probability that the prediction set makes an error conditionally on being selected (and on the past). It is given by 
\begin{align}
\IER_t(\mathcal{R}^{\OSCI})=\Pi_t(q_t^{\OSCI}),\:\:t\geq 1\label{IEROSCI} \; ,
\end{align}
for $q_t^{\OSCI}\in [0,B]$, 
for the functional
\begin{align}
\label{PiCtSt}
    \Pi_t:=\Pi_{\cC_t,S_t} : q\in [0,B] \mapsto \:&\E[\Err_t(Y_t, \cC_t(X_t,q)) \:|\:  S_t(X_t, q)=1,  \mathcal{F}_{t-1}] \; ,
\end{align}
with $ \Pi_t (q)=0$ by convention if $\P(S_t(X_t, q)=1 \:|\:  \mathcal{F}_{t-1})=0$.  When $X_t$ is independent of the past $\mathcal{F}_{t-1}$, the expectation in \eqref{PiCtSt} is taken only with respect to the distribution of $(X_t,Y_t)$. Under dependence, however, the distribution of $(X_t,Y_t)$ is conditional on $\mathcal{F}_{t-1}$. 
\et{Note that $\Pi_t(q)$ is not defined for $q\notin [0,B]$. If $q_t^{\OSCI}\geq B$, we have either $S_t(X_t, q_t^{\OSCI})=0$ or $\cC_t(X_t, q_t^{\OSCI})=\mathcal{Y}$ (Definition~\ref{def:thr_proc}) and thus $\IER_t(\mathcal{R}^{\OSCI})=0$, so that IER control is trivial in that case. By contrast, since $\mathcal{C}_t(X_t,q_t^{\OSCI})$ is arbitrary for $q_t^{\OSCI}<0$, we do not seek an IER bound in this case, see \S~\ref{sec:method}. } 

\et{The error function $\Pi_t$ may take various forms depending on the setting; these  will be explored in depth in \S~\ref{sec:baoselect}, \S~\ref{sec:selectiveclassifiid} and \S~\ref{sec:selectiveNDiid}.
 Typically, under specific assumptions (e.g., iid and regularity),  we show therein that the (possibly random) $\Pi_t$ gets close to a (deterministic) benchmark function $\Pi^\bench$ when $t$ gets large. }
More generally, let us introduce the distance 
\begin{equation}
D_t:=\abs{\Pi_t(q_t^{\OSCI})-\Pi^\bench(q_t^{\OSCI})}\ind{q^\OSCI_t\in [0,B]}\leq \|\Pi_t-\Pi^\bench\|_\infty ,\:\:t\geq 1 \; , \label{equDtT}
\end{equation}
for some deterministic function $\Pi^\bench:[0,B]\rightarrow \R$ called the {\it benchmark IER function}. 

In the particular case where $\Pi_t$ is a deterministic function not depending on $t$ (iid and time-constant prediction rule, see Section~\ref{sec:baoselect}), we choose $\Pi^\bench:=\Pi_t$ to obtain $D_t=0$. \et{We refer to this case as the 'non-adaptive' case, as it generally  corresponds to a prediction rule that is fixed and only depends on the hold-out initial sample.} 
However, in general, $\Pi_t$ is a random function (e.g., 'adaptive' case where an updated prediction rule is used) so that the error $D_t>0$ should be taken into account in the analysis, and $\Pi^\bench$ is chosen to make this error as small as possible. \et{Typically, $\Pi^\bench$ is chosen as the limit of $\Pi_t$ so that $D_t$ measures the error from using the current estimated score/selection rule instead of the limiting oracle rule at the current threshold.}

\et{Since we aim to prove that $\IER_t(\mathcal{R}^{\OSCI})=\Pi_t(q_t^{\OSCI})$ is close to $\alpha$,  a natural candidate $q^\bench_{\alpha}$ for the limit of $q_t^{\OSCI}$ when $t$ is large is to solve $\Pi^{\bench} (q)=\alpha$ with respect to $q$. To this end, we introduce the following assumption.}

\begin{assumption}\label{ass:easy}
The benchmark function  $\Pi^\bench : q \in [0,B] \mapsto \Pi^\bench(q)\in [0,1]$ is a (deterministic) function with $\Pi^\bench(q^\bench_{\alpha})=\alpha$ for some $q^\bench_{\alpha}\in (0,B)$. In addition, $\Pi^\bench$
is nonincreasing on $[0,B]$ and continuously differentiable on some compact $[\ell_\alpha,u_\alpha]\subset (0,B)$ with $\ell_\alpha<q^\bench_{\alpha}<u_\alpha$ and $(\Pi^\bench)'(q^\bench_{\alpha})<0$.
\end{assumption}

As we will see further on, this regularity assumption is  mild and will be satisfied in our application cases (restricting the $\alpha$ range in some cases). 
It can also be relaxed (in particular the global non-increasingness can be dropped) and we refer to \S~\ref{sec:genapprox} for more details. \et{As we will see in the application cases (\S~\ref{sec:regression},\S~\ref{sec:selectiveclassifiid},\S~\ref{sec:selectiveNDiid}), the adaptation stage is necessary to ensure that $q_t^{\OSCI}$ converges to some  $q^\bench_{\alpha}$ that correctly incorporates the true model parameters.}

\subsection{Main results}\label{sec:mainresults}

We start by recalling the following classical conditions (see, e.g., \cite{chen2001stochastic}) on the step size sequence of the Robbins-Monro algorithm:
\begin{equation}
    \label{equ-stepseq}
    \sum_{j\geq 1}  \gamma_{j} = \infty \:\:\mbox{ and }\:\: \sum_{j\geq 1}  \gamma^2_{j} < \infty \; .
\end{equation}
The following result holds (see \S~\ref{proof:th:consistency} for a proof).
\begin{theorem}
    \label{th:consistency}[Strong consistency]
   Consider the \method\  threshold sequence $(q_t^{\OSCI})_{t\geq 1}$  with a nonincreasing step size sequence $(\gamma_j)_{j \geq 1}$  satisfying \eqref{equ-stepseq} and starting point $q_1^{\OSCI}\in [0,B]$. 
Let $\Pi^{\bench}$ be a function satisfying Assumption~\ref{ass:easy} with some $q^\bench_{\alpha}\in (0,B)$.   Recall that $D_t$ and $J(t)$ are defined in  \eqref{equDtT} and  \eqref{def:Jt}, respectively. 
Then,  almost surely on the event $\{\lim_{t\to \infty}J(t)=+\infty\}\cap \{\lim_{t\to \infty}D_t=0\}$,  we have 
$\lim_{t\to \infty}\qt_t^{\OSCI}=q^{\bench}_{\alpha}.$
\end{theorem}

Theorem~\ref{th:consistency} is well known in the case of no selective inference, with the slightly relaxed assumption that $\Pi^\bench$ is continuous and only crosses $\alpha$ in $q^\bench_{\alpha}$ with $\Pi^\bench(0)>\alpha$ and $\Pi^\bench(B)<\alpha$ \citep{chen2001stochastic}\footnote{In light of our proof and our relaxed assumption Assumption~\ref{ass:crossing}, Theorem~\ref{th:consistency} also extends to that assumption.}. Hence, Theorem~\ref{th:consistency}  provides initial insights into the conditions under  which  convergence can be expected to hold under selection. First, the number of selections should be large enough. To quantify this, 
let us introduce 
$(\underline{J}(t))_{t\geq 1}$ and $(\bar{J}(t))_{t\geq 1}$ two increasing integer {\it deterministic} sequences with $\underline{J}(t)\leq t$ for all $t\geq 1$ and let
\begin{equation}
    \label{equprobaJtsympa}
    \omega(t) =1-\P\big(\underline{J}(t)\leq J(t) \leq \bar{J}(t)\big), \:\:t\geq 1 \; .
\end{equation}
Second, the (random) quantity $D_t$ \eqref{equDtT} should be small enough. 
This is quantified via a {\it deterministic} sequence $(\bar{D}_t)_{t\geq 1}$, tending to zero and eventually nonincreasing, such that there exists a constant $C_1>0$ such that 
\begin{equation}
    \label{equ:Dtcondition}
    \P( D_t> C_1 \bar{D}_t)\leq \frac{6\eta_t}{\pi^2t^2},    \mbox{ for all $t\geq 1 \; ,$}
\end{equation}
 for $(\eta_t)_{t\geq 1}\in (0,1)$ some nonincreasing sequence with $\log(1/\eta_t)\leq \log(t)$. In the non-adaptive case, $D_t=0$ and we can always choose $\bar{D}_t=0$, $\eta_t=1/t$ to ensure \eqref{equ:Dtcondition}.
\et{In the adaptive case, a typical choice of $\eta_t$ will be $C_2/t$ for some constant $C_2\geq 1$.}

\begin{theorem}
    \label{th:rateconvfull}[In-probability rate]
Consider the \method\  procedure $\mathcal{R}^{\OSCI}$ with any selection/prediction set rule and a threshold sequence $(q_t^{\OSCI})_{t\geq 1}$ \eqref{equACIquantile} with step size sequence $\gamma_j= c j^{-\beta}$, \et{$c>0$,} $\beta\in (1/2,1)$ and starting point $q_1^{\OSCI}\in [0,B]$.  
  Let $\Pi^{\bench}$ be a function satisfying Assumption~\ref{ass:easy} with some $q^\bench_{\alpha}\in (0,B)$ and the corresponding $D_t$ being defined by \eqref{equDtT}, which is assumed to be bounded by $(\bar{D}_t)_{t\geq 1}$ as in \eqref{equ:Dtcondition}, with a rate $(\eta_t)_{t\geq 1}\in (0,1)$.
 Then  there exist constants\footnote{In all the results of the main text, the constants may depend on everything but the time $t$. We refer the reader to the non-asymptotic bounds of \S~\ref{sec:genapprox} to track the exact dependence of these constants with respect to the other quantities.} $T_0\geq 1$ and $C_0>0$ such that if $t\geq T_0$,
\begin{equation}\label{equ-conviid-specialcases-full}
\P(|\IER_t(\mathcal{R}^{\OSCI}) - \alpha|+ \abs{\qt_t^{\OSCI}-q^{\bench}_{\alpha}} \leq C_0 (\rho_t+ \bar{D}_{k(t)}))\geq 1-3\eta_{k(t)} -\omega(k(t))-\omega(t+1) \; ,
\end{equation}
where $(\omega(t))_{t\geq 1}$ is given by \eqref{equprobaJtsympa} and $\rho_t$, $k(t)$, $\underline{J}(t)$, $\bar{J}(t)$ are given by one of the following:
\begin{itemize}
    \item[(i)] \et{Loose $J(t)$ concentration bounds:} $\rho_t= (\log t)^{1/2} t^{1/2-\beta} + t^{\beta-1} $ with  \et{$\underline{J}(t)\sim \underline{r} t$, $\underline{r}\in(0,1)$,} $\bar{J}(t)=t$ and $k(t)=\lfloor c't/2\rfloor$, $c'\in(0,1)$;
    \item[(ii)] \et{Sharper $J(t)$ concentration bounds:}
    $\rho_t= (\log t)^{1/2} \Delta(t)^{1/2} \bar{J}(t)^{-\beta} 
 + \Delta(t)^{-1}  \bar{J}(t)^{\beta} $ with $\Delta(t):=\underline{J}(t+1)-\bar{J}(k(t))$ and any choice of $k(t)$ with $1\ll k(t)< t$, $\underline{J}(t)$, $\bar{J}(t)$ such that 
    $\underline{J}(t+1)\sim \underline{J}(k(t))\sim\bar{J}(t+1)\sim \bar{J}(k(t))$, $\Delta(t)\sim \bar{J}(t+1)-\underline{J}(k(t))$, $\Delta(t)\to \infty$ and $\bar{J}(t)^{\beta}\ll \Delta(t) \ll \bar{J}(t)$.
\end{itemize}
\end{theorem}
Theorem~\ref{th:rateconvfull} is proved in \S~\ref{sec:proofcor:rateconvfull}. 
\et{There are two ways to use Theorem~\ref{th:rateconvfull}, either with (i) for linear-order lower bound on $J(t)$ or (ii) for sharper concentration of $J(t)$. Both are valid, but (ii) leads to a better rate than (i), while (i) applies to a wider range of situations (see Remark~\ref{rem:ratetypical} below). Hence, choosing either of these two cases will depend on the considered type of selection: informative selection rules will mostly rely on (i), while other -- easier -- selection rules will mostly rely on (ii).}
\et{Theorem~\ref{th:rateconvfull} extends the results in \cite{angelopoulos2024online} to the selective inference case. In fact, even in case of no selective inference ($S_t(X,q)\equiv 1$) and no adaptation ($D_t=0$), Theorem~\ref{th:rateconvfull} shows an in-probability bound (rather than an in-expectation bound in \cite{angelopoulos2024online}), which is new in the literature to the best of our knowledge.}

\begin{remark}[Typical rate under concentration]\label{rem:ratetypical}
Under a concentration property for $J(t)$ when $t$ grows, a typical choice for the bounds $\underline{J}(t)$ and $\bar{J}(t)$ are given by $\lfloor rt \pm C(\log t)^{\lambda'}t^{\lambda} \rfloor$ for some selection rate $r\in (0,1)$ and concentration rate exponent $\lambda\in [1/2,1)$ (with a constant $C>0$, $\lambda'>0$ and $\omega(t)\to 0$). Then Theorem~\ref{th:rateconvfull} (case (ii)) with $k(t)=t-(1\vee \lfloor t^{4\beta/3}\rfloor)$ leads to  $\rho_t(\beta)= (\log t)^{1/2} t^{-\beta/3}$ when $\beta\in (1/2, 3/4)$, with $\beta>3\lambda/4$. It is clear that this rate is better than the one produced by Theorem~\ref{th:rateconvfull} in the case (i) on the range $3\lambda/4<\beta<3/4$. Theorem~\ref{th:rateconvfull} then leads to the optimized rate $\rho_t(\beta)+ \bar{D}_{k(t)} $ with
\begin{align}\label{equ:raterho-typicalconc}
\rho_t(\beta)=\rho_t(\beta;\beta_0)&:= \left\{\begin{array}{cc}
   (\log t)^{1/2} t^{-\beta/3}    & \mbox{ if $\beta\in (\beta_0,3/4)$;} \\
   (\log t)^{1/2} t^{1/2-\beta} + t^{\beta-1}   & \mbox{ if $\beta\in (1/2,\beta_0] \cup [3/4,1)$,} \\
\end{array}\right.
\end{align}
with $\beta_0:=(1/2)\vee (3\lambda/4)$ (hence $\beta_0=1/2$ and $(1/2,\beta_0]=\emptyset$ if the concentration holds with $\lambda\leq 2/3$).
\end{remark}

\begin{remark}[Choosing $\beta$]\label{rem:choosebeta}
    Our in-probability rate provides guidelines for choosing $\beta$ in the typical case of Remark~\ref{rem:ratetypical}:
    \begin{itemize}
        \item[-] simple choice: minimizing the rate $\rho_t(\beta)$ in \eqref{equ:raterho-typicalconc} with respect to $\beta$ suggests to choose $\beta=3/4$, to obtain $\rho_t(\beta)=(\log t)^{1/2} t^{-1/4}$. This is our default recommendation.
        \item[-] adaptive choice: recall the smaller $\beta$, the better the bound \eqref{FCPbound_decre} on $\FP-\alpha$ (for the same selection number). Hence, if the order of $\bar{D}_{k(t)}$ is known, a choice is to take $\beta\in (1/2,3/4]$ the smallest such that $\rho_t(\beta)\lesssim \bar{D}_{k(t)}$. 
    \end{itemize}
    Note that both choices coincide if $\bar{D}_{k(t)}\ll (\log t)^{1/2}t^{-1/4}$. In this case, we obtain the same rate $t^{-1/4}$ (up to a log term) both for the $\FP-\alpha$ and $\abs{\IER-\alpha}$ bounds. 
\end{remark}

\begin{theorem}
\label{th:genL2-first}[In-expectation rate]
Consider the \method\  procedure $\mathcal{R}^{\OSCI}$ with any selection/prediction set rule and a threshold sequence $(q_t^{\OSCI})_{t\geq 1}$ \eqref{equACIquantile} with step size sequence $\gamma_j= c j^{-\beta}$, $\beta\in (0,1)$ and starting point $q_1^{\OSCI}\in [0,B]$. 
Assume that for all $t\geq 1$, the selection indicator $S_t(X_t,q^{\OSCI}_{t})$ at time $t$ is independent of the past $\mathcal{F}_{t-1}$. 
  Consider some function $\Pi^{\bench}$ satisfying Assumption~\ref{ass:easy} with some $q^\bench_{\alpha}\in (0,B)$ and the corresponding $D_t$ being defined by \eqref{equDtT} and consider $(\underline{J}(t))_{t\geq 1}$ and $(\bar{J}(t))_{t\geq 1}$ two sequences tending to infinity, and $(\omega(t))_{t\geq 1}$ as in \eqref{equprobaJtsympa}.
\et{Then, there exists $C_0>0$ and $T_0\geq 1$ such that for any nondecreasing sequence $k(t)<t$, we have for any $t\geq 1$, with $k(t)\geq T_0$,
    \begin{align}
&\e{\abs{\IER_t(\mathcal{R}^{\OSCI}) - \alpha}\ind{q_t^\OSCI\in [0,B]}}^2+\e{\paren{q^{\OSCI}_{t}-q^{\bench}_{\alpha}}^2 }\nonumber\\
&\leq  C_0 \omega(t)+ \frac{C_0 (k(t))^{\beta}}{(\underline{J}(t))^{\beta}}\paren{1+\sum_{k=1}^{t-1} \widetilde{D}_k} + C_0\probp{J(k(t))< T_0}
 ,\label{equ:genL2withoutt:rate}
\end{align}}
where $(\widetilde{D}_t)_{t\geq 1}$ is a (deterministic) sequence tending to zero and eventually decreasing such that $\E[D_t]\leq \widetilde{D}_t$ for all $t\geq 1$.
\end{theorem}
Theorem~\ref{th:genL2-first} is proved in \S~\ref{sec:cor:genL2}. In case of no selective inference ($S_t(X,q)\equiv 1$) and no adaptation ($D_t=0$) the bound is of order $t^{-\beta}$ \et{(take $k(t)=T_0$, so that $J(k(t))=k(t)=T_0$ almost surely)} and we recover Theorem 2.2 of  \cite{costa2020non} and Proposition~2 of \cite{angelopoulos2024online}. \et{Theorem~\ref{th:genL2-first} thus extends these works to the case of selection and adaptive scores, which added significant technicalities in the analysis and corresponds  in \eqref{equ:genL2withoutt:rate} to the remainder terms $\widetilde{D}_t$ and $\underline{J}(t),\omega(t)$, respectively. The additional probability $\probp{J(k(t))< T_0}$ can in general be  bounded with a concentration argument --  similarly to $\omega(t)$ -- and by suitably choosing $k(t)$ (ideally  of order logarithmic with $t$).} 
\et{Note that the assumption on the selection rule is stronger in Theorem~\ref{th:genL2-first} than in Theorem~\ref{th:rateconvfull}:  $S_t(X_t,q^{\OSCI}_{t})$ should be independent of the past $\mathcal{F}_{t-1}$, thus  it excludes any adaptive history-dependent selection rule. 
This assumption holds for time-constant selection rules and under independence of the stream.
}

\et{
\begin{remark}[Price of adaptation]
  \label{rem:adaptation}  
  The price of adaptation can be identified easily in the remainder terms of each bound. For Theorem~\ref{th:rateconvfull} (high-probability, both in cases (i) and (ii)), the adaptation stage is quantified by $D_t>0$. Markedly, it deteriorates the rate only if $\bar{D}_{k(t)}$ dominates $\rho_t$, because it comes in an additive manner. 
  By contrast, for Theorem~\ref{th:genL2-first} (in-expectation), 
  the term $1+\sum_{k=1}^{t-1} {\widetilde{D}_k}$ introduced by adaptation is multiplicative and not additive to the rate $(\underline{J}(t))^{-\beta}$, typically of order $t^{-\beta}$. Hence, the adaptation stage can deteriorate the rate  even if $\sum_{k=1}^{t-1} {\widetilde{D}_k}$ is tending slowly to infinity. This generally requires taking $\beta$ large enough to compensate for the adaptation error.
\end{remark}
}

\et{
\begin{remark}[Threshold convergence]
  \label{rem:thconv}  
In Theorems~\ref{th:rateconvfull}~and~\ref{th:genL2-first}, the two quantities $|\IER_t(\mathcal{R}^{\OSCI}) - \alpha|$ and $\abs{\qt_t^{\OSCI}-q^{\bench}_{\alpha}}$ are both shown to be small. However, they do not have a symmetric role: while our first concern is to show that $|\IER_t(\mathcal{R}^{\OSCI}) - \alpha|$ is small, the fact that $\qt_t^{\OSCI}$ stabilizes around $q^{\bench}_{\alpha}$ is more technical and can be seen as a restriction since these results will {\it de facto} provide IER control only in those situations. Nevertheless, we note that the most general version of our bound can involve a dynamic choice of $q^{\bench}_{\alpha}$, that is, $q^{\bench}_{t,\alpha}$ depending on the time point $t$, see Theorems~\ref{th:almostsureconv}~and~\ref{th:genL2} in \S~\ref{sec:genapproxbounds}. This leaves the door open for deriving IER control in more dynamic regimes.
\end{remark}
}

In summary, Theorems~\ref{th:rateconvfull}~and~\ref{th:genL2-first}  showed that, provided that suitable bounding sequences 
$\underline{J}(t)$, $\bar{J}(t)$ can be built for $J(t)$, and provided that $D_t$ is appropriately vanishing,  we can obtain explicit convergence rates of the instantaneous error rate to $\alpha$ and of the threshold sequence to some $q^\bench_\alpha$.

\subsection{\et{Result in the informative case}}\label{sec:incompresscritic}

\et{We consider an informative selection rule of the form $\St_t(x, \qt) = \ind{W_t(x) >q}$, for some test statistic $W_t(x)$, see \eqref{eq:selectionstat}. 
In that case, the selection $S_t(X_t,q_t^\OSCI)$  depends on the actual threshold $q_t^\OSCI$ at time $t$.  Since $S_t(X,q_t^\OSCI)$ necessarily depends on $\mathcal{F}_{t-1}$ for an informative selection, the in-expectation result (Theorem~\ref{th:genL2-first}) is out of reach, and we should make use of the in-probability bound (Theorem~\ref{th:rateconvfull}).
}

\et{To quantify theoretically the performance  of the algorithm \method\  in the informative case we introduce an incompressible  error probability term, that we will denote by $\delta$. 
More formally, assume that, when $t$ is large, the test statistic $W_t(x)$ in \eqref{eq:selectionstat} gets close to $W(x)$, a random variable satisfying the following assumption, which also involves a new quantity $B^*\leq B$:
\begin{assumption}
    \label{ass:support}
$W(X)$ has a continuous density supported in $[\kappa,B^*]$, which is positive on $(\kappa,B^*)$, for some $0\leq \kappa<B^*\leq B$. 
\end{assumption}
The failure probability of \method\  is then defined by
\begin{equation}
    \label{deltaprobaJtfinite}
    \delta:= \inf_{q\in [0,B^*)}\inf_{t\geq 1}\probp{\sup_{T\geq t} q_T^\OSCI>q}\in [0,1],
\end{equation}
which corresponds to the probability that the \method\:threshold sequence fails to eventually remain below a critical selection threshold.  
This failure probability $\delta$ is possibly nonzero (Lemma~\ref{lem:deltapositive}) and 'incompressible' in the sense that the only type of results which is achievable is that, $q_t^\OSCI$ converges appropriately with probability (roughly) greater than $ 1-\delta$ (Lemma~\ref{lem:impossible}).  
Intuitively, this comes from the fact that if $q_T^\OSCI>B^*$ when $T$ is large, then the procedure can only converge to a limiting point above $B^*$, which violates Assumption~\ref{ass:support}. 
}

\et{Similarly, given the benchmark function $\Pi^\bench$ and $q^{\bench}_{\alpha}$ such that $\Pi^\bench(q^{\bench}_{\alpha})=\alpha$, 
\method\: could not reach the optimal threshold if $q^{\bench}_{\alpha}\geq B^*$. This means that our analysis will be done in the `super-critical' regime where $\alpha$ is above the {\it critical level} defined by 
\begin{equation}
\label{equalphastargen}
\alpha^*:=\Pi^\bench(B^*)\in [0,1].    
\end{equation}
This quantity is similar to the criticality phenomenon in batch multiple testing \citep{Chi2007}. We show here that it can also occur in the online setting when $\alpha^*=\Pi^\bench(B^*)>0$. The critical value $\alpha^*$ can be roughly interpreted as the smallest level that \method\:can achieve.
}

\et{Our analysis is able to show a convenient behavior of \method\: with probability (roughly) at least $1-\delta$, as the following result shows.
\begin{theorem}
[In probability rate for informative selection]
    \label{rem:coroonanEvent}
    Let $\alpha>\alpha^*$ for $\alpha^*$ the critical level \eqref{equalphastargen}. 
    Considering $\delta$ the incompressible probability of failure \eqref{deltaprobaJtfinite}, for any $\delta'>\delta$ there exists $T_0=T_0(\delta')$ and $\bar q=\bar q(\delta')\in (0,B^*)$ such that the following holds.
In the same setting of Theorem~\ref{th:rateconvfull} except that $\omega(t)$ in \eqref{equprobaJtsympa} and $\bar{D}_t$ in \eqref{equ:Dtcondition} are computed in restriction to respective events $\Omega_t=\set{\forall \ell\in [T_0,t], q^{\OSCI}_\ell\leq \bar q}\in \mathcal{F}_{t-1}$, that is for $t\geq 1$, 
\begin{equation}\label{equrestrictioncond}
    \omega(t) =\P(\Omega_t)-\P\big(\Omega_t ,\underline{J}(t)\leq J(t) \leq \bar{J}(t)\big); \:\:\P(\Omega_t, D_t> C_1 \bar{D}_t)\leq \frac{6\eta_t}{\pi^2t^2},
\end{equation}
we have for $t\geq T_0$, 
     \begin{equation*}
\P(|\IER_t(\mathcal{R}^{\OSCI}) - \alpha|+ \abs{\qt_t^{\OSCI}-q^{\bench}_{\alpha} } \leq C_0 (\rho_t+ \bar{D}_{k(t)}))\geq 1-C_0(\eta_{k(t)} +\omega(k(t))+\omega(t+1))-\delta',
\end{equation*} 
for $\rho_t$, $k(t)$ given either by case (i) or case (ii) in the statement of Theorem~\ref{th:rateconvfull}.
\end{theorem}
Theorem~\ref{rem:coroonanEvent} is proved in \S~\ref{sec:proof:rem:coroonanEvent} and will be used for informative selections to obtain convergence rates, up to the incompressible term $\delta$.
For brevity, we will sometimes write that the IER “does converge” in this informative setting. This is understood to mean that the convergence holds outside an event of probability approximately $\delta$.
}

\et{
\begin{remark}[Diminishing $\delta$]
A finite-horizon proxy of $\delta$ is computed in \S~\ref{sec:deltaprobaJtfinite} for a grid of hyperparameters.
   A way to make $\delta$ as small as possible is to optimize the hyperparameter $c,\beta,q_1$ of \method, as we propose in \S~\ref{sec:select_process}, possibly in combination with the restarting technique presented in Remark~\ref{rem:rerun}. 
\end{remark}
}

\section{\et{Application to regression}}\label{sec:regression}
\label{sec:baoselect}

\subsection{Setting and assumptions}

Throughout this section, we consider the regression case $\cY=\R$, $\cX=[0,1]^d$, as already exemplified in \S~\ref{subsec:prelim:setting} with $B=1$. We consider a prediction set rule $\cC_t(\cdot,\cdot)$ as in \eqref{eq:set_pred_funct} with 
a two-sided $V_t(x,y)$ as in \eqref{regressionVadaptive} (a one-sided version is also considered in \S~\ref{rk:plb}), for some estimators $\hat{\mu}_t(x)$ (resp. $\hat{\sigma}_t^{2}(x)>0$) of the conditional expectation (resp. variance) of $Y_{t}$ given $X_{t}=x$, so that the prediction set reduces to a prediction interval \eqref{equCtregression}. 

Except in \S~\ref{sec:AR} (autoregressive case), we focus in this section on the following iid online regression model 
\begin{equation}\label{regression:model}
    Y_t = \mu(X_t) + \sigma(X_t) \xi_t, \:\:t\geq 1,
\end{equation}
for iid couples $(X_t,\xi_t)$, $t\geq 1$, following an unknown distribution, and distributed as $(X,\xi)$, for which $X$ is independent of $\xi$. 
\begin{assumption}\label{ass:densityf}
In the iid regression model \eqref{regression:model}, the distribution of the noise $\xi$ has a symmetric continuous density $f$ with respect to the Lebesgue measure which is positive on its support with $M_1:=\sup_{y\in \R} f(y)<\infty$,   $M_2:=\sup_{y\in \R} \{|y|f(y)\}<\infty$,   $M_3:=\sup_{x\in \R^d}(1/\sigma(x))<\infty$ and $\E \xi^2 =1$.
\end{assumption}
Note that these assumptions imply in particular $\mu(x)=\E[Y_t\:|\: X_t=x]$ and $\sigma(x)=\V^{1/2}[Y_t\:|\: X_t=x]$ and we further assume that $\hat{\mu}_t$ and $\hat{\sigma}_t$ converge to these respective quantities at some rate: for  $t$ large enough and some constant $C_1>0$,
\begin{align}
&\P[\E [|\hat{\mu}_t(X_t)-{\mu}(X_t)| + |\hat{\sigma}_t(X_t)-{\sigma}(X_t)|\mid \mathcal{F}_{t-1}] \geq C_1 (\log t)^{\nu_0} t^{-\nu}] \leq \frac{3}{\pi^2t^3};\label{ass:esti}\\
& \E[(\hat{\mu}_t(X_t)-{\mu}(X_t))^2] + \E[(\hat{\sigma}_t(X_t)-{\sigma}(X_t))^2] \leq C_1 (\log t)^{2\tau_0}t^{-2\tau},\label{ass:esti:average}
\end{align}
for some $\nu_0,\tau_0\geq 0$ and $\nu,\tau\in (0,1/2]$.
\et{Guarantee \eqref{ass:esti} is a high-probability conditional L1-type bound and \eqref{ass:esti:average} is an unconditional $L^2$ bound.}
As a typical instance, if  $\mu(\cdot)$ and $\sigma^2(\cdot)$ are Lipschitz functions, there exist kernel-based estimators $\hat{\mu}_t$ and $\hat{\sigma}_t$ (based upon the $t-1$ sample $\mathcal{F}_{t-1}$) such that \eqref{ass:esti} and \eqref{ass:esti:average} holds with $\nu=\tau=1/(2+d)$ \citep{tsybakov2008nonparametric, wang2008effect}. \et{Similarly, neural-network estimators can be used with other regularity classes for $\mu(\cdot)$ and $\sigma^2(\cdot)$ \citep{schmidt2020nonparametric}.}

The IER function $ \Pi_t $ \eqref{PiCtSt} gets close in this case to the IER benchmark function
\begin{align}
    \label{PizeroClassicalscore}
    \Pi^\bench ( q)&=2\bar F (\bar \Phi^{-1}((1-q)/2)); \quad q^\bench_{\alpha}=1-2\bar{\Phi}(\ol{F}^{-1}(\alpha/2)),
\end{align}
respectively, 
where  $\bar{F}$ denotes the upper tail of the noise $\xi$ distribution. Clearly, $\Pi^\bench$ satisfies Assumption~\ref{ass:easy}. The distance $D_t$ \eqref{equDtT} between $\Pi_t(q)$ and 
$\Pi^\bench$ can be bounded from \eqref{ass:esti}-\eqref{ass:esti:average} (see Proposition~\ref{prop:Dtiideasyselect}). 

We apply our theoretical analysis in the next subsections for 
various selection rules: 
 time-constant, adaptive and informative.
 
\et{
By \eqref{PizeroClassicalscore}, the oracle threshold is common to all selection rules that we will consider in the next subsections (except in \S~\ref{rk:plb} for which a one-sided score is considered), which may seem surprising at first. This is specific to this regression setting and to the form of the score function used. This is not the case  in the classification/testing framework, see \S~\ref{sec:selectiveclassifiid} and \S~\ref{sec:selectiveNDiid}).  Nevertheless, the threshold $\qt_t^{\OSCI}$ still uses the specificity of the selection rule to ensure the adversarial FCP control \eqref{FCPbound_decre}.  Algorithm~\ref{alg:PredictXoriented} implements the selective regression procedure, denoted by  $\mathcal R^{SCI}$. 
}

\subsection{Time-constant selection}\label{sec:appliiidregression_timeconstant}\label{sec:xoriented:optimality}

We now apply our main results (Theorems~\ref{th:rateconvfull}~and~\ref{th:genL2-first}) in case of a time-constant $X$-oriented selection rule,  \et{where $S_t(x,q)$ is of the form $S(x)$}.

\begin{corollary}\label{cor:rateconveasyselect}
    Consider the online iid regression model \eqref{regression:model} satisfying Assumption~\ref{ass:densityf} and consider $q^{\bench}_\alpha$ given by \eqref{PizeroClassicalscore}.
Consider $\mathcal{R}^{\OSCI}$ 
using  a time-constant $X$-oriented selection (Algorithm~\ref{alg:PredictXoriented}) with $\P(S(X_1)=1)$ fixed in $(0,1]$, $\mathcal{F}_{t-1}$-measurable estimators $\hat{\mu}_t, \hat{\sigma}_t$, and a threshold sequence $(q_t^{\OSCI})_{t\geq 1}$ \eqref{equACIquantile}  with a step size sequence $\gamma_j= c j^{-\beta}$, $\beta\in (1/2,1)$, and  starting point $q_1^{\OSCI}\in [0,B]$. Then 
there exist
    $T_0\geq 1$ and $C_0>0$ such that for $t\geq T_0$, the following holds:
    \begin{itemize}
        \item In-probability rate: 
        assuming that the estimators satisfy \eqref{ass:esti},  we have
\begin{equation}\label{equ-conviid-specialcases-easyselect}
\probp{|\IER_t(\mathcal{R}^{\OSCI}) - \alpha|+ |\qt_t^{\OSCI}-q^{\bench}_\alpha| \leq C_0 (\rho_t(\beta)+(\log t)^{\nu_0} t^{-\nu})}\geq 1-1/t,
\end{equation}
for $\rho_t(\beta)=\rho_t(\beta;1/2)$ given by \eqref{equ:raterho-typicalconc}.
\item In-expectation rate: choosing $\beta>1- \tau$ and assuming that the estimators satisfy \eqref{ass:esti} and \eqref{ass:esti:average}, we have 
\begin{equation}\label{equ:MSEconv-regressioneasy}
\paren{\e{\abs{\IER_t(\mathcal{R}^{\OSCI}) - \alpha}}}^2 + \e{\paren{q_{t}^{\OSCI}-q^{\bench}_\alpha}^2}
\leq C_0 (\log t)^{\tau_0\et{+\beta}} t^{1- \tau-\beta }.    
\end{equation}
    \end{itemize}
\end{corollary}

Corollary~\ref{cor:rateconveasyselect} is proved in \S~\ref{proof:cor:rateconveasyselect}. The proof is based on more general results valid in an iid model (not necessarily a regression model), see \S~\ref{sec:appli_timeconstant}. 
According to Remark~\ref{rem:choosebeta}, if $\nu$ is known, an adaptive choice of $\beta$ as a function of $\nu$ can be done in \eqref{equ-conviid-specialcases-easyselect}, by letting $\beta=(3\nu)\wedge (3/4)\vee (1/2+\epsilon)$, for $\epsilon$ arbitrary small. For instance,  when $\nu=1/(2+d)$ for a dimension $d\geq 1$, this gives $\beta=3/4$ for $d\in \{1,2\}$, $\beta=3/5$ for $d=3$ and $\beta=1/2+\epsilon$ for $d\geq 4$.

\begin{remark}\label{rem:condcovopt}
    The oracle threshold $q^\bench_\alpha$ in \eqref{PizeroClassicalscore} does not depend on the selection, so in addition to the convergence property of Corollary~\ref{cor:rateconveasyselect}, $\mathcal{R}^{\OSCI}$ can produce prediction sets with asymptotic {\it conditional} coverage, that is, conditionally on $X=x$ for any $x$. In addition, any other method satisfying this asymptotic conditional coverage should produce prediction sets with a larger averaged length, up to some vanishing terms. This conditional coverage and optimality result are postponed into \S~\ref{sec:condcovopt} for space reasons. 
\end{remark}

\subsection{\et{Selection for large predicted outcomes}}\label{sec:adaptiveXorientedSelection}

We consider here the case where the user wants to select examples whose predicted outcome exceeds $y_0$, that is: 
\begin{equation}
    \label{equ:adapselectextreme}
    S_t(x)=\ind{\hat{\mu}_t(x)\geq y_0},
\end{equation}
for a value $y_0\in \R$.
To use Theorem~\ref{th:rateconvfull} in this context, we should obtain high probability bounds for the selection number $J(t)$ \eqref{def:Jt}. Since $\P(\hat{\mu}_t(X)\geq y_0\:|\: \mathcal{F}_{t-1})$ gets close to $\P(\mu(X)\geq y_0)$ when $t$ gets large, this leads to the following assumption.
\begin{assumption}\label{ass:densitymu}
In the iid regression model~\eqref{regression:model}, $\mu(X)$ has a positive density in a neighborhood of $y_0$ which is upper bounded by $M_4>0$.
\end{assumption}

\begin{corollary}\label{cor:rateconvadaptselect}
    Consider the online iid regression model \eqref{regression:model} satisfying Assumptions~\ref{ass:densityf}~and~\ref{ass:densitymu}.
   Consider $\mathcal{R}^{\OSCI}$ 
   (Algorithm~\ref{alg:PredictXoriented}), 
   using   
    the adaptive $X$-oriented selection \eqref{equ:adapselectextreme}, $\mathcal{F}_{t-1}$-measurable estimators $\hat{\mu}_t, \hat{\sigma}_t$ satisfying \eqref{ass:esti}, and a threshold sequence $(q_t^{\OSCI})_{t\geq 1}$ \eqref{equACIquantile} with a step size sequence $\gamma_j= c j^{-\beta}$, $\beta\in (1/2,1)$, and  starting point $q_1^{\OSCI}\in [0,B]$. Then 
      there exist
     $T_0\geq 1$ and $C_0>0$ such that for $t\geq T_0$, 
  \begin{align}\label{equ-conviid-specialcases-adaptselect}
\probp{|\IER_t(\mathcal{R}^{\OSCI}) - \alpha|+ |\qt_t^{\OSCI}-q^{\bench}_\alpha| \leq C_0 (\rho_t(\beta)+(\log t)^{\nu_0} t^{-\nu})}\geq 1-C_0/t^{1/2},
\end{align}
\et{with $\rho_t(\beta)=\rho_t(\beta;\beta_0)$ as in \eqref{equ:raterho-typicalconc} with $\beta_0=(3/4)(1-\nu/2)\in [9/16,3/4]$.}
\end{corollary}

Corollary~\ref{cor:rateconvadaptselect} is proved in \S~\ref{sec:proofadaptiveXorientedSelection}.
In a nutshell, the difference with the time-constant selections considered in \S~\ref{sec:appliiidregression_timeconstant} is that the estimation error slows down the concentration of $J(t)$ here, which impacts the convergence rate when $\beta\leq \beta_0$.

\subsection{\et{Length restricted prediction interval}}\label{sec:lengthrestriction}

\et{We now consider the informative rule  where the user reports the prediction interval \eqref{equCtregression} only if it is of Lebesgue measure smaller than some $2\lambda_0>0$, that is, $ S_t(x,q)=\ind{\hat{\sigma}_t(x) \Phi^{-1}(q/2+1/2)<\lambda_0}=\ind{W_t(x)>q}$ is of the form \eqref{eq:selectionstat} for
\begin{equation}
    \label{equ:Wtselectlength}
  W_t(x)=2\Phi (\lambda_0/\hat{\sigma}_t(x))-1.
\end{equation}
The theoretical properties are derived from the general methodology developed in \S~\ref{sec:incompresscritic}. First, we define the limit of $W_t(x)$ as
\begin{equation}
    \label{equ:Wselectlength}
  W(x)=2\Phi (\lambda_0/\sigma(x))-1,
\end{equation}
and we note that Assumption~\ref{ass:support} is satisfied with $\kappa=2\Phi (\lambda_0/\sigma_{\max})-1$ and $B^*=2\Phi (\lambda_0/\sigma_{\min})-1$  if the following assumption holds: 
\begin{assumption}
    \label{ass:supportsigma}
$\sigma(X)$ has a continuous positive density supported in $[\sigma_{\min},\sigma_{\max}]$ for some $0<\sigma_{\min}< \sigma_{\max}<\infty$.
\end{assumption}
The critical value $\alpha^*$ in \eqref{equalphastargen} is thus $\alpha^*=\Pi^\bench(B^*)=2\bar F(\lambda_0/\sigma_{\min})$. %
Note that $\alpha^*$ increases as $\lambda_0$ decreases or $\sigma_{\min}$ increases, as expected because the problem is more demanding when the user asks for more accuracy or data are more variable. 
Applying Theorem~\ref{rem:coroonanEvent} leads to the following result.
\begin{corollary} \label{cor:restrictedlength}
Consider the online iid regression model \eqref{regression:model} satisfying Assumptions~\ref{ass:densityf}~and~\ref{ass:supportsigma}.
       Let $\alpha\in (\alpha^*,1)$ and consider $\delta$ the incompressible rate \eqref{deltaprobaJtfinite}, for any $\delta'>\delta$ there exists $T_0=T_0(\delta')$ and $\bar q=\bar q(\delta')\in (0,B^*)$ such that the following holds.
The procedure $\mathcal{R}^{\OSCI}$ (Algorithm~\ref{alg:PredictXoriented}), 
   using   
    the informative selection \eqref{eq:selectionstat} with $W_t(x)$ in \eqref{equ:Wtselectlength}, $\mathcal{F}_{t-1}$-measurable estimators $\hat{\mu}_t, \hat{\sigma}_t$ satisfying \eqref{ass:esti}, and a threshold sequence $(q_t^{\OSCI})_{t\geq 1}$ \eqref{equACIquantile} with a step size sequence $\gamma_j= c j^{-\beta}$, $\beta\in (1/2,1)$, and  starting point $q_1^{\OSCI}\in [0,B]$ is such that there exists
      $C_0>0$ such that for all $t\geq T_0$, 
  \begin{align}\label{equ-convlength}
\probp{|\IER_t(\mathcal{R}^{\OSCI}) - \alpha|+ |\qt_t^{\OSCI}-q^{\bench}_\alpha| \leq C_0 (\rho_t(\beta)+(\log t)^{\nu_0} t^{-\nu})}\geq 1-\delta'-C_0/t, 
\end{align}
for $\rho_t(\beta) = \rho_t(\beta;3/4)=(\log t)^{1/2}t^{1/2-\beta} + t^{\beta-1} $ as in \eqref{equ:raterho-typicalconc}.
\end{corollary}
The proof of Corollary~\ref{cor:restrictedlength} is provided in \S~\ref{proof:cor:restrictedlength}.
Recall that since the selection is informative, $\qt_t^{\OSCI}$ may not converge to the oracle threshold on an 'unfortunate event' of probability (roughly) at most $\delta$ (see \S~\ref{sec:incompresscritic}). To this extent, the term $1-\delta'$ in \eqref{equ-convlength} is unavoidable. The convergence rate $\rho_t(\beta;3/4)$ is also slightly slower than in the results of the previous subsections.  
}

\subsection{\et{Informative prediction lower bound}}\label{rk:plb}

\et{The selection considered in this subsection corresponds to the situation where the user is only interested in examples with predicted intervals above some benchmark value $y_0\in \R$. We consider here the one-sided predicted intervals $\cC_t(X_t,q_t^\OSCI)$ given by 
\begin{equation}
\mathcal{C}_t(x,q) = \big[\hat{\mu}_t(x)-\sigma_0\Phi^{-1}(q),+\infty\big)\label{equCtregressiononesided},
\end{equation}
if $q\in ( 0,1)$, with $\mathcal{C}_t(x,q) =[\hat{\mu}_t(x),+\infty)$ if $q\leq 0$ and $\mathcal{C}_t(x,q) =\R$ if $q\geq 1$. 
Let us consider the homoscedastic iid regression model with ${\sigma}(x)=\sigma_0>0$ for simplicity.
This corresponds to the informative selection 
$ S_t(x,q)=\ind{\hat{\mu}_t(x)-\sigma_0\Phi^{-1}(q)>y_0}=\ind{W_t(x)>q}$
of the form \eqref{eq:selectionstat} for
\begin{equation}
    \label{equ:Wtlb}
  W_t(x)= \Phi\bigg(\frac{\hat{\mu}_t(x)-y_0}{\sigma_0}\bigg).
\end{equation}
The limit of $W_t(x)$ is defined as
\begin{equation}
    \label{equ:Wlb}
 W(x)= \Phi\bigg(\frac{{\mu}(x)-y_0}{\sigma_0}\bigg),
\end{equation}
and we note that Assumption~\ref{ass:support} in \S~\ref{sec:incompresscritic} is satisfied with 
$\kappa=\Phi ((\mu_{\min}-y_0)/\sigma_0)$
and 
$B^*=\Phi ((\mu_{\max}-y_0)/\sigma_0)$
if the following assumption holds: 
\begin{assumption}
    \label{ass:supportmu}
 $\mu(X)$ has a continuous positive density supported in $[\mu_{\min},\mu_{\max}]$ for some $-\infty<\mu_{\min}<\mu_{\max}<\infty$ with $\mu_{\max}>y_0$. 
\end{assumption}
The benchmark function $\Pi^\bench$ is given here by 
$$\Pi^\bench(q)=\P\bigg( \Phi\bigg(\frac{{\mu}(X)-Y}{\sigma_0}\bigg)>q\:\bigg|\: \Phi\bigg(\frac{{\mu}(X)-y_0}{\sigma_0}\bigg)>q\bigg)= \P(\Phi(-\xi)>q)=F(-\Phi^{-1}(q)).$$
The critical value $\alpha^*$ in \eqref{equalphastargen} is thus $\alpha^*=\Pi^\bench(B^*)= F((y_0-\mu_{\max})/\sigma_0)$.
Note that $\alpha^*$ increases as $y_0$ increases or $\mu_{\max}$ decreases, which is well expected because the problem gets increasingly more difficult. 
Applying Theorem~\ref{rem:coroonanEvent} leads to the following result.
\begin{corollary} \label{cor:lb}
Consider the online homoscedastic iid regression model \eqref{regression:model} with ${\sigma}(x)=\sigma_0>0$ and satisfying Assumptions~\ref{ass:densityf}~and~\ref{ass:supportmu}.
       Let $\alpha\in (\alpha^*,1)$ and consider $\delta$ the incompressible rate \eqref{deltaprobaJtfinite}, for any $\delta'>\delta$ there exists $T_0=T_0(\delta')$ and $\bar q=\bar q(\delta')\in (0,B^*)$ such that the following holds.
The procedure $\mathcal{R}^{\OSCI}$ (Algorithm~\ref{alg:PredictXoriented}), 
   using   
    the informative selection \eqref{eq:selectionstat} with $W_t(x)$ in \eqref{equ:Wtlb},    a $\mathcal{F}_{t-1}$-measurable estimator $\hat{\mu}_t$ satisfying \eqref{ass:esti} (with $\hat{\sigma}_t\equiv \sigma_0$), and a threshold sequence $(q_t^{\OSCI})_{t\geq 1}$ \eqref{equACIquantile} with a step size sequence $\gamma_j= c j^{-\beta}$, $\beta\in (1/2,1)$, and  starting point $q_1^{\OSCI}\in [0,B]$ is such that there exist
      $C_0>0$ such that for all $t\geq T_0$, 
  \begin{align}\label{equ-convlb}
\probp{|\IER_t(\mathcal{R}^{\OSCI}) - \alpha|+ |\qt_t^{\OSCI}-q^{\bench}_\alpha| \leq C_0 (\rho_t(\beta)+(\log t)^{\nu_0} t^{-\nu})}\geq 1-\delta'-C_0/t,
\end{align}
for $\rho_t(\beta) = \rho_t(\beta;3/4)=(\log t)^{1/2}t^{1/2-\beta} + t^{\beta-1} $ as in \eqref{equ:raterho-typicalconc}.
\end{corollary}
The proof of Corollary~\ref{cor:lb} is provided in \S~\ref{proof:cor:lb}.
}

\subsection{Autoregressive model}\label{sec:AR}

Theorems~\ref{th:rateconvfull} and \ref{th:genL2-first} are valid beyond the iid model and are applied here in a 
specific dependent context, by considering an autoregressive time series. 
Let us consider an  $\mbox{AR}(d)$ process with Gaussian innovations, parameterized by $\varphi\in\R^d$ ($d\geq 1$):
\begin{equation}\label{regressionAR:model}
Y_t=\sum_{k=1}^d \varphi_k Y_{t-k}+Z_t,\:\:\:  t\in\mbz,
\end{equation}
where the variables $(Z_t)_{t\in\mbz}$ are iid $\mathcal{N}(0,\sigma^2)$ (with $\sigma^2=1$), and $(Y_t)_{t\in\mbz}$ is a stationary process. We assume  that the polynomial $P(z)=1-\sum_{k=1}^d \varphi_k z^k$ has no root of modulus smaller or equal to $1$, which ensures that the model \eqref{regressionAR:model} is well defined \citep{anderson2011statistical}. Also, for simplicity, we assume that $\E(Y_t)=0$ for all $t\in\mbz$.

We consider the setting where we predict $Y_t$ from the covariate $X_t=(Y_{t-1},\dots,Y_{t-d})\in\R^d$, which has the specific property to be $\mathcal{F}_{t-1}$-measurable, because we have $\mathcal{F}_{t-1}=\sigma(Y_{-(n-1)}, \dots, Y_{t-1})$ here, recall \eqref{filtration}. In this setting, the holdout set $\cD_0$ is not independent from $(X_t,Y_t)_{t\geq 1}$ and we assume that the whole observed process $(X_t,Y_t)_{t\geq -(n-1)}$ is as in \eqref{regressionAR:model}. For adaptation, we consider an $\mathcal{F}_{t-1}$-measurable estimator $\wh{\varphi}_t\in\R^d$  of $\varphi$ and the corresponding estimator $\hat{\mu}_t(x)=\scal[0]{\wh{\varphi}_t}{x}$  of $\mu(x)=\scal[0]{\varphi}{x}$, where $\scal[0]{\cdot}{\cdot}$ denotes the inner product in $\R^d$. With this estimator and $\hat{\sigma}_t(x)\equiv 1$, we consider the adaptive score $V_t(x,y)$  \eqref{regressionVadaptive}, the prediction interval $\mathcal{C}_t(x,q)$  \eqref{equCtregression}, which leads to a benchmark function $\Pi^\bench ( q)=1-q$ and benchmark threshold $q^\bench_{\alpha}=1-\alpha$ ($F=\Phi$ here). 
In addition, there exists $\wh{\varphi}_t$ such that \eqref{ass:esti} and \eqref{ass:esti:average} are satisfied in the $\mbox{AR}(d)$ model with $\nu_0=3/2$, $\nu=1/2$, $\tau_0=0$ and $\tau=1/2$. For instance, this is the case by using the regularized least-squares method  \citep{goldenshluger2001nonasymptotic}.
We thus obtain a result similar to Corollary~\ref{cor:rateconveasyselect} in the (non-iid) autoregressive case.

\begin{corollary}\label{cor:rateconveAR}
    Consider the online stationary $\mbox{AR}(d)$ model \eqref{regressionAR:model} with Gaussian innovations, parameterized by $\varphi\in\R^d$. 
Consider $\mathcal{R}^{\OSCI}$ using   
    a time-constant $X$-oriented selection (Algorithm~\ref{alg:PredictXoriented}) with a probability $\P(S(X_1)=1)$ fixed in $(0,1]$, estimators $\hat{\mu}_t(x)=\scal[0]{\wh{\varphi}_t}{x}$, $\hat{\sigma}_t\equiv 1$, and a threshold sequence $(q_t^{\OSCI})_{t\geq 1}$ \eqref{equACIquantile}  with a step size sequence $\gamma_j= c j^{-\beta}$, $\beta\in (1/2,1)$, and  starting point $q_1^{\OSCI}\in [0,B]$. Then there exist $T_0\geq 1$ and $C_0>0$ such that for $t\geq T_0$, 
    \begin{itemize}
        \item
    In-probability rate: assuming that the estimator $\hat{\mu}_t(x)=\scal[0]{\wh{\varphi}_t}{x}$ satisfies
        \eqref{ass:esti} with $\nu_0=3/2$, $\nu=1/2$,  we have
\begin{equation}\label{equ:ARConstSelectProba}
\probp{|\IER_t(\mathcal{R}^{\OSCI}) - \alpha|+ |\qt_t^{\OSCI}-(1-\alpha)| \leq C_0 (\rho_t(\beta)+(\log t)^{3/2} t^{-1/2})}\geq 1-1/t,
\end{equation}
for $\rho_t(\beta)=\rho_t(\beta;1/2)$ given by \eqref{equ:raterho-typicalconc}.
\item  
In-expectation rate: if in addition the constant selection is  $S_t(x,q)\equiv 1$ (no selective inference) and assuming that the estimator $\hat{\mu}_t(x)=\scal[0]{\wh{\varphi}_t}{x}$ satisfies \eqref{ass:esti:average} with $\tau_0=0$, $\tau=1/2$, we have 
\begin{equation}\label{equ:ARNoSelectExpect}
\paren{\e{\abs{\IER_t(\mathcal{R}^{\OSCI}) - \alpha}}}^2 + \e{\paren{q_{t}^{\OSCI}-(1-\alpha))}^2}
\leq C_0 t^{-\beta+1/2}.   
\end{equation}
    \end{itemize}
\end{corollary}

Corollary~\ref{cor:rateconveAR} is proved in \S~\ref{proof:cor:rateconveAR}. 
This shows that our theory can be applied in this non-iid setting as well. It is to our knowledge the first result that shows a convergence for the threshold/IER of an ACI-type method in the autoregressive model. This model was also considered in \cite{zaffran2022adaptive} (Example A.3) with $d=1$, but the problem and criterion therein was markedly different (convergence of the average length of a partly oracle, non-adaptive, ACI prediction set). 

In the cases considered, the convergence rates obtained in Corollary~\ref{cor:rateconveAR} are better than those of Corollary~\ref{cor:rateconveasyselect}, because the estimation task is easier in this parametric (dependent) model. 

Note that Corollary~\ref{cor:rateconveAR} encompasses the selection case $S(X_t)=\ind{Y_{t-1}\geq y_0}$ for some benchmark $y_0$ (which was not covered by our previous iid results). We apply this for the time series of electricity prices where a selection corresponds to the occurrence of an extreme event in \S~\ref{sec:elec}.

\section{\et{Application to selective classification}}\label{sec:selectiveclassifiid}

We consider in this section the classification task with $K$ classes, for which $\cY=\range{K}$.

\subsection{Setting and notation}

Consider the  iid classification model given by
\begin{equation}
    \label{model:classif}
    (X_t,Y_t) \in [0,1]^d \times [K], t\geq 1, \mbox{ iid  }\sim (X,Y),
\end{equation}
with an arbitrary distribution. We denote the true class probabilities by
\begin{equation}
  \label{equ:pit}
  \pi(y|x):= \P(Y_{t}= y\:|\: X_t=x), \:y\in [K],\:x\in \mathcal{X},
\end{equation}
and these probabilities are estimated with $\mathcal{F}_{t-1}$-measurable estimators 
$\hat{\pi}_{t}(y|x)\in [0,1]$. We assume that they satisfy, for a $t$ large enough and $C_1>0$, 
\begin{align}
\probp{\e{\max_{y\in \range{K}}|\hat{\pi}_{t}(y|X_t)-{\pi}(y|X_t)|\mid \mathcal{F}_{t-1}} \geq C_1 (\log t)^{\nu_0} t^{-\nu}}\leq \frac{3}{\pi^2t^3}\label{ass:esticlassif},
\end{align}
for $\nu_0>0$ and $\nu\in (0,1/2]$. This property is classical if the functions $x\mapsto {\pi}(y|x)$, $y\in [K]$, belong to some regular classes \citep{tsybakov2007fast,bos2022convergence}.

\et{A prediction set consists only  of the maximum estimated class probability, that is,  $\{\hat{Y}_t(x)\}=\arg\max_{k \in \range{K}}\hat{\pi}_{t}(k|x)$ (considering any tie-breaking rule and assuming $B=1$), the prediction set rule\footnote{It is of the form \eqref{eq:set_pred_funct} by taking $V_t(x,y)=\ind{y\neq \hat{Y}_t(x)}$.} is $C_t(x,q)=\{\hat{Y}_t(x)\}$ for $q<1$ and $C_t(x,q)=[K]$ otherwise. The selection rule $\St_t(x, q)$ is the informative selection \eqref{eq:selectionstat} with the statistics
\begin{align}
 W_t(x)&=\max_{y\in \range{K}} \hat{\pi}_{t}(y|x).\label{equ:Wtselecclassif}
\end{align}
This corresponds to selecting the examples for which the estimated probability of being in a class is large enough. 
Algorithm~\ref{alg:BinaryClassification} implements the selective classification procedure.}

\subsection{Convergence rate}

\et{
A convergence result is obtained by applying the pipeline presented in \S~\ref{sec:incompresscritic} to this classification context: we consider the limit of $W_t(x)$ given by 
\begin{equation}
    \label{equWClassif}
    W(x):=\max_{y\in \range{K}} {\pi}(y|x)\in [1/K,1].
\end{equation}
Assumption~\ref{ass:support} then entails that the maximum class probability is supported in $[\kappa,B^*]$. 
In addition, the IER function \eqref{PiCtSt} is given by 
$\Pi_{\cC_t, S_t}(q)=\P\big[Y_t\neq \hat{Y}_t(X_t) \mid W_t(X_t)>q, \mathcal{F}_{t-1}\big] = 1- \E\big[\pi(\hat{Y}_t(X_t)|X_t)\mid W_t(X_t)>q, \mathcal{F}_{t-1}\big]$, 
which is intended to get close to 
$\Pi^\bench(q)=1-\e{W(X)\mid W(X)>q}$
when $t$ grows, 
which is a special form studied in \S~\ref{sec:benchclasstest}. In particular, it is proved therein that the critical value \eqref{equalphastargen} is $\alpha^*=1-B^*$ and if $\alpha\in (\alpha^*,1-\e{W(X)})$, the benchmark function $\Pi^\bench$ satisfies Assumption~\ref{ass:easy}  with $q_\alpha^\bench\in (\kappa,1-\alpha)$.
We can thus apply Theorem~\ref{rem:coroonanEvent} to obtain the following result.
}

\begin{theorem}\label{th:rate:classif}
Consider the iid classification model \eqref{model:classif} satisfying Assumption~\ref{ass:support} with $W(x)$ given by \eqref{equWClassif} and support values  $1/K\leq \kappa\leq B^*\leq 1$, and respective estimators $ \hat{\pi}_{t}(y|x), x\in [0,1]^d,y\in [K],t\geq 1,$ of the class probabilities $\pi(y|x), x\in [0,1]^d,y\in [K]$ in \eqref{equ:pit} satisfying \eqref{ass:esticlassif}. 
Let $\alpha\in (0,1)$ be an error level  such that $1-B^*<\alpha<1-\e{W(X)}$. Consider the \method\  selective classification method (Algorithm~\ref{alg:BinaryClassification}) 
 for the threshold sequence $(q_t^{\OSCI})_{t\geq 1}$  with step size sequence $\gamma_j= c j^{-\beta}$, $\beta\in (1/2,1)$ and starting point $q_1^{\OSCI}\in [0,1-\alpha]$. 
Consider $\delta$ given in \eqref{deltaprobaJtfinite}. Then for any $\delta'>\delta$, there exist $T_0=T_0(\delta')\geq 1$ and $C_0=C_0(\delta')>0$ such that for all $t\geq T_0$,
\begin{equation}\label{equ-convclassif}
\probp{|\IER_t(\mathcal{R}^{\OSCI}) - \alpha|+ \abs{\qt_t^{\OSCI}-q^{\bench}_{\alpha}} \leq C_0(\rho_t(\beta)+(\log t)^{\nu_0/2} t^{-\nu/2})}\geq 1-\delta'- C_0/t,
\end{equation}
for $\rho_t(\beta) = \rho_t(\beta;3/4)=(\log t)^{1/2}t^{1/2-\beta} + t^{\beta-1} $ as in \eqref{equ:raterho-typicalconc}, $q^{\bench}_{\alpha}$ the only point of $(\kappa,1-\alpha)$ satisfying \eqref{equqalphainfo} and for $\IER_t(\mathcal{R}^{\OSCI})=\P(Y_t\neq \hat{Y}_t(X_t)\mid S^{\OSCI}_t(X_t)=1, \mathcal{F}_{t-1})$. 
\end{theorem}
Theorem~\ref{th:rate:classif} is proved in \S~\ref{proof:th:rate:classif}.  Note that the condition $\alpha<1-\e{W(X)}$ excludes 
targets above the Bayes classification error $1-\e{W(X)}=\P(Y\neq {Y}^*(X))$ (see \eqref{equ:classiforaclepredict}). It is mainly technical and only excludes favorable situations for which selecting all examples already achieves error below $\alpha$.  
By contrast, the condition $\alpha>1-B^*$ is a limitation when $B^*<1$, that is, when the variable ${\pi}(\hat{Y}_t(X)|X)=\max_{y\in \range{K}} {\pi}(y|X)$ cannot be arbitrarily close to $1$.

In Theorem~\ref{th:rate:classif}, optimizing $\rho_t(\beta)$ with respect to $\beta$ suggests to choose $\beta=3/4$. We can also opt for $\beta=(1-\nu/2)\wedge (1/2+\nu/2)\leq 3/4$ when $\nu$ is known, which diminishes $\beta$ (thus improving the FCP bound \eqref{FCPbound_decre}) and provides the same global rate $\rho_t(\beta)+(\log t)^{\nu_0/2} t^{-\nu/2} \propto t^{-\nu/2}$ (up to log terms).

\subsection{Optimality}

Theorem~\ref{th:rate:classif} states that the \method\  procedure gets close to the procedure selecting at time $t$ if $S^*(X_t)=1$ and predicting the class $Y^*(X_t)$, where
\begin{align}
   Y^*(x)&=\arg\max_{y\in \range{K}}\{\pi(y|x)\};\label{equ:classiforaclepredict}\\
    S^*(x)&=\ind{W(x)>q_\alpha^\bench}.\label{equ:classiforacleselect}
\end{align}
Reformulating Theorem~1 of \cite{zhao2023controlling} in our context, the latter procedure is optimal in the sense that it maximizes the probability of correct selection among all selection rules and classifiers with IER control at level $\alpha$ (see Proposition~\ref{prop:optimalclassif} for a slightly more general statement).  
Combined with Theorem~\ref{th:rate:classif}, this implies the following result.

\et{
\begin{corollary}
    \label{optimalityclassifiid}
 In the setting of Theorem~\ref{th:rate:classif}, consider the corresponding online selective classification method $\mathcal{R}^{\OSCI}=(S^{\OSCI}_t(X_t),\hat{Y}_t(X_t))_{t\geq 1}$. Then for any $\delta'\in(\delta,1)$, there exists $C_0=C_0(\delta')$ and $T_0=T_0(\delta')$, such that for $\rho_t(\beta)=\rho_t(\beta;3/4) =(\log t)^{1/2}( t^{1/2-\beta}+ t^{\beta-1} )$ as in \eqref{equ:raterho-typicalconc}, the following holds:
 \begin{itemize}
     \item[(i)] Marginal selection-conditional error rate:  for $t\geq T_0$, 
     $$\P(Y_t\neq \hat{Y}_t(X_t)\mid S^{\OSCI}_t(X_t)=1)\leq (1-\delta')^{-1}(\alpha +C \delta')+C_0 (\rho_t(\beta) +(\log t)^{\nu_0/2} t^{-\nu/2});$$
     \item[(ii)] Mimicking the oracle: for $t\geq T_0$,
     \begin{align*}
     \e{\ind{Y_t= {Y}^*(X_t)} S^*(X_t)}
&\leq (1-\delta')^{-1}\E\big[ \ind{Y_t= \hat{Y}_t(X_t)}S^{\OSCI}_t(X_t)\big] \\
&\:\:\:\:\:+  C_0(\rho_t(\beta) +(\log t)^{\nu_0/2} t^{-\nu/2});
\end{align*}
     \item[(iii)] Approximately dominating all other valid procedures:
      any selective classification procedure $(S_t(X_t),\tilde{Y}_t(X_t))_{t\geq 1}$ (for which each function $S_t(\cdot)$ and $\tilde{Y}_t(\cdot)$ is $\mathcal{F}_{t-1}$-measurable) is dominated by \method\  in the following sense:  if for all $t\geq T_0$, $\P(Y_t\neq \tilde{Y}_t(X_t)\mid S_t(X_t)=1)\leq \alpha+C\delta'+\upsilon_t$, for some constant $C>0$ and rate $\upsilon_t\to 0$, then for all $t\geq T_0$
    \begin{align}\label{eq-OSCIpower-class}
     \E\big[ \ind{Y_t=\tilde{Y}_t(X_t)}S_t(X_t)\big] \leq\:& (1-\delta')^{-1} \E\big[ \ind{Y_t= \hat{Y}_t(X_t)}S^{\OSCI}_t(X_t)\big]\\
     &+ C_0(\rho_t(\beta) +(\log t)^{\nu_0/2} t^{-\nu/2})  + C \delta'+ C\upsilon_t,\nonumber
    \end{align}
 \end{itemize}
  where in these inequalities, $C>0$ is a constant not depending on $\delta'$.
\end{corollary}
}

Corollary~\ref{optimalityclassifiid} is proved in \S~\ref{proof:optimalityclassifiid}. It says that in the selective classification setting, \method\  is optimal up to remainder terms that are small if $\delta$ \eqref{deltaprobaJtfinite} is small and $t$ is large.

\section{\et{Application to conformal testing}}\label{sec:selectiveNDiid}

\et{We consider in this section the online  testing paradigm where $\mathcal{Y}=\{0,1\}$, for which we wish to test ``$Y_t=0$'' versus ``$Y_t=1$'' at each time point $t$.
The label $Y_t$ is revealed immediately after  the decision  has been made.}

\subsection{Setting and notation}

Consider the iid model
\begin{equation}
    \label{model:ND}
    (X_t,Y_t) \in [0,1]^d \times \{0,1\}, t\geq 1, \mbox{ iid  }\sim (X,Y),
\end{equation}
with an arbitrary distribution.
To discriminate between the null ``$Y_t=0$'' and the alternative ``$Y_t=1$'', a key quantity is the local false discovery rate \citep{ETST2001} 
\begin{align}
{\lfdr}(x) &:= \mP(Y=0\mid X=x)\label{equ:lfdr}.
\end{align}
Suppose an estimator sequence $(\widehat{\lfdr}_t(x))_{t\geq 1}$ is available, that satisfies the following property: for $C_1>0$ and $t$ large enough,
\begin{align}
\probp{\e{|\widehat{\lfdr}_t(X_t)-\lfdr(X_t)|\mid \mathcal{F}_{t-1}} \geq C_1 (\log t)^{\nu_0} t^{-\nu}}\leq \frac{3}{\pi^2t^3}.\label{ass:estiND}
\end{align}
This assumption is met under regularity conditions for $x\mapsto \lfdr(x)$ \citep{tsybakov2007fast,bos2022convergence}, possibly after a dimension reduction stage \cite{ni2025sequential}.

\et{In conformal testing, the decision is whether to reject the null, which means that the prediction set rule\footnote{It is of the form \eqref{eq:set_pred_funct} by taking $V_t(x,y)=\ind{y=0}$.} is $\cC_t(x,q)=\set{1}$ for $q<1$ and $\cC_t(x,q)=\set{0,1}$ for $q\geq 1$ ($B=1$ here). The selection rule $\St_t(x, q)$ is the informative selection \eqref{eq:selectionstat} with the statistics
\begin{align}
    W_t(x)= 1-\widehat{\lfdr}_t(x).
    \label{equ:WtselecND}
\end{align}
This corresponds to selecting the examples for which the estimated null probability is small enough. Algorithm~\ref{alg:onlineBH} implements the conformal testing procedure.
}

\subsection{Convergence rate}

\et{
A convergence result is obtained by applying the pipeline presented in \S~\ref{sec:incompresscritic} to this testing framework: we consider the limit of $W_t(x)$ given by 
\begin{equation}
    \label{equW}
    W(x):=1-{\lfdr}(x)\in [0,1].
\end{equation}
Assumption~\ref{ass:support} then entails that the alternative class probability is supported in $[\kappa,B^*]$. 
In addition, the IER function \eqref{PiCtSt} is given by $\Pi_{\cC_t, S_t}(q)=\P(Y_t=0 \mid 1-\widehat{\lfdr}_t(X_t) >q, \mathcal{F}_{t-1})$. This corresponds to an instantaneous false discovery rate, conditional both on rejection and the history. For a large time horizon, it is intended to approach the benchmark function
\begin{equation}\label{equ:qvaluefunction}
\Pi^\bench(q)= \probp{Y=0 \mid 1-{\lfdr}(X) >q}=\e{{\lfdr}(X)\mid {\lfdr}(X)<1-q},
\end{equation}
which is  termed the positive FDR \citep{Storey2003, Efron2008} in the multiple testing literature, for the (optimal) case that the rejection region is  $\lfdr(X)<1-q$.
The function $\Pi^\bench$ is studied in \S~\ref{sec:benchclasstest}. In particular, it is proved therein that the critical value \eqref{equalphastargen} is $\alpha^*=1-B^*$ and if $\alpha\in (\alpha^*,1-\e{W(X)})$, the benchmark function $\Pi^\bench$ satisfies Assumption~\ref{ass:easy}  with $q_\alpha^\bench\in (\kappa,1-\alpha)$.
We can thus apply Theorem~\ref{rem:coroonanEvent} to obtain the following result.
}

\begin{theorem}\label{th:rate:ND}
Consider the iid testing model \eqref{model:ND} satisfying Assumption~\ref{ass:support} with the test statistics $W(x)=1-{\lfdr}(x)$ \eqref{equW} and support values $0\leq \kappa<B^*\leq 1$, and respective estimators $( \widehat{\lfdr}_t(x))_{x\in [0,1]^d,t\geq 1}$ of $(\lfdr(x))_{x\in [0,1]^d}$ satisfying \eqref{ass:estiND}. 
Let an error level $\alpha\in (0,1)$ be such that $1-B^*<\alpha<\P(Y=0)$. Consider the online conformal testing  method (Algorithm~\ref{alg:onlineBH}) 
for the threshold sequence $(q_t^{\OSCI})_{t\geq 1}$  with step size sequence $\gamma_j= c j^{-\beta}$, $\beta\in (1/2,1)$ and starting point $q_1^{\OSCI}\in [0,1-\alpha]$. 
Consider $\delta$ given in \eqref{deltaprobaJtfinite}. Then for any $\delta'>\delta$, there exist $T_0=T_0(\delta')\geq 1$ and $C_0=C_0(\delta')>0$ such that for all $t\geq T_0$,
\begin{equation}\label{equ-convND}
\P\paren{|\IER_t(\mathcal{R}^{\OSCI}) - \alpha|+ \abs{\qt_t^{\OSCI}-q^{\bench}_{\alpha}} \leq C_0(\rho_t(\beta)+(\log t)^{\nu_0/2} t^{-\nu/2})}\geq 1-\delta'- C_0/t,
\end{equation}
for $\rho_t(\beta)=\rho_t(\beta;3/4)=(\log t)^{1/2}t^{1/2-\beta} + t^{\beta-1} $ as in \eqref{equ:raterho-typicalconc}, $q^{\bench}_{\alpha}$ the only point of $(\kappa,1-\alpha)$ satisfying \eqref{equqalphainfo} and where $\IER_t(\mathcal{R}^{\OSCI})=\P(Y_t=0\mid S^{\OSCI}_t(X_t)=1, \mathcal{F}_{t-1})$ is the instantaneous false discovery rate. 
\end{theorem}
Theorem~\ref{th:rate:ND} is proved in \S~\ref{proof:th:rate:ND}. 
Note that Assumption~\ref{ass:support} with the test statistics $W(x):=1-{\lfdr}(x)$ means that the distribution of $\lfdr(X)$ has a continuous density on $[\alpha^*,1-\kappa]$, which is positive on $(\alpha^*,1-\kappa)$. In particular, there is no criticality ($\alpha^*=0$) if $\lfdr(X)$ has a distribution that can be arbitrary close to zero, which means that null observations and alternative observations can be arbitrarily well separated with a positive probability. If $\alpha^*>0$, however, they cannot be arbitrarily well separated and Theorem~\ref{th:rate:ND} imposes $\alpha>\alpha^*$ to have the desired high probability convergence rate.

\subsection{Optimality}

Theorem~\ref{th:rate:ND} ensures that the \method\ method is close to the oracle testing procedure 
\begin{align}
    S^*(x)&=\ind{\lfdr(x)<1-q_\alpha^\bench}.\label{equ:NDoracleselect}
\end{align}
This procedure corresponds to the oracle procedure of \cite{SC2007}.
Reformulating Theorem~1 of \cite{SC2007} in our context, the latter procedure is optimal in the sense that it maximizes the probability of correct detection among all detection rules with IER control at level $\alpha$ (see Proposition~\ref{prop:optimalND} for a slightly more general and detailed statement). 

\et{
\begin{corollary}
    \label{optimalityNDiid}
 In the setting of Theorem~\ref{th:rate:ND}, consider the corresponding conformal testing method $\mathcal{R}^{\OSCI}=((S^{\OSCI}_t(X_t))_{t\geq 1},(\{1\})_{t\geq 1})$. Then for any $\delta'\in(\delta,1)$, there exists $C_0=C_0(\delta')$ and $T_0=T_0(\delta')$, such that we have for $\rho_t(\beta)=\rho_t(\beta;3/4)=(\log t)^{1/2} t^{1/2-\beta} + t^{\beta-1}$ as in \eqref{equ:raterho-typicalconc},
 \begin{itemize}
     \item[(i)] One-step marginal false discovery rate control: for all $t\geq T_0$,
     $$\P(Y_t=0\mid S^{\OSCI}_t(X_t)=1)\leq (1-\delta')^{-1}(\alpha +C\delta')+C_0\rho_t(\beta)+C_0(\log t)^{\nu_0/2} t^{-\nu/2};$$
     \item[(ii)] Mimicking the oracle: for $t\geq T_0$,
     $$\e{Y_t S^*(X_t)}\leq (1-\delta')^{-1} \E\big[ Y_t S^{\OSCI}_t(X_t)\big] +  C_0\rho_t(\beta)+C_0(\log t)^{\nu_0/2} t^{-\nu/2} ; $$
     \item[(iii)] Approximately dominating all other valid procedures:
     for any online testing procedure $((S_t(X_t))_{t\geq 1},(\{1\})_{t\geq 1})$ (for which $S_t(\cdot)$ is $\mathcal{F}_{t-1}$-measurable) is dominated by \method\  in the following sense:  if for all $t\geq T_0$, $\P(Y_t=0\mid S_t(X_t)=1)\leq \alpha+C\delta'+\upsilon_t$, for some rate $\upsilon_t\to 0$, then for all $t\geq T_0$
    \begin{equation}\label{eq-OSCIpower-ND}
     \E\big[ Y_t S_t(X_t)\big] \leq (1-\delta')^{-1} \E\big[ Y_t S^{\OSCI}_t(X_t)\big]+ C_0\rho_t(\beta)+C_0(\log t)^{\nu_0/2} t^{-\nu/2}  + C \delta'+ C\upsilon_t,
    \end{equation}
 \end{itemize}
   where in the above inequalities $C>0$ is a constant not depending on $\delta'$.
\end{corollary}
}

Corollary~\ref{optimalityNDiid} is proved in \S~\ref{proof:optimalityNDiid}. 
According to Corollary~\ref{optimalityNDiid}, the conformal testing procedure of Algorithm~\ref{alg:onlineBH} can be seen as a counterpart of the Sun and Cai procedure \citep{SC2007} in an online testing context, that mimics asymptotically the oracle (with a rate here)\et{, up to a term involving the failure probability $\delta$.} In addition, \et{our procedure provides an FDP guarantee in an adversarial setting where the number of rejections is large. This also contrasts with the procedure SAST of \cite{gang2023structure}, for which no such adversarial FDP bound
is provided, see Table~\ref{tab:novelties} and the discussion in \S~\ref{sec:comparison}.}

Finally, let us mention that our theory also covers the case of online selection by prediction (considered in \cite{Jin2023selection} in the batch setting), see \S~\ref{sec:selectpredict} for more details.

\section{Numerical experiments}\label{sec:num}

In this section, we evaluate our methods\footnote{Our code  is available at \url{https://github.com/pierreHmbt/onlineSCI}.} on synthetic data for the different applications presented above. 
Additional real data examples are given in \S~\ref{sec:realdata} (\texttt{CIFAR-10} data set) and \S~\ref{sec:elec} (electricity prices), for which the results are analogous. \et{Guidelines for hyperparameter tuning are also provided in \S~\ref{sec:select_process}. 
}

We show that the ACI method of \cite{angelopoulos2024online} fails to control the FCP over the selected examples and that the inflation in FCP can be large when selection is unaccounted for; we demonstrate the  actual convergence of \method\  in a manner which is in accordance with the proved theoretical results;  we show the large potential advantage of using adaptive scores over non-adaptive scores with regard to overall and instantaneous power.   
\subsection{Online conformal testing}\label{subsec-simND}
For the online testing described in  \S~\ref{sec:selectiveNDiid}, we consider the following iid data generation.  Examples are sampled from the mixture model
\begin{align*}
    Y \sim \text{Bernoulli}(0.2) \: \text{and} \: (X \mid Y=y) \sim \mathcal{N}((0, 0), I_2) \cdot \ind{y=0} + \mathcal{N}((3, 3), I_2) \cdot \ind{y=1} \; ,
\end{align*}
where $I_2$ is the identity matrix of size $2 \times 2$, $Y=0$ corresponds to a null and $Y=1$ to an alternative. In this specific application, our goal is to identify the observations $X_t$ associated to $Y_t=1$. 

We apply Algorithm \ref{alg:onlineBH} with $\alpha=0.1$, $q_1=0.5$ and $\gamma_j = j^{-3/4} $ (as suggested by Remark~\ref{rem:choosebeta}), and compute an estimator $\widehat{\lfdr}_t$ of $\lfdr$ \eqref{equ:lfdr} 
by using a Random Forest (RF) classifier with the default parameters in \texttt{scikit-learn} \citep{scikit-learn}. At each time, the estimation is performed in two different ways:
\begin{enumerate}
    \item {Non-adaptive:} the RF estimator is trained with an independent data set $\cD_0$ of size $100$;
    \item {Adaptive:} the RF estimator is updated at each time $t$ with the new available pair  $(X_{t-1}, Y_{t-1})$.
\end{enumerate}

\begin{figure}[h!]
		\centering
		\includegraphics[width=.45\linewidth]{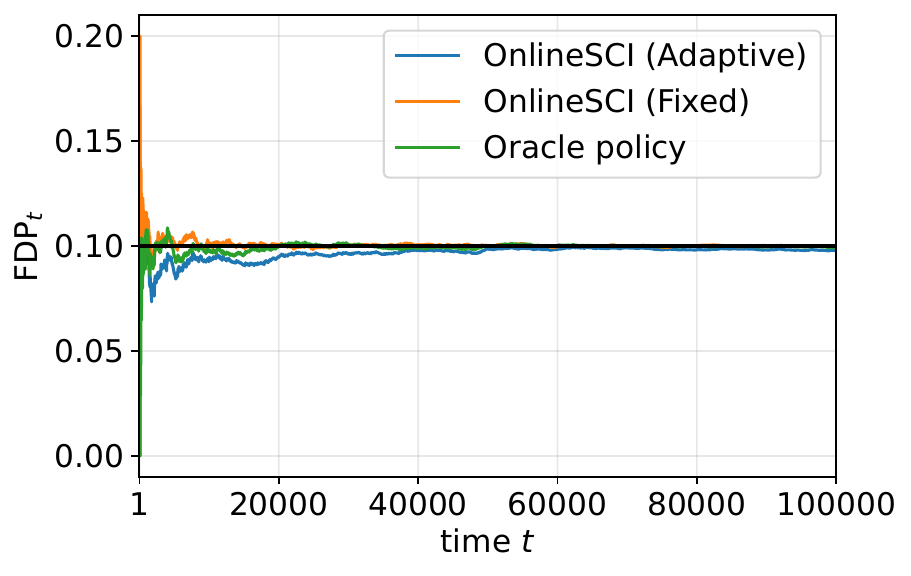}
        \includegraphics[width=.45\linewidth]{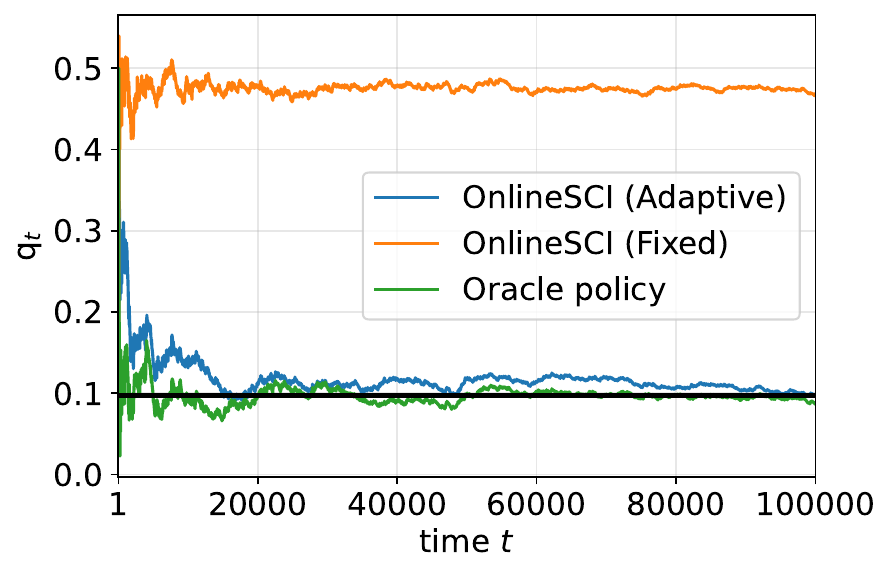} \\
        \includegraphics[width=.45\linewidth]{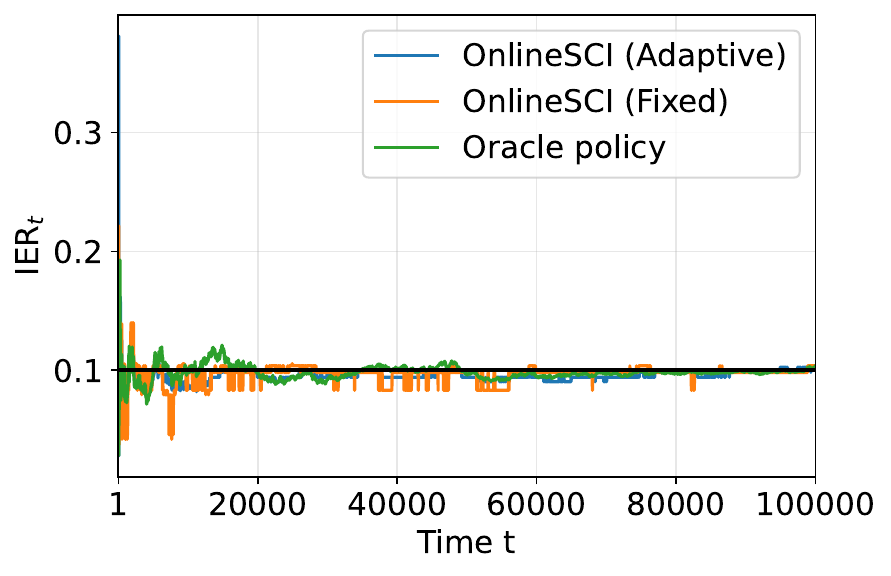}
        \includegraphics[width=.45\linewidth]{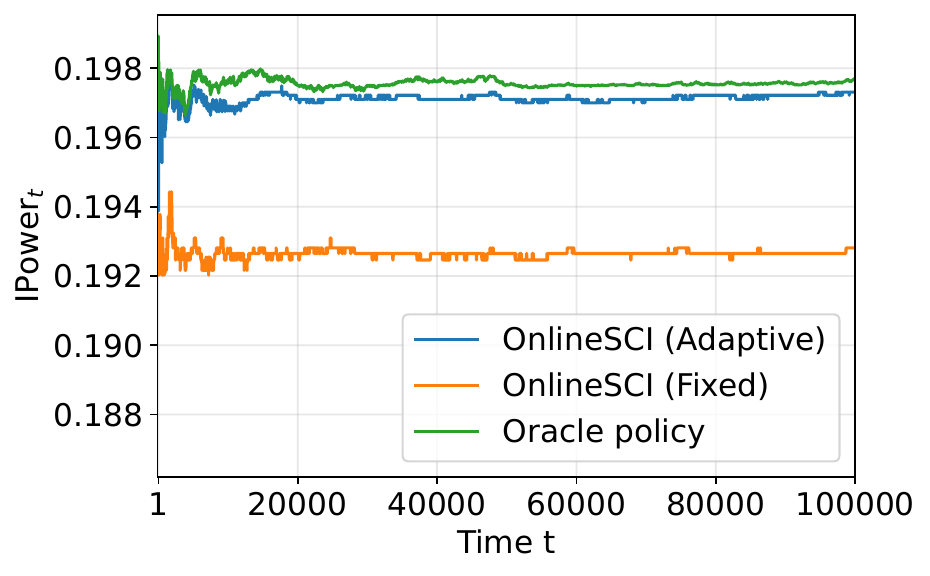}
        \caption{Online conformal testing: results for \method~and the oracle policy on a particular data set, see details in \S~\ref{subsec-simND}. Top left panel: $\FDP_t$ versus time step $t$. Top right panel: method threshold $q_t$ versus time step $t$. Bottom row: $\IER_t$ and instantaneous power at each time step. 
        }
        \label{fig:ND_xp}
\end{figure}

Finally, the following benchmark procedures are also considered:
\begin{itemize}
    \item Oracle policy: Algorithm \ref{alg:onlineBH} with $W_t(x)=1-\lfdr(x) = \probp{Y=1 \mid X=x}$. It is considered as the best benchmark procedure, because it achieves both FCP control and threshold convergence to the optimal solution of \cite{SC2007}. Recall that this is not the case for the procedure with constant threshold $q^\bench_\alpha$, which has no adversarial finite sample FCP bound guarantee.  
	\item 
    Levels based On Recent Discovery++ (LORD++) \citep{Javanmard18, ramdas2017online}, the adaptive version SAFFRON \citep{ramdas2018saffron} and the label-aware version FeedBack LORD \citep{lu2025feedback}: At each time step $t$, we construct a conformal $p$-value \cite{vovk2009line} with the data-driven score $V_t(x, 0) = 1 - \widehat{\lfdr}_t(x)$ and where only the $(X_s)_{1 \leq s \leq t}$ associated to $Y_s = 0$ are considered in the calibration sample. 
    The $p$-values are not independent hence these procedures are not guaranteed to control the online FDR.
\end{itemize}

Finally, for each procedure, we evaluate the $\FCP$ (which is the FDP in that case), $\IER$, and the instantaneous power, that is,  
\begin{equation}
    \label{Ipower}
\mbox{IPower}_t := \probp{S_t(X_t)=1, Y_t=1 \mid \mathcal{F}_{t-1}} = \e{Y_t S_t(X_t)\mid \mathcal{F}_{t-1}}\leq \P(Y_1=1)=0.2.
\end{equation}
as in  our optimality result, see  \S~\ref{sec:selectiveNDiid}.

Figure \ref{fig:ND_xp} displays the results obtained using (a) \method~with both adaptive and fixed statistics, and (b) the oracle policy. We see that the FDP approaches the nominal level $\alpha=0.1$ for all methods. In addition, the sequence $(q_t)_{t \geq 1}$   converges to $q^0_{\alpha}$ for \method~with adaptive statistics but not  for fixed statistics. This is  expected because the lfdr estimate does not approach the true one when the statistics are fixed. 
We also clearly see that \method~with adaptive statistics returns comparable results as the oracle policy with $\IER_t$ close to $\alpha$ (bottom left panel) and a similar power (bottom right panel). In addition, the convergence of the FDP is faster than the convergence of $q_t$, which is reasonable because the setting is iid and we expect that the FDP concentrates better in this favorable case. Figure~\ref{fig:ND_xp_multiple} shows that this is not specific to one data generation. 
Finally, LORD makes no discovery in this particular example. This is not surprising, as LORD is designed to operate under much more difficult conditions than those considered in this paper (in particular, LORD does not use the previous labels).

\subsection{Online selective classification with $2$ classes}\label{subsec-simClassification}

For online selective classification (see \S~\ref{sec:selectiveclassifiid}), we consider the following iid data generation: examples are sampled from the mixture model
\begin{align*}
    &Y = \text{Bernoulli}(0.5)
     \text{ and } (X \mid Y=y) \sim \mathcal{N}((0, 0), I_2) \cdot \ind{y=0} + \mathcal{N}((1, 1), I_2) \cdot \ind{y=1}\; .
\end{align*}
We apply Algorithm \ref{alg:BinaryClassification} with $\alpha=0.1$, $q_1=0.8$ and $\gamma_j = 0.5 \cdot j^{-3/4} $ and compute the estimator $\hat{\pi}_{t}(y|x)$ of $\pi_t(y|x)$ \eqref{equ:pit} using a Gaussian  Naive Bayes (GNB) classifier with the default parameters in \texttt{scikit-learn} \citep{scikit-learn}. As in the previous subsection, the learning  is performed in two different ways:
\begin{enumerate}
  \item {Non-adaptive:} the GNB estimator is trained with an independent data set $\cD_0$ of size $10$;
   \item {Adaptive:} the GNB estimator is updated at each time step $t$ with the new available pair  $(X_{t-1}, Y_{t-1})$.
  \end{enumerate}
The benchmark procedures are:
\begin{itemize}
    \item Oracle policy: Algorithm~\ref{alg:BinaryClassification} for which the true probabilities $\pi_t(y|x)$ \eqref{equ:pit} are used instead of their estimates $\hat{\pi}_{t}(y|x)$;
    \item ACI with decaying step sizes \citep{angelopoulos2024online} where the non-conformity score is $1-\widehat{\pi}_t(x,y)$.
\end{itemize}
\et{The instantaneous power is defined here by
\begin{equation}
    \label{Ipowerclassif}
\mbox{IPower}_t := \probp{S_t(X_t)=1,Y_t= \hat{Y}_t(X_t)\mid \mathcal{F}_{t-1}} ,
\end{equation}
as in our optimality result, see \S~\ref{sec:selectiveclassifiid}.
}

\begin{figure}[h!]
		\centering
		\includegraphics[width=.45\linewidth]{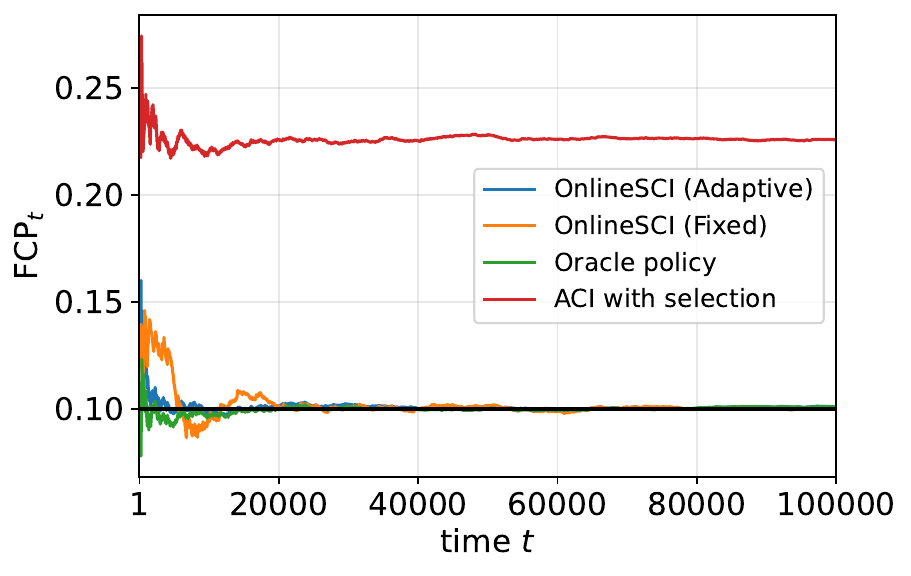}
        \includegraphics[width=.45\linewidth]{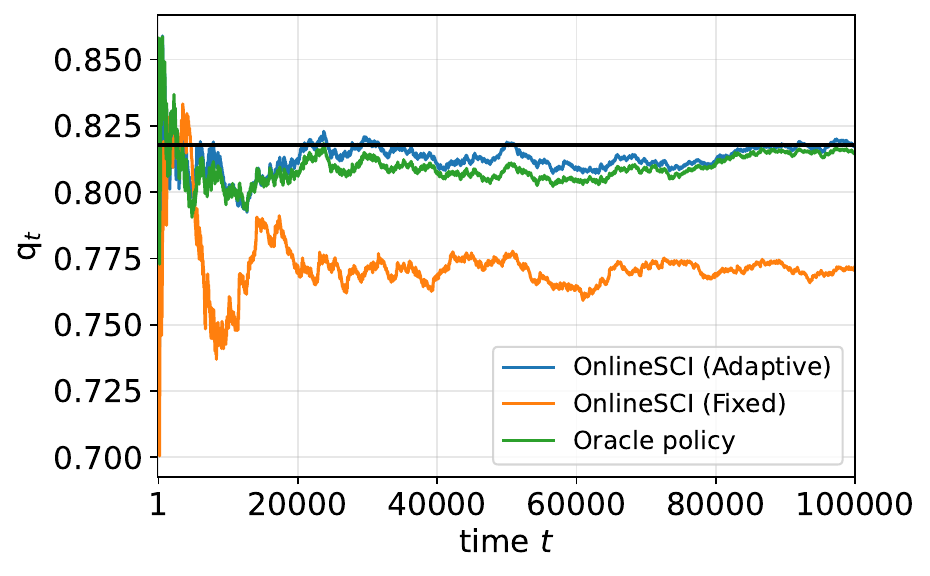} \\
        \includegraphics[width=.45\linewidth]{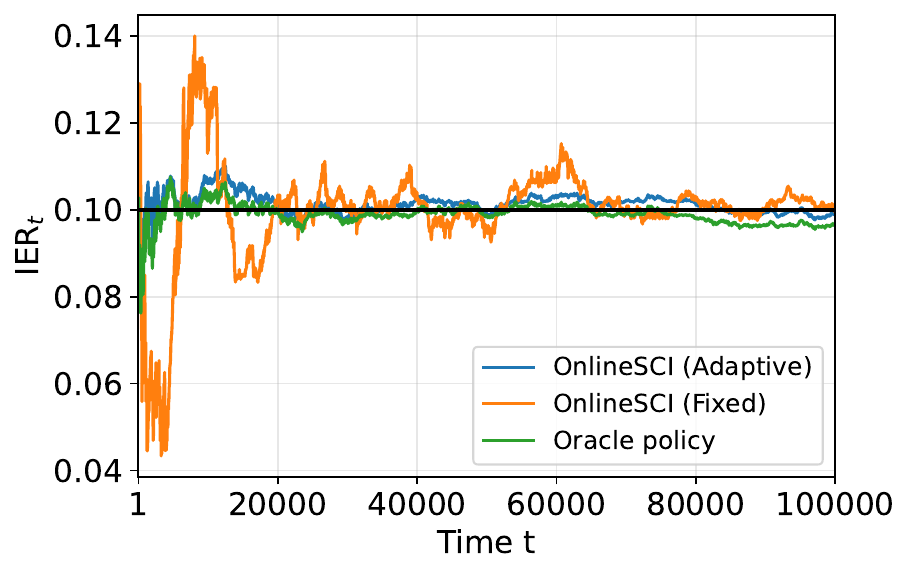}
        \includegraphics[width=.45\linewidth]{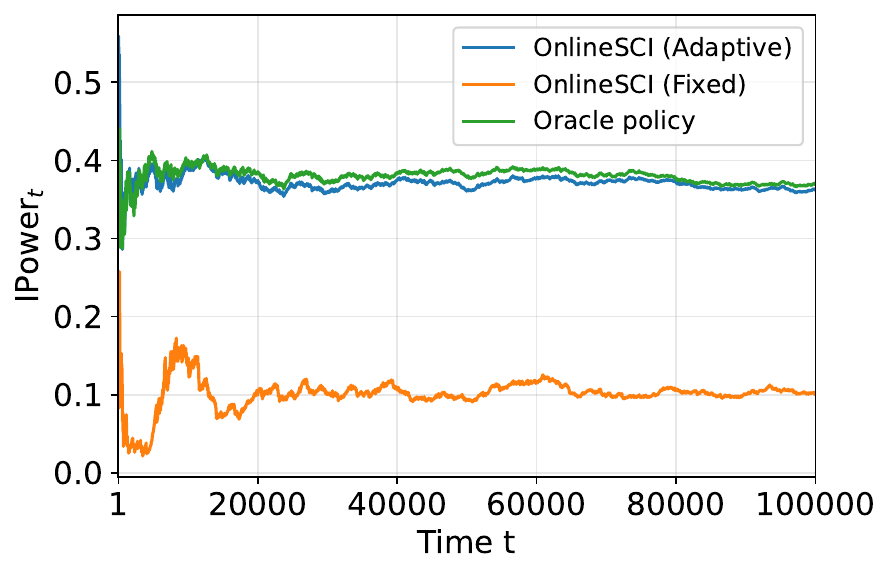}
        \caption{
        Online selective classification: results for \method~and the oracle policy on a particular data set, see details in  \S~\ref{subsec-simClassification}. Top left panel: $\FCP_t$ versus  $t$. Top right panel: method threshold $q_t$ versus  $t$. Bottom row: $\IER_t$ and instantaneous power versus $t$.
        }
        \label{fig:SC_xp}
\end{figure}

Figure \ref{fig:SC_xp} illustrates the results obtained with the different algorithms. As expected, the $\FCP$ is controlled at level $\alpha=0.1$ for all methods except ACI of \cite{angelopoulos2024online} \pie{with a naive selection; meaning that a selection is made if the size of the set is $1$} (top left panel), which has a much larger FCP \et{(the same conclusion is obtained for DtACI \citep{gibbs2024conformal}, see \S~\ref{sec:compaDTACI})}. Furthermore, the threshold sequences $(q_t)_{t \geq 1}$ are all converging, but with different limits; only the \method~with adaptive scores (blue curve) and the oracle policy (green curve) have thresholds converging to $q^0_\alpha$. Figure \ref{fig:SC_xp} (bottom right panel) also shows that the adaptive method strongly dominates the non-adaptive one in terms of instantaneous power.

\section{Discussion}\label{sec:discuss}

In this work, we introduced the method \method\  which supports various online selective conformal inferences. We show both theoretically and empirically that it makes the best of two worlds: while it comes with an FCP bound guarantee over all possible sequences of observations, it controls the history-conditional instantaneous error rate at a level close to $\alpha$ when the time gets large in the favorable case where the data are identically distributed (iid or autoregressive) and under mild assumptions. At this moment, we are not aware of any other method achieving both tasks together: on the one hand, the (generally unknown) method using the constant oracle threshold $q^\bench_\alpha$ does not provide the FCP bound guarantee in the adversarial setting. On the other hand, the selective DtACI method of \cite{bao2025cap} does not control the instantaneous error rate in the iid setting (\S~\ref{sec:compaDTACI}). Finally, in our leading cases (iid online regression, selective classification, conformal testing), we establish that the method \method\  satisfies an optimality property,  up to remainder terms that vanish at specific rates. 

However, for informative selections, we have shown that a theoretical failure probability $\delta\geq 0$ comes into play in our IER bounds, which records the unfavorable trajectories of \method. \et{This failure may be due to a bad choice of the initial parameters of the algorithm that makes it `frozen'. The error $\delta$ is possibly nonzero (Lemma~\ref{lem:deltapositive}), but it can be reduced by restarting the algorithm and adjusting the parameter choice if necessary, see Remark~\ref{rem:rerun}.}

This work opens several interesting avenues for future research.
For a non-adaptive score function, our method relies only on the observed outcomes at the selected time points. This is particularly desirable in applications where measuring the outcome is costly. For instance, in health screening, individuals might undergo further diagnostic testing only if selected by an online conformal testing procedure, where the null hypothesis is that the person is healthy. Online conformal testing remains applicable even when the health status of unselected individuals is unknown. An interesting direction for future work is to move beyond a fixed score function by developing adaptive scores that incorporate evidence from the covariates, which are observed at all time points.

Our convergence rates also provided guidelines for choosing the step size sequence $\gamma_t$ in \method\   (Remark~\ref{rem:choosebeta}). While these practical recommendations seem to be effective in practice, the picture is far from being complete, and it is interesting to explore other ways to set this step size sequence, for instance by using part of the past data. This is a very active area of research in the corresponding machine learning literature 
which promises interesting future developments and recommendations for \method.

\et{
Another direction for future work is to establish convergence results for \method\: under distribution shift. Our most general bounds (Section~E) essentially require no modeling assumptions and provide a starting point for accommodating such settings. However, the nature of the distribution shift would need to be specified so that the order of the remainder terms can be made explicit.
}

\section*{Acknowledgments}
We thank two anonymous referees, the Associate Editor, and the Editor for their constructive comments and valuable suggestions, which have helped considerably improve the manuscript.
We would like to warmly acknowledge Badr-Eddine Ch\'erief-Abdellatif for his involvement in the Emergence project which greatly inspired and supported the work on this manuscript. We thank Félix Laplante for his help with the Streamlit app. We also thank Antoine Godichon-Baggioni, Jean-Marc Bardet, Olivier Wintenberger and Aaditya Ramdas for very helpful discussions. 
The authors acknowledge grants ANR-21-CE23-0035 (ASCAI) and ANR-23-CE40-0018-01 (BACKUP) of the French National Research Agency ANR, the Emergence project MARS of Sorbonne Universit\'e, and Israeli Science Foundation grant no. 406/24.

\bibliographystyle{plain}
\bibliography{biblio}

\appendix

\section{Terminology and relation to existing work}\label{sec:comparison}

  Online conformal inference refers to ACI type procedures, i.e., gradient descent algorithms with the pinball loss \citep{gibbs2021adaptive}. Our method is also of the ACI type, but accounts for inference  only at selected time points; hence it is called {\it online selective conformal inference}. 
  
  \et{ACI-type methods are 'adaptive' in the sense of adapting with respect to a covariate shift. This is achieved by updating the empirical quantile over time using the gradient descent algorithm. In this work, we use 'adaptive' in a different sense: we say that the (ACI-type) method is {\it adaptive} if it updates over time the non-conformity score, and more generally the inference. So an adaptive \method\ will use a selection rule $S_t$ that depends on $X_t$ and on the past $\mathcal F_{t-1}$ (all data before time $t$). 
  There are other ACI-works using an adaptive score function \cite{bao2025cap,gibbs2024conformal} but without quantifying precisely the impact of this adaptation stage on the error rate. A contribution of our work is to identify and explicitly quantify the adaptation cost $D_t$  \eqref{equDtT} in our remainder terms. 
  }

\et{For online selective interval estimation, \citep{weinstein2020online}  proposed controlling the expected FCP when confidence intervals are reported sequentially  for selected parameters. Their methods can be translated to reporting prediction intervals over time for selected examples in settings where  the outcome $Y_t$ is never revealed. The setting in which selection takes place, but past outcomes   $(Y_s)_{s<t}$ are available at time $t$ (i.e., before  constructing the prediction interval for $Y_t$), has only recently been studied by  \cite{bao2025cap, sale2025online,zheng2026online}.  
}

\et{
Unlike in the works of \cite{bao2025cap, sale2025online,zheng2026online}, we consider the case that the selection can depend on the resulting prediction set:  after construction, the analyst may decide not to select the example if the prediction set is not sufficiently informative. In the special case where the examples are iid, these works provide finite-sample IER control (specifically, they focus  on controlling the probability of non-coverage conditional on being selected, $\P(Y_t\notin \cC_t \:|\:  S_t=1)$), but only for more restrictive classes of selection rules that do not allow selection to be based on the  informativeness of the resulting prediction sets.  We obtain  only asymptotic IER control, up to remainder terms, but for a much more general class of selection rules, including rules that depend on  the resulting prediction sets. 
}

\et{
In \cite{bao2025cap}, for the non-iid setting, another algorithm is suggested that provides FCP control  on the selected. 
First, they showed that by constructing prediction sets using ACI or a variant of ACI and reporting them only for selected examples, we may have $\FCP_t(\cR) \gg \alpha$ for all $t$ even in the iid setting (see also  \S~\ref{subsec-simClassification}). 
Then, they provided a simple fix to ACI for asymptotic $\FCP_t$ control in the adversarial setting, but it  does not bound the history-conditional instantaneous miscoverage probability among the selected. 
They also proposed a selective version of online DtACI \cite{gibbs2024conformal}, but it does not solve the above problem:  while OSCI asymptotically controls the IER, the selective DtACI does not, and has a threshold sequence not converging for long time horizons, see \S~\ref{sec:compaDTACI}. 
}

\et{
In \cite{angelopoulos2024online},  a procedure is suggested that has good properties both in the adversarial and iid settings, but only when  no selection takes place. Arguably, (implicit) selection always takes place when multiple inferences are considered over time: analysts do not care about all inferences produced, and would like the error guarantee to be on the inferences that are of interested to them. Moreover, as argued in  \cite{angelopoulos2024online}, both error control of $\FCP_t$ in the adversarial setting and control of the instantaneous miscoverage probability in the iid setting are of interest. 
We adopt  the approach in \cite{angelopoulos2024online} of suggesting a procedure that has good properties  both in the adversarial and iid setting, but unlike in \cite{angelopoulos2024online}, we provide the error control on the selected examples. Importantly, our approach is  adaptive: we incorporate  all past information in order to improve the score and the selection rule, so  we aim to control the history-conditional instantaneous error rate
$\IER=\P(Y_t\notin \cC_t \:|\:  S_t=1, \mathcal F_{t-1})$ (in online conformal testing,  the IER is $\P(Y_t=0 \:|\:  S_t=1, \mathcal F_{t-1})$). 
}

\et{ In the machine learning literature, gradient descent type algorithms with a selection stage have been  developed to improve the convergence rate \citep{namkoong2017adaptive,lu2022general}, but not for approaching probabilities/quantiles conditionally on a prescribed user-given rule (that is, selection is not included in the criterion at hand), to the best of our knowledge. One innovation of this manuscript is to derive explicit rates for IER convergence of \method\ in the context of selection. A closely related work is \cite{labopin2019conditional} that aims to build an online approximation of the quantile of the distribution of the score conditionally on $X=x_0$, by selecting at time $t$ only if $X_t$ is in a neighborhood of $x_0$ (with some specific vanishing size).
The problem therein is to obtain an asymptotic conditional coverage, so it is markedly different than ours where the selection is driven by the user (although we can also recover the conditional coverage with optimality in the specific iid online regression model, see \S~\ref{sec:xoriented:optimality}). 
}

\et{Regarding our specific online selective classification application, we note that it differs from the classification-with-abstention and classification-with-rejection literature in its inferential goal. The online branch of this literature typically treats abstention as an action with a cost and studies regret or excess-abstention guarantees \citep{gangrade2021online}, rather than controlling  FCP/IER among selected classifications.
}

\et{For online testing, a very recent work \citep{lu2025feedback} suggests modifying classic online algorithms (e.g., LORD or SAFFRON) that do not have access to the true labels $(Y_t)_{t\geq 1}$, to the setting where $Y_t$ is observed after the decision $S_t$ is made. Their methods, unlike ours,  guarantee  finite-sample error control, but at the cost of lower power. For the example in Figure~\ref{fig:intro_ND}, their method LF makes no discoveries. See \S~\ref{subsec-simND} for additional power comparisons. 
}

\et{Another recent work closely related to our online testing application is SAST \citep{gang2023structure}, which differs from both our and \cite{lu2025feedback} in that it  is built for the classical online testing setting in which the hypothesis states remain unobserved. It uses conditional local false discovery rates, estimated under a two-groups mixture model with smoothly varying sparsity and signal distributions, and adaptively allocates alpha-wealth to improve power.   SAST provides asymptotic online FDR control under its modeling assumptions and is designed to improve power in standard online multiple testing, specifically in settings where the two-groups mixture model can be well approximated. Unlike our method, it does not provide  adversarial FDP bounds or use label feedback to update the procedure.  Moreover, our framework is not limited to online testing: the same OnlineSCI recursion also covers online selective prediction intervals and selective classification, including informative selection rules that may depend on the resulting prediction set.
}

\et{Finally, let us note that while explicit rates for \method\: are provided in our work, we do not think these are optimal in general. It is likely that a method only focusing on the iid framework improves the convergence rate simply by inverting an empirical counterpart of the IER quantity, in the spirit of \cite{gang2023structure} for the testing setting. For this inversion method, the rate would be non-parametric and depending on the regularity of the IER functional. By contrast, rates of \method\: always involve an additional $\beta$-term coming from the step sequence $(\gamma_t)_t$, which may deteriorate the convergence rate. However,  again,  the merit of \method\: is to provide some convergence rate of IER (in 'favorable' situations) {\it in addition to} an adversarial FCP control, which would not be obtained for such an inversion procedure.}

\begin{table}[h!]
\centering
\scriptsize
\setlength{\tabcolsep}{3pt}
\renewcommand{\arraystretch}{1.05}

{\scriptsize (a) Online prediction}

\vspace{1mm}

\begin{tabularx}{\textwidth}{
>{\raggedright\arraybackslash}p{0.13\textwidth}
>{\raggedright\arraybackslash}p{0.12\textwidth}
>{\raggedright\arraybackslash}p{0.18\textwidth}
>{\raggedright\arraybackslash}p{0.22\textwidth}
>{\raggedright\arraybackslash}X
}
\hline
\rowcolor{blue!15}
\textbf{Methods} &
\textbf{Selection} &
\textbf{Cumulative target} &
\textbf{One-step target and conditioning} &
\textbf{Properties and optimality} \\
\hline

ACI \cite{gibbs2021adaptive} and DtACI \cite{gibbs2024conformal} &
No selection &
Nonselective pathwise FCP bound &
No selected-error target &
Distribution-shift accommodation; no selective result \\

\hline
Decaying-step ACI \cite{angelopoulos2024online}&
No selection  &
Nonselective pathwise FCP bound &
One-step unconditional miscoverage asymptotically under regularity conditions &
Threshold convergence to the oracle in the nonselective setting \\

\hline
CAP \cite{bao2025cap}  EXPRESS \cite{sale2025online} PEMI \cite{zheng2026online}&
Restricted selection &
Finite-sample expected FCR &
Marginal selection-conditional coverage, not history-conditional
 &
Exchangeability/i.i.d. assumptions; no conditional rate, optimality or online recursion \\

\hline
Selective DtACI \cite{bao2025cap} &
General selection; not investigating the informative case &
Selective pathwise FCP bound &
No proved IER guarantee &
Dynamic step-size adaptation; no conditional rate or optimality \\

\hline
OnlineSCI (this paper) &
General selection; including the informative case  &
Selective pathwise FCP bound &
History-conditional IER approximation with explicit rate under concentration and regularity assumptions &
Model-specific oracle convergence and optimality. Informative selection handled with remainders depending on a failure probability $\delta$. \\

\hline
\end{tabularx}

\vspace{2mm}

{\scriptsize (b) Online testing}

\vspace{1mm}

\begin{tabularx}{\textwidth}{
>{\raggedright\arraybackslash}p{0.1\textwidth}
>{\raggedright\arraybackslash}p{0.15\textwidth}
>{\raggedright\arraybackslash}p{0.18\textwidth}
>{\raggedright\arraybackslash}p{0.18\textwidth}
>{\raggedright\arraybackslash}X
}
\hline
\rowcolor{blue!15}
\textbf{Methods} &
\textbf{Label / feedback} &
\textbf{Cumulative target} &
\textbf{One-step target and conditioning} &
\textbf{Properties and optimality} \\
\hline

LORD++ \cite{ramdas2017online} &
Label not revealed; $p$-values available online &
Online finite sample unconditional FDR or mFDR &
No history-conditional instantaneous false discovery guarantee &
$\alpha$-investing framework; no pathwise FDP bound or optimality \\

\hline
GAIF / Feedback LORD \cite{lu2025feedback} &
Label revealed after decision &
Online finite sample unconditional FDR or mFDR &
No history-conditional instantaneous false discovery guarantee &
Feedback-enhanced online testing; information regime closer to OnlineSCI but no pathwise FDP guarantee or optimality \\

\hline
SAST \cite{gang2023structure} &
Label not revealed; local-fdr inversion &
Asymptotic online FDR under the model assumptions &
No history-conditional instantaneous false discovery guarantee  &
Local-fdr estimation and model-based convergence to the oracle; no pathwise FDP bound \\

\hline
OnlineSCI &
Label revealed after decision; adaptive test statistics may use past labels &
Pathwise FDP bound &
History-conditional instantaneous false-discovery rate asymptotic control &
Convergence toward the oracle lfdr rule; one-step marginal  false discovery rate asymptotic control / optimality in super-critical regimes, with a remainder term including a failure probability $\delta$ \\

\hline
\end{tabularx}

\vspace{2mm}

\noindent

\caption{\et{
Comparison of error targets, information regimes, and theoretical guarantees in related online methods.
 The main novelty of the method is the combination of 
pathwise selected-error bounds, history-conditional IER rates, adaptive scores, and threshold-coupled selection in one recursion. 
Our algorithms cover both the 'prediction set' (panel a) and 'testing' (panel b) frameworks. They also cover the online selective classification setting, which, to the best of our knowledge, has not previously been considered in the online literature with FCP or IER control (so this case is omitted from this comparison table). 
FCR/FDR are expectations of cumulative error proportions; pathwise FCP/FDP bounds concern the realized stream.
A guarantee in one column does not generally imply a guarantee in another.
}
    \label{tab:novelties}}
\end{table}

\section{Complementary results}\label{sec:otherresults}

We gather here some results that have been postponed to the appendix for space reasons.

\subsection{Rates for iid model and a time-constant selection}\label{sec:appli_timeconstant}

Theorem~\ref{th:rateconvfull} leads to the following result.

\begin{corollary}\label{cor:rateconv}
    Consider an iid model and the procedure $\mathcal{R}^{\OSCI}$ using   
    a time-constant $X$-oriented selection with $\P(S(X_1)=1)\in(0,1]$ and some prediction rule, with a corresponding threshold sequence $(q_t^{\OSCI})_{t\geq 1}$ \eqref{equACIquantile} with a step size sequence $\gamma_j= c j^{-\beta}$, $\beta\in (1/2,1)$, and  starting point $q_1^{\OSCI}\in [0,B]$. Then  
    there exist $q^{\bench}_{\alpha}\in (0,B)$, $T_0\geq 1$ 
and $C_0>0$ such that for $t\geq T_0$, we have
\begin{equation}\label{equ-conviid-specialcases}
\probp{|\IER_t(\mathcal{R}^{\OSCI}) - \alpha|+ \abs{\qt_t^{\OSCI}-q^{\bench}_{\alpha}} \leq C_0 (\rho_t(\beta)+\bar{D}_t) }\geq 1-C_0/t,
\end{equation}
where $\rho_t(\beta)=\rho_t(\beta;1/2)$ as in \eqref{equ:raterho-typicalconc}, $(\bar{D}_t)_{t\geq 1}$ satisfies \eqref{equ:Dtcondition} for $\eta_t=1/t$ and a sequence $(D_t)_{t\geq 1}$ in \eqref{equDtT} \et{such that $D_{k(t)}\propto D_t$ for $k(t)\propto t$ as $t$ is large,} and based on a function $\Pi^{\bench}$ which satisfies Assumption~\ref{ass:easy} with the above $q^\bench_{\alpha}$. 
\end{corollary}

Corollary~\ref{cor:rateconv} is proved in \S~\ref{sec:proofcor:rateconv}. 
In addition to the previous result, Theorem~\ref{th:genL2-first} also provides the following in-expectation convergence result.

\begin{corollary}
    \label{th:MSEconv}
  Consider an iid model and a time-constant $X$-oriented selection with $\P(S(X_1)=1)\in(0,1]$. 
     Then the \method\  threshold sequence $(q_t^{\OSCI})_{t\geq 1}$ with this selection rule, a step size sequence $\gamma_j= c j^{-\beta}$, $\beta\in (0,1)$ and  starting point $q_1^{\OSCI}\in [0,B]$ is such that  there exist $q^{\bench}_{\alpha}\in (0,B)$, $T_0\geq 1$ 
and $C_0>0$ such that for $t\geq T_0$, we have
\et{
\begin{equation}\label{equ:MSEconv}
\paren{\e{\abs{\IER_t(\mathcal{R}^{\OSCI}) - \alpha}\ind{q_t^\OSCI\in [0,B]}}}^2 + \e{\paren{q_{t}^{\OSCI}-q^{\bench}_\alpha}^2}
\leq  C_0 (t/\log(t))^{-\beta} \paren{1 + \sum_{k=1}^{t-1}\widetilde{D}_k},    
\end{equation}
}
where $(\widetilde{D}_t)_{t\geq 1}$ is a (deterministic) sequence tending to zero and eventually decreasing such that $\E[D_t]\leq \widetilde{D}_t$ for all $t\geq 1$.
 The sequence $(D_t)_{t\geq 1}$ in \eqref{equDtT} is based on a function $\Pi^{\bench}$ which is arbitrary but must satisfy Assumption~\ref{ass:easy} with the above $q^\bench_{\alpha}$.
\end{corollary}

The proof of Corollary~\ref{th:MSEconv} is provided in \S~\ref{proofth:MSEconv}. 
Note that Corollaries~\ref{cor:rateconv}~and~\ref{th:MSEconv} are used in \S~\ref{sec:appliiidregression_timeconstant} in the regression case, but can be applied in the classification case as well.

\subsection{Conditional coverage and optimality}\label{sec:condcovopt}

\et{In addition to the marginal selected coverage property of Corollary~\ref{cor:rateconveasyselect}, we show here that $\mathcal{R}^{\OSCI}$ can produce prediction sets with asymptotic {\it pointwise conditional} coverage, that is, conditionally on $X=x$ for any $x$.} In addition, any other method satisfying this asymptotic conditional coverage should produce prediction sets with a larger averaged length, up to some vanishing terms.
For this, we use the following assumption (made in addition to Assumption~\ref{ass:densityf}):
\begin{assumption}\label{ass:densityf:decreasingonly}
In the regression model \eqref{regression:model}, the density $f$ of $\xi$ is decreasing on $[0,+\infty)$. Also, $\sup_{x\in \mathcal{X}}\sigma(x)<+\infty$.
\end{assumption}
With Assumption~\ref{ass:densityf:decreasingonly}, 
the prediction set of \method\  approaches 
\begin{equation}\label{optpredictionset}
\cC^*(x)=\brac{\mu(x)\pm \sigma(x) \ol{F}^{-1}(\alpha/2) }=\{y\in \mathcal{Y}\::\:  f_{Y|X=x}(y) \geq c^*(x)\},  
\end{equation}
with $f_{Y|X=x}(y)=f(|y-\mu(x)|/\sigma(x))/\sigma(x)$ being the density of the distribution of $Y$ given $X=x$ and $c^*(x)=f(\ol{F}^{-1}(\alpha/2))/\sigma(x)$. Hence, $\cC^*(x)$ is a level set of the conditional density that satisfies the conditional coverage $\probp{Y\in \cC^*(x)\mid X=x} = 1-\alpha$, and it is well known that its length $2\sigma(x) \ol{F}^{-1}(\alpha/2)$ is minimal among all prediction sets with  conditional coverage at least $1-\alpha$, see \cite{clayton2006learning} among others. 
Combining this with Corollary~\ref{cor:rateconveasyselect} shows the optimality result.

Formally, we should consider the slightly stronger pointwise version of convergence condition \eqref{ass:esti:average} for $\hat{\mu}_t$ and $\hat{\sigma}_t$: for $t$ large enough and some constant $C_1>0$,
\begin{align}
\sup_{x\in \mathcal{X}} \E[(\hat{\mu}_t(x)-{\mu}(x))^2] + \sup_{x\in \mathcal{X}}\E[(\hat{\sigma}_t(x)-{\sigma}(x))^2] &\leq C_1 (\log t)^{2\tau_0} t^{-2\tau},\label{ass:esti:averagepoint}
\end{align}
for some $\tau_0>0$ and $\tau\in (0,1/2]$.

\begin{corollary}\label{cor:optrateiidregress}    
Consider the online iid regression model \eqref{regression:model} satisfying Assumptions~\ref{ass:densityf}~and~\ref{ass:densityf:decreasingonly}.
Consider $\mathcal{R}^{\OSCI}$ as given in Algorithm~\ref{alg:PredictXoriented}, using   
    a time-constant $X$-oriented selection with $\P(S(X_1)=1)$ fixed in $(0,1]$, $\mathcal{F}_{t-1}$-measurable estimators $\hat{\mu}_t, \hat{\sigma}_t$ satisfying \eqref{ass:esti}-\eqref{ass:esti:averagepoint}, and a threshold sequence $(q_t^{\OSCI})_{t\geq 1}$ \eqref{equACIquantile} with a step size sequence $\gamma_j= c j^{-\beta}$, $\beta\in (1/2,1)$, and  starting point $q_1^{\OSCI}\in [0,B]$.
    Then $\mathcal{C}^{\OSCI}_t(x) = \big[ \hat{\mu}_t(x) \pm \hat{\sigma}_t(x) \Phi^{-1}\big(\frac{q_t^{\OSCI}+1}{2}\big)\big]$ is such that:
    \begin{itemize}
        \item[(i)] it has approximate pointwise conditional coverage: for all  $x\in \mathcal{X}$,
        \begin{equation}
            \label{equcondconv}
            \probp{Y_t\in \cC^{\OSCI}_t(X_t)\mid X_t=x}\geq 1-\alpha - C_0 [\rho_t(\beta)+(\log t)^{\tau_0}t^{-\tau}+(\log t)^{\nu_0} t^{-\nu}], 
        \end{equation}
        for some constant $C_0>0$ and $\rho_t(\beta)=\rho_t(\beta;1/2)$ as in \eqref{equ:raterho-typicalconc};
        \item[(ii)] 
it mimics the oracle procedure:  choosing $\beta>1- \tau$, 
for all $x\in \mathcal{X}$, and $M>0$,
\begin{align}\label{equmimicoracleiidregress}
\e{|\cC^{\OSCI}_t(x)|\wedge M}\leq |\cC^*(x)|  + C_0(\log t)^{\tau_0/2} t^{1/2- \tau/2-\beta/2 }+ M/t ,
\end{align}
for some constant $C_0>0$ (depending on all other constants but not on $M$).
        \item[(iii)] it approximately dominates any other method with approximate pointwise conditional coverage: choosing $\beta>1- \tau$,
        for any method producing a sequence $\cC'_t(x),t\geq 1,x\in \mathcal{X}$ of ($\mathcal{F}_{t-1}$-measurable) prediction set rules with the conditional coverage property $\forall x\in \mathcal{X}$, $\probp{Y_t\in \cC'_t(x)\mid X_t=x}\geq 1-\alpha - \eta'_t$, for some rate $\eta'_t>0$ tending to zero,  we have for $t$ large enough, for all $x\in \mathcal{X}$, and $M>0$,
        \begin{equation}
            \label{equcondoptimality}
            \e{|\cC_t^{\OSCI}(x)|\wedge M} \leq \e{|\cC'_t(x)|}+C_0 (\log t)^{\tau_0/2} t^{1/2- \tau/2-\beta/2 }+C_0 \eta'_t+M/t,
        \end{equation}
        for some constant $C_0>0$ (depending on all other constants but not on $M$).
    \end{itemize}
\end{corollary}

Corollary~\ref{cor:optrateiidregress} is proved in \S~\ref{proof:cor:optrateiidregress}. It relies both on the in-probability\footnote{The conditional coverage (i) relies only on in-probability bounds and thus can be extended to more general selection rules.} and  in-expectation bounds obtained in Corollary~\ref{cor:rateconveasyselect}. A byproduct of \eqref{equcondoptimality} is that, choosing $\beta>1- \tau$, and $M=\sqrt{t}$, we have $\limsup_{t\to \infty}\paren{\e{|\cC_t^{\OSCI}(x)|\wedge \sqrt{t}}-\e{|\cC'_t(x)|}}\leq 0.$
Hence, Corollary~\ref{cor:optrateiidregress} states that no method can outperform \method\  by a large amount, which is an optimality property. While other methods could also satisfy a similar optimality property (as the oracle one if $q^\bench_\alpha$ is known), recall that \method\  also markedly has a valid FCP bound \eqref{FCPbound_decre} in the adversarial model. Hence, it makes the best of two worlds.

\subsection{Feasibility in the informative case}\label{sec:commentdelta}

The next result shows that $\delta$ is incompressible: it is impossible to provide that $\P(\lim_t q_t^\OSCI = q_\alpha^\bench)> 1-\delta$. 

\begin{lemma}[Incompressibility of $\delta$]
    \label{lem:impossible}
    Consider the \method\  threshold sequence $(q_t^{\OSCI})_{t\geq 1}$ \eqref{equACIquantile}, $\delta$ given by \eqref{deltaprobaJtfinite} and any $q^\bench\in (0,B^*)$ with $B^*$ as in Assumption~\ref{ass:support}. Then $\P(\lim_t q_t^\OSCI = q^\bench)\leq 1-\delta$.
\end{lemma}

\begin{proof}
Assume $\delta>0$ (otherwise the result is clearly true). 
By definition, we have for all $q\in [0,B^*)$ and $t\geq 1$, $\probp{\sup_{T\geq t} q_T^\OSCI>q}\geq \delta$. Hence taking any $q\in (q^\bench,B^*)$ (which is possible because $q^\bench\in (0,B^*)$), we have 
$$\probp{q_t^\OSCI \nrightarrow q^\bench} \geq \probp{\lim_t\sup_{T\geq t} q_T^\OSCI\geq q}=\lim_t\probp{\sup_{T\geq t} q_T^\OSCI\geq q}\geq \delta,$$  
which shows the result.
\end{proof}

By contrast, $q_t^\OSCI$ is uniformly away from $B^*$ by below with probability $1-\delta'$ for $\delta'>\delta$:

\begin{lemma}[Feasibility without the $\delta$ term]
    \label{lem:possible}
    Consider the \method\  threshold sequence $(q_t^{\OSCI})_{t\geq 1}$ \eqref{equACIquantile} and $\delta$ given by \eqref{deltaprobaJtfinite}. Then for all $\delta'>\delta$, there exists $T_0=T_0(\delta')$ and $\bar q=\bar q(\delta')\in (0,B^*)$ such that $\probp{\sup_{T\geq T_0} q_T^\OSCI\leq \bar q}>1-\delta'$.
\end{lemma}

\begin{proof}
   This comes from \eqref{deltaprobaJtfinite} and the basic property of an infimum.
\end{proof}

\subsection{\et{Benchmark function for classification and testing}}\label{sec:benchclasstest}

\et{We consider both the classification case and testing cases, for which the benchmark function can be expressed as ($B=1$)
\begin{align}
    \Pi^\bench(q)&= 1-\e{W(X)\mid W(X)>q},\:\:\:q\in [0,1],\label{PizeroInfo}
\end{align}
with the convention $\Pi^\bench(q)=1-B^*$ when  $\P(W(X)>q)=0$, for some random variable $W(X)$ satisfying Assumption~\ref{ass:support}.
Applying Lemma~\ref{lem:decreasing}, if $\alpha\in (1-B^*,1-\e{W(X)})$, the benchmark function $\Pi^\bench$ satisfies Assumption~\ref{ass:easy}  with $q_\alpha^\bench\in (\kappa,1-\alpha)$ such that 
\begin{equation}
\label{equqalphainfo}
    \e{W(X)\mid W(X)>q_\alpha^\bench}=1-\alpha.
\end{equation} 
Indeed, Lemma~\ref{lem:decreasing} ensures that 
$\Pi^\bench$ is decreasing on $[\kappa,B^*]$ and continuous on $[0,1]$ with $\Pi^\bench(0)=\Pi^\bench(\kappa)=1-\E[W(X)]$ and $\Pi^\bench(B^*)=\Pi^\bench(1)=1-B^*$. In addition, 
 $\Pi^\bench$ is differentiable on $(\kappa,B^*)$, with a continuous derivative on any compact set of $(\kappa,B^*)$ which is away from $0$.  
As a result, $\Pi^\bench$ satisfies Assumption~\ref{ass:easy}. 
}
\et{As defined in \eqref{sec:incompresscritic}, the critical level \eqref{equalphastargen} is thus given by 
\begin{equation}
\label{equalphastar}
\alpha^*:=1-B^*.    
\end{equation}
}

\subsection{Rate and optimality for selection by prediction}\label{sec:selectpredict}

Consider the iid regression model \eqref{regression:model} with $\sigma(x)=1$ for simplicity. The authors of 
\cite{Jin2023selection} are concerned with testing the null hypothesis $Y_t\leq y_0$ (for prespecified $y_0\in \R$) in order to discover $Y_t>y_0$. We consider this task in an online fashion, by applying the testing example (\S~\ref{sec:selectiveNDiid}) with $\tilde{Y}_t = \ind{Y_t>y_0}$ and specific estimates in the adaptive score. 

More specifically, we use \method\  with the prediction set rule $\mathcal{C}_t(x,q)=(y_0,+\infty)$ (so that the non-coverage error is $\Err_t(y, C)=\ind{y\leq y_0}$) and the selection rule is \eqref{eq:selectionstat} with the test statistics $W_t(x)=1-\widehat{\lfdr}_t(x)$ \eqref{equ:WtselecND}, where $\widehat{\lfdr}_t(x)$ is an estimate of the local-fdr \eqref{equ:lfdr}, defined in this framework by
\begin{align}
{\lfdr}(x) &= \mP(Y\leq y_0\mid X=x) = F\paren{y_0-{\mu}(x)}\label{equ:lfdr:JC},
\end{align}
where $F$ is the cdf of the noise $\xi$ of the regression model.
\method\  and Algorithm~\ref{alg:onlineBH} with $\tilde{Y}_t = \ind{Y_t>y_0}$ gives rise to a new online algorithm for selection by prediction. 

We can derive convergence properties of this new procedure by applying Theorem~\ref{th:rate:ND} and Corollary~\ref{optimalityNDiid} in this particular setting.
For this, we should check Assumption~\ref{ass:support} with the test statistics $W(x):=1-{\lfdr}(x)$. 

\begin{assumption}
\label{ass:selecpredict}
The online iid regression model \eqref{regression:model} is such that $\mtc{X}= [0,1]^d$, $\sigma(x)=1$, $x\in \mtc{X}$, the regression function $\mu(\cdot)$ is a Lipschitz function. 
In addition, Assumption~\ref{ass:supportmu} holds and 
the noise $\xi$ has a continuous density which is positive on $[y_0-\mu_{\max},y_0-\mu_{\min}]$. 
\end{assumption}
Under Assumption~\ref{ass:selecpredict}, it is clear that Assumption~\ref{ass:support} holds with $\kappa=1-F(y_0-\mu_{\min})\geq 0$ and $B^*=1-F(y_0-\mu_{\max})\leq 1$ (hence $\alpha^*=F(y_0-\mu_{\max})$). In addition, by considering a classical kernel-based lfdr estimate $\widehat{\lfdr}_t(x)$ of $\lfdr(x)$ (not using the explicit expression \eqref{equ:lfdr:JC} because $F$ is unknown), we can meet assumption 
 \eqref{ass:estiND} for $\nu_0=1$ and $\nu=1/(2+d)$ \citep{tsybakov2008nonparametric}. Applying Theorem~\ref{th:rate:ND} and Corollary~\ref{optimalityNDiid} gives the convergence rate $(\log t)^{1/2} t^{1/2-\beta} + t^{\beta-1} +t^{-1/(4+2d)}$ with optimality. Here, the IER can be interpreted as an instantaneous false discovery rate control, while the optimality in \eqref{eq-OSCIpower-ND} is in terms of instantaneous true discovery rate. Here, note that a criticality phenomenon typically occurs, that is, $\alpha^*>0$, for instance when the noise $\xi$ has a density supported in $\R$, as in the Gaussian case. Hence, not all error level $\alpha>0$ can be reached for the selection by prediction problem.

\section{Additional numerical experiments}\label{secaddxp}

\subsection{Additional figures for synthetic data}

Figures~\ref{fig:ND_xp_multiple}~and~\ref{fig:SC_xp_multiple} are given below.

\begin{figure}[h!]
		\centering
		\includegraphics[width=.45\linewidth]{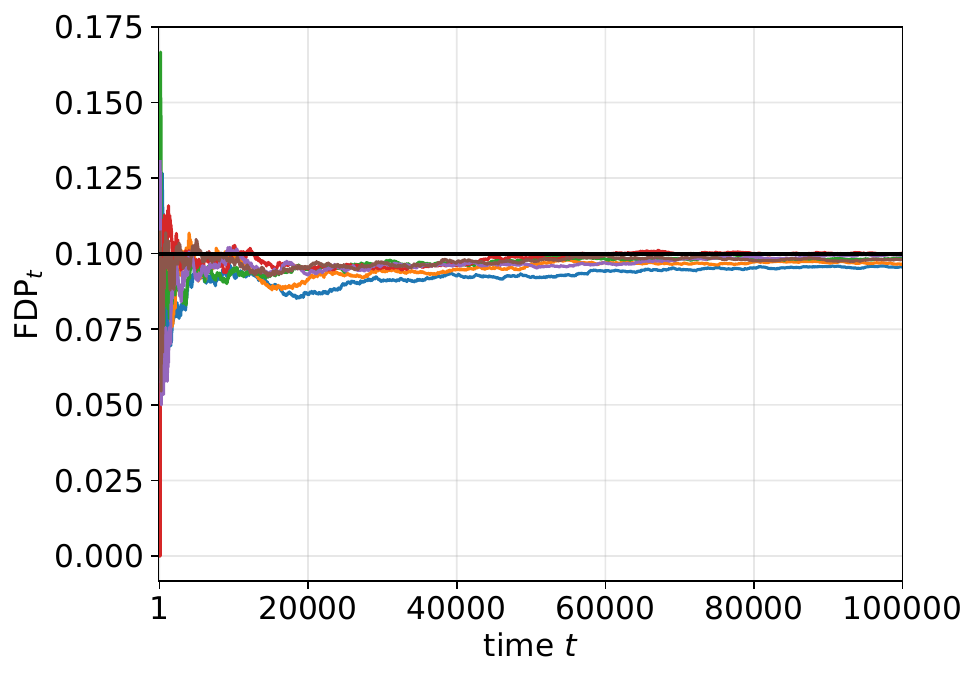}
        \includegraphics[width=.45\linewidth]{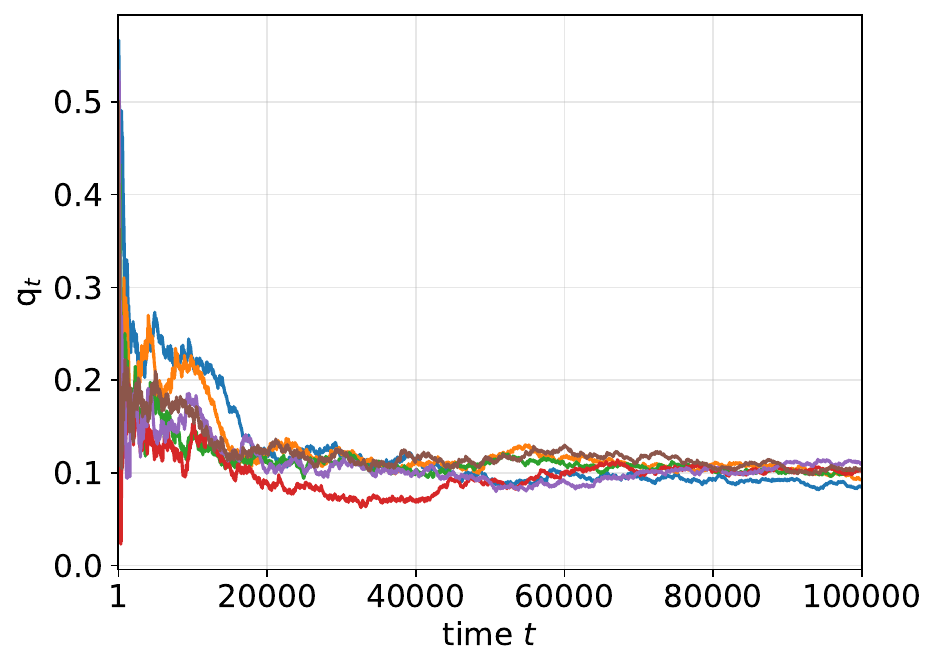} \\
        \caption{Online conformal testing with adaptive \method: same as Figure~\ref{fig:ND_xp} top panels with 
        $5$ different data set generations from the same model. \et{Black line corresponds to $\alpha=0.1$.}}
        \label{fig:ND_xp_multiple}
\end{figure}

\begin{figure}[h!]
		\centering
		\includegraphics[width=.45\linewidth]{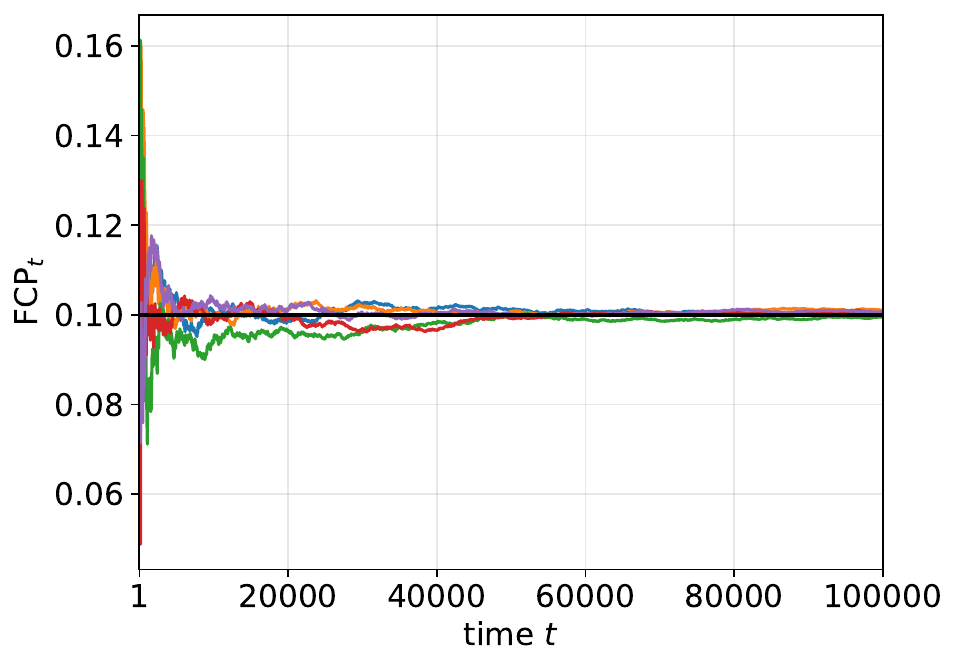}
        \includegraphics[width=.45\linewidth]{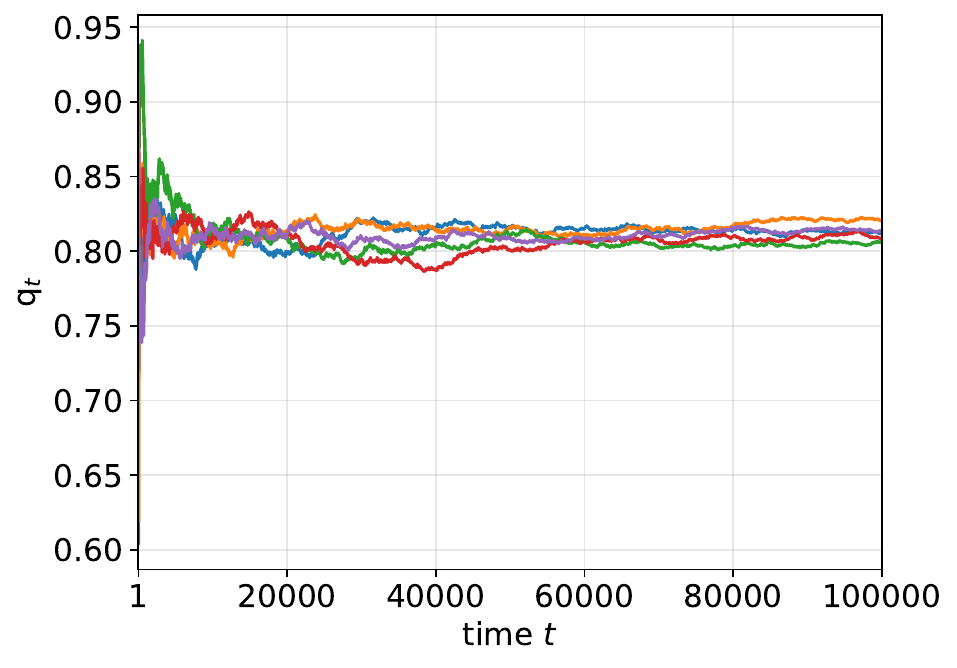} \\
        \caption{
        Online selective classification with adaptive \method: same as Figure~\ref{fig:SC_xp} top panels with 
        $5$ different data set generations from the same model. \et{Black line corresponds to $\alpha=0.1$.}}
        \label{fig:SC_xp_multiple}
\end{figure}

\subsection{\et{Comparison with DtACI}}\label{sec:compaDTACI}
\pie{
We compare our method with DtACI \citep{gibbs2024conformal} applied to the selective classification problem (see section \ref{subsec-simClassification}). We evaluate two versions of DtACI:
\begin{itemize}
    \item DtACI with naive selection, where a prediction is selected only if the size of the returned prediction set is equal to $1$;
    \item 
An informative version of selective DtACI \citep{bao2025cap}, where      the threshold update of the DtACI algorithm is only performed if a selection occurs, that is, $\max_{y\in \range{K}} \hat{\pi}_{t}(y|X_t)>q_t$, with a non-coverage error occurring if and only if $Y_t\neq \hat{Y}_t(X_t)$ (by using the notation of \S~\ref{sec:selectiveclassifiid}).
\end{itemize}
The hyperparameters $\eta$ and $\sigma$ of DtACI are set as specified in \citep{gibbs2024conformal}. Additionally, the set of candidate $\gamma$ values is defined as $\left\{\frac{i}{k} \cdot 0.1, i = 1, \dots, k\right\}$ with $k=30$.\\
}

\begin{figure}[h!]
		\centering
		\includegraphics[width=.45\linewidth]{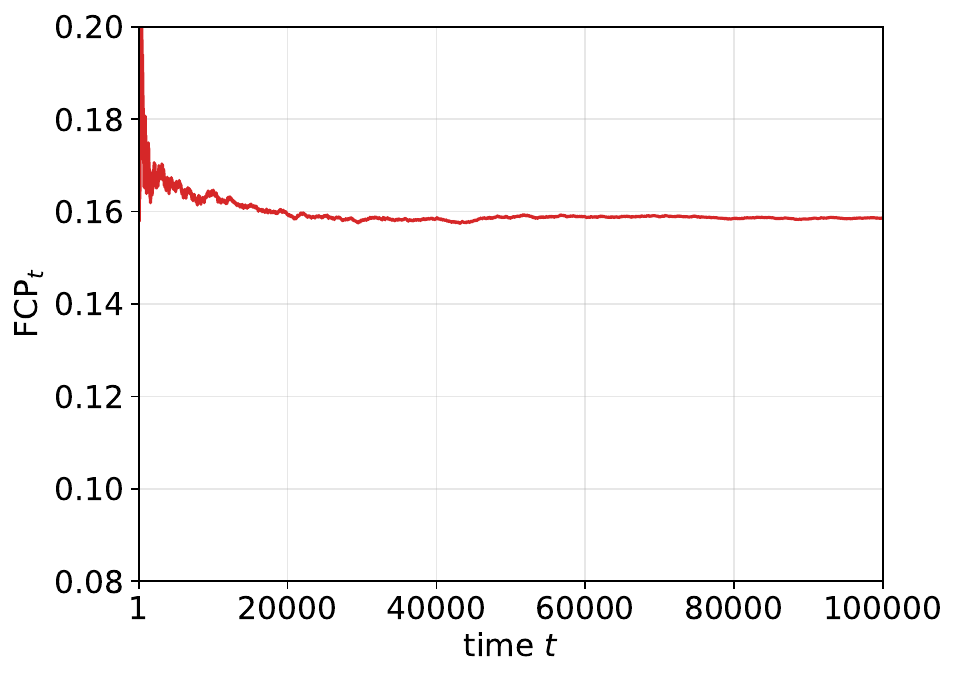}
        \includegraphics[width=.45\linewidth]{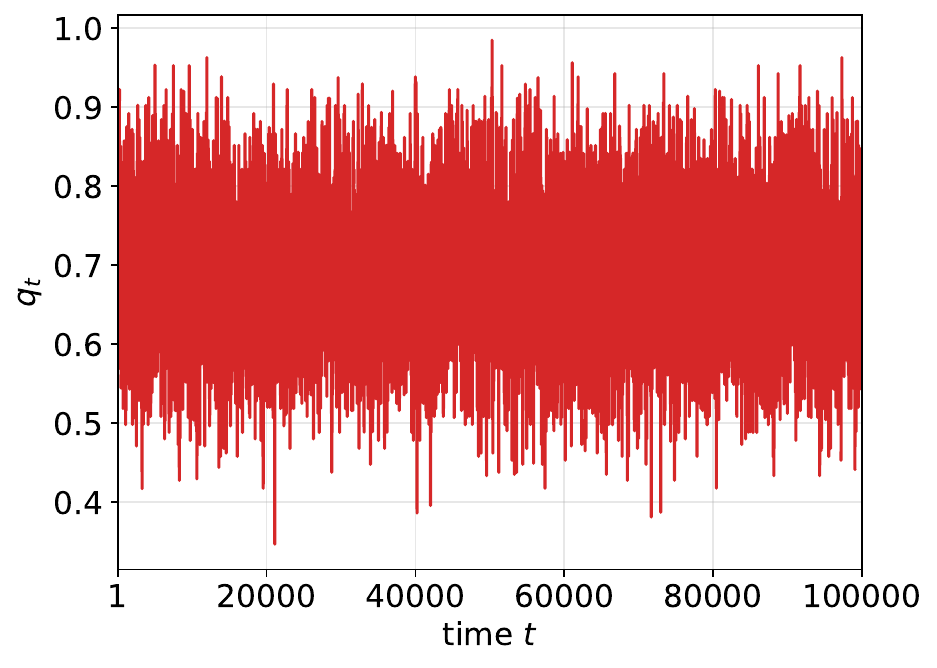} \\
        \caption{
        \pie{
        False classification rate $(\FCP_t)_t$ and threshold sequence $(q_t)_t$ of DtACI with naive selection (see text) in an online selective classification framework.
        }
        \label{fig:DtACI_naive}
        }
\end{figure}

\begin{figure}[h!]
		\centering
		\includegraphics[width=.45\linewidth]{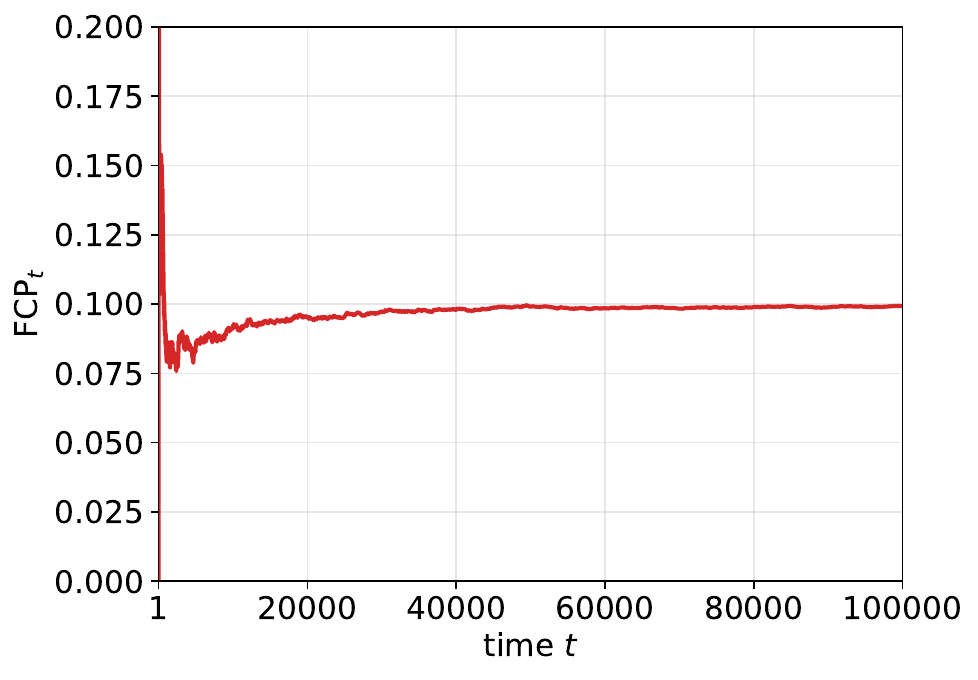}
        \includegraphics[width=.45\linewidth]{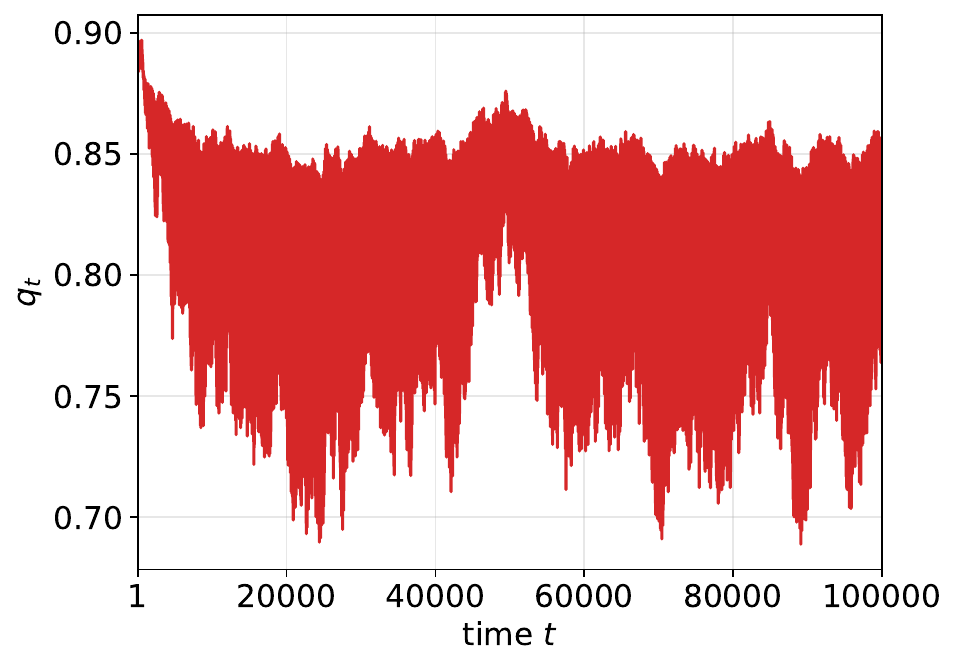} \\
        \caption{
        \pie{
         False classification rate $(\FCP_t)_t$ and threshold sequence $(q_t)_t$ of the new informative version of DtACI (see text)  in an online selective classification framework.
        }
        \label{fig:DtACI_select}
        }
\end{figure}

\begin{figure}[h!]
		\centering
		\includegraphics[width=.7\linewidth]{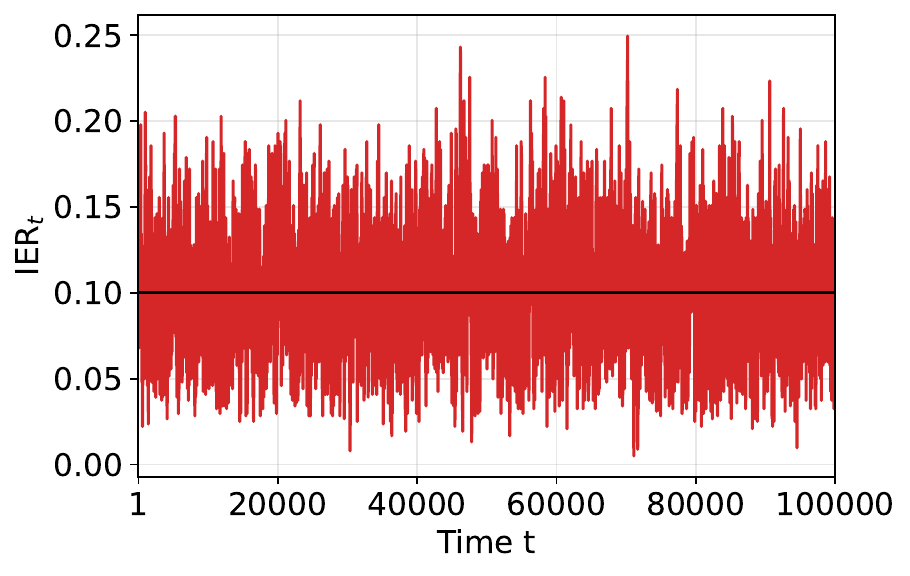} \\
        \caption{
        \pie{
        IER of the new informative version of DtACI (see text) in an online selective classification framework.
        }
        \label{fig:DtACI_IER}
        }
\end{figure}

\pie{
As expected, DtACI with naive selection fails to control the FCP (false classification rate here) at level $\alpha=0.1$ (Figure \ref{fig:DtACI_naive}, left panel). Additionally, the sequence $(q_t)_t$ does not converge (Figure \ref{fig:DtACI_naive}, right panel).
In contrast, our extended version of DtACI successfully controls the FCP at level $\alpha=0.1$ (Figure \ref{fig:DtACI_select}, left panel). However, similar to the naive version, and unlike \method, the sequence $(q_t)_t$ remains unstable (Figure \ref{fig:DtACI_select}, right panel). In addition, 
its estimated IER does not approach $\alpha$ (Figure \ref{fig:DtACI_IER}).
}

\subsection{Real data: online conformal testing and selective classification for \texttt{CIFAR-10} data set}\label{sec:realdata}

We use the image data set \texttt{CIFAR-10} (\url{https://www.cs.toronto.edu/~kriz/cifar.html}), which consists of $60000$ $32 \times 32$ colour images in $10$ classes, with $6000$ images per class.

We consider two applications: online conformal testing and online selective classification. 
\begin{enumerate}
    \item For online conformal testing, the null classes are `\textit{car}' and `\textit{truck}', and the alternative is `\textit{cat}'. The total number of images in the hold-out set $\mathcal{D}_0$ and in the online stream of data are equal to $15000$ and $3000$, respectively. Using the hold-out set, we compute $\widehat{\lfdr}_t$ with a convolutional Neural Network (CNN) with 2 convolutional layers, one pooling layer, and 3 fully-connected trained using a cross-entropy loss and a stochastic gradient descent with a learning rate of $0.001$ and a momentum equal to $0.9$. \et{Finally, $\alpha=0.1$, $q_1=0.4$, and $\gamma_t = t^{-3/4}/2$}. As in the synthetic section, we compare our method to \et{LORD++}.
    \item For online selective classification, we consider the classes `\textit{deer}' and `\textit{cat}'. The total number of images in the hold-out and online stream of data are equal to $10000$ and $2000$, respectively.  We estimate $\probp{Y=y \mid X=x}$ with a CNN with 2 convolutional layers, one pooling layer, and 3 fully-connected trained using a cross-entropy loss and a stochastic gradient descent with a learning rate of $0.001$ and a momentum equal to $0.9$. \et{Finally, $\alpha=0.1$, $q_1=0.8$, and $\gamma_t = t^{-3/4}/2$}. As in the synthetic section, we compare our method to ACI \citep{angelopoulos2024online} where we select if the returned prediction set at time $t$ is of size one.
\end{enumerate}

\begin{figure}[h!]
		\centering
        \begin{subfigure}[b]{0.9\linewidth}
        \includegraphics[width=.45\linewidth, height=.3\linewidth]{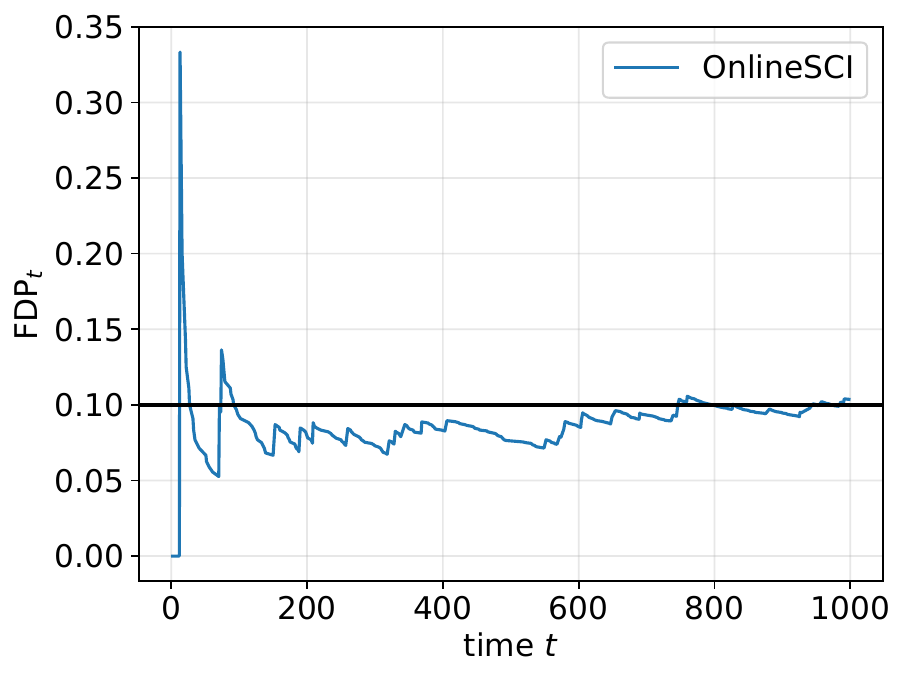}
        \includegraphics[width=.45\linewidth, height=.3\linewidth]{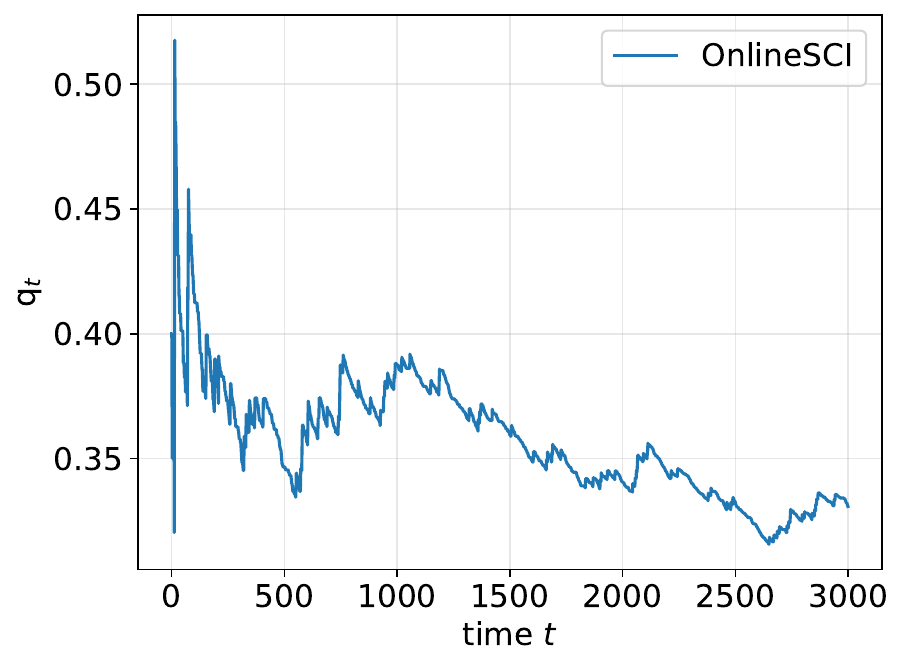}
        \caption{Online conformal testing}
        \label{fig:SC_xp_cifar2_multiple_1}
        \end{subfigure}\\\vspace{1em}
        \begin{subfigure}[b]{0.9\linewidth}
        \includegraphics[width=.45\linewidth, height=.3\linewidth]{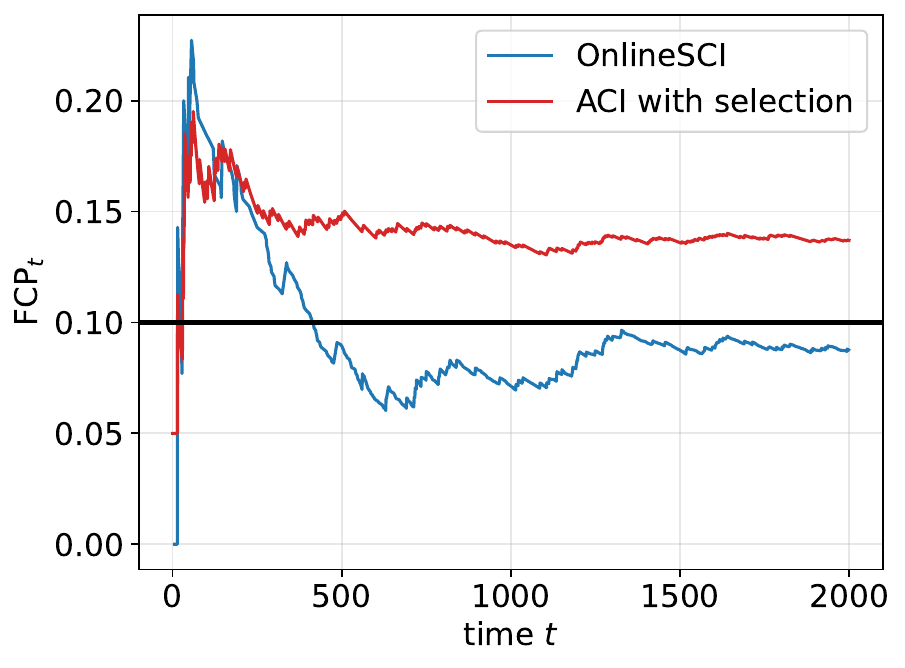}
        \includegraphics[width=.45\linewidth, height=.3\linewidth]{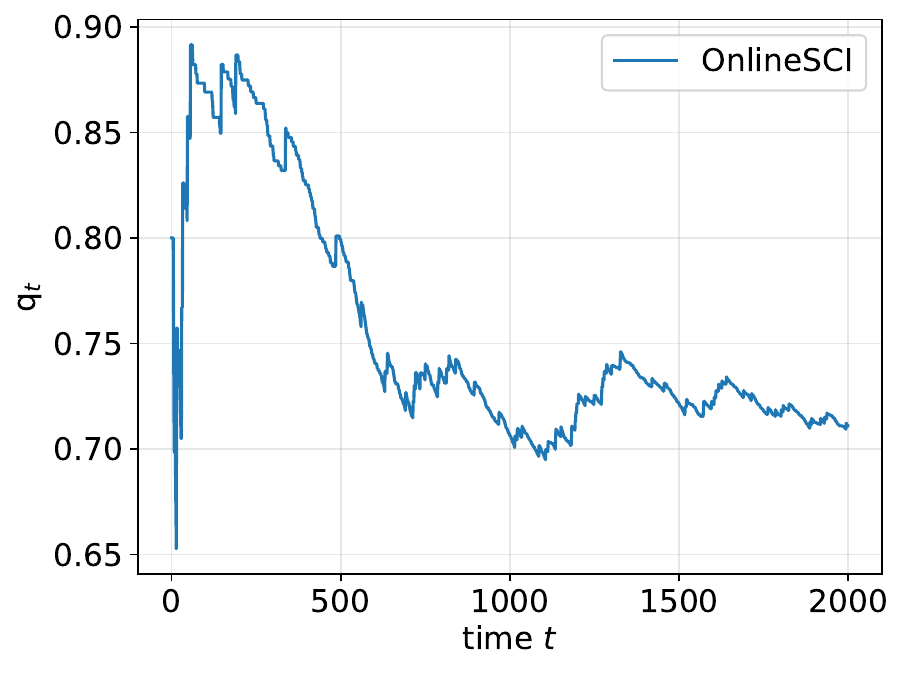}
        \caption{Online selective classification}
        \label{fig:SC_xp_cifar2_multiple_2}
        \end{subfigure}
        \caption{Top panel: Results for \method~on \texttt{CIFAR-10} for online conformal testing. Bottom panel: Results for \method~on \texttt{CIFAR-10} for online selective classification. Left panels: $\FCP_t$ versus time step $t$. Right panels: method threshold $q_t$ versus time step $t$.}
        \label{fig:SC_xp_cifar2_multiple}
\end{figure}

\begin{figure}[htbp]
  \centering
  \begin{subfigure}[b]{1.\linewidth}
    \includegraphics[width=\linewidth]{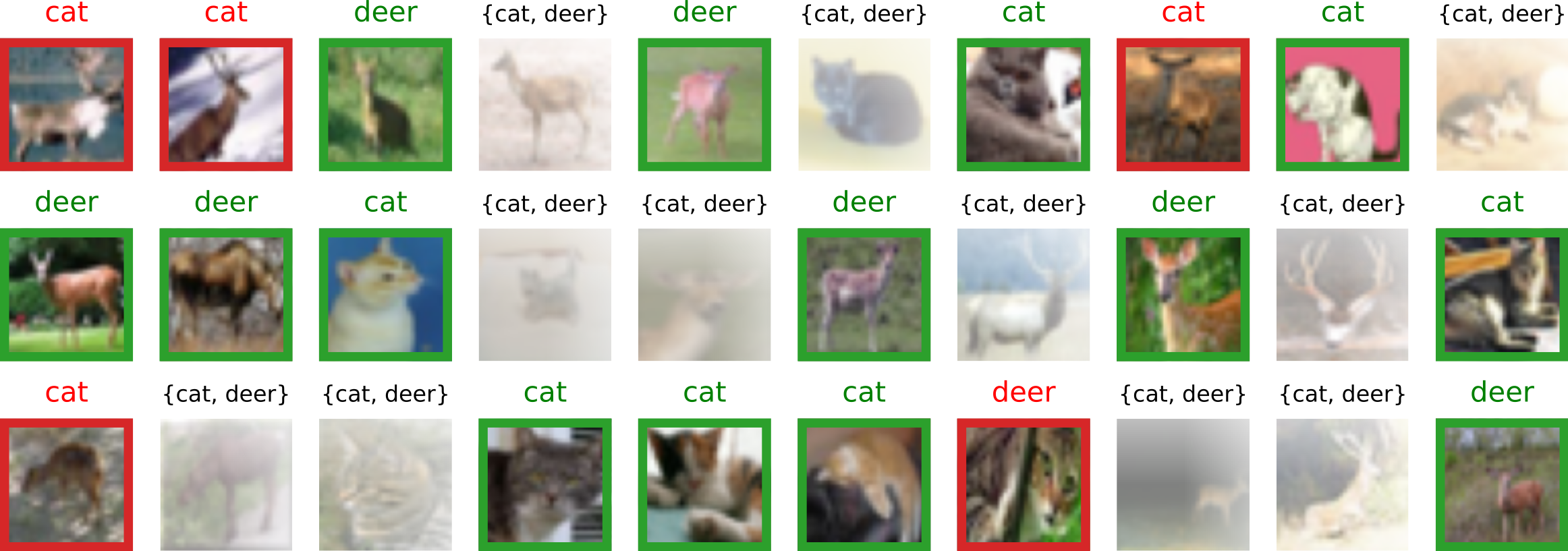}
    \caption{Classical conformal prediction output at $\alpha = 0.1$.}
    \label{fig:intro-selective-classif-top_panel}
  \end{subfigure}

  \vskip\baselineskip 

  \begin{subfigure}[b]{1.\linewidth}
    \includegraphics[width=\linewidth]{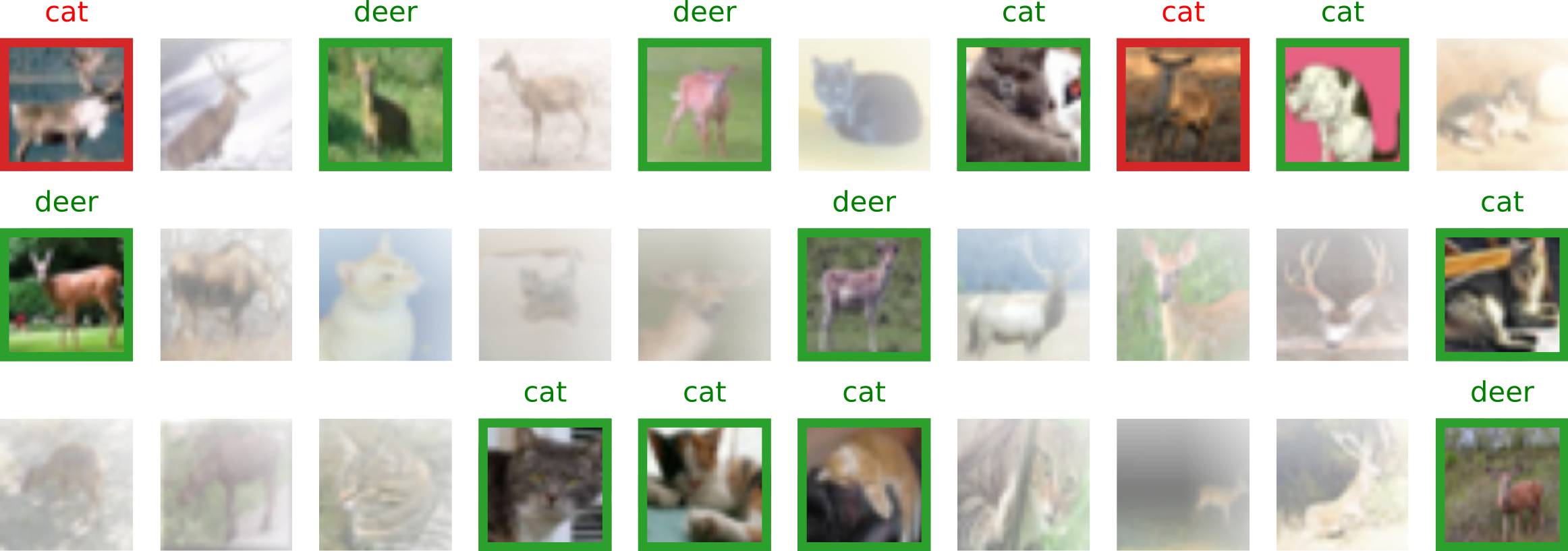}
    \caption{New method \method. }
    \label{fig:intro-selective-classif-bottom_panel}
  \end{subfigure}

  \caption{Illustration of the online selective classification application for \texttt{CIFAR-10} with classes `\textit{cat}' and `\textit{deer}'. The square color indicates whether a classification error occurred: red if it did, green otherwise.  When there is no color, this means that we do not make a selection. The classes on top of each selected image corresponds to the prediction set made by the convolutional neural network classifier.}
  \label{fig:intro-selective-classif-combined_figure}
\end{figure}

Figure \ref{fig:SC_xp_cifar2_multiple} displays, for the two applications, the $\FCP_t$ (left panel) and the sequence $(q_t)_{t\geq1}$ (right panel) returned by \method\  (and by ACI for selective classification). We see that \method~controls the $\FCP$ at the nominal level $\alpha$ (blue curves) in both cases, which is not the case of ACI in the case of online selective classification (red curve). As in the synthetic data, \et{LORD++} is not able to make discovery; the first $p$-values are not accurate enough, as we do not have enough data at these time points. An additional illustration is given in Figure~\ref{fig:intro-selective-classif-combined_figure} for online selective classification.

\subsection{Real data: electricity prices}\label{sec:elec}

To assess the flexibility of \method, we end the numerical experiments with an example where the data stream is not meant to be iid . This time series data set consists of French electricity spot prices from 2016 to 2019. It is an hourly data set, containing $(3 \cdot 365 + 366) \cdot 24 = 35064$ observations. The goal is to predict at day $D$ the $24$ prices of day $D+1$ with the following explanatory variables: day-ahead forecast consumption, day-of-the-week, $24$ prices of the day $D-1$ and $24$ prices of the day $D-7$. A complete presentation of this data set is given in \cite{zaffran2022adaptive}. For simplicity, we use the price predictions of year $2019$ already computed in \url{https://github.com/mzaffran/AdaptiveConformalPredictionsTimeSeries} (non-adaptive estimator).

\begin{figure}[h!]
		\centering
        \includegraphics[width=.47\linewidth, height=.3\linewidth]{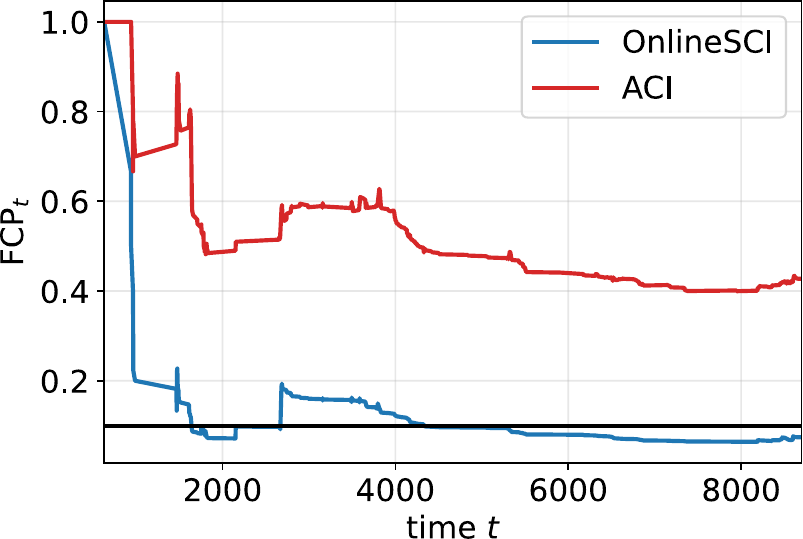}
        \includegraphics[width=.47\linewidth, height=.3\linewidth]{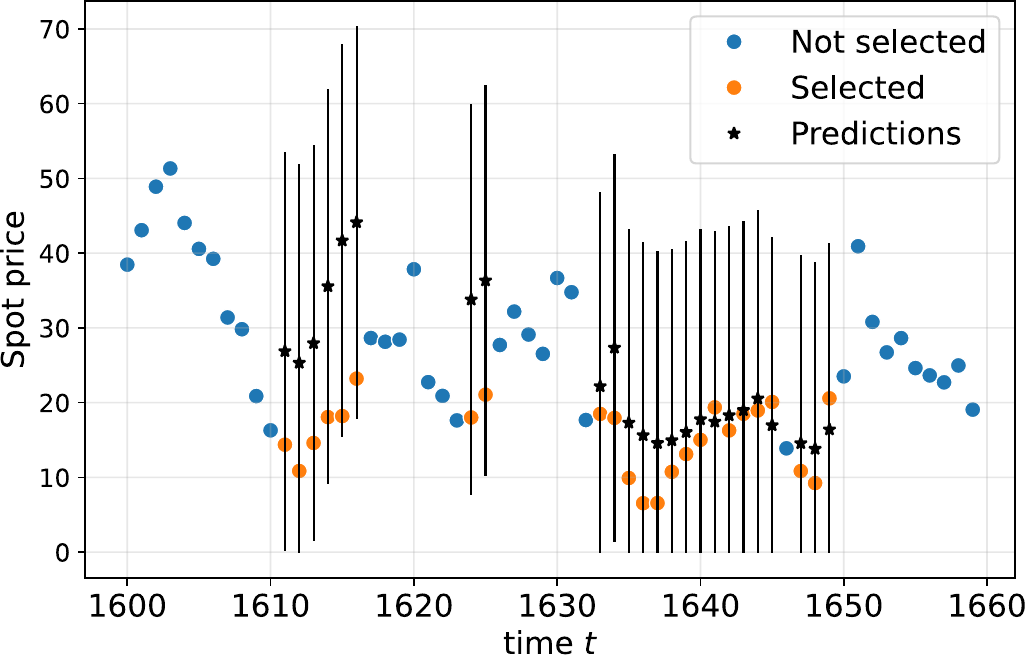}
        \caption{Result for the edf data set, see detail in \S~\ref{sec:elec}.}
        \label{fig:EDF}
\end{figure}

We consider the online selective conformal prediction intervals of \S~\ref{sec:AR}. Specifically, the problem at hand is to construct an interval around the predicted value at time $t$ only if the electricity spot price at time $t-1$ was lower than $20$, that is, $S_t(X_t)=\ind{Y_{t-1}\leq 20}$ for $X_t=(Y_{t-1},\dots,Y_{t-d})\in\R^d$.  The $\FDP_t$ obtained with \method~and \texttt{ACI} when $\alpha=0.1$, $\gamma_t = 10 \cdot t^{-3/4}$ and $q_1=0.8$ are given in Figure \ref{fig:EDF} (Left panel). Examples of prediction sets constructed by \method~on the selected points are displayed in Figure \ref{fig:EDF} (Right panel). 
The conclusion is qualitatively similar to above, and is in accordance with the theoretical findings of \S~\ref{sec:AR} in the autoregressive model.

\subsection{\et{Guidelines and procedure for hyperparameter tuning}}\label{sec:select_process}

\et{We give some insight on how the choice of the hyperparameters of our method $(q_1, (\gamma_t)_t)$ influences its performances. We also propose a systematic approach to select them. \\
}

\et{Let us consider the online conformal testing problem  of Section \ref{subsec-simND}. In this iid setting, the result of \S~\ref{sec:selectiveNDiid} shows that $q_t$ should converge to $q^*$ with high probability. Hence, intuitively, the closer $q_1$ is of $q^*$ the better it is if the step sizes are small enough (the ideal but unrealistic case being that $q_1 = q^*$ and $\gamma_t = 0$). Conversely, if $q_1$ is far from $q^*$, large step sizes are needed to rapidly approach it. Hence, we see that there is an interplay between $q_1$ and $(\gamma_t)_t$. }\\

\paragraph*{\method\:with automatic hyperparameter tuning}
\et{Regarding the previous discussion, we now propose an automatic procedure for tuning the hyperparameters of \method\:
from the holdout set $\cD_0 := (X_{i}, Y_{i})_{-(n-1)\leq i\leq 0}$ of size $n$ available at time $t=1$ .  
First, compute the standard empirical-quantile batch estimator $\widehat{q}^0$ of $q^*$ on $\cD_0$, which will be used as the starting  threshold $q_{-(n-1)}$ of \method (by shifting the index of the recursion). Then, considering a sequence $\gamma^0_t = c \cdot (n+t)^{-\beta}$, $t=-(n-1),\dots,0$, run \method\: for multiple values of $c$ and $\beta$ on a specific grid. Finally, select a pair $(c, \beta)$ from the resulting trajectories using a constrained power criterion, e.g., the pair maximizing power while returning an FCP $< 1.2 \cdot \alpha$ (other criteria could be used depending on the application).} 

\paragraph*{Methods}
\et{We consider the online conformal testing version of \method, assuming for simplicity that the predictor (local-fdr) is the oracle (true) one. We compare the performance of \method\: using  (i) hyperparameters selected via our automated procedure described above when the burn-in dataset $\cD_0$ is of size $n=500$ (in the iid model, our theoretical guarantees are maintained conditionally on $\cD_0$); (ii) hyperparameters in specific ranges: for the initial point $q_1$, we consider a grid of $10$ values defined as $q_1\in \{0.1 + i \cdot 0.8/9 , i = 0, \dots, 9\}$. We also consider the oracle threshold $q^*$ and $\widehat{q}_0$ which is an estimator of $q^*$ calculated on an independent dataset of size $500$. For the step-size, we use $\gamma_t = c \cdot t^{-\beta}$, where the parameters $(c, \beta)$ are chosen in a double grid of size $10 \times 10$ taking values from $0.2$ to $1$ for $c$ and from $0.5$ to $0.9$ for $\beta$; (iii) we also test the specific ideal pair $(q_0, \gamma_t) = (q^*, 0)$.\\
}

\begin{figure}[h!]
		\centering
		\includegraphics[width=.48\linewidth]{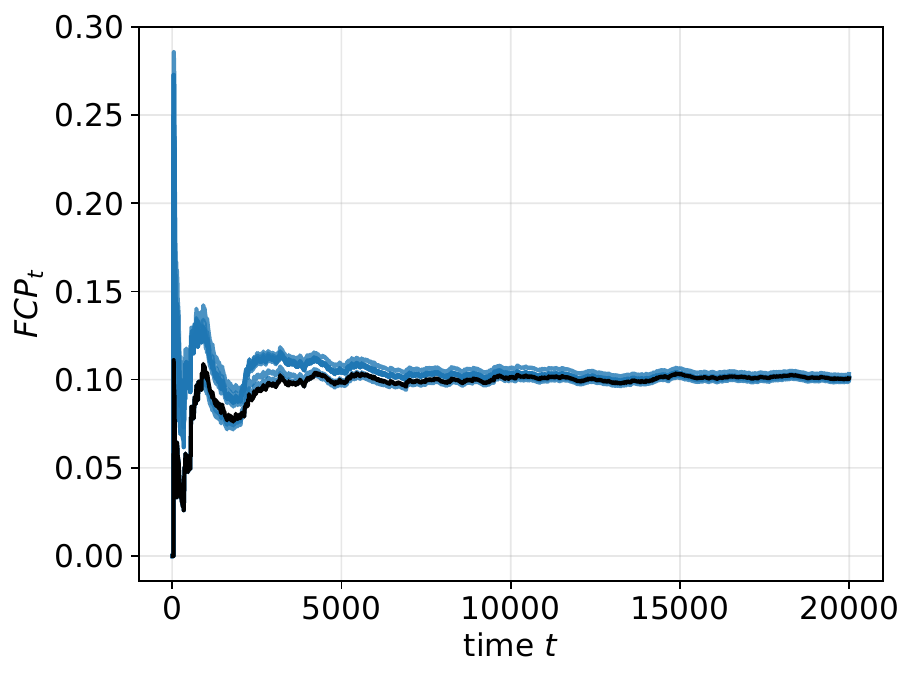}
        \includegraphics[width=.48\linewidth]{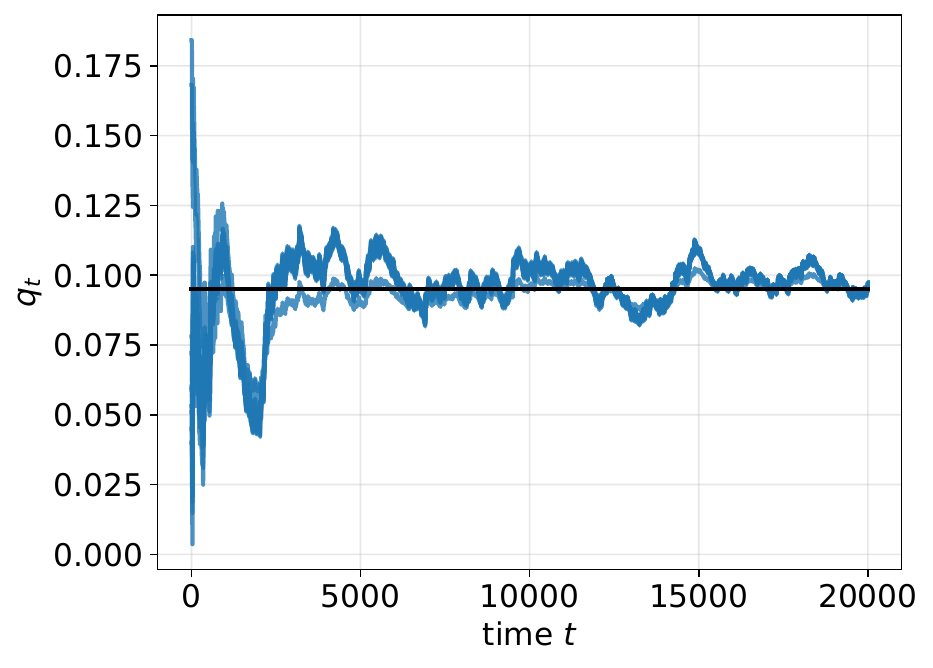} \\
        \caption{\et{False discovery proportion (left) and threshold sequence (right) for different versions of \method\: in a conformal testing setting. Trajectories obtained using $10$ different burn-in datasets for selecting the hyperparameters (blue) or the ideal hyperparameter $(q^*, 0)$ (black). }}%
        \label{fig:ND_procedure}
\end{figure}

\begin{figure}[h!]
		\centering
		\includegraphics[width=.48\linewidth]{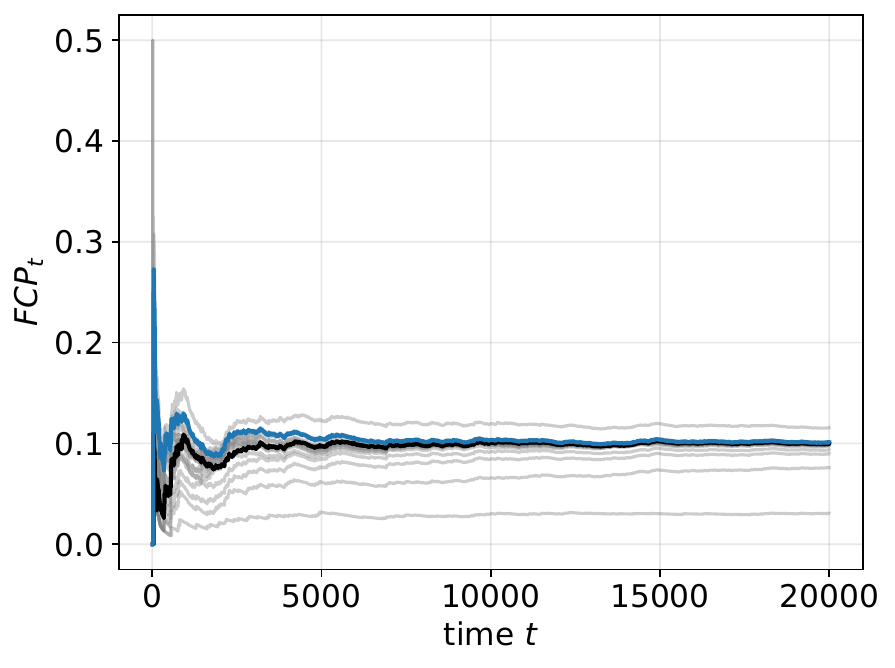}
        \includegraphics[width=.48\linewidth]{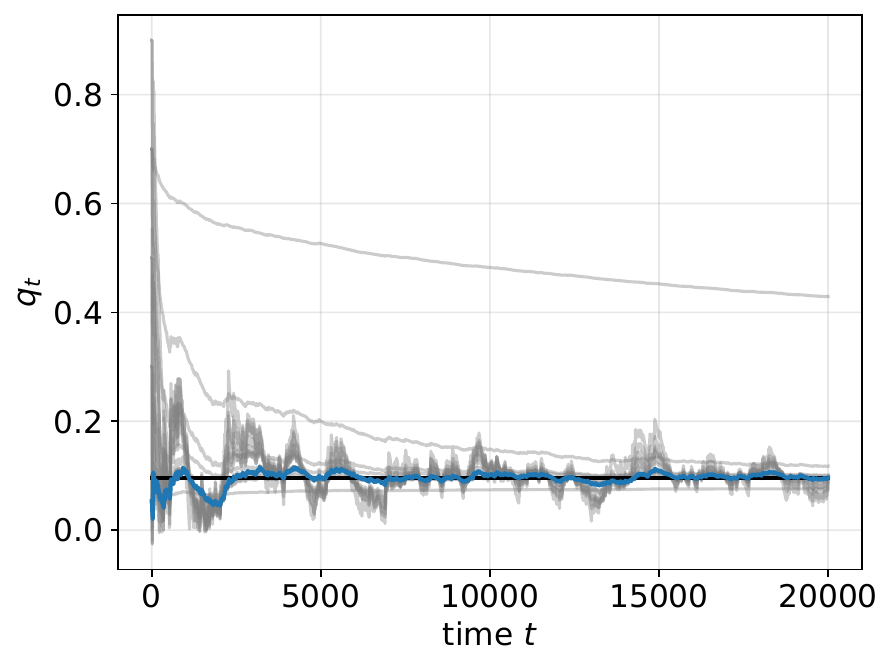} \\
        \caption{
        \et{
        False discovery proportion (left) and threshold sequence (right) for different versions of \method\: in a conformal testing setting. Trajectories obtained using our procedure for selecting the hyperparameters (blue), the ideal hyperparameter $(q^*, 0)$ (black), hyperparameters exploring a grid (grey curves).
        }}
        \label{fig:ND_procedure_compare}
\end{figure}

\paragraph*{Results}
\et{
Figure \ref{fig:ND_procedure} presents the trajectories of $\FDP_t$ and $q_t$ for the automated method (i) above (blue), across $10$ burn-in dataset $\cD_0$ generations and for the ideal method (iii) (black). 
We add in Figure \ref{fig:ND_procedure_compare} the grid methods (ii) in grey. 
This shows that the automated method produces better results than most of those obtained from the grid and closely matches the oracle. 
These results empirically validate our automatic hyperparameter tuning in this setting. Other values in the grid can be tried out via the interactive website \url{https://onlinesciapp-srjzthcznjslvzjiavpahm.streamlit.app/}.
}

\subsection{\et{Empirical results related to $\delta$ from Equation \eqref{deltaprobaJtfinite}}}\label{sec:deltaprobaJtfinite}

\pie{
We empirically evaluate the probability that our method stops, for a given set of hyperparameters. This probability is as a finite-horizon proxy for $\delta$ in Equation~\eqref{deltaprobaJtfinite}. To do so, we place ourselves in the online conformal testing version of \method, assuming for simplicity that the predictor (local-fdr) is the oracle (true) one. For the following grid of size $10 \times 10 \times 10$,
$q_0 \in \{0.01, 0.09, 0.18, \dots, 0.8\}$, $c \in \{1, 1.44, 1.89, \dots, 5\}$, and $\beta \in \{0.5, 0.54, 0.59, \dots, 0.9\},$ we estimate the probability that the method stops before $T = 1000$ by generating $1000$ datasets for each combination of hyperparameters.
}

\begin{figure}[h!]
		\centering
		\includegraphics[width=1.\linewidth]{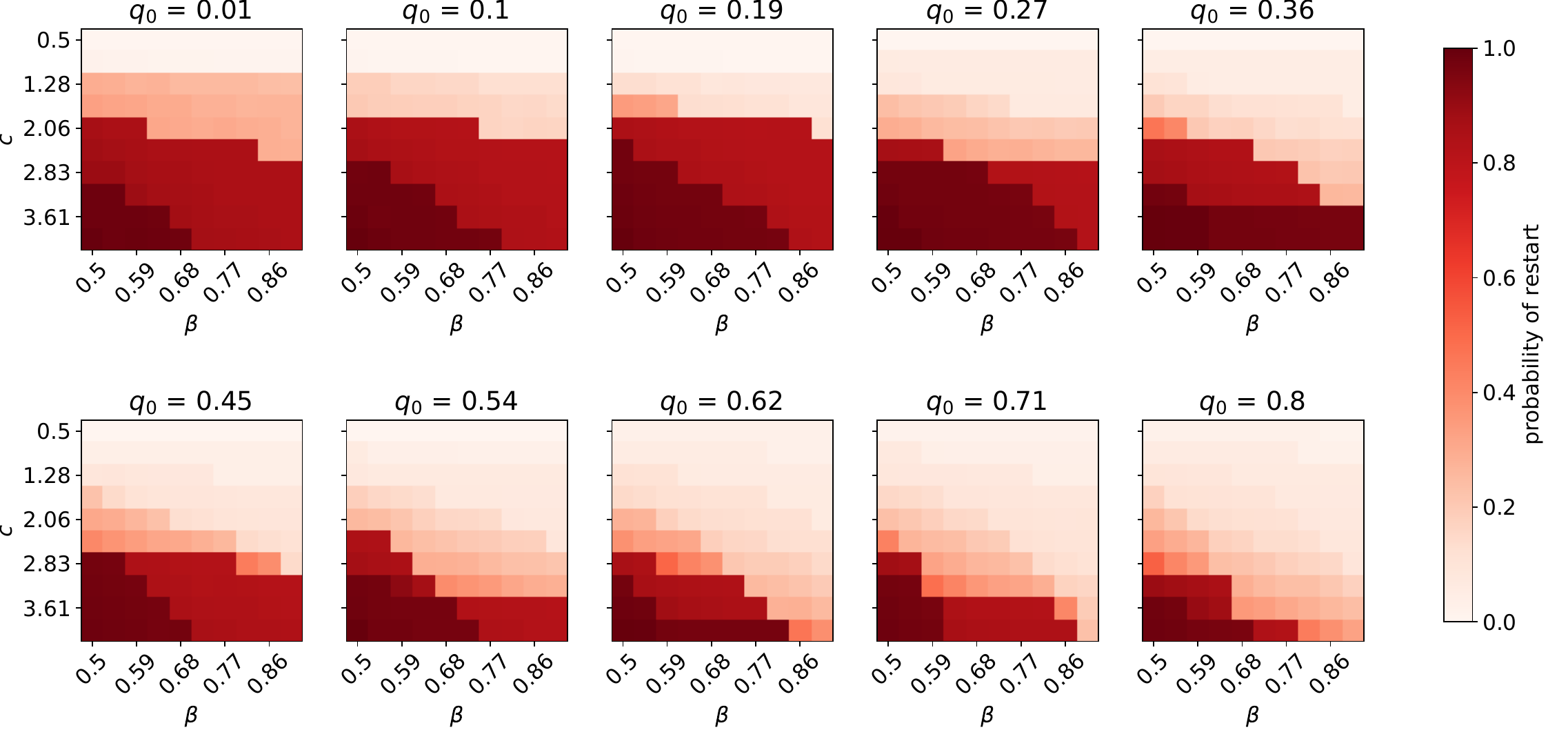}
        \caption{Heatmaps of the probability to restart in the online testing problem for different values of $q_0$, $c$, and $\beta$.}
        \label{fig:delta_estim}
\end{figure}

\pie{
Figure~\ref{fig:delta_estim} displays the heatmaps associated with this experiment. For a given $q_0$, we observe that the probability of restarting increases as $c$ grows and $\beta$ decreases. This is expected, since $\gamma_t = c \cdot t^{-\beta}$, meaning that the first steps are very large, making it plausible that the initial sequence of $q_t$ values exceeds $1$ (the restarting threshold in the online testing problem) if we make an error. Interestingly, the probability of restarting globally decreases as $q_0$ grows. One possible explanation is that, in the online testing problem, the selection criterion is $\ind{1 - \lfdr(x) > q}$. Hence, for a large $q_t$, selection only occurs when the prediction is very certain, resulting in no errors during the early steps. Since the sequence $\gamma_t$ decays rapidly, this results in a lower chance of restarting (see Figure~\ref{fig:delta_estim}, bottom-right panel, for instance).
}

\section{General approximation results}\label{sec:genapprox}\label{sec:genapproxbounds}

First, we relax Assumption~\ref{ass:easy}.

\begin{assumption}\label{ass:crossing}
The benchmark function $\Pi^\bench : q \in [0,B] \mapsto \Pi^\bench(q)\in [0,1]$ is a deterministic function such that there exists a value $q^\bench_{\alpha}\in (0,B)$ with $\Pi^\bench(q^\bench_{\alpha})=\alpha$ and for all $q\in [0,q^\bench_{\alpha})$, $\Pi^\bench(q)>\alpha$, for all $q\in (q^\bench_{\alpha},B]$, $\Pi^\bench(q)<\alpha$. 
In addition, there exists a compact $[\ell_\alpha,u_\alpha]\subset (0,B)$ such that $\ell_\alpha<q^\bench_{\alpha}<u_\alpha$, $\inf_{q\in (0,\ell_\alpha)}\Pi^\bench(q)>\alpha$ and $\sup_{q\in (u_\alpha,B)}\Pi^\bench(q)<\alpha$.
\end{assumption}

Typically, the first statement of Assumption~\ref{ass:crossing} is satisfied if $\Pi^\bench$ is continuous decreasing on $[\ell_\alpha,u_\alpha]$ and $\Pi^\bench(q^\bench_{\alpha})=\alpha$.
The last statement of Assumption~\ref{ass:crossing} is satisfied if $\Pi^\bench$ is either non-increasing on $[0,1]$, or continuous over $(0,1)$ with $\lim_{0^+} \Pi^\bench>\alpha$ and $\lim_{B^-} \Pi^\bench<\alpha$. In these cases, we can choose $\ell_\alpha=q^\bench_{\alpha}/2$ and $u_\alpha=(q^\bench_{\alpha}+B)/2$. 
Finally, the slope of $\Pi^\bench$ on the vicinity of $q^\bench_{\alpha}$ is a useful quantity to obtain convergence rate, and we let
\begin{align}
a&:=\inf\Big\{\frac{|\Pi^\bench(q)-\alpha|}{|q-q^{\bench}_{\alpha}|}\::\: q\in [\ell_\alpha,u_\alpha], \Pi^\bench(q)\neq \alpha\Big\} \geq 0\label{equ:a};\\
b&:=\sup\Big\{\frac{|\Pi^\bench(q)-\alpha|}{|q-q^{\bench}_{\alpha}|}\::\: q\in [\ell_\alpha,u_\alpha], q\neq q^{\bench}_{\alpha}\Big\} \leq +\infty\label{equ:b}.
\end{align}
 Then, the minimum/maximum slope of $\Pi^\bench$ with respect to $q^{\bench}_{\alpha}$ can be quantified by introducing
\begin{align}
 \tilde{a}&:=a \wedge \frac{\inf_{q\in (0,\ell_\alpha)}\Pi^\bench(q)-\alpha}{B}\wedge \frac{\alpha-\sup_{q\in (u_\alpha,B)}\Pi^\bench(q)}{B} \nonumber\\
 &\:\:\:\:\:\:\wedge \paren{\frac{1-\alpha}{B+\alpha \sup_{j\geq 1}\gamma_j}}\wedge \paren{\frac{\alpha}{B+(1-\alpha)\sup_{j\geq 1}\gamma_j }} \label{equ:atildewithoutt}
;\\
 \tilde{b}&:=b \vee \frac{1-\alpha}{q^{\bench}_{\alpha}-\ell_\alpha}\vee \frac{\alpha}{u_\alpha-q^{\bench}_{\alpha}}.\label{equ:btildewithoutt}
\end{align}

To derive rates we also consider the following mild assumption:
\begin{assumption}\label{ass:slope}
\mbox{The constant $a$ \eqref{equ:a} is positive and the constant $b$ \eqref{equ:b} is finite.}
\end{assumption}

The following result is obvious.
\begin{lemma}
    Assumptions~\ref{ass:crossing}~and~\ref{ass:slope} are both satisfied under Assumption~\ref{ass:easy}.
\end{lemma}

Our first nonasymptotic bound is as follows. 

\begin{theorem}\label{th:almostsureconv:withoutt}[In-probability bound]
Consider the \method\  threshold sequence $(q_t^{\OSCI})_{t\geq 1}$ \eqref{equACIquantile} with any nonincreasing step size sequence $(\gamma_j)_{j \geq 1}$ and starting point $q_1^{\OSCI}\in [0,B]$. 
Consider some function $\Pi^{\bench}$ satisfying Assumption~\ref{ass:crossing} with some $q^\bench_{\alpha}\in (0,B)$ and the corresponding $D_t$ being defined by \eqref{equDtT}. Recall that $J(t)=\sum_{k=1}^{t-1} S_{k}(X_{k}, \qt_{k}^{\OSCI}) + 1$ is given by \eqref{def:Jt}.
Let for all $T>t\geq 0$,
\begin{align}
\label{equmart2}
\epsilon_{t,T} &:=  \sum_{k=t+1}^{T}\gamma_{J(k)}  (\Err_k(Y_k, \cC_k(X_k,q_k^{\OSCI}))-  \Pi_{\cC_k,S_k}(q_k^{\OSCI}))  \cdot \ind{\qt_k^{\OSCI} \geq 0}  \cdot S_k(X_k, \qt_k^{\OSCI}) .
\end{align}
Then 
for all $T>t\geq 1$, we have pointwise
\begin{align}
&\abs{\qt_T^{\OSCI}-q^{\bench}_{\alpha}}\label{eq:QuantileBoundsLearningwithoutt}\\\leq \:&
 (\max_{t\leq k\leq T} D_k)/\tilde{a} +\max_{t\leq k\leq T-1} | \epsilon_{k-1,T-1} |  + \gamma_{J(t)} + \paren{\tilde{a}\sum_{j=J(t)}^{J(T+1)-1} \gamma_{j}}^{-1} \paren{\abs{\epsilon_{t-1,T}} +B+\gamma_1} ,\nonumber
\end{align}
with the convention $1/0=+\infty$ if $\sum_{j=J(t)}^{J(T+1)-1} \gamma_{j}=0$ or $\tilde{a}=0$. 
In addition, the process $(\epsilon_{t,T})_{T\geq t}$ vanishes in the following sense: for all $x,y>0$, $T\geq t$,
\begin{align}\label{epssmall}
\P\Big(\epsilon_{t,T}>x \Big) \leq \P\bigg(\sum_{j=J(t)}^{J(T+1)-1} \gamma_{j}^2 > y\bigg)+  e^{-2x^2/y}.
\end{align}
\end{theorem}

The proof is given in \S~\ref{proof:th:almostsureconvwithoutt}, in a slightly more general context where $\Pi^\bench$ is allowed to depend on $t$. 

Let us provide some intuition for the bound given in \eqref{eq:QuantileBoundsLearningwithoutt}. 
First, while the bound is on the term $\abs{\qt_T^{\OSCI}-q^{\bench}_{\alpha}}$, it is valid for any $t<T$ and the quality of the bound depend on the choice of $t$. For instance, we can take $t=t(T)=T/2$, but other choices can lead to better bound (see, e.g., Theorem~\ref{th:rateconvfull} case (ii)). Second, the process $\epsilon_{t-1,T}$ is guaranteed to be small by \eqref{epssmall}, when $\sum_{j=J(t)}^{J(T+1)-1} \gamma_{j}^2$ is not too large, which is true if we choose a step size sequence such that $\sum_{j\geq 1} \gamma_{j}^2 <\infty$ and if the procedure makes enough selections when $t,T$ are large. 
Third, $\sum_{j=J(t)}^{J(T+1)-1} \gamma_{j}$ should be large, which is true if $\sum_{j\geq 1}  \gamma_{j} = \infty$ and $J(T+1)$ is large. Hence, it is natural to choose a step size sequence with \eqref{equ-stepseq}.
Fourth,  $D_k$ should be small in order to obtain a small bound.

Our second nonasymptotic bound is as follows. 
\begin{theorem}\label{th:genL2withoutt}[In-expectation bound]
Consider the \method\  threshold sequence $(q_t^{\OSCI})_{t\geq 1}$  \eqref{equACIquantile} with any starting point $q_1^{\OSCI}\in [0,B]$ and a nonincreasing step size sequence $(\gamma_j)_{j \geq 1}$. Consider some function $\Pi^{\bench}$ satisfying Assumption~\ref{ass:crossing} with some $q^\bench_{\alpha}\in (0,B)$ and the corresponding $D_t$ being defined by \eqref{equDtT}. Consider  $\tilde{a}$ given by \eqref{equ:atildewithoutt} and $J(\cdot)$  given by \eqref{def:Jt}.
Assume that  for some $t\geq j_0\geq 1$, 
\begin{itemize}
    \item[(i)] $1-2\tilde{a} \gamma_{j_0}\geq 0$ and  $\gamma_j/\gamma_{j+1}\leq 1/(1-\tilde{a}\gamma_j)$ for $j_0\leq j\leq t$;
    \item[(ii)] $S_k$ is independent of the past $\mathcal{F}_{k-1}$, for $j_0\leq k\leq t$.
\end{itemize} Then, for any $t_0$ with $j_0\leq t_0 \leq t$, on the event $J(t_0)\geq j_0$, we have pointwise
    \begin{align}
&\e{\paren{q^{\OSCI}_{t}-q^{\bench}_{\alpha}}^2\mid S_{1}\dots S_{t-1} }\nonumber\\
&\leq 
\frac{\gamma_{J(t)}}{\gamma_{J(t_0)}} \paren{(B+\sup_{j\geq 1}\gamma_{j})^2 \vee \frac{4}{\tilde{a}}} + 4(B +\sup_{j\geq 1}\gamma_{j})\gamma_{J(t)}\sum_{k=t_0}^{t-1} S_k \e{D_k\mid S_{1}\dots S_{k-1}},\label{equ:genL2withoutt}
\end{align}
where we denoted $S_t=S_t(X_t,q^{\OSCI}_{t})$ for all $t\geq 1$.
\end{theorem}

Theorem~\ref{th:genL2withoutt} is proved in \S~\ref{proofth:genL2withoutt}. 
This bound is simpler than the one obtained in Theorem~\ref{th:almostsureconv:withoutt}, but relies on a stronger independence assumption (ii). Assumption (ii) is satisfied typically when the $(X_t,Y_t)$ are independent over time and the selection rule only involve $X_t$ ($X$-oriented selection). To make the bound in \eqref{equ:genL2withoutt} small, the sequence $\gamma_j$ should vanishes while $J(t)$ should grows to infinity when $t$ grows. The other remainder term is small when $D_t$ is small enough. Finally, the assumption (i) is relatively light, and is satisfied for the sequence $\gamma_j=c j^{-\beta}$ for $\beta\in (0,1)$. Note that the two conditions in \eqref{equ-stepseq} are not required here.

Finally, bounding $\abs{\qt_t^{\OSCI}-q^{\bench}_{\alpha}}$ directly leads to a bound on the $\IER$.

\begin{lemma}\label{lem:boundIER}
Consider the \method\  threshold sequence $(q_t^{\OSCI})_{t\geq 1}$  \eqref{equACIquantile} with any nonincreasing step size sequence $(\gamma_j)_{j \geq 1}$ and starting point $q_1^{\OSCI}\in [0,B]$. 
Consider some function $\Pi^{\bench}$ satisfying Assumption~\ref{ass:crossing} with some $q^\bench_{t,\alpha}\in (0,B)$ and the corresponding $D_t$ being defined by \eqref{equDtT}.
Considering in addition $\tilde{b}$ as in \eqref{equ:btildewithoutt}. Then we have for all $t\geq 1$, provided that $q_t^\OSCI\in [0,B]$,
\begin{equation}\label{equ:boundIERgen}
    |\IER_t(\mathcal{R}^{\OSCI}) - \alpha| \leq  D_t + \tilde{b} |q_t^{\OSCI} -q_\alpha^\bench|.
\end{equation}
\end{lemma}

Lemma~\ref{lem:boundIER} is proved in \S~\ref{sec:proofboundIER}. It is typically combined with the bounds above to obtain bounds on $\abs{\IER_t(\mathcal{R}^{\OSCI})-\alpha}$.

\section{Proofs for the general bounds in \S~\ref{sec:genapproxbounds}}

\subsection{Proof of Theorem~\ref{th:almostsureconv:withoutt}}\label{proof:th:almostsureconvwithoutt}

We prove Theorem~\ref{th:almostsureconv:withoutt} by proving a slightly more general result, for which $\Pi^\bench_t$ may depend on time.

\subsubsection{Extending Theorem~\ref{th:almostsureconv:withoutt}}

We extend Theorem~\ref{th:almostsureconv:withoutt} to a more general situation where the benchmark function $\Pi^\bench_t$ possibly depends on the time $t$. Hence, $q^\bench_\alpha$, $a$ and $\tilde{a}$ now depend on $t$, and $D_t$ reads $\abs{\Pi_t(q_t^{\OSCI})-\Pi^\bench_t(q_t^{\OSCI})}$. 
For simplicity, and with some abuse of notation, the numbers of the equations are unchanged when we refer to that quantities. For instance, we still refer to \eqref{equDtT} for $D_t$, even if $\Pi^\bench_t$ depends on $t$.

\begin{theorem}\label{th:almostsureconv}
Consider the \method\  threshold sequence $(q_t^{\OSCI})_{t\geq 1}$ \eqref{equACIquantile} with any step size sequence $(\gamma_j)_{j \geq 1}$ and starting point $q_1^{\OSCI}\in [0,B]$. For each $t\geq 1$, consider some function $\Pi^{\bench}_t$ satisfying Assumption~\ref{ass:crossing} with some $q^\bench_{t,\alpha}\in (0,B)$ and the corresponding $D_t$ in \eqref{equDtT}. Let $a_t\geq 0$ and $\tilde{a}_t\geq 0$ being defined by \eqref{equ:a} and \eqref{equ:atildewithoutt}, respectively, in which $\Pi^\bench(q)$ and $q^{\bench}_{\alpha}$ have been replaced by $\Pi^\bench_t(q)$ and $q^{\bench}_{t,\alpha}$, respectively.
Then for all $T>t\geq 1$, we have pointwise
\begin{align}
\abs{\qt_T^{\OSCI}-q^{\bench}_{T,\alpha}}&\leq 
 \max_{t\leq k\leq T} (D_k/\tilde{a}_k) +\max_{t\leq k\leq T-1} | \epsilon_{k-1,T-1} | +\max_{t\leq k\leq T-1}|q^\bench_{T,\alpha} -q^\bench_{k,\alpha}| + \max_{J(t)\leq j\leq J(T-1)}\gamma_{j} \nonumber\\
&+
\frac{1}{\min_{t\leq k\leq T} \tilde{a}_k}\brac{ \left(\sum_{j=J(t)}^{J(T+1)-1} \gamma_{j}\right)^{-1} \paren{\abs{\epsilon_{t-1,T}} +B+\sup_{j\geq 1}\gamma_j} },\label{eq:QuantileBoundsLearning}
\end{align}
with the convention $1/0=+\infty$ if $\sum_{j=J(t)}^{J(T+1)-1} \gamma_{j}=0$ and the process defined by $(\epsilon_{t,T})_{T\geq t}$ \eqref{equmart2} satisfying \eqref{epssmall}.
\end{theorem}

Theorem~\ref{th:almostsureconv} is proved in \S~\ref{proof:th:almostsureconv}.
Compared to the original Theorem~\ref{th:almostsureconv:withoutt}, the bound in Theorem~\ref{th:almostsureconv} essentially contains the additional term $\max_{t\leq k\leq T-1}|q^\bench_{T,\alpha} -q^\bench_{k,\alpha}|$ measuring how $\Pi^\bench_k$ (hence $q^\bench_{k,\alpha}$) moves when $k$ varies from $t$ to $T-1$. It is zero when  $\Pi^\bench_k$ does not depend on $k$, and we recover Theorem~\ref{th:almostsureconv:withoutt}.
Another difference is that the sequence $(\gamma_{j})_{j\geq 1}$ is non necessarily assumed to be nonincreasing for the sake of generality.

\subsubsection{Proof of Theorem~\ref{th:almostsureconv}}\label{proof:th:almostsureconv}

Let us denote $q_t^{\OSCI}$ by $q_t$ for short. 
First consider for a given $t\geq 1$ the functions
\begin{align}
\widetilde{\Pi}_t(q)&:=\E[\reaction_t(Y_t, \cC_t(X_t,q)) \:|\:  S_t(X_t, q)=1,  \mathcal{F}_{t-1}] = \Pi_t(q) \ind{\qt \in [0,B]} + \ind{\qt < 0} , \:\:q\in \R;\label{equ:pitilde}\\
    \widetilde{\Pi}_t^{\bench}(q)&:=\Pi^\bench_t(q)\ind{q\in [0,B]}+\ind{q<0}, \:\:q\in \R,\label{equ:pizerotilde}
\end{align}
where $\Pi_t(q) $ is given by \eqref{PiCtSt}.
Note that $D_t=\abs{\Pi_{t}(q_t^{\OSCI})-\Pi^\bench_t(q_t^{\OSCI})}\ind{q^\OSCI_t\in [0,B]}$ given by \eqref{equDtT} is also such that
\begin{equation}
    \label{equDtTtilde}
D_t=\abs{\widetilde{\Pi}_t(q_t^{\OSCI})-\widetilde{\Pi}^\bench_t(q_t^{\OSCI})},\:\:t\geq 1.
\end{equation}
Also, by Assumption~\ref{ass:crossing} for $\Pi^\bench_t$ and with $\tilde{a}_t$ corresponding to definition \eqref{equ:atildewithoutt}, we have 
\begin{equation}\label{equ:usefulslope}
\inf_{t\geq 1} \frac{\abs{\widetilde{\Pi}_t^{\bench}(q_t)-\alpha}}{\abs{q_t-q^{\bench}_{t,\alpha}}}\geq \inf\set{\frac{\abs{\widetilde{\Pi}_t^{\bench}(q)-\alpha}}{\abs{q-q^{\bench}_{t,\alpha}}}, -\alpha \sup_{j\geq 1}\gamma_{j} \leq q\leq  B+(1-\alpha)\sup_{j\geq 1}\gamma_{j} }
 \geq \tilde{a}_t.
\end{equation}
Lastly, our proof rely on the following quantities, for $\epsilon\in (0,\alpha\wedge (1-\alpha))$,
\begin{align}
    q^\bench_{t,\alpha+\epsilon}&:=\inf\set{q\in \R\::\: \widetilde{\Pi}^\bench_{t}(q)\leq \alpha+\epsilon}\label{qzeroPluseps};\\
    q^\bench_{t,\alpha-\epsilon}&:=\sup\set{q\in \R\::\: \widetilde{\Pi}^\bench_{t}(q)\geq \alpha-\epsilon}\label{qzeroMinuseps}.
\end{align}
Note that by Assumption~\ref{ass:crossing} for $\Pi^\bench_t$, these quantities always exist and we have $-\alpha \sup_j \gamma_j /2\leq q^\bench_{t,\alpha+\epsilon}\leq q^\bench_{t,\alpha}$ and $B\geq q^\bench_{t,\alpha-\epsilon}\geq q^\bench_{t,\alpha}$. In addition, we also have the relations:
\begin{align}
     q^\bench_{t,\alpha-\epsilon}-q^\bench_{t,\alpha}&\leq \epsilon/\tilde{a}_t\label{qzeroMinuseps-a};\\
   q^\bench_{t,\alpha} - q^\bench_{t,\alpha+\epsilon}&\leq \epsilon/\tilde{a}_t\label{qzeroPluseps-a}.
\end{align}
Indeed, let us establish \eqref{qzeroPluseps-a} (the proof for \eqref{qzeroMinuseps-a} is similar) by assuming $q^\bench_{t,\alpha} > q^\bench_{t,\alpha+\epsilon}$ (otherwise there is nothing to prove). By \eqref{qzeroPluseps}, there exists a sequence $(\eta_s)_{s\geq 1}>0$ tending to zero, such that for large $s$,  ${\Pi}^\bench_{t}(q^\bench_{t,\alpha+\epsilon}+\eta_s)\leq \alpha+\epsilon$. Hence, by \eqref{equ:usefulslope}, we have
$$
q^\bench_{t,\alpha}-\paren{q^\bench_{t,\alpha+\epsilon}+\eta_s}  = \frac{q^\bench_{t,\alpha}-\paren{q^\bench_{t,\alpha+\epsilon}+\eta_s}}{{\Pi}^\bench_{t}(q^\bench_{t,\alpha+\epsilon}+\eta_s)-\alpha} \paren{{\Pi}^\bench_{t}(q^\bench_{t,\alpha+\epsilon}+\eta_s)-\alpha}\leq \epsilon/\tilde{a}_t .
$$
Hence we obtain \eqref{qzeroPluseps-a} by taking $s$ tending to infinity.\\

Now that these notation have been introduced, we prove Theorem~\ref{th:almostsureconv}. We obtain by \eqref{equACIquantile} that for $T\geq 1$,
\begin{align}
q_{T+1} - q_1&= \sum_{t=1}^{T} (q_{t+1}-q_t)= \sum_{t=1}^{T}\gamma_{J(t)} \paren{\reaction_{t}(Y_t, \cC_t(X_t, \qt_t),q_t)-\alpha} \cdot S_t(X_t, \qt_t)\label{equalphatinfo}.
\end{align}
Considering the process
\begin{align}
Z_{T}&:= \sum_{t=1}^{T}\gamma_{J(t)}  (\reaction_t(Y_t, \cC_t(X_t,q_t),q_t)- \widetilde{\Pi}_t(q_t)) \cdot S_t(X_t, \qt_t) ,\:\:\:T\geq 1,\label{equmartinfo}
\end{align}
we  have 
$
q_{T+1} - q_1
= Z_{T} +   \sum_{t=1}^{T}\gamma_{J(t)} (\widetilde{\Pi}_t(q_t)-\alpha) \cdot S_t(X_t, \qt_t)$.
This entails, for $1\leq t<T$, 
\begin{align*}
q_{T+1} - q_{t}
= Z_{T} - Z_{t-1} 
&+  \sum_{k=t}^{T}\gamma_{J(k)} (\widetilde{\Pi}_k(q_k)-\alpha) \cdot S_k(X_k, \qt_k)
\end{align*}

As a result, by \eqref{equmart2}, $Z_{T-1} - Z_{t-1}=\epsilon_{t-1,T-1}$ and we have for any $1\leq t<T$,
\begin{align}\label{equmart3}
\qt_T-\qt_t &=\epsilon_{t-1,T-1}   + \sum_{k=t}^{T-1}\gamma_{J(k)}\paren{\widetilde{\Pi}_k(\qt_k) -\alpha}\cdot S_k(X_k,\qt_k).
\end{align}

Now fix any $1\leq t<T$ and consider the set of ``crossing-point times''
$$
\mathcal{C}=\{ t\leq k\leq T-1 \::\: \widetilde{\Pi}_k(\qt_k)\leq \alpha<\widetilde{\Pi}_{k+1}(\qt_{k+1}) \mbox{ or }\widetilde{\Pi}_k(\qt_k)> \alpha\geq \widetilde{\Pi}_{k+1}(\qt_{k+1})\}
$$

\paragraph*{Case 1:  $\mathcal{C}$ is not empty}
In this case, we may let $t^*= \max \mathcal{C}$.  
By definition, we have either of the two following cases:
\begin{itemize}
\item[(i)] $\widetilde{\Pi}_{t^*}(\qt_{t^*})\leq \alpha \mbox{ and } \forall s\in \{t^*+1,\dots,T\}, \widetilde{\Pi}_s(\qt_{s})> \alpha$;
\item[(ii)] $\widetilde{\Pi}_{t^*}(\qt_{t^*})>\alpha \mbox{ and } \forall s\in \{t^*+1,\dots,T\}, \widetilde{\Pi}_s(\qt_{s})\leq \alpha$.
\end{itemize}
In case (i), we have by \eqref{equmart3} that 
$\qt_T-\qt_{t^*} \geq \epsilon_{t^*-1,T-1} +\gamma_{J(t^*)}(\widetilde{\Pi}_{t^*}(\qt_{t^*}) -\alpha)S_{t^*}(X_{t^*},\qt_{t^*})$, that is, 
$\gamma_{J(t^*)}( \alpha-\widetilde{\Pi}_{t^*}(\qt_{t^*}) )\geq  \gamma_{J(t^*)}( \alpha-\Pi^\bench_{t^*}(\qt_{t^*}) )S_{t^*}(X_{t^*},\qt_{t^*})\geq \epsilon_{t^*-1,T-1}  -(\qt_T-\qt_{t^*})$ and since 
$\widetilde{\Pi}_T(\qt_{T}) >\alpha\geq \widetilde{\Pi}_{t^*}(\qt_{t^*})$ and by \eqref{equDtTtilde} for $D_T$ and $D_{t^*}$, we then have ${\Pi}^{\bench}_T(\qt_{T})+D_T >\alpha\geq {\Pi}^\bench_{t^*}(\qt_{t^*})-D_{t^*}$.
By \eqref{qzeroPluseps}, this means that $\qt_{t^*}\geq q^\bench_{t^*,\alpha+D_{t^*}}$  with $q^\bench_{t^*,\alpha+D_{t^*}}\leq q^\bench_{t^*,\alpha}$.
Thus 
\begin{align*}
\qt_{t^*}-\qt_T \geq q^\bench_{t^*,\alpha+D_{t^*}}-\qt_T&=q^\bench_{t^*,\alpha+D_{t^*}}-q^\bench_{t^*,\alpha} + q^\bench_{t^*,\alpha}-q^\bench_{T,\alpha} + q^\bench_{T,\alpha}-\qt_T \\
&\geq q^\bench_{T,\alpha}-\qt_T -\paren{q^\bench_{t^*,\alpha}-q^\bench_{t^*,\alpha+D_{t^*}}} - \abs{q^\bench_{t^*,\alpha}-q^\bench_{T,\alpha}}.
\end{align*}
Similarly, since by \eqref{qzeroMinuseps}, $q_T\leq q^\bench_{T,\alpha-D_{T}}$  with $q^\bench_{T,\alpha-D_{T}}\geq q^\bench_{T,\alpha}$, we have
$
\qt_T-q^\bench_{T,\alpha} \leq q^\bench_{T,\alpha-D_{T}}-q^\bench_{T,\alpha}.
$
Combining the above bounds leads to 
\begin{align*}
|\qt_T-q^\bench_{T,\alpha}| \leq \gamma_{J(t^*)} - \epsilon_{t^*-1,T-1} + \abs{q^\bench_{t^*,\alpha}-q^\bench_{T,\alpha}}+ \paren{ q^\bench_{t^*,\alpha}-q^\bench_{t^*,\alpha+D_{t^*}}} \vee \paren{q^\bench_{T,\alpha-D_{T}}-q^\bench_{T,\alpha}}.
\end{align*}
Now, by \eqref{qzeroMinuseps-a} and \eqref{qzeroPluseps-a}, we have $q^\bench_{T,\alpha-D_{T}}-q^\bench_{T,\alpha} \leq D_T/\tilde{a}_T$ and $q^\bench_{t^*,\alpha}-q^\bench_{t^*,\alpha+D_{t^*}}\leq D_{t^*}/\tilde{a}_{t^*}$ respectively. 
This provides 
\begin{align}\label{firstcaseproof2}
|\qt_T-q^\bench_{T,\alpha}| \leq  \max_{J(t)\leq j\leq J(T-1)}\gamma_{j}+\max_{t\leq k\leq T-1} | \epsilon_{k-1,T-1} | + \max_{t\leq k\leq T-1}|q^\bench_{T,\alpha} -q^\bench_{k,\alpha}|+ \max_{t\leq k\leq T} \paren{D_{k}/\tilde{a}_k} .
\end{align}
In case (ii), by \eqref{equmart3}, we have similarly
$\qt_T-\qt_{t^*} \leq \epsilon_{t^*-1,T-1}   +\gamma_{J(t^*)}( \widetilde{\Pi}_{t^*}(\qt_{t^*})-\alpha)S_{t^*}(X_{t^*},\qt_{t^*})\leq \epsilon_{t^*-1,T-1}  +\gamma_{J(t^*)}( \widetilde{\Pi}_{t^*}(\qt_{t^*})-\alpha)$. Also, by  \eqref{equDtTtilde}, and since $\widetilde{\Pi}_T(\qt_{T})-D_T \leq \alpha <\widetilde{\Pi}_{t^*}(\qt_{t^*})+D_{t^*}$ in this case, we have both 
$q_T\geq q^\bench_{T,\alpha+D_{T}}$ and $q_{t^*}\leq q^\bench_{t^*,\alpha-D_{t^*}}$ and thus
$$\qt_T-\qt_{t^*}\geq \qt_T-q^\bench_{T,\alpha}+q^\bench_{T,\alpha}-q^\bench_{t^*,\alpha}+ q^\bench_{t^*,\alpha} -q^\bench_{t^*,\alpha-D_{t^*}}\geq  \qt_T-q^\bench_{T,\alpha} -\abs{q^\bench_{T,\alpha}-q^\bench_{t^*,\alpha}}-\paren{q^\bench_{t^*,\alpha-D_{t^*}}-q^\bench_{t^*,\alpha}}.$$
 In particular, $\qt_T-q^\bench_{T,\alpha} \leq \epsilon_{t^*-1,T-1} + \abs{q^\bench_{T,\alpha}-q^\bench_{t^*,\alpha}} + \paren{q^\bench_{t^*,\alpha-D_{t^*}}-q^\bench_{t^*,\alpha}}+\gamma_{J(t^*)}$.
 In addition, since $q^\bench_{T,\alpha}-\qt_T\leq q^\bench_{T,\alpha}-q^\bench_{T,\alpha+D_{T}}$, this case also gives \eqref{firstcaseproof2} by using \eqref{qzeroMinuseps-a} and \eqref{qzeroPluseps-a}.
 
\paragraph*{Case 2:  $\mathcal{C}$ is empty.}
In this case, we have either of the two following cases:
\begin{itemize}
\item[(i)] $\widetilde{\Pi}_k(\qt_k)\leq \alpha$ for all $t\leq k\leq T$;
\item[(ii)] $\widetilde{\Pi}_k(\qt_k)\geq \alpha$ for all $t\leq k\leq T$.
\end{itemize}
In case (i), we have by \eqref{equmart3} that 
$
\qt_{T+1}-\qt_t \leq \epsilon_{t-1,T} - \min_{t\leq k\leq T}|\widetilde{\Pi}_k(\qt_k)-\alpha|\sum_{k=t}^{T}\gamma_{J(k)}S_k(X_k,\qt_k) .
$
In addition $|\qt_{T+1}-\qt_t| \leq B+\max_{1 \leq j \leq J(T+1)}(\gamma_{j})$ by Lemma~\ref{lem:boundv}, hence we have
$$
\min_{t\leq k\leq T}|\widetilde{\Pi}_k(\qt_k)-\alpha| 
\leq 
  \paren{\sum_{k=t}^{T}\gamma_{J(k)}S_k(X_k,\qt_k)}^{-1} \paren{|\epsilon_{t-1,T}| +B+\max_{1 \leq j \leq J(T+1)}\gamma_{j}},
$$
 with the convention $1/0=+\infty$ if $\sum_{k=t}^{T}\gamma_{J(k)}S_k(X_k,\qt_k)=0$. 
As a result, by \eqref{equDtTtilde} of $D_k$, and since $\abs{\Pi^\bench_k(\qt_k)-\alpha}\geq \tilde{a}_k \abs{\qt_k-q^{\bench}_{k,\alpha}}$ by \eqref{equ:a}, we have
\begin{align}\label{equinterm1}
&\min_{t\leq k\leq T}\abs{\qt_k-q^{\bench}_{k,\alpha}}\\
&\leq \max_{t\leq k\leq T}(D_k/\tilde{a}_k) +
\max_{t\leq k\leq T}(1/\tilde{a}_k)  \paren{\sum_{k=t}^{T}\gamma_{J(k)}S_k(X_k,\qt_k)}^{-1} \paren{|\epsilon_{t-1,T}| +B+\max_{1 \leq j \leq J(T+1)}\gamma_{j}}.\nonumber
 \end{align}
Now, by \eqref{equmart3}, we have for all $k\in [t, T]$ that $\qt_T-\qt_k\leq \epsilon_{k-1,T-1}$, hence
\begin{align*}
    \qt_{T}-q^{\bench}_{T,\alpha} &= q_T - q_k + q_k - q^{\bench}_{k,\alpha} +q^{\bench}_{k,\alpha} -q^{\bench}_{T,\alpha}\\
    &\leq \epsilon_{k-1,T-1} + |q^{\bench}_{T,\alpha}-q^{\bench}_{k,\alpha}| + |\qt_k-q^{\bench}_{k,\alpha}|.
\end{align*}
In addition, since $\widetilde{\Pi}_T(\qt_T)\leq \alpha$, we have ${\Pi}^{\bench}_T(\qt_T)\leq \alpha+D_T$ and $\qt_{T}\geq q^\bench_{T,\alpha+D_{T}}$, which gives $q^{\bench}_{T,\alpha}-\qt_{T}\leq q^{\bench}_{T,\alpha}-q^\bench_{T,\alpha+D_{T}}\leq D_T/\tilde{a}_T$ by \eqref{qzeroPluseps-a}.
We hence obtain
\begin{align}
|\qt_{T}-q^{\bench}_{T,\alpha}|-\max_{t\leq k\leq T-1} \abs{\epsilon_{k-1,T-1}} 
- \max_{t\leq k\leq T-1}|q^{\bench}_{T,\alpha}-q^{\bench}_{k,\alpha}| 
\leq \paren{\min_{t\leq k\leq T}|\qt_k-q^{\bench}_{k,\alpha}|} \wedge \paren{D_T/\tilde{a}_T} . \label{equinterm2}
\end{align}
Combining \eqref{equinterm1} and \eqref{equinterm2} leads to
\begin{align}
&|\qt_{T}-q^{\bench}_{T,\alpha}|-\max_{t\leq k\leq T-1} \abs{\epsilon_{k-1,T-1}} - \max_{t\leq k\leq T-1}|q^{\bench}_{T,\alpha}-q^{\bench}_{k,\alpha}|\nonumber\\
&\leq 
\max_{t\leq k\leq T}(D_k/\tilde{a}_k) +
\max_{t\leq k\leq T}(1/\tilde{a}_k)  \paren{\sum_{k=t}^{T}\gamma_{J(k)}S_k(X_k,\qt_k)}^{-1} \paren{|\epsilon_{t-1,T}| +B+\max_{1 \leq j \leq J(T+1)}\gamma_{j}}.\label{equ-truc-info}
\end{align}
In case (ii), we have similarly 
\begin{align*}
\qt_{T+1}-\qt_{t} &\geq \epsilon_{t-1,T} +\min_{t\leq k\leq T}|\alpha-\widetilde{\Pi}_k(\qt_k)|\sum_{k=t}^{T}\gamma_{J(k)}S_k(X_k,\qt_k) \\
&\geq \epsilon_{t-1,T} +\paren{\min_{t\leq k\leq T} \set{\tilde{a}_k\abs{\qt_k-q^{\bench}_{T,\alpha}}} -\max_{t\leq k\leq T}  D_k}\sum_{k=t}^{T}\gamma_{J(k)}S_k(X_k,\qt_k),
\end{align*}
which also yields \eqref{equinterm1}.
Furthermore, for all $k\in [t, T]$, 
 $\qt_T-\qt_k\geq \epsilon_{k-1,T-1}$, hence $$q^{\bench}_{T,\alpha}-\qt_T\leq \max_{t\leq k\leq T-1} |\epsilon_{k-1,T-1}|  + \min_{t\leq k\leq T}|q^{\bench}_{k,\alpha}-\qt_k|+\max_{t\leq k\leq T-1}|q^{\bench}_{T,\alpha}-q^{\bench}_{k,\alpha}|.$$
On the other hand, since $\widetilde{\Pi}_T(\qt_T)\geq \alpha$, we have ${\Pi}^{\bench}_T(\qt_T)\geq \alpha-D_T$ and $\qt_{T}\leq q^\bench_{T,\alpha-D_{T}}$, which gives $\qt_{T}-q^{\bench}_{T,\alpha}\leq q^\bench_{T,\alpha-D_{T}}-q^{\bench}_{T,\alpha}\leq D_T/\tilde{a}_T$ by \eqref{qzeroMinuseps-a}.
This also yields \eqref{equinterm2} and thus also \eqref{equ-truc-info}.
Combining the two cases above gives the bound (by Lemma~\ref{lem:computing}):
\begin{align*}
&\abs{\qt_T-q^{\bench}_{T,\alpha}}\\
&\leq \max_{t\leq k\leq T}(D_k/\tilde{a}_k) 
 +\max_{t\leq k\leq T-1} | \epsilon_{k-1,T-1} | +\max_{t\leq k\leq T-1}|q^\bench_{T,\alpha} -q^\bench_{k,\alpha}| \\
&+\big(\max_{J(t)\leq j\leq J(T-1)}\gamma_{j}\big)\vee \paren{\max_{t\leq k\leq T}\paren{1/\tilde{a}_k}\brac{ \left(\sum_{j=J(t)}^{J(T+1)-1} \gamma_{j}\right)^{-1} \paren{\abs{\epsilon_{t-1,T}} +B+\sup_{j\geq 1}\gamma_j} }}.
\end{align*}
This gives the final bound \eqref{eq:QuantileBoundsLearning}. 

Let us now prove the last statement and show that the process defined by $(\epsilon_{t,T})_{T\geq t}$ \eqref{equmart2} satisfies \eqref{epssmall}. Define the filtration: 
$$
\cG_t=\sigma(\mathcal{F}_t,S_{t+1}(X_{t+1}, \qt_{t+1})),\:\: t\geq 1.
$$
We first prove that the process $(Z_t)_{t\geq 1}$ defined in \eqref{equmartinfo} is a centered martingale with respect to the filtration $(\cG_t)_{t\geq 1}$ with $\E(Z_T^2)\leq \sum_{t=1}^{J(T)} \gamma_{t}^2$. 
First, for $t\geq 1$, $Z_t$ is $\cF_{t}$ measurable and thus $\cG_{t}$ measurable, and for $t\geq 2$,
\begin{align*}
&\Ec{Z_{t}-Z_{t-1}}{\cG_{t-1}} \\&=\Ec{\gamma_{J(t)}  (\reaction_t(Y_t, \cC_t(X_t,\qt_t),q_t)- \widetilde{\Pi}_t(\qt_t)) \cdot S_t(X_t, \qt_t)}{\cG_{t-1}}\\
&=0+\gamma_{J(t)} S_t(X_t, \qt_t)\Ec{\reaction_t(Y_t, \cC_t(X_t,\qt_t),q_t)- \widetilde{\Pi}_t(\qt_t) }{S_t(X_t,\qt_t)=1,\cF_{t-1}}.
\end{align*}
Now, the latter is equal to zero by definition of $\widetilde{\Pi}_t(q)$ and because $q_t$ is $\cF_{t-1}$ measurable.
This means that $\Ec{Z_{t}}{\cF_{t-1}}=Z_{t-1}$ for $t\geq 2$ and thus $(Z_t)_{t\geq 1}$ is a martingale with respect to the filtration $(\cG_t)_{t\geq 1}$.
Now, this implies that $(\epsilon_{t,T})_{T\geq t}$ is also a centered martingale with respect to the filtration $(\cG_{T})_{T\geq t}$. In addition, we have for $T>t$,
$$
\Delta\epsilon_T:=\epsilon_{t,T}-\epsilon_{t,T-1} = \gamma_{J(T)}(\reaction_T(Y_T, \cC_T(X_T,\qt_T))- \widetilde{\Pi}_T(\qt_T)) \cdot S_T(X_T, \qt_T) \in [A_T, B_T],
$$
with $A_T=-\widetilde{\Pi}_T(\qt_T)\gamma_{J(T)} S_T(X_T, \qt_T)$ and $B_T=(1-\widetilde{\Pi}_T(\qt_T))\gamma_{J(T)} S_T(X_T, \qt_T)$ which are both $\cG_{T-1}$ measurable.
In addition, the square compensator is
\begin{align*}
\langle \epsilon\rangle_T&:=\sum_{k=t}^T   \E[(\Delta\epsilon_k)^2\:|\: \cG_{T-1}] \\
&=\sum_{k=t}^T   \E[ (\gamma_{J(k)}(\reaction_k(Y_k, \cC_k(X_k,\qt_k))- \widetilde{\Pi}_k(\qt_k)) \cdot S_k(X_k, \qt_k))  )^2 \:|\: \cG_{T-1}] \\
&\leq  (1/4)\sum_{k=t}^{T}  \gamma_{J(k)}^2 \cdot S_k(X_k, \qt_k) =(1/4)\sum_{j=J(t)}^{J(T+1)-1}  \gamma_{j}^2,
\end{align*}
by using Lemma~\ref{lem:computing}.
Since 
$$
\mathcal{D}_T = \sum_{k=t}^T (B_k-A_k)^2  \leq \sum_{k=t}^T \gamma^2_{J(k)} \cdot S_k(X_k, \qt_k)=\sum_{j=J(t)}^{J(T+1)-1}  \gamma^2_{j},
$$
by using again Lemma~\ref{lem:computing}.
we have by Lemma~\ref{lem:concentrationMart} that
for all $x,y>0$, $T\geq t$,
$$
\P\Big(\epsilon_{t,T}>x, \sum_{j=J(t)}^{J(T+1)-1}  \gamma_{j}^2\leq 2y/3 \Big) \leq  e^{-3x^2/y}.
$$

\begin{lemma}\label{lem:computing}
For any function $\Psi(\cdot)$, for all $T\geq 1$, we have
\begin{equation}\label{equ:computing}
    \sum_{k=1}^{T}\Psi(\gamma_{J(k)}) S_k(X_k,\qt_k) = \sum_{j=1}^{J(T+1)-1} \Psi(\gamma_{j}) .
\end{equation}
\end{lemma}

\begin{proof}
Denote for all $j\geq 1$,
$$
\tau_j = \min\{k\geq 1 \::\: J(k+1)-1=j\}  = \min\Big\{k\geq 1 \::\: \sum_{i=1}^{k} \St_{i}(X_i, \qt_i) =j\Big\}
$$
the first time $k$ for which there have been $j$ selections at times $\leq k$ (by convention $\tau_j=+\infty$ if the $j$-th selection never occurs).
Note that when $\tau_j<+\infty$, $J(\tau_j+1)-1=j$, $J(\tau_j)=j$ and $\St_{k}(X_k, \qt_k)=1$ for $k=\tau_j$.
Hence, for any function $\Psi(\cdot)$, for all $T\geq 1$, denoting $J_0=J(T+1)-1$, so that $j\leq J_0$ if and only if $\tau_j\leq T$
\begin{align*}
    \sum_{k=1}^{T}\Psi(\gamma_{J(k)}) S_k(X_k,\qt_k) &=  \sum_{j=1}^{J_0} \Psi(\gamma_{J(\tau_j)}) \sum_{k=1}^{T} \ind{k=\tau_j} S_k(X_k,\qt_k) \\
    &=  \sum_{j=1}^{J_0} \Psi(\gamma_{j}) \sum_{k=1}^{T} \ind{k=\tau_j}  \\
&=  \sum_{j=1}^{J_0} \Psi(\gamma_{j}) . 
\end{align*}
This gives \eqref{equ:computing}.
\end{proof}

\subsection{Proof of Theorem~\ref{th:genL2withoutt}}\label{proofth:genL2withoutt}

Similarly to \S~\ref{proof:th:almostsureconv}, we show a version of Theorem~\ref{th:genL2withoutt} which is slightly more general, accommodating the case where $\Pi^\bench_t$ may depend on time.

\subsubsection{Extending Theorem~\ref{th:genL2withoutt}}

\begin{theorem}\label{th:genL2}
Consider the \method\  threshold sequence $(q_t^{\OSCI})_{t\geq 1}$ \eqref{equACIquantile} with any starting point $q_1^{\OSCI}\in [0,B]$ and any step size sequence $(\gamma_j)_{j \geq 1}$. For each $t\geq 1$, consider some function $\Pi^{\bench}_t$ satisfying Assumption~\ref{ass:crossing} with some $q^\bench_{t,\alpha}\in (0,B)$ and the corresponding $D_t$ in \eqref{equDtT}.  
Let $a_t$ and $\tilde{a}_t$ being defined by \eqref{equ:a} and \eqref{equ:atildewithoutt}, respectively, in which $\Pi^\bench(q)$ and $q^{\bench}_{\alpha}$ have been replaced by $\Pi^\bench_t(q)$ and $q^{\bench}_{t,\alpha}$, respectively.
Let 
$\tilde{a}=\inf_{t\geq 1} \tilde{a}_t \geq 0$. 
Assume that  for some $t\geq j_0\geq 1$, 
\begin{itemize}
    \item $(\gamma_j)_{j_0\leq j\leq t}$ is nonincreasing, $1-2\tilde{a} \gamma_{j_0}\geq 0$ and  $\gamma_j/\gamma_{j+1}\leq 1/(1-\tilde{a}\gamma_j)$ for $j_0\leq j\leq t$;
    \item $S_k$ is independent of any $\mathcal{F}_{k-1}$-measurable variable, for $j_0\leq k\leq t$.
\end{itemize} Then, for any $t_0$ with $j_0\leq t_0 \leq t$ with $J(t_0)\geq j_0$, we have pointwise
    \begin{align}
\e{\paren{q^{\OSCI}_{t}-q^{\bench}_{t,\alpha}}^2\mid S_{1}\dots S_{t-1} }
&\leq 
\frac{\gamma_{J(t)}}{\gamma_{J(t_0)}} \paren{(B+\sup_{j\geq 1}\gamma_{j})^2 \vee \frac{4}{\tilde{a}}} \nonumber\\
&+4(B+\sup_{j\geq 1}\gamma_{j})\gamma_{J(t)}\sum_{k=t_0}^{t-1}\abs{q^{\bench}_{k+1,\alpha}-q^{\bench}_{k,\alpha}} \nonumber\\
&+ 4 \gamma_{J(t)}\sum_{k=t_0}^{t-1} S_k\paren{ (B +\sup_{j\geq 1}\gamma_{j})\e{D_k\mid S_{1}\dots S_{k-1}} - \gamma_{J(k)}}_+  \label{boundL2withminus},
\end{align}
where we denoted $S_t=S_t(X_t,q^{\OSCI}_{t})$. Note that the bound \eqref{boundL2withminus} is a bit better than \eqref{equ:genL2withoutt}, because of the ``$- \gamma_{J(k)}$'' term in the last sum of the right-hand side.
\end{theorem}
Theorem~\ref{th:genL2} is proved in \S~\ref{proofth:genL2}.

\subsubsection{Proof of Theorem~\ref{th:genL2}} \label{proofth:genL2}

In this proof, we denote $q^{\OSCI}_{t}$ by $q_{t}$ and $S_t(X_t,q_t)$ by $S_t$ for short. 
We also denote $S_{1:i}$ for the sequence $S_1,\dots,S_i$.  
We have for all $t\geq t_0$, 
\begin{align*}
&\frac{\e{\paren{q_{t}-q^{\bench}_{t,\alpha}}^2\mid S_{1:{t-1}} }}{\gamma_{J(t)}}\vee \frac{4}{\tilde{a}} - \frac{\e{\paren{q_{t_0}-q^{\bench}_{t_0,\alpha}}^2\mid S_{1:{t_0-1}} }}{\gamma_{J(t_0)}}\vee \frac{4}{\tilde{a}} \\
  &=\sum_{k=t_0}^{t-1} \paren{ \frac{\e{\paren{q_{k+1}-q^{\bench}_{k+1,\alpha}}^2\mid S_{1:{k}} }}{\gamma_{J(k+1)}}\vee \frac{4}{\tilde{a}} -\frac{\e{\paren{q_{k}-q^{\bench}_{k,\alpha}}^2\mid S_{1:{k-1}} }}{\gamma_{J(k)}}\vee \frac{4}{\tilde{a}}} .
  \end{align*}
We now use the following relation (see below for a proof): for $k\geq t_0$ such that $J(k)\geq j_0$,
\begin{align}
&  \frac{\e{\paren{q_{k+1}-q^{\bench}_{k+1,\alpha}}^2\mid S_{1:k}}}{\gamma_{J(k+1)}} \leq \frac{\e{\paren{q_{k}-q^{\bench}_{k,\alpha}}^2\mid S_{1:k-1}}}{\gamma_{J(k)}}\vee \frac{4}{\tilde{a}}\label{equ:intermnewproof_main}\\
& \:\:\:\:\:\:\:\:\:\:\:\:\:\:+ 4(B+\sup_{j\geq 1}\gamma_{j})\abs{q^{\bench}_{k+1,\alpha}-q^{\bench}_{k,\alpha}}
+ 4S_k \paren{(B +\sup_{j\geq 1}\gamma_{j}) \e{D_k\mid S_{1:k-1}} - \gamma_{J(k)}}_+\nonumber.
\end{align}
Putting this relation into the above display provides that, provided that $J(t_0)\geq j_0$, for all $t\geq t_0$,
\begin{align*}
&\frac{\e{\paren{q_{t}-q^{\bench}_{t,\alpha}}^2\mid S_{1:{t-1}} }}{\gamma_{J(t)}}\vee \frac{4}{\tilde{a}} 
\leq 4(B+\sup_{j\geq 1}\gamma_{j})\sum_{k=t_0}^{t-1}\abs{q^{\bench}_{k+1,\alpha}-q^{\bench}_{k,\alpha}}\\
&\:\:\:\:\:\:\:\:\:\:\:\:\:\:+\frac{(B+\sup_{j\geq 1}\gamma_{j})^2}{\gamma_{J(t_0)}} \vee \frac{4}{\tilde{a}}  + 4\sum_{k=t_0}^{t-1} S_k\paren{(B +\sup_{j\geq 1}\gamma_{j}) \e{D_k\mid S_{1:k-1}} - \gamma_{J(k)}}_+,
\end{align*}
which proves the result. 

Let us now finally prove \eqref{equ:intermnewproof_main}. 
To show this, first observe that by definition
\begin{align*}
    \paren{q_{k+1}-q^{\bench}_{k+1,\alpha}}^2 &= \paren{q_{k}-q^{\bench}_{k,\alpha} + q^{\bench}_{k,\alpha}-q^{\bench}_{k+1,\alpha} + \gamma_{J(k)}\paren{\reaction_k\paren{Y_k, \cC_k(X_k, q_k), q_k} - \alpha}S_k}^2 \\
    &\leq \paren{q_{k}-q^{\bench}_{k,\alpha}}^2    + 2 \gamma_{J(k)} \paren{q_{k}-q^{\bench}_{k,\alpha}}\paren{\reaction_k\paren{Y_k, \cC_k(X_k, q_k), q_k} - \alpha}S_k \\
    &+ 2\paren{q_{k}-q^{\bench}_{k,\alpha}}\paren{q^{\bench}_{k,\alpha}-q^{\bench}_{k+1,\alpha}} + 2\gamma_{J(k)}^2S_k+ 2\paren{q^{\bench}_{k,\alpha}-q^{\bench}_{k+1,\alpha}}^2\\
    &\leq \paren{q_{k}-q^{\bench}_{k,\alpha}}^2    + 2 \gamma_{J(k)} \paren{q_{k}-q^{\bench}_{k,\alpha}}\paren{\reaction_k\paren{Y_k, \cC_k(X_k, q_k), q_k} - \alpha}S_k\\
    &+ 2\gamma_{J(k)}^2S_k+ 4(B+\sup_{j\geq 1}\gamma_{j})\abs{q^{\bench}_{k+1,\alpha}-q^{\bench}_{k,\alpha}}.
\end{align*}
Denote for $q\in \R$, $\widetilde{\Pi}_k(q)=\E[\reaction_k(Y_k, \cC_k(X_k,q)) \:|\:   \mathcal{F}_{k-1},S_k=1 ] = \Pi_k(q)  \cdot \ind{q \in [0,B]} + \ind{q < 0} $, where $\Pi_k(q):=\Pi_{\cC_k,S_k}(q) $ is given by \eqref{PiCtSt}. 
We thus obtain
 \begin{align}
  &\frac{\e{\paren{q_{k+1}-q^{\bench}_{k+1,\alpha}}^2\mid S_{1:k}}}{\gamma_{J(k+1)}} -4(B+\sup_{j\geq 1}\gamma_{j})\abs{q^{\bench}_{k+1,\alpha}-q^{\bench}_{k,\alpha}}\nonumber\\
  &\leq  \frac{\e{\paren{q_{k}-q^{\bench}_{k,\alpha}}^2\mid S_{1:k}}}{\gamma_{J(k+1)}}  +\frac{2\gamma_{J(k)}^2}{\gamma_{J(k+1)}}S_k + \frac{2\gamma_{J(k)}}{\gamma_{J(k)+1}}\e{ \paren{q_{k}-q^{\bench}_\alpha}^2 \frac{\widetilde{\Pi}_k(q_k)-\alpha}{q_{k}-q^{\bench}_{k,\alpha}} \mid S_1,\dots,S_k}S_k\nonumber\\
   &= \frac{\e{\paren{q_{k}-q^{\bench}_{k,\alpha}}^2\mid S_{1:k}}}{\gamma_{J(k+1)}}  +\frac{2\gamma_{J(k)}^2}{\gamma_{J(k+1)}}S_k + \frac{2\gamma_{J(k)}}{\gamma_{J(k)+1}}\e{ \paren{q_{k}-q^{\bench}_\alpha}^2 \frac{\widetilde{\Pi}_k(q_k)-\alpha}{q_{k}-q^{\bench}_{k,\alpha}} \mid S_{1:k-1}}S_k,
  \label{equintermnewproof1}
\end{align}
where we used in the last step that by assumption $S_k$ is independent of the $\mathcal{F}_{k-1}$-measurable variable $q_k$.
Similarly, we have 
\begin{align*}
&   \frac{\e{\paren{q_{k}-q^{\bench}_{k,\alpha}}^2\mid S_{1:k}}}{\gamma_{J(k+1)}} \\
&=\frac{\e{\paren{q_{k}-q^{\bench}_{k,\alpha}}^2\mid S_{1:k}}}{\gamma_{J(k)}}(1-S_k) + \frac{\e{\paren{q_{k}-q^{\bench}_{k,\alpha}}^2\mid S_{1:k}}}{\gamma_{J(k)+1}}S_k\\
&=\frac{\e{\paren{q_{k}-q^{\bench}_{k,\alpha}}^2\mid S_{1:k-1}}}{\gamma_{J(k)}}(1-S_k) + \frac{\e{\paren{q_{k}-q^{\bench}_{k,\alpha}}^2\mid S_{1:k-1}}}{\gamma_{J(k)+1}}S_k,
\end{align*}
 In addition, by Assumptions~\ref{ass:crossing}~and~\ref{ass:slope}, and
since $-\alpha \sup_{j\geq 1}\gamma_{j}\leq q_k\leq B+(1-\alpha)\sup_{j\geq 1}\gamma_{j}$ and $\abs{q_{k}-q^{\bench}_{k,\alpha}}\leq B+\sup_{j\geq 1}\gamma_{j}$ almost surely by Lemma~\ref{lem:boundv}, we derive by using \eqref{equ:usefulslope} and \eqref{equDtTtilde} that
\begin{align*}
\paren{q_{k}-q^{\bench}_{k,\alpha}}^2 \frac{\widetilde{\Pi}_k(q_k)-\alpha}{q_{k}-q^{\bench}_{k,\alpha}}&\leq 
- \paren{q_{k}-q^{\bench}_{k,\alpha}}^2 \abs{\frac{\widetilde{\Pi}_k^{\bench}(q_k)-\alpha}{q_{k}-q^{\bench}_{k,\alpha}}} + \abs{q_{k}-q^{\bench}_{k,\alpha}} \abs{\widetilde{\Pi}_t(q_t^{\OSCI})-\widetilde{\Pi}^\bench_t(q_t^{\OSCI})}\\
&\leq - \tilde{a}\paren{q_{k}-q^{\bench}_{k,\alpha}}^2 + (B +\sup_{j\geq 1}\gamma_{j}) D_k.
\end{align*}
Putting the latter back into \eqref{equintermnewproof1} leads to 
{\small
\begin{align*}
  &\frac{\e{\paren{q_{k+1}-q^{\bench}_{k+1,\alpha}}^2\mid S_{1:k}}}{\gamma_{J(k+1)}} -4(B+\sup_{j\geq 1}\gamma_{j})\abs{q^{\bench}_{k+1,\alpha}-q^{\bench}_{k,\alpha}}-\frac{\e{\paren{q_{k}-q^{\bench}_{k,\alpha}}^2\mid S_{1:k-1}}}{\gamma_{J(k)}}(1-S_k) \nonumber\\
  &\leq  S_k\paren{\frac{2\gamma_{J(k)}^2}{\gamma_{J(k)+1}} +\frac{\e{\paren{q_{k}-q^{\bench}_{k,\alpha}}^2\mid S_{1:k-1}}}{\gamma_{J(k)+1}}  + \frac{2\gamma_{J(k)}}{\gamma_{J(k)+1}}\e{ - \tilde{a}\paren{q_{k}-q^{\bench}_{k,\alpha}}^2 + (B +\sup_{j\geq 1}\gamma_{j}) D_k \mid S_{1:k-1}}}\\
& =  S_k\frac{\gamma_{J(k)}}{\gamma_{J(k)+1}} \paren{ 2\gamma_{J(k)}+(1-2\tilde{a}\gamma_{J(k)} )\frac{\e{\paren{q_{k}-q^{\bench}_{k,\alpha}}^2\mid S_{1:k-1}}}{\gamma_{J(k)}}  + 2(B +\sup_{j\geq 1}\gamma_{j}) \e{D_k\mid S_{1:k-1}}} \\
& \leq  S_k\frac{\gamma_{J(k)}}{\gamma_{J(k)+1}} \paren{ 2\gamma_{J(k)}+(1-2\tilde{a}\gamma_{J(k)} )\paren{\frac{\e{\paren{q_{k}-q^{\bench}_{k,\alpha}}^2\mid S_{1:k-1}}}{\gamma_{J(k)}}\vee \frac{4}{\tilde{a}}}  + 2(B +\sup_{j\geq 1}\gamma_{j}) \e{D_k\mid S_{1:k-1}}} 
\end{align*}
}
Now using $\gamma_j/\gamma_{j+1}\leq 1/(1-\tilde{a}\gamma_j)$ and $2(1-\tilde{a} \gamma_{j})\geq 1$ for $j=J(k)\geq j_0$, the latter is at most $S_k$ times
\begin{align*}
 & \frac{\e{\paren{q_{k}-q^{\bench}_{k,\alpha}}^2\mid S_{1:k-1}}}{\gamma_{J(k)}}\vee \frac{4}{\tilde{a}}
 + 4(B +\sup_{j\geq 1}\gamma_{j}) \e{D_k\mid S_{1:k-1}}\\
 &+ \frac{\gamma_{J(k)}}{1-\tilde{a}\gamma_{J(k)}} \paren{2-  \tilde{a}  \paren{\frac{\e{\paren{q_{k}-q^{\bench}_{k,\alpha}}^2\mid S_{1:k-1}}}{\gamma_{J(k)}}\vee \frac{4}{\tilde{a}}} } \\
 \leq \:&\frac{\e{\paren{q_{k}-q^{\bench}_{k,\alpha}}^2\mid S_{1:k-1}}}{\gamma_{J(k)}}\vee \frac{4}{\tilde{a}}
 + 4(B +\sup_{j\geq 1}\gamma_{j}) \e{D_k\mid S_{1:k-1}} - 4\gamma_{J(k)}.
 \end{align*}
This establishes \eqref{equ:intermnewproof_main}.

 \subsection{Proof of Lemma~\ref{lem:boundIER}}\label{sec:proofboundIER}
 
Note that, with the notation \eqref{equ:pizerotilde}, the following inequality is true:
\begin{equation}\label{equ:usefulslopebis}
\sup_{t\geq 1} \frac{\abs{\widetilde{\Pi}^{\bench}(q^{\OSCI}_t)-\alpha}}{\abs{q^{\OSCI}_t-q^{\bench}_{\alpha}}}\leq \sup\set{\frac{\abs{\widetilde{\Pi}^{\bench}(q)-\alpha}}{\abs{q-q^{\bench}_{\alpha}}}, -\alpha \sup_{j\geq 1}\gamma_{j} \leq q\leq  B+(1-\alpha)\sup_{j\geq 1}\gamma_{j} } 
 \leq \tilde{b} .
\end{equation}

If $q_t^{\OSCI}\in [0,B]$, we have by the triangle inequality,
 \begin{align*}
 |\IER_t(\mathcal{R}^{\OSCI}) - \alpha|&=|\Pi_t(q_t^{\OSCI})-\alpha|\\
  &\leq |\Pi_t(q_t^{\OSCI})-\Pi^\bench(q_t^{\OSCI})| + |\Pi^\bench(q_t^{\OSCI}) -\Pi^\bench(q_\alpha^\bench)|\\
 &\leq D_t + \tilde{b} |q_t^{\OSCI} -q_\alpha^\bench|,
 \end{align*}
because $\frac{|\Pi^\bench(q_t^{\OSCI}) -\Pi^\bench(q_\alpha^\bench)|}{|q_t^{\OSCI} -q_\alpha^\bench|} =\frac{\abs{\widetilde{\Pi}^{\bench}(q^{\OSCI}_t)-\alpha}}{\abs{q^{\OSCI}_t-q^{\bench}_{\alpha}}} \leq \tilde{b} $ since $q_t^{\OSCI}\in [0,B]$ and by \eqref{equ:usefulslopebis}.

\section{Proofs of main results}

\subsection{Proof of Theorem~\ref{thangelinfo}}\label{proof:thangelinfo}

We first state and prove the following lemma.
\begin{lemma}\label{lem:boundv}
For all $t\geq 1$ and any sequence $(X_{t},Y_{t})_{t\geq 1}$ the sequence of thresholds $(\qt_t^{\OSCI})_{t\geq 1}$ of \method\  is such that $\qt_t^{\OSCI}\in [-\alpha \max_{1 \leq j \leq J(t)}\gamma_{j}, B+(1-\alpha)\max_{1 \leq j \leq J(t)}\gamma_{j}]$.
\end{lemma}
\begin{proof}
    \et{We prove the slightly stronger result $\qt_t^{\OSCI}\in [-\alpha \max_{1 \leq j \leq J(t-1)}\gamma_{j}, B+(1-\alpha)\max_{1 \leq j \leq J(t-1)}\gamma_{j}]$ by recursion with respect to $t$. The result is clear for $t=2$ because $q_1\in [0,B]$ (note that $J(1)=1$). Assume it is true at time $t$ and that there is a selection at time $t$ (otherwise the result is trivial). The result is clear if $q_t\in [0,B]$. 
    If $\qt_{t} < 0$, then $\qt_{t+1} \geq \qt_t \geq -\alpha \max_{1 \leq j \leq J(t-1)}\gamma_{j} \geq -\alpha \max_{1 \leq j \leq J(t)}\gamma_{j}$. Also, $\qt_{t+1} = q_t  + \gamma_{J(t)}(1-\alpha) \leq \gamma_{J(t)}(1-\alpha)\leq B+(1-\alpha)\max_{1 \leq j \leq J(t)}\gamma_{j}.$ 
    If $\qt_{t} > B$, either there is no selection, and in that case $q_{t+1}=q_t$, or there is a selection and thus no mistake, and in that case $\qt_{t+1} \leq \qt_t$, which means $\qt_{t+1} \leq B+(1-\alpha)\max_{1 \leq j \leq J(t-1)} \gamma_{j}$ and $\qt_{t+1} \geq B-\alpha \gamma_{J(t)} \geq -\alpha \max_{1 \leq j \leq J(t)} \gamma_{j}$. In any case, $\qt_{t+1}^{\OSCI}\in [-\alpha \max_{1 \leq j \leq J(t)}\gamma_{j}, B+(1-\alpha)\max_{1 \leq j \leq J(t)}\gamma_{j}]$ and the recursion is proved.
    }
\end{proof}

We now turn to prove Theorem~\ref{thangelinfo}, by following \cite{angelopoulos2024online} (Theorem~2 therein), but restricting only to selected indices. 
We have by definition of $\FP_t$ that, for all $t\geq 1$, letting $N:=N(t)$ (assumed positive otherwise the result is clear),
\begin{align*}
\FP_t(\mathcal{R})-\alpha 
&= N^{-1}\sum_{k=1}^N (\Err_{\tau_k}(X_{\tau_k}, Y_{\tau_k}, \qt_{\tau_k}) - \alpha)\\
&\leq N^{-1}\sum_{k=1}^N (\reaction_{\tau_k}(X_{\tau_k}, Y_{\tau_k}, \qt_{\tau_k}) - \alpha),
\end{align*}
where we used that $\Err_t(x, y, q) \leq \reaction_t(x, y, q)$ for all $x,y,q$ and $\tau_k$ denotes the time at which we have made the $k$-th selection. Then the last display is
\begin{align*}
&N^{-1}\sum_{k=1}^N \gamma^{-1}_{k} [\gamma_{k}(\reaction_{\tau_k}(X_{\tau_k}, Y_{\tau_k}, \qt_{\tau_k}) - \alpha)]\\
&= N^{-1}\sum_{k=1}^N \sum_{j=1}^{k}(\gamma^{-1}_{j}-\gamma^{-1}_{j-1})[\gamma_{k}(\reaction_{\tau_k}(X_{\tau_k}, Y_{\tau_k}, \qt_{\tau_k}) - \alpha)]\\
&= N^{-1}\sum_{j=1}^{N} (\gamma^{-1}_{j}-\gamma^{-1}_{j-1}) \sum_{k=j}^N  [\gamma_{k}(\reaction_{\tau_k}(X_{\tau_k}, Y_{\tau_k}, \qt_{\tau_k}) - \alpha)] \\
&= N^{-1}\sum_{j=1}^{N} (\gamma^{-1}_{j}-\gamma^{-1}_{j-1}) (q_{\tau_{N+1}}-q_{\tau_j}), 
\end{align*}
where $\gamma^{-1}_0 = 0$ by convention \et{and because $q_{\tau_{k+1}}= q_{\tau_{k}+1}=q_{\tau_k}+\gamma_{J(\tau_k)}(\reaction_{\tau_k}(X_{\tau_k}, Y_{\tau_k}, \qt_{\tau_k}) - \alpha)$, with $J(\tau_k)=N(\tau_k-1)+1=k$. }
Hence we obtain
\begin{align*}
\FP_t(\mathcal{R})-\alpha &\leq \big\lvert N^{-1}\sum_{j=1}^{N} (\gamma^{-1}_{j}-\gamma^{-1}_{j-1}) (\qt_{\tau_{N+1}} - \qt_{\tau_j}) \big\rvert \\
&\leq N^{-1}\sum_{j=1}^{N} \lvert \gamma^{-1}_{j}-\gamma^{-1}_{j-1}\rvert \lvert \qt_{\tau_{N+1}} - \qt_{\tau_j}\rvert \; ,
\end{align*}
Now, using Lemma~\ref{lem:boundv}, we have $\lvert \qt_{\tau_{N+1}} - \qt_{\tau_j}\rvert \leq (B + \max_{1 \leq j \leq J}(\gamma_j))$ with $J:=J(t)$ and thus
\begin{align*}
\FP_t(\mathcal{R})-\alpha &\leq \dfrac{B + \max_{1 \leq j \leq J}(\gamma_j)}{N}\sum_{j=1}^{N} \lvert \gamma^{-1}_{j}-\gamma^{-1}_{j-1} \rvert \\
&= \dfrac{B + \max_{1 \leq j \leq J}(\gamma_j)}{N} \left(\gamma^{-1}_{1} + \sum_{j=2}^{N} \lvert \gamma^{-1}_{j}-\gamma^{-1}_{j-1} \rvert \right) \; .
\end{align*}
This gives the first FCP upper-bound. When $(\gamma_t)_{t \geq 1}$ is nonincreasing, $\sum_{j=1}^{N} \lvert \gamma^{-1}_{j}-\gamma^{-1}_{j-1} \rvert = \sum_{j=1}^{N} \gamma^{-1}_{j}-\gamma^{-1}_{j-1} = \gamma^{-1}_N$ and we obtain the second FCP upper-bound.

\subsection{Proof of Theorem~\ref{th:consistency}}\label{proof:th:consistency}

As proven in \S~\ref{proof:th:almostsureconv}, the process $(Z_t)_{t\geq 1}$ defined in \eqref{equmartinfo} is a centered martingale with respect to the filtration $(\cG_t)_{t\geq 1}$ with $\sup_{T\geq 1}\E(Z_T^2)\leq \sum_{t=1}^{+\infty} \gamma_{t}^2<+\infty$ by assumption.  
Therefore, by Doob's martingale convergence theorem, there exists some random variable $Z$ such that $Z_t\to Z$ when $t$ tends to infinity, almost surely. This means that, almost surely, $Z_t$ is a Cauchy sequence. 

Recall now the bound \eqref{eq:QuantileBoundsLearningwithoutt}, proved in Theorem~\ref{th:almostsureconv:withoutt}: for all $T>t\geq 1$, we have pointwise
\begin{align*}
\abs{\qt_T^{\OSCI}-q^{\bench}_{\alpha}}\leq \:&
 (\max_{t\leq k\leq T} D_k)/\tilde{a} +\max_{t\leq k\leq T-1} | \epsilon_{k-1,T-1} | \nonumber\\
 &+ \gamma_{J(t)} + \paren{\tilde{a}\sum_{j=J(t)}^{J(T+1)-1} \gamma_{j}}^{-1} \paren{\abs{\epsilon_{t-1,T}} +B+\gamma_1}.
\end{align*}

Fix now $\varepsilon>0$ and assume that we are on the event $\{\lim_{t\to \infty}J(t)=+\infty\}\cap \{\lim_{t\to \infty}D_t=0\}$. Since $Z_t$ is a Cauchy sequence (almost surely) then there exist $T_1=T_1(\varepsilon)$ such that for all $T\geq t\geq T_1$ we have
$$\abs{\epsilon_{t-1,T}}\vee \max_{t\leq k\leq T-1}|\epsilon_{k-1,T-1}|=|Z_T-Z_{t-1}|\vee \max_{t\leq k\leq T-1}|Z_{T-1} - Z_{k-1}|\leq \varepsilon.$$

Since $\lim_{j\to \infty}\gamma_{j}=0$, there exists $J_1=J_1(\varepsilon)$ such that for $J\geq J_1$, $\gamma_{J}\leq \varepsilon$. Also, since $\lim_{t\to \infty} J(t)=+\infty$, there exists $T_2=T_2(\varepsilon)$ such that for all $t\geq T_2$, $J(t)\geq J_1$ and thus $\gamma_{J(t)}\leq \varepsilon$. 
Furthermore, since $\lim_{t\to \infty}D_t=0$, there exists $T_3=T_3(\varepsilon)$ such that for all $T\geq t\geq T_3$,  $\max_{t\leq k\leq T} D_k\leq \tilde{a}\:\varepsilon$. 

Now consider $t=\max(T_1,T_2,T_3)$. Since  $\sum_{j=J(t)}^{J(T+1)-1}\gamma_{j}$ tends to infinity when $T$ tends to infinity, there exists $T_4=T_4(\varepsilon)\geq t$ such that for $T\geq T_4$, $$\paren{\sum_{j=J(t)}^{J(T+1)-1} \gamma_{j}}^{-1}\leq \frac{\tilde{a}}{\varepsilon+B+\gamma_1} \varepsilon.$$
This provides that for any $T\geq T_4$, we have $\abs{\qt_T^{\OSCI}-q^{\bench}_{\alpha}}\leq 4\varepsilon$ and concludes the proof.

\subsection{Proof of Theorem~\ref{th:rateconvfull}}\label{sec:proofcor:rateconvfull}

We establish the following more general result. 

\begin{proposition}
    \label{cor:rateconvfull_append}
  Consider  the \method\  threshold sequence $(q_t^{\OSCI})_{t\geq 1}$ \eqref{equACIquantile} with step size sequence $\gamma_j= c j^{-\beta}$, $\beta\in (1/2,1)$ and starting point $q_1^{\OSCI}\in [0,B]$. Consider some function $\Pi^{\bench}$ satisfying Assumptions~\ref{ass:crossing}~and~\ref{ass:slope} with some $q^\bench_{\alpha}\in (0,B)$ and the corresponding $D_t$ being defined by \eqref{equDtT}, which is assumed to be bounded by $(\bar{D}_t)_{t\geq 1}$ as in \eqref{equ:Dtcondition}, with a rate $(\eta_t)_{t\geq 1}\in (0,1)$. 
In addition, we consider $(\underline{J}(t))_{t\geq 1}$ and $(\bar{J}(t))_{t\geq 1}$ two increasing integer (deterministic) sequences and consider $\omega(t)$ given by \eqref{equprobaJtsympa}.
  Then $(q_t^{\OSCI})_{t\geq 1}$ is such that there exist $T_0=T_0(\beta)\geq 1$ and $C_0=C_0(\beta,a,c)>0$ such that if $t>k\geq T_0$, we have with probability at least $1-3\eta_k-\omega(k)-\omega(t+1)$, 
\begin{align}
\abs{\qt_t^{\OSCI}-q^{\bench}_{\alpha}} &\leq 
 x(k,t,c)+ \tilde{a}^{-1} C_1\bar{D}_k +c (\underline{J}(k))^{-\beta}
 + \left(\tilde{a}c\sum_{j=\bar{J}(k)}^{\underline{J}(t+1)-1} j^{-\beta} \right)^{-1} \paren{x(k,t,c) +B+c },\label{eq:simpleboundbao}
\end{align}
where 
\begin{align}
x(k,t,c)&:=\brac{c^2(\log(t-k)+\log(1/\eta_t)) \sum_{j=\underline{J}(k)}^{\bar{J}(t+1)-1} j^{-2\beta}}^{1/2}\label{equ:yeta}.
\end{align}
\end{proposition}

Proposition~\ref{cor:rateconvfull_append} implies Theorem~\ref{th:rateconvfull} (case (ii)): indeed, considering the regime of $k(t)$, $\Delta(t)$ given in Theorem~\ref{th:rateconvfull} (ii), and apply the two first statements of Lemma~\ref{lemsuminte}:
 \begin{align*}
 \sum_{j=\underline{J}(k(t))}^{\bar{J}(t+1)-1} j^{-\beta} &\sim    \sum_{j=\bar{J}(k(t))}^{\underline{J}(t+1)-1} j^{-\beta} \sim \Delta(t) \bar{J}(t)^{-\beta};\\
   \sum_{j=\underline{J}(k(t))}^{\bar{J}(t+1)-1} j^{-2\beta} &\sim \Delta(t)  \bar{J}(t)^{-2\beta}.
\end{align*}
This means that there exists $T_0=T_0(\beta)\geq 1$ such that for all $t\geq T_0$, we have $x(k,t,c)\leq 2c (\log(t-k(t))+\log(1/\eta_t))^{1/2} \Delta(t)^{1/2} \bar{J}(t)^{-\beta}$ and thus the following holds
\begin{align*}
&\P(|\IER_t(\mathcal{R}^{\OSCI}) - \alpha|\ind{\qt_t^{\OSCI}\in [0,B]}+ \abs{\qt_t^{\OSCI}-q^{\bench}_{\alpha}} \leq C_0 (\rho_t(\beta)+ \bar{D}_{k(t)}))\\
&\geq 1-3\eta_{k(t)} -\omega(k(t))-\omega(t+1),
\end{align*}
where we have used Lemma~\ref{lem:boundIER}. Since $\rho_t(\beta)+ \bar{D}_{k(t)}$ tends to zero when $t$ tends to infinity, we have that $\qt_t^{\OSCI}\in [0,B]$ provided that $\abs{\qt_t^{\OSCI}-q^{\bench}_{\alpha}} \leq C_0 (\rho_t(\beta)+ \bar{D}_{k(t)})$ when $t$ is large enough. This leads to \eqref{equ-conviid-specialcases-full}.
Case (i) is completely similar by applying the two last statements of Lemma~\ref{lemsuminte} (by taking $k(t)=\lfloor c't/2\rfloor$ and $\Delta(t)\sim c't/2$).
\\

Let us now prove Proposition~\ref{cor:rateconvfull_append} from the general bound \eqref{eq:QuantileBoundsLearningwithoutt}. 
By assumption, for all $T>t\geq 2$, there is an event $\Omega$ with probability at least $1-\omega(k)-\omega(t+1)$, for which
$$
\underline{J}(t+1)\leq J(t+1) \leq \bar{J}(t+1)  \mbox{ and } \underline{J}(k)\leq J(k) \leq \bar{J}(k),
$$
and thus on $\Omega$, we have $J(k)^{-\beta} \leq \underline{J}(k)^{-\beta}$ and both
\begin{align*}
&\sum_{j=\underline{J}(k)}^{\bar{J}(t+1)-1} j^{-\beta}\geq \sum_{j=J(k)}^{J(t+1)-1} j^{-\beta} \geq \sum_{j=\bar{J}(k)}^{\underline{J}(t+1)-1} j^{-\beta}; \\
&\sum_{j=J(k+1)}^{J(t+1)-1} j^{-2\beta} \leq \sum_{j=\underline{J}(k+1)}^{\bar{J}(t+1)-1} j^{-2\beta}.
\end{align*}
In addition, we have by \eqref{equ:Dtcondition},
\begin{align}
    \P(\max_{k\leq \ell\leq t} D_\ell \geq C_1 \bar{D}_k )&\leq \sum_{\ell=k}^t \P( D_\ell \geq C_1 \bar{D}_k )\leq \sum_{\ell=k}^t \P( D_\ell \geq C_1 \bar{D}_\ell )\nonumber\\
    &\leq \eta_k\frac{6}{\pi^2} \sum_{\ell\geq 1} \ell^{-2} =\eta_k\label{equ:intermDt}.
\end{align}
Also, we have
\begin{align*}
\P\Big(\Omega \cap \{\max_{k\leq \ell\leq t-1} | \epsilon_{\ell-1,t-1} | >x\} \Big) &\leq \sum_{\ell = k}^{ t-1}\P\Big( \Omega \cap \{| \epsilon_{\ell-1,t-1} | >x\} \Big) \\
&\leq 2\sum_{\ell = k}^{ t-1} \P\bigg(\Omega \cap \Big\{c^2\sum_{j=J(\ell)}^{J(t)-1} j^{-2\beta} > y\Big\}\bigg)+  2\sum_{\ell = k}^{ t-1} e^{-x^2/y}\\
&\leq 2(t-k) \P\bigg(c^2\sum_{j=\underline{J}(k)}^{\bar{J}(t+1)-1} j^{-2\beta} > y\bigg) + 2(t-k)e^{-x^2/y} \\&\leq 2\eta_t \leq 2\eta_{k},
\end{align*}
by choosing $y=c^2\sum_{j=\underline{J}(k)}^{\bar{J}(t+1)-1} j^{-2\beta}$ and $x=x(k,t,c)$ given by \eqref{equ:yeta}.
Combining this with \eqref{eq:QuantileBoundsLearningwithoutt} and \eqref{equ:intermDt}, we obtain the desired result.

\subsection{Proof of Theorem~\ref{th:genL2-first}}\label{sec:cor:genL2}

We easily check that the considered step size sequence $(\gamma_j)_{j\geq 1}$ satisfies the assumption of Theorem~\ref{th:genL2withoutt}. The latter hence entails that for any $t_0 \leq t$ with $t_0\geq j_0$,
\et{
\begin{align*}
\e{\paren{q_{t}-q^{\bench}_\alpha}^2}
&\leq \e{ \e{\paren{q_{t}-q^{\bench}_\alpha}^2\mid S_{1}\dots S_{t-1}}\ind{\underline{J}(t)\leq J(t), j_0\leq J(t_0)}} \\
&\:\:\:\:+ (B+\gamma_1)^2 (\omega(t) +  \probp{J(t_0) < j_0})\\
&\leq  \paren{(B+\gamma_1)^2 \vee \frac{4}{\tilde{a}}} \frac{\gamma_{\underline{J}(t)}}{\gamma_{t_0}}+ 4((B +\gamma_1)\vee 1)\gamma_{\underline{J}(t)}\sum_{k=t_0}^{t-1} \e{S_k}\e{D_k} \\
&\:\:\:\:+ (B+\gamma_1)^2  (\omega(t) +  \probp{J(t_0) < j_0}),
\end{align*}
}
where we used that by the independence assumption
$$\e{S_k \e{D_k\mid S_{1}\dots S_{k-1}}}=\e{S_k}\e{ \e{D_k\mid S_{1}\dots S_{k-1}}} = \e{S_k} \e{D_k}.$$
\et{We apply this with $T_0=j_0$, $k(t)=t_0$ and use Lemma~\ref{lemsuminte} to obtain the desired bound for  $\e{\paren{q_{t}-q^{\bench}_\alpha}^2}$.
Next, we use Lemma~\ref{lem:boundIER} with Jensen's inequality to derive that
\begin{align*}
 \paren{\e{|\IER_t(\mathcal{R}^{\OSCI}) - \alpha|\ind{q_t^\OSCI\in [0,B]}}}^2&\leq  \paren{\E[D_t] + \tilde{b} \e{\paren{q_t^{\OSCI} -q_\alpha^\bench}^2}^{1/2}}^2\\&\leq 2 \paren{\E[D_t]}^2 + 2\tilde{b}^2 \e{\paren{q_t^{\OSCI} -q_\alpha^\bench}^2} ,
\end{align*}
which gives 
\begin{align*}
&\e{\abs{\IER_t(\mathcal{R}^{\OSCI}) - \alpha}\ind{q_t^\OSCI\in [0,B]}}^2\\
&\leq 2(\E[D_t])^2 +C_0 \omega(t)+ C_0(\underline{J}(t))^{-\beta}\paren{1+\sum_{k=1}^{t-1} \e{D_k}} +  C_0\probp{J(t_0) < j_0}
  .
\end{align*}
}
We finish the proof by using $\E[D_t]\leq \widetilde{D}_t$ and that $(\underline{J}(t))^{-\beta} \paren{1 + \sum_{k=1}^{t-1}\widetilde{D}_k} \gtrsim \widetilde{D}_t\gg (\widetilde{D}_t)^2$.

\subsection{Proof of Theorem~\ref{rem:coroonanEvent}}\label{sec:proof:rem:coroonanEvent}

First note that by the proof of Theorem~\ref{th:rateconvfull} the following result holds.
\begin{proposition}
[Restriction to an event]
    \label{rem:coroonanEvent0}
  Consider the same setting of Theorem~\ref{th:rateconvfull} except that $\omega(t)$ in \eqref{equprobaJtsympa} and $\bar{D}_t$ in \eqref{equ:Dtcondition} are computed in restriction to respective events $\Omega_t\in \mathcal{F}_{t-1}$, for $t\geq 1$, as in \eqref{equrestrictioncond}. The conclusion is that for $t\geq T_0$, 
     \begin{align}
&\P(|\IER_t(\mathcal{R}^{\OSCI}) - \alpha|+ \abs{\qt_t^{\OSCI}-q^{\bench}_{\alpha} } \leq C_0 (\rho_t+ \bar{D}_{k(t)}), \Omega_t)
\nonumber\\
&\geq \P(\cap_{t\geq T_0}\Omega_t)-C_0(\eta_{k(t)} +\omega(k(t))+\omega(t+1)).\nonumber
\end{align}
\end{proposition}
We now prove Theorem~\ref{rem:coroonanEvent} by applying Proposition~\ref{rem:coroonanEvent0}: we have for $t\geq T_0$, 
    \begin{align*}
        &\P(|\IER_t(\mathcal{R}^{\OSCI}) - \alpha|+ \abs{\qt_t^{\OSCI}-q^{\bench}_{\alpha} } \leq C_0 (\rho_t+ \bar{D}_{k(t)}))\\
        &\geq 
        \P(|\IER_t(\mathcal{R}^{\OSCI}) - \alpha|+ \abs{\qt_t^{\OSCI}-q^{\bench}_{\alpha} } \leq C_0 (\rho_t+ \bar{D}_{k(t)}), \Omega_t)\\
        & \geq  \P(\forall t\geq T_0, q_t\leq \bar q)-C_0(\eta_{k(t)} +\omega(k(t))+\omega(t+1)).
    \end{align*}
     We conclude by applying Lemma~\ref{lem:possible}.

\section{Proofs for regression cases}

We provide here proofs of our results in the regression cases (except for the informative selections, that are postponed to  \S~\ref{sec:proofsinforegress}). 

\subsection{Bounding $D_t$ in the online iid regression model}\label{sec:Dtboundiidregression}

Lemma~\ref{lem:forDT} leads to the following result by an integration.

\et{
\begin{proposition}\label{prop:Dtiideasyselect}
    For the online iid regression model \eqref{regression:model} with Assumption~\ref{ass:densityf}, consider the  distance $D_t\leq \|\Pi_t-\Pi^\bench\|_\infty$ \eqref{equDtT} for any $(\mathcal{F}_{t-1},X_t)$-measurable selection rule $S_t$ and the IER function \eqref{PiCtSt} and benchmark function $ \Pi^\bench$ \eqref{PizeroClassicalscore} of \S~\ref{sec:regression}. 
    Then, for all $t\geq 1$, we have almost surely 
   \begin{align}
&       \|\Pi_t-\Pi^\bench\|_\infty  \label{boundDtregression}\\
&       \leq   M_3\:\E[M_1 |\hat{\mu}_t(X_t)-{\mu}(X_t)|+ 2(M_2\vee 1)|\hat{\sigma}_t(X_t)-{\sigma}(X_t)| \:|\: S_t=1, \mathcal{F}_{t-1}].
\nonumber
    \end{align}
    In particular, for all $t\geq 1$, we have almost surely
    \begin{align}
      & \|\Pi_t-\Pi^\bench\|_\infty  \label{boundDtregression2}\\
&       \leq  \probp{S_t=1\:|\:  \mathcal{F}_{t-1}}^{-1}\paren{ M_3\:\E[M_1 |\hat{\mu}_t(X_t)-{\mu}(X_t)|+ 2(M_2\vee 1)|\hat{\sigma}_t(X_t)-{\sigma}(X_t)| \:|\:  \mathcal{F}_{t-1}]}.
\nonumber
    \end{align}
\end{proposition}
}

This shows that $D_t$ converges appropriately to zero when choosing appropriate estimators $\hat{\mu}_t$ and $\hat{\sigma}_t$. 
Assuming \eqref{ass:esti} and \eqref{ass:esti:average}, for a time-constant $X$-oriented selection rule, the bound \eqref{boundDtregression2} leads to condition \eqref{equ:Dtcondition} with $\bar{D}_t=O((\log t)^{\nu_0} t^{-\nu})$, and to $\widetilde{D}_t=O((\log t)^{\tau_0}t^{-\tau})$. 

\subsection{Proof of Corollary~\ref{cor:rateconv}}\label{sec:proofcor:rateconv}

We apply Theorem~\ref{th:rateconvfull} as in Remark~\ref{rem:ratetypical}, in each time-constant $X$-oriented case for $r:=\P(S(X_1)=1)>0$:
\begin{itemize}
    \item $r=1$, that is, $S_k(x,q)\equiv 1$: 
    in that case, $J(t)=t$ and thus $\omega(t)=0$ with $\underline{J}(t)= t = \bar{J}(t)$ and $\eta_t=1/t$. 
    We can apply Theorem~\ref{th:rateconvfull} both in case (i) if $\beta\geq 3/4$ and case (ii) with $k(t)=t- \Delta(t)$ and $\Delta(t)=1\vee \lfloor t^{4\beta/3} \rfloor$ if $\beta< 3/4$. This leads to the desired conclusion with $\rho_t(\beta)$ as in \eqref{equ:raterho-typicalconc}.
\item $r<1$ and $S(x)=S_k(x,q)$ does not depend on $q$ and $k$: then 
$J(t+1)-1=\sum_{i=1}^{t} S(X_{i})$ follows a binomial distribution of parameter $t,r:=\P(S(X)=1)$. 
Thus Bernstein's inequality gives that, if $\log(t)/t \leq 2 r(1-r)$, 
$$\P\Big(|\mathcal{B}(t,r)-tr|\leq 2\sqrt{r(1-r)t\log(t)}\Big)\geq 1-2/t.$$
As a result, we may choose $\bar{J}(t)\geq 1+\lceil (t-1)r + \{(t-1)\log(t-1)\}^{1/2}\rceil $ and $\underline{J}(t)\leq 1+\lfloor (t-1)r - \{(t-1)\log(t-1)\}^{1/2}\rfloor$ to obtain the required assumption with $\omega(t)=2/t$. 
Hence, $\underline{J}(t)$ and $\bar{J}(t)$ can be chosen as $\lfloor rt \pm C(\log t)^{1/2} t^{\lambda} \rfloor$ for $\lambda=1/2$ and some constant $C>0$. Hence Remark~\ref{rem:ratetypical} applies with $\lambda=1/2$, which gives the result.
\end{itemize}

\subsection{Proof of Corollary~\ref{th:MSEconv}} \label{proofth:MSEconv}

\et{Let $r:=\P(S(X_1)=1)$
and apply Theorem~\ref{th:genL2-first} which is possible because $S_k(X_k,q_k)=S(X_k)$ is always independent of $\mathcal{F}_{k-1}$ for a time-constant selection and in the iid model. To deal with $\omega(t)$, we choose $\omega(t)=2/t$ by taking  $\bar{J}(t):=1+\lceil (t-1)r + \{(t-1)\log(t-1)\}^{1/2}\rceil $, $\underline{J}(t):=1+\lfloor (t-1)r - \{(t-1)\log(t-1)\}^{1/2}\rfloor$ and $\omega(t)=2/t$, as justified in \S~\ref{sec:proofcor:rateconv} (which is also a valid choice if $r=1$). 
The additional term $ \probp{J(k(t))< T_0}$ is bounded as follows: considering $k(t)=\lceil C \log(2t)\rceil $ for a constant $C>0$ to be chosen further on, we have 
$$
\e{J(k(t))} \geq \e{J(k(t)) \ind{J(k(t))\geq \underline{J}(k(t))}}\geq \underline{J}(k(t)) \cdot (1-2/k(t))\geq T_0,
$$
for $t$ large enough. Hence, still using Bernstein's inequality, 
$$
\probp{J(k(t))< T_0} \leq \probp{J(k(t))\leq \E(J(k(t)))/2} \leq e^{-C_0 (\underline{J}(k(t)))^2/k(t)}\leq e^{-C'_0 k(t)}
$$
for some constant $C_0,C'_0>0$. Hence,  taking $k(t)=\log(t)/C'_0$ gives $ \probp{J(k(t))< T_0}\leq 1/t$. The result follows.
}

\subsection{Proof of Corollary~\ref{cor:rateconveasyselect}}\label{proof:cor:rateconveasyselect}

Let us first check that Assumption~\ref{ass:easy} is satisfied with $B=1$ and $\Pi^\bench$ and $q^\bench_{\alpha}$ given in \eqref{PizeroClassicalscore}. Clearly, $\Pi^\bench$ is decreasing and 
$
(\Pi^{\bench})'(q)=- \frac{f}{\phi}\big( \Phi^{-1}\big(\frac{1-q}{2}\big)\big), \:q\in (0,1),
$
 is a continuous function on the compact interval $[q^{\bench}_{\alpha}/2,(1+q^{\bench}_{\alpha})/2]$. In addition,  $$(\Pi^{\bench})'(q^\bench_\alpha)=- \frac{f}{\phi}\big( \ol{F}^{-1}\big(\alpha/2\big)\big)<0.$$ 
Hence, Assumption~\ref{ass:easy} holds.

Now, we apply Corollary~\ref{cor:rateconv} and Corollary~\ref{th:MSEconv} with the bounds on $\bar{D}_t=O((\log t)^{\nu_0} t^{-\nu})$ and $\widetilde{D}_t=O((\log t)^{\tau_0}t^{-\tau})$ obtained in \S~\ref{sec:Dtboundiidregression}. This gives \eqref{equ-conviid-specialcases-easyselect}.

Finally, note that from \eqref{equ-conviid-specialcases-easyselect} we have $\P(q_t^\OSCI\notin [0,B])\leq 1/t$. This can in turn be used into Corollary~\ref{th:MSEconv} to remove the indicator $\ind{q_t^\OSCI\in [0,B]}$ in the IER expectation (by assuming both \eqref{ass:esti} and \eqref{ass:esti:average}) and gives \eqref{equ:MSEconv-regressioneasy}.

\subsection{Proof of Corollary~\ref{cor:optrateiidregress}}\label{proof:cor:optrateiidregress}

Let us fix $x\in \mathcal{X}$ and proves point (i). By \eqref{equ:controlcondx}, we have
\begin{align*}
\probp{Y_t\notin \cC^{\OSCI}_t(x)\mid X_t=x}&=
\e{\probp{V_t(x,Y_t) >q_t^{\OSCI} \mid X_t=x, \mathcal{F}_{t-1}}} \\
&\leq \e{\Pi^\bench(q_t^{\OSCI})} + M_1  (1/{\sigma}(x)) \e{\paren{\hat{\mu}_t(x)-{\mu}(x)}^2}^{1/2}\\
&\:\:+ M_2 \e{\paren{\hat{\sigma}_t(x)  / \sigma(x)-1}^2}^{1/2}\\
&\leq \alpha +\e{\abs{\Pi^\bench(q_t^{\OSCI})-\Pi^\bench(q_\alpha^{\bench})}} + C_0 (\log t)^{\tau_0}t^{-\tau},
\end{align*}
by using \eqref{ass:esti:averagepoint}, and for some constant $C_0>0$ (for short, we shall still denote the constant by $C_0$ in the lines below although this constant should in principle be updated).  Now, we use that
\begin{align*}
\e{\abs{\Pi^\bench(q_t^{\OSCI})-\Pi^\bench(q_\alpha^{\bench})}} &\leq \e{\abs{\Pi^\bench(q_t^{\OSCI})-\Pi^\bench(q_\alpha^{\bench})}\ind{\abs{q_t^{\OSCI}-q_\alpha^{\bench}}\leq \eta_t}} + 1/t,
\end{align*}
by using Corollary~\ref{cor:rateconveasyselect} and letting $\eta_t:=C_0 [\rho_t(\beta)+(\log t)^{\nu_0} t^{-\nu}]$. By the mean value theorem, we have that, for $t$ large enough, when $\abs{q_t^{\OSCI}-q_\alpha^{\bench}}\leq\eta_t$,
\begin{align*}
\abs{\Pi^\bench(q_t^{\OSCI})-\Pi^\bench(q_\alpha^{\bench})}&\leq \eta_t \max_{q\in [q_\alpha^{\bench}-\eta_t, q_\alpha^{\bench}+\eta_t]}
\abs{(\Pi^\bench)'(q)} \\&= \eta_t \max_{q\in [q_\alpha^{\bench}/2, (B+q_\alpha^{\bench})/2]}\abs{\frac{f}{\phi}\paren{ \Phi^{-1}\paren{\frac{1-q}{2}}}}\leq C_0 \eta_t. \end{align*}
This gives \eqref{equcondconv} of point (i).

Now, let us prove point (ii). We relate $|\cC^{\OSCI}_t(x)|=2\hat{\sigma}_t(x) \Phi^{-1}\Big(\frac{q_t^{\OSCI}+1}{2}\Big) $ to $|\cC^*(x)|=2\sigma(x) \ol{F}^{-1}(\alpha/2)$ by using Corollary~\ref{cor:rateconveasyselect}: first note that applying the in-probability bound \eqref{equ-conviid-specialcases-easyselect}, we have for $t$ large enough, for all $M>0$,
\begin{align*}
\e{|\cC^{\OSCI}_t(x)|\wedge M }
&\leq 
2\e{\hat{\sigma}_t(x) \Phi^{-1}\Big(\frac{q_t^{\OSCI}+1}{2}\Big) \ind{\abs{q_t^{\OSCI}-q_\alpha^{\bench}}\leq \eta_t}} + M/t\\
&\leq 
2\e{\hat{\sigma}_t(x) \Phi^{-1}\Big(\frac{q_t^{\OSCI}+1}{2}\Big) \ind{q_t^{\OSCI}\in[q_\alpha^{\bench}/2, (1+q_\alpha^{\bench})/2]}} + M/t.
\end{align*}
Now, by the triangle inequality and the mean-value theorem,
\begin{align*}
&\e{\hat{\sigma}_t(x) \Phi^{-1}\Big(\frac{q_t^{\OSCI}+1}{2}\Big)\ind{q_t^{\OSCI}\in[q_\alpha^{\bench}/2, (1+q_\alpha^{\bench})/2]} } - \sigma(x) \ol{F}^{-1}(\alpha/2)\\
&\leq \e{(\hat{\sigma}_t(x)-\sigma(x)) \Phi^{-1}\Big(\frac{q_t^{\OSCI}+1}{2}\Big)\ind{q_t^{\OSCI}\in[q_\alpha^{\bench}/2, (1+q_\alpha^{\bench})/2]} } \\
&+\sigma(x) \e{\paren{\Phi^{-1}\Big(\frac{q_t^{\OSCI}+1}{2}\Big)- \Phi^{-1}\Big(\frac{q_\alpha^{\bench}+1}{2}\Big)}\ind{q_t^{\OSCI}\in[q_\alpha^{\bench}/2, (1+q_\alpha^{\bench})/2]} } \\
&\leq C_0 \paren{\e{\abs{\hat{\sigma}_t(x)-\sigma(x)} }+ \e{\abs{q_t^{\OSCI}-q_\alpha^{\bench}}}},
\end{align*}
because $\sigma(x)$ is uniformly upper-bounded by Assumption~\ref{ass:densityf:decreasingonly}. 
Combining the latter with the previous display, we obtain
\begin{align*}
\e{|\cC^{\OSCI}_t(x)|\wedge M}&\leq |\cC^*(x)| + M/t +  C_0 \e{\paren{\hat{\sigma}_t(x)-\sigma(x)}^2}^{1/2} + C_0\e{\paren{q_t^{\OSCI}-q_\alpha^{\bench}}^2}^{1/2}  \\
&\leq |\cC^*(x)| + M/t +  C_0 (\log t)^{\tau_0}t^{-\tau} + C_0 [(\log t)^{\tau_0} t^{1- \tau-\beta }]^{1/2},
\end{align*}
by using the in-expectation bound \eqref{equ:MSEconv-regressioneasy} of Corollary~\ref{cor:rateconveasyselect}.
This gives \eqref{equmimicoracleiidregress}.

Finally, we prove (iii). Consider any procedure $\cC'_t(x),t\geq 1,x\in \mathcal{X}$ with the conditional coverage property $\probp{Y_t\in \cC'_t(x)\mid X_t=x}\geq 1-\alpha - \eta'_t$. Since $\cC^*(x), x\in \mathcal{X}$ satisfies \eqref{optpredictionset}, we have 
 almost surely, for all $y\in \mathcal{Y}$,
$$( f_{Y|X=x}(y) - c^*(x))(\ind{y\in \cC^*(x)}-\ind{y\in \cC'_t(x)}) \geq 0.$$ 
Hence, by integrating with respect to the Lebesgue measure,
 \begin{align*}
      0&\leq \int_{\R} (f_{Y|X=x}(y) - c^*(x))(\ind{y\in \cC^*(x)}-\ind{y\in \cC'_t(x)}) dy\\
      &=\int_{\R} f_{Y|X=x}(y) \ind{y\in \cC^*(x)}dy-\int_{\R} f_{Y|X=x}(y) \ind{y\in \cC'_t(x)}) dy - c^*(x)(|\cC^*(x)|-|\cC'_t(x)|).
\end{align*}
      Now, since $\probp{Y_t\in \cC^*(x)\mid X_t=x}=1-\alpha$ (because the distribution of $Y_t$ given $X_t=x$ is continuous) and by assumption $\probp{Y_t\in \cC'_t(x)\mid X_t=x}\geq 1-\alpha - \eta'_t$, we integrate with respect to $\mathcal{F}_{t-1}$, to obtain
      \begin{align*}
      0&\leq \eta'_t - c^*(x)(|\cC^*(x)|-\e{|\cC'_t(x)|}).
      \end{align*}
Since $c^*(x)>0$ and $\sigma(x)$ is uniformly lower bounded, this gives
\begin{align}\label{equintermoptimal}
|\cC^*(x)| \leq \e{|\cC'_t(x)|}+C_0 \eta'_t,
  \end{align}
which means that, up to $C_0\eta'_t$, the prediction set $\cC^*(x)$ is pointwise better than $\cC_t(x)$. 
 Combining the latter with \eqref{equmimicoracleiidregress} concludes the proof.

\subsection{Proof of Corollary~\ref{cor:rateconvadaptselect}}\label{sec:proofadaptiveXorientedSelection}

First note that $\P\paren{\mu(X_t)\geq y_0}>0$ by Assumption~\ref{ass:densitymu}.
Here we have $\Pi_t(q)=\probp{V_t(X_t,Y_t)>q\mid \hat{\mu}_t(X_t)\geq y_0,\cF_{t-1}}$ with $V_t(x,y)$ as in \eqref{regressionVadaptive}
and with the benchmark function $\Pi^\bench$ as in \eqref{PizeroClassicalscore}.

Consider for $r_t:=\sqrt{C_1 (\log t)^{\nu_0} t^{-\nu}}$  
the probability
\begin{equation}
    \label{epsilonbaocase}
    \varepsilon_t(\mathcal{F}_{t-1}):=\probp{\abs{\hat{\mu}_t(X_t)-\mu(X_t)}> r_t \mid  \mathcal{F}_{t-1}}.
\end{equation}
By \eqref{ass:esti} and since by Markov's inequality $\varepsilon_t(\mathcal{F}_{t-1})\leq \e{\abs{\hat{\mu}_t(X_t)-\mu(X_t)}\mid \mathcal{F}_{t-1}}/r_t$, we have
\begin{align}
\probp{\varepsilon_t(\mathcal{F}_{t-1})> r_t}\leq \probp{\e{\abs{\hat{\mu}_t(X_t)-\mu(X_t)}\mid \mathcal{F}_{t-1}} > r^2_t} \leq \frac{3}{\pi^2t^3}\label{equ:concvarepilonbaocase}.
\end{align}
In addition, we have 
\begin{align}
\abs{\Pc{\hat{\mu}_t(X_t)\geq y_0}{\cF_{t-1}}-\P\paren{\mu(X_t)\geq y_0}}&\leq \Pc{\abs{\hat{\mu}_t(X_t)-\mu(X_t)}>r_t}{\cF_{t-1}}+\P\paren{\abs{\mu(X_t)-y_0}\leq r_t}\nonumber\\
&\leq \varepsilon_t(\mathcal{F}_{t-1})+ 2M_4 r_t \leq (1+2M_4) r_t\label{usefulboundbaocase},
\end{align}
for $t$ larger than some $T_0$ (only depending on the  neighborhood of Assumption~\ref{ass:densitymu}) and on the event where $\varepsilon_t(\mathcal{F}_{t-1})\leq r_t$.

Now, we would like to apply Theorem~\ref{th:rateconvfull}, which requires to establish \eqref{equ:Dtcondition} and to find sequences bounding $J(t)$ as in \eqref{equprobaJtsympa}.

\paragraph*{Bounding $D_t$}

Let us first show for $t\geq 1$,
\begin{equation}
    \label{equ:Dtconditionbaocase}
    \P( D_t> C_0 r^2_t)\leq \frac{6}{\pi^2t^3},
\end{equation}
for some constant $C_0>0$. 
From \eqref{usefulboundbaocase}, we have for $t\geq T_0$ (with an updated $T_0$), that $\Pc{\hat{\mu}_t(X_t)\geq y_0}{\cF_{t-1}}\geq r_0:=\P\paren{\mu(X_t)\geq y_0}/2>0$, provided that $\varepsilon_t(\mathcal{F}_{t-1})\leq r_t$. Hence, by applying \eqref{boundDtregression} in Proposition~\ref{prop:Dtiideasyselect}, we obtain
\begin{align*}
     D_t\leq  \|\Pi_t-\Pi^\bench\|_\infty 
       &\leq   M_3\:\E[M_1 |\hat{\mu}_t(X_t)-{\mu}(X_t)|+ (M_2\vee 1)|\hat{\sigma}_t(X_t)-{\sigma}(X_t)| \:|\: S_t(X_t)=1, \mathcal{F}_{t-1}]\\
       &\leq r_0^{-1} M_3\:\E[M_1 |\hat{\mu}_t(X_t)-{\mu}(X_t)|+ (M_2\vee 1)|\hat{\sigma}_t(X_t)-{\sigma}(X_t)| \:|\: \mathcal{F}_{t-1}],
    \end{align*}
for $t\geq T_0$ and provided that $\varepsilon_t(\mathcal{F}_{t-1})\leq r_t$. 
As a result, we obtain (by adjusting appropriately the constant $C_0>0$)
\begin{align*}
    \P( D_t> C_0 r^2_t) &\leq \probp{\E [|\hat{\mu}_t(X_t)-{\mu}(X_t)| + |\hat{\sigma}_t(X_t)-{\sigma}(X_t)|\mid \mathcal{F}_{t-1}] \geq C_1 (\log t)^{\nu_0} t^{-\nu}} + \probp{\varepsilon_t(\mathcal{F}_{t-1})>r_t}\\
    &\leq \frac{3}{\pi^2t^3}+\frac{3}{\pi^2t^3},
\end{align*}
by using \eqref{ass:esti} and \eqref{equ:concvarepilonbaocase}, which shows \eqref{equ:Dtconditionbaocase}.

\paragraph*{Bounding $J(t)$}

Let us now prove for $t\geq T_0$ (with an updated constant $T_0$) 
\begin{equation}
    \label{equprobaJtsympabaoselect}
   \P\big( J(t) \notin [\underline{J}(t),\bar{J}(t)]\big)\leq 1/t,
\end{equation}
for functions $\underline{J}(t)<\bar{J}(t)$ to be suitably chosen.
By applying Lemma~\ref{lem:JtconcentrMart} \eqref{equJtconcentrMart}-\eqref{equJtconcentrMartotherside}, \et{we have that 
$$
    \probp{\abs{J(t+1)-1 - \sum_{k=1}^t\probp{\hat{\mu}_k(X_k)\geq y_0 \mid \mathcal{F}_{k-1}}} \geq \paren{\frac{t\log(4t)}{2}}^{1/2}}\leq 1/(2t).
$$
This means that for all $t\geq T_0$ (for some $T_0\geq  1$),
$$
    \probp{\abs{\et{J(t)} - \sum_{k=\lfloor \sqrt{t}\rfloor}^t \probp{\hat{\mu}_k(X_k)\geq y_0 \mid \mathcal{F}_{k-1}}} \geq \paren{t\log(4t)}^{1/2}}\leq 1/(2t).
$$
In addition, we have by \eqref{equ:concvarepilonbaocase} and \eqref{usefulboundbaocase} that
\begin{align*}
    &\probp{\sum_{k=\lfloor \sqrt{t}\rfloor}^t \abs{\probp{\hat{\mu}_k(X_k)\geq y_0 \mid \mathcal{F}_{k-1}}- \probp{{\mu}(X_k)\geq y_0 }}> (1+2M_4) \sum_{k=\lfloor \sqrt{t}\rfloor}^t  r_k} \\
    &\leq \probp{\exists k\in [\lfloor \sqrt{t}\rfloor,t] \::\: \varepsilon_k(\mathcal{F}_{k-1})> r_k}\leq \sum_{k=\lfloor \sqrt{t}\rfloor}^t \frac{3}{\pi^2 k^3} \leq C_0/t^{1/2},
\end{align*}
by adjusting again the constant $T_0$. }
Note that 
$$
(1+2M_4) \sum_{k=T_0}^t  r_k = (1+2M_4)(C_1)^{1/2}\sum_{k=T_0}^t   C_1 (\log k)^{\nu_0/2} k^{-\nu/2} \leq C_0 (\log t)^{\nu_0/2} t^{1-\nu/2},
$$
for some constant well chosen constant $C_0>0$. By adjusting $C_0>0$, this gives \eqref{equprobaJtsympabaoselect} for 
\begin{align*}
    \underline{J}(t) &= \lfloor t\:\probp{{\mu}(X)\geq y_0 } - C_0(\log t)^{\nu_0/2} t^{1-\nu/2}\rfloor;\\
    \bar{J}(t) &=\lfloor t\:\probp{{\mu}(X)\geq y_0 } + C_0(\log t)^{\nu_0/2} t^{1-\nu/2}\rfloor,
\end{align*}
because $(\log t)^{\nu_0/2} t^{1-\nu/2} \gg \paren{t\log(4t)}^{1/2}$.

\paragraph*{Conclusion} 

We apply Theorem~\ref{th:rateconvfull} with Remark~\ref{rem:ratetypical} and $\lambda=1-\nu/2$.

\subsection{Proof of Corollary~\ref{cor:rateconveAR}}\label{proof:cor:rateconveAR}

Let us first derive the in-probability bound. For this, we apply Theorem~\ref{th:rateconvfull}, which is possible because $\Pi^\bench ( q)=1-q$ satisfies Assumption~\ref{ass:easy} with $q_\alpha^\bench=1-\alpha$. 

\paragraph*{Bounding $D_t$}

Recall \eqref{equDtT}, so that $D_t\leq \|\Pi_t-\Pi^\bench\|_\infty$. 
Since $Y_t$ conditionally on $X_t $ follows a Gaussian distribution with mean $\mu(X_t)=\scal[0]{\varphi}{X_t}$ and variance $1$, we apply Lemma~\ref{lem:forDT}, which leads to ($\hat{\sigma}_t(x)=\sigma(x)=1$,  $g=\phi$ the density of the standard Gaussian distribution):
\begin{align}
D_t &\leq \E[ \sup_{q\in [0,B]}\abs{\probp{V_t(X_t,Y_t) >q \:|\: X_t, \mathcal{F}_{t-1}} - \Pi^0(q)} \:|\:  S_t(X_t, q)=1,  \mathcal{F}_{t-1}] \nonumber\\
&\leq M_1  {\sigma}^{-1} |\hat{\mu}_t(X_t)-{\mu}(X_t)|,  \label{equboundDtAR}
\end{align}
with $M_1=\phi(0)=1/\sqrt{2\pi}$ (because $\sup_{x\in \R} \{|x|\phi(x)\}<\infty$).
Hence, \eqref{ass:esti} with $\nu_0=3/2$, $\nu=1/2$ provides \eqref{equ:Dtcondition} with $\bar D_t =O( (\log t)^{3/2} t^{-1/2})$ ($\eta_t=1/t$). 

\paragraph*{Bounding $J(t)$}

We have $J(t)
=1+ \sum_{k=1}^{t-1}\ind{S(X_k)=1}
=:K(X_1,\dots,X_{t-1})
$. The function $K(x_1,\dots,x_{t-1})$ satisfies the following separately
bounded condition: for $x_1,\dots,x_{t-1}\in \R^d$ and $x'_1,\dots,x'_{t-1}\in \R^d$ with $x'_k=x_k$ for all $k$ except for one $k_0\in [t-1]$, we have
$$
\abs{K(x_1,\dots,x_{t-1})-K(x'_1,\dots,x'_{t-1})}\leq 1.
$$
By Lemma~\ref{lem:concentrationJtAR}, denoting $r=\P(S(X_1)=1)>0$ and choosing $\bar{J}(t)=1+\lceil ((t-1)r+ \sqrt{M_0\log(2t)(t-1)}\rceil$ and $\underline{J}(t)=1+\lfloor ((t-1)r- \sqrt{M_0\log(2t)(t-1)}\rfloor$ we have \eqref{equprobaJtsympa} with $w(t)=t^{-1}$.

Applying Theorem~\ref{th:rateconvfull} with Remark~\ref{rem:ratetypical} and $\lambda=1/2$ thus gives the in-probability bound \eqref{equ:ARConstSelectProba}.

\et{Let us now prove the in-expectation bound under the no-selective inference assumption. We apply Theorem~\ref{th:genL2-first}, for which the selection indicator $S_t(X_t,q^{\OSCI}_{t})\equiv 1$ is deterministic so is independent of the past $\mathcal{F}_{t-1}$. Since $J(t)=t$ for all $t\geq 1$ in this case, the concentration terms vanishes (with $\bar{J}(t)=\underline{J}(t)=J(t)=t$ and $\omega(t)=0$).
Also, $\probp{J(k(t))< T_0}=0$ for $k(t)\equiv T_0$. Finally, \eqref{equboundDtAR} and \eqref{ass:esti:average} with $\tau_0=0$, $\tau=1/2$ provides $\E[D_t]\leq \widetilde{D}_t$ with $\widetilde{D}_t=O(  t^{-1/2})$. This gives $\sum_{k=1}^{t-1} \widetilde{D}_k = O(t^{1/2})$ and thus the desired bound follows.
}

\section{Proofs for informative selections}

In these proofs, we denote $\Pi_t$ for $\Pi_{\cC^{\OSCI}_t, S^\OSCI_t}$ and $q_t$ for $q^{\OSCI}_t$ for short.

\subsection{General lower bound for $J(t)$}

\et{
\begin{lemma}[Lower bounding $J(t)$ in the informative case]\label{lem:Jtinfo}
Consider the iid model for $(X_t,Y_t)_{t\geq 1}$. 
Consider the threshold $(q_t^{\OSCI})_{t\geq 1}$ \eqref{equACIquantile}, used with any informative selection rule of the form $S_t(x,q)=\ind{W_t(x)>q}$ for some $\mathcal{F}_{t-1}$-measurable function $W_t:x\in \mathcal{X}\mapsto W_t(x)\in \R$. Consider any (deterministic) function $W:x\in \mathcal{X}\mapsto W(x)\in \R$ and assume that the random variable $W(X)$ has a density with respect to the Lebesgue measure bounded uniformly by $C>0$. 
Assume also that $\P(W(X)>\bar q)>0$ for some $\bar q \in (0,B)$. 
Let for all $t\geq 1$ and some $r_t>0$ vanishing as $t$ tends to infinity, 
\begin{equation}
    \label{equ-epsiFtlemma}
    \varepsilon_t(\mathcal{F}_{t-1}) :=  \probp{|W_t(X_t)-W(X_t)|>r_t\mid \mathcal{F}_{t-1}},
\end{equation}
for which we assume that for $t\geq T_0$ (for some $T_0\geq 1$), we have  $ r_t\leq \probp{W(X)>\bar q  }/(2(C+1))$.
Then, for any integer sequence $m_t<t$ tending to infinity as $t$ tends to infinity, we have for all $t\geq T_0$, 
\begin{equation}
    \label{equprobaJtsympaclassif}
   \P\big( J(t) <\underline{J}(t), \cap_{k=m_t}^t \Omega_k,\cap_{k=m_t}^t \{\varepsilon_k(\mathcal{F}_{k-1})\leq  r_k\}\big)\leq 1/t,
\end{equation}
where $\Omega_t$ is such that 
$\Omega_t\subset \set{\forall k\in [m_t,t], q^{\OSCI}_k\leq \bar q}\in \mathcal{F}_{t-1}$
and 
\begin{align}
    \label{equJunderlem}
    \underline{J}(t)&:=0\vee \big\lfloor (t-m_t)\:\probp{W(X)>\bar q  }/2 - \sqrt{(t-1) \log(t-1)/2}\big\rfloor .
\end{align}
In addition, for $t\geq T_0$, on the events $\Omega_t$ and $\varepsilon_t(\mathcal{F}_{t-1})\leq r_t$, we have pointwise
\begin{equation}
    \label{equ:probarejectlb}
    \probp{W_t(X_t)>q_t\mid \mathcal{F}_{t-1}}\geq \probp{W(X)>\bar q  }/2 .
\end{equation}
\end{lemma}
}

\et{
\begin{proof}
 Denote $q_t^{\OSCI}$ by $q_t$ for short. By applying Lemma~\ref{lem:JtconcentrMart} \eqref{equJtconcentrMartotherside}, we have that for all $t\geq T_0$, since $m_t\leq t-1$,
\begin{align*}
    &\probp{J(t)-1 \leq \sum_{k=m_{t}}^{t-1}\probp{W_k(X_k)>q_k\mid \mathcal{F}_{k-1}}- x, \cap_{k=m_{t}}^{t} \Omega_k,\cap_{k=m_{t}}^{t} \{\varepsilon_k(\mathcal{F}_{k-1})\leq  r_k\}}\\&\leq e^{-2x^2/(t-1)}.
\end{align*}
In addition, by successively applying Lemma~\ref{lem:concentrespJt} \eqref{eqeintermlemmaJtconc}, and \eqref{equ-epsiFtlemma}, we have almost surely
\begin{align*}
\abs{\probp{W_t(X_t)>q_t\mid \mathcal{F}_{t-1}}- \probp{W(X_t)>q_t  \mid \mathcal{F}_{t-1}}}&\leq C r_t + \probp{|W_t(X_t)-W(X_t)|>r_t\mid \mathcal{F}_{t-1}} \\
&= C r_t + \varepsilon_t(\mathcal{F}_{t-1}).
\end{align*}
This means that, for all $t\geq T_0$, on the event $\Omega_t$ and $\varepsilon_t(\mathcal{F}_{t-1})\leq  r_t$, 
\begin{align}
\probp{W_t(X_t)>q_t\mid \mathcal{F}_{t-1}}&\geq  \probp{W(X_t)>q_t  \mid \mathcal{F}_{t-1}}- C r_t - \varepsilon_t(\mathcal{F}_{t-1})\label{equintermusefulforcorollary}\\
&\geq \probp{W(X)>\bar q  }- C r_t - \varepsilon_t(\mathcal{F}_{t-1})\nonumber\\
&\geq \probp{W(X)>\bar q  }- (C+1) r_t \geq \probp{W(X)>\bar q  }/2, 
\nonumber
\end{align}
for which we used in addition that $2(C+1) r_t\leq \probp{W(X)>\bar q  }$ in the last inequality. Combining this with the first display of the proof provides the result.
\end{proof}
}

\subsection{Proofs for informative regression cases}\label{sec:proofsinforegress}

\subsubsection{Proof of Corollary~\ref{cor:restrictedlength}}\label{proof:cor:restrictedlength}

\et{By Theorem~\ref{rem:coroonanEvent}, we have to identify $\omega(t)$ in \eqref{equprobaJtsympa} and $\bar{D}_t$ in \eqref{equ:Dtcondition}, which are computed in restriction to respective events $\Omega_t=\set{\forall k\in [T_0,t], q_k\leq \bar q}\in \mathcal{F}_{t-1}$. 
\paragraph*{Bounding $J(t)$}
First note that by Assumption~\ref{ass:supportsigma}, the density of $W(X)$ is upper-bounded by some $C_0>0$ and we have $\P(W(X)>\bar q)>0$ because $\bar q<B^*$. Hence, by applying Lemma~\ref{lem:Jtinfo}, we obtain that for $t\geq T_0$ (up to update $T_0$), we have 
\begin{equation*}
   \P\big( J(t) <\underline{J}(t), \Omega_t,\cap_{k=m_t}^t \{\varepsilon_k(\mathcal{F}_{k-1})\leq  r_k\}\big)\leq 1/t,
\end{equation*}
where $\underline{J}(t)=\lfloor t \:\probp{W(X)>\bar q  }/5\rfloor$, which is of order $t$ and where we choose $m_t=\lfloor t/2\rfloor$. 
Consider in the sequel $r_t=\sqrt{C'_1(\log t)^{\nu_0} t^{-\nu}}$ for $C'_1>0$ a large enough constant.
\paragraph*{Bounding $D_t$}
First note that by applying Lemma~\ref{lem:concentrespJt}, on the event $\Omega_t$ and provided that $\varepsilon_t(\mathcal{F}_{t-1})\leq r_t$, we have 
\begin{align*}
    \probp{W_t(X_t)>q_t\mid \mathcal{F}_{t-1}}&\geq  \probp{W(X)>\bar q} -C_0 r_t - \varepsilon_t(\mathcal{F}_{t-1}) \\
    &\geq \probp{W(X)>\bar q}/2>0.
\end{align*}
Hence, we can apply Proposition~\ref{prop:Dtiideasyselect} to get that for all $t\geq T_0$ (possibly by increasing $T_0$), 
\begin{align*}
       \P( D_t > r^2_t, \Omega_t,\cap_{k=m_t}^t \{\varepsilon_k(\mathcal{F}_{k-1})\leq  r_k\})\leq \frac{6}{\pi^2t^3}
\end{align*}
for a constant $C'_1>0$ chosen large enough.
\paragraph*{Conclusion}
We apply now Theorem~\ref{rem:coroonanEvent} (i) (or more precisely Proposition~\ref{rem:coroonanEvent0} on the event including $\cap_{k=m_t}^t \{\varepsilon_k(\mathcal{F}_{k-1})\leq  r_k\}$ and $\Omega_t$) to deduce that for $t\ge T_0$, we have with $k(t)=c't/2$;
  \begin{align*}
&\P(|\IER_t(\mathcal{R}^{\OSCI}) - \alpha|+ \abs{\qt_t^{\OSCI}-q^{\bench}_{\alpha} } \leq C_0 ((\log t)^{1/2} t^{1/2-\beta} + t^{\beta-1}+ r^2_t), \Omega_t,\cap_{k=m_t}^t \{\varepsilon_k(\mathcal{F}_{k-1})\leq  r_k\})
\nonumber\\
&\geq 1-\delta' - \P\big(\big(\cap_{k=m_t}^t \{\varepsilon_k(\mathcal{F}_{k-1})\leq  r_k\}\big)^c\big)-C_0/t ,
\end{align*}
up to increase the constant $C_0>0$. 
Let us now show that 
\begin{equation}
    \label{equ-interlength}
    \P\big(\big(\cap_{k=m_t}^t \{\varepsilon_k(\mathcal{F}_{k-1})\leq  r_k\}\big)^c\big)\leq \frac{5}{\pi^2t^2},
\end{equation}
(up to increase the constant $C_0$), which leads to \eqref{equ-convlength}, where we recall that $m_t=\lfloor t/2\rfloor$.
First, the quantity $\varepsilon_t(\mathcal{F}_{t-1}) $ in \eqref{equ-epsiFtlemma} can be upper-bounded as follows: we have for $r_t>0$ small enough whenever $t\geq T_0$, that 
\begin{align*}
    \varepsilon_t(\mathcal{F}_{t-1}) &= \probp{|W_t(X_t)-W(X_t)|>r_t\mid \mathcal{F}_{t-1}}\\
    &\leq \probp{|W_t(X_t)-W(X_t)|>r_t, \hat{\sigma}_t(X_t)\geq \sigma_{\min}/2 \mid \mathcal{F}_{t-1}} + \probp{\hat{\sigma}_t(X_t)< \sigma_{\min}/2\mid \mathcal{F}_{t-1}}\\
    &\leq \probp{2|\Phi (\lambda_0/\hat{\sigma}_t(X_t))- \Phi (\lambda_0/\sigma(X_t))|>r_t, \hat{\sigma}_t(X_t)\geq \sigma_{\min}/2\mid \mathcal{F}_{t-1}} \\
    &+ \probp{\hat{\sigma}_t(X_t)< \sigma_{\min}/2\mid \mathcal{F}_{t-1}}\\
    &\leq \probp{ (4/\sigma^2_{\min})(\lambda_0/\sqrt{2\pi}) |\hat{\sigma}_t(X_t))-\sigma(X_t)|>r_t\mid \mathcal{F}_{t-1}} + \probp{\hat{\sigma}_t(X_t)< \sigma_{\min}/2\mid \mathcal{F}_{t-1}}\\
    &\leq 2 \probp{ (4/\sigma^2_{\min})(\lambda_0/\sqrt{2\pi}) |\hat{\sigma}_t(X_t))-\sigma(X_t)|>r_t\mid \mathcal{F}_{t-1}}\\
    &\leq \frac{(8/\sigma^2_{\min})(\lambda_0/\sqrt{2\pi})}{r_t} \e{  |\hat{\sigma}_t(X_t))-\sigma(X_t)|\mid \mathcal{F}_{t-1}},
\end{align*}
where we apply the Markov inequality and the fact that $\Phi$ is $(1/\sqrt{2\pi})$-Lipschitz. 
Considering 
 $r_t$ with $C_1'\geq C_1(8/\sigma^2_{\min} (\lambda_0/\sqrt{2\pi})$, this gives   
\begin{align*}
    &\P(\varepsilon_t(\mathcal{F}_{t-1})\geq r_t) \leq \P(\e{  |\hat{\sigma}_t(X_t))-\sigma(X_t)|\mid \mathcal{F}_{t-1}} \geq C_1(\log t)^{\nu_0} t^{-\nu})
    \leq \frac{3}{\pi^2t^3},
\end{align*}
and thus 
\begin{align}\label{equ:intermunionboundnotenough}
   \P\big(\big(\cap_{k=m_t}^t \{\varepsilon_k(\mathcal{F}_{k-1})\leq  r_k\}\big)^c\big)&\leq \frac{3}{\pi^2} \sum_{k=m_t}^t k^{-3} \sim \frac{9}{2\pi^2}  t^{-2}
\end{align}
by applying Lemma~\ref{lemsuminte}. This leads to \eqref{equ-interlength}.
}

\subsubsection{Proof of Corollary~\ref{cor:lb}}\label{proof:cor:lb}

\et{The proof is analogous to \S~\ref{proof:cor:restrictedlength} by adjusting the upper bound of
\begin{align*}
    \varepsilon_t(\mathcal{F}_{t-1}) &= \probp{|W_t(X_t)-W(X_t)|>r_t\mid \mathcal{F}_{t-1}}
\end{align*}
to the specific $W_t(x)$ \eqref{equ:Wtlb} and $W(x)$ \eqref{equ:Wlb} at hand here. This upper bound can be obtained as follows:
\begin{align*}
    \varepsilon_t(\mathcal{F}_{t-1}) &\leq \probp{\bigg|\Phi\bigg(\frac{\hat{\mu}_t(X_t)-y_0}{\sigma_0}\bigg) -\Phi\bigg(\frac{{\mu}(X_t)-y_0}{\sigma_0}\bigg) \bigg|>r_t\mid \mathcal{F}_{t-1}} \\
    &\leq \probp{ [1/(\sigma_{0}\sqrt{2\pi})] |\hat{\mu}_t(X_t))-\mu(X_t)|>r_t\mid \mathcal{F}_{t-1}} \\
    &\leq \frac{[1/(\sigma_{0}\sqrt{2\pi})]}{r_t} \e{  |\hat{\mu}_t(X_t))-\mu(X_t)|\mid \mathcal{F}_{t-1}}.
    \end{align*}
    Hence, by considering $r_t=\sqrt{C'_1(\log t)^{\nu_0} t^{-\nu}}$ with $C_1'\geq C_1[1/(\sigma_{0}\sqrt{2\pi})]$, we obtain
    by \eqref{ass:esti}
    \begin{equation}
    \label{equ-interlb}
    \P\big(\big(\cap_{k=m_t}^t \{\varepsilon_k(\mathcal{F}_{k-1})\leq  r_k\}\big)^c\big)\leq \frac{5}{\pi^2t^2}. 
\end{equation}
This proves \eqref{equ-convlb}.
}

\subsection{Proofs for selective classification results in \S~\ref{sec:selectiveclassifiid}}

\subsubsection{Proof of Theorem~\ref{th:rate:classif}}\label{proof:th:rate:classif}

\et{We apply Theorem~\ref{rem:coroonanEvent} (case (ii)) with an appropriate choice of $\omega(t)$ in \eqref{equprobaJtsympa} and $\bar{D}_t$ in \eqref{equ:Dtcondition}, recalling that the latter two conditions are computed in restriction to  $\Omega_t=\set{\forall k\in [T_0,t], q_k\leq \bar q}\in \mathcal{F}_{t-1}$.}
\et{Also, by taking $r_t:=\sqrt{C_1(\log t)^{\nu_0} t^{-\nu}}$, we have 
\begin{align*}
    \varepsilon_t(\mathcal{F}_{t-1}) &:=  \probp{|W_t(X_t)-W(X_t)|>r_t\mid \mathcal{F}_{t-1}}\\
    &\leq \probp{\max_{y\in \range{K}}|\hat{\pi}_t(y|X_t)-\pi(y|X_t)|> r_t \mid  \mathcal{F}_{t-1}}
\end{align*}
and thus by Markov's inequality, for $t\geq T_0$ (for some $T_0\geq 1$),
\begin{align}
\probp{\varepsilon_t(\mathcal{F}_{t-1})> r_t}\leq \probp{\e{\max_{y\in \range{K}}|\hat{\pi}_t(y|X_t)-{\pi}(y|X_t)|\mid \mathcal{F}_{t-1}} > r^2_t} \leq \frac{3}{\pi^2t^3}\label{equ:concvarepilonclassif},
\end{align}
where we used \eqref{ass:esticlassif} in the last inequality.
}

\paragraph*{Bounding $J(t)$}
\et{
We have $\P(W(X)>\bar q)>0$ (Assumption~\ref{ass:support} and $\bar q<B^*$) and the density of $W(X)$ is continuous on a compact and thus is bounded. Lemma~\ref{lem:Jtinfo} thus gives for $t\geq T_0$ (for some $T_0\geq 1$),
\begin{equation*}
   \P\big( J(t) <\underline{J}(t), \Omega_t, \cap_{k=m_t}^t \{\varepsilon_k(\mathcal{F}_{k-1})\leq  r_k\}\big)\leq 1/t,
\end{equation*}
 with $\underline{J}(t)=\lfloor t \:\probp{W(X)>\bar q  }/5\rfloor$, which is of order $t$ and $m_t=\lfloor t/2\rfloor$. Hence, with the help of \eqref{equ:concvarepilonclassif}, we obtain for $t\geq T_0$ (for some $T_0\geq 1$),
 \begin{equation*}
   \P\big( J(t) <\underline{J}(t), \Omega_t\big)\leq 2/t.
\end{equation*}
}

\paragraph*{Bounding $D_t$}

Let us show that for all $t\geq T_0$ (with an updated constant $T_0$),
\begin{equation}
    \label{boundDtclassif}
    \probp{D_t > C_0 r_t, \Omega_t} \leq  \frac{3}{\pi^2t^3}.
\end{equation}
For this, we have for all $q\in [0,\bar q]$,  
\begin{align*}
    \Pi_t(q)&= 1-\e{  \pi(\hat{Y}_t|X_t) \mid W_t(X_t)>q, \mathcal{F}_{t-1}},
\end{align*}
where we recall that $\{\hat{Y}_t\}=\arg\max \hat{\pi}_t(\cdot|X_t)$. 
Since $\pi(\hat{Y}_t|X_t) \leq \max_{y\in [K]}\pi(y|X_t)=W(X_t)$ by definition, we have on the event where $\varepsilon_t(\mathcal{F}_{t-1})\leq r_t$,
\begin{align*}
    1-\Pi_t(q) 
    &=\e{  \pi(\hat{Y}_t|X_t) \mid W_t(X_t)>q, \mathcal{F}_{t-1}}\\
    &\leq  \e{ W(X_t) \mid W_t(X_t)>q, \mathcal{F}_{t-1}}\\
    &=\e{ W(X_t) \frac{\ind{W_t(X_t)>q}}{\probp{W_t(X_t)>q\mid \mathcal{F}_{t-1}}}\mid \mathcal{F}_{t-1}} \\
    &\leq \e{ W(X_t) \frac{\ind{W(X_t)>q - r_t}}{\probp{W(X_t)>q+ r_t}-r_t}\mid \mathcal{F}_{t-1}} + r_t\\
    &\leq \e{ W(X_t) \frac{\ind{W(X_t)>q - r_t}}{\probp{W(X_t)>q- r_t}-C_0r_t}} + r_t,
\end{align*}
because since the density of $W(X_t)$ is bounded we have $\probp{q-r_t<W(X_t)\leq q + r_t}\leq C_0 r_t$, 
and by assuming that $C_0,T_0$ are large enough so that $r_t\leq (1-\bar q) /2$ and $C_0r_t<\probp{W(X_t)>\bar q}/2$. 
Since $1/(1-u)\leq 1+2u$ for all $u\in (0,1/2]$, it entails that on the event where $\varepsilon_t(\mathcal{F}_{t-1})\leq r_t$,
\begin{align*}
    1-\Pi_t(q)     &\leq \e{ W(X_t) \frac{\ind{W(X_t)>q - r_t}}{\probp{W(X_t)>q- r_t}}\paren{1+2\frac{C_0r_t}{\probp{W(X_t)>q- r_t}}}} + r_t\\
    &\leq \paren{1+2\frac{C_0r_t}{\probp{W(X_t)>\bar q}}}\e{ W(X_t) \frac{\ind{W(X_t)>q - r_t}}{\probp{W(X_t)>q- r_t}}} + \varepsilon_t\\
    &=\paren{1+2\frac{C_0r_t}{\probp{W(X_t)>\bar q}}}(1-\Pi^\bench(q-r_t)) + r_t.
\end{align*}
Since $\sup_{q\in [0,\bar q]}\abs{\Pi^\bench(q-r_t)-\Pi^\bench(q+r_t)} \leq C_0 r_t$ by Lemma~\ref{lem:decreasing} (by upper bounding the negative term in the derivative by $0$), we have finally (still by adjusting the constant $C_0$), 
$$
\sup_{q\in [0,\bar q]} (\Pi^\bench(q)-\Pi_t(q)) \ind{\varepsilon_t(\mathcal{F}_{t-1})\leq r_t}\leq C_0 r_t
$$
Similarly, we obtain $
\sup_{q\in [0,\bar q]} (\Pi_t(q)-\Pi^\bench(q)) \ind{\varepsilon_t(\mathcal{F}_{t-1})\leq r_t}\leq C_0 r_t
$
which gives finally, for $t\geq T_0$,  on the event $\Omega_t$
\begin{equation}
    \label{equusefulforoptim}
  |\Pi_t(q_t)-\Pi^\bench(q_t)| \ind{q_t\in [0,B]} \ind{\varepsilon_t(\mathcal{F}_{t-1})\leq r_t} \leq C_0 r_t,
\end{equation}
which gives \eqref{boundDtclassif} with \eqref{equ:concvarepilonclassif}.
 
\paragraph*{Conclusion}

We apply Theorem~\ref{rem:coroonanEvent} (case (i)), to obtain \eqref{equ-convclassif}.

\subsubsection{Optimality of the oracle procedure}

\begin{proposition}\label{prop:optimalclassif}
Consider the iid classification model \eqref{model:classif} satisfying Assumption~\ref{ass:support} with $W(x):=\max_{y\in \range{K}} \P(Y= y\mid X=x)$ and support values  $1/K\leq \kappa\leq B^*\leq 1$, and an error level $\alpha\in (0,1)$ such that $1-B^*<\alpha<1-\e{W(X)}$. Then any selective classification procedure $(S_t(X_t),\tilde{Y}_t(X_t))_{t\geq 1}$ (for which each function $S_t(\cdot)$ and $\tilde{Y}_t(\cdot)$ is $\mathcal{F}_{t-1}$-measurable) is dominated by the oracle procedure $(S^*(X_t),Y^*(X_t))_{t\geq 1}$ given by \eqref{equ:classiforaclepredict}-\eqref{equ:classiforacleselect} in the following sense:  if for all $t\geq T_0$, $\P(Y_t\neq \tilde{Y}_t(X_t)\mid S_t(X_t)=1)\leq \alpha+z$, for $z\geq 0$, then for all $t\geq T_0$
    \begin{equation}\label{eq-optimalpower-class}
        \E\big[ \ind{Y_t= \tilde{Y}_t(X_t)}S_t(X_t)\big]\leq \E\big[ \ind{Y_t=Y^*(X_t)}S^*(X_t)\big]+z q^\bench_\alpha/(1-\alpha-q^\bench_\alpha).
    \end{equation}
\end{proposition}
For $z=0$, this result is (essentially) proved in Theorem~1 of \cite{zhao2023controlling}. Hence, Proposition~\ref{prop:optimalclassif} is an online extension, accommodating the case where $z>0$.

\begin{proof}
First note that by \eqref{equqalphainfo}, we have 
\begin{align*}
    \P(Y_t\neq {Y}^*_t(X_t)\mid S^*(X_t)=1)&=\e{1-W(X_t)\mid S^*(X_t)=1}\\
    &= 1- \e{W(X_t)\mid W(X_t)>q^\bench_\alpha} = \alpha
\end{align*}
and that $q^\bench_\alpha<1-\alpha$.
Also, by denoting 
\begin{equation}
    \label{equfSC}f:x\in (0,1)\mapsto f(x)=(x-\alpha)/(1-x),
\end{equation} which is continuous and increasing, we have
\begin{align*}
S^*(X_t)&=\ind{1-W(X_t)<1-q^\bench_\alpha} =\ind{f(1-W(X_t))<c^*}\\
&=\ind{1-W(X_t)-\alpha - c^* W(X_t)<0},
\end{align*}
for $c^*:=f(1-q^\bench_\alpha)>0$. 
Now, let us consider any selective classification procedure $(S_t(X_t),\tilde{Y}_t(X_t))_{t\geq 1}$, for which each function $S_t(\cdot)$ and $\tilde{Y}_t(\cdot)$ are $\mathcal{F}_{t-1}$-measurable.
Hence, almost surely, we have
$$
(1-W(X_t)-\alpha - c^* W(X_t))(S^*(X_t)-S_t(X_t))\leq 0.
$$
Rearranging this inequality, we obtain
\begin{align*}
&c^* (W(X_t) S^*(X_t)- W(X_t)S_t(X_t))\\
&\geq (1-W(X_t)-\alpha)S^*(X_t) - (1-W(X_t)-\alpha)S_t(X_t),
\end{align*}
and thus because $W(X_t)\geq \pi(\tilde{Y}_t(X_t)|X_t)$, 
\begin{align*}
&c^* (W(X_t) S^*(X_t)- \pi(\tilde{Y}_t(X_t)|X_t)) S_t(X_t))\\
&\geq (1-W(X_t)-\alpha)S^*(X_t) - (1-\pi(\tilde{Y}_t(X_t)|X_t)-\alpha)S_t(X_t).
\end{align*}
By integrating the latter in $X_t$ and $\mathcal{F}_{t-1}$, and by assuming $\P(Y_t\neq \tilde{Y}_t(X_t)\mid S_t(X_t)=1)\leq \alpha+z$, we obtain
\begin{align*}
&c^* (\E(W(X_t) S^*(X_t))- \E(\pi(\tilde{Y}_t(X_t)|X_t) S_t(X_t))\\
&\geq \E((1-W(X_t)-\alpha)S^*(X_t)) - \E((1-\pi(\tilde{Y}_t(X_t)|X_t)-\alpha)S_t(X_t)),
\end{align*}
that is,
\begin{align*}
c^* (\E(\ind{Y_t= {Y}^*_t(X_t)} S^*(X_t) )- \E(\ind{Y_t= \tilde{Y}_t(X_t)} S_t(X_t)))\geq 0  - z \P(S_t(X_t)=1)\geq -z,
\end{align*}
because
$
\P(Y_t\neq \tilde{Y}_t(X_t)\mid X_t, \mathcal{F}_{t-1}) = 1- \pi(\tilde{Y}_t(X_t)|X_t)
$
and  
$
\P(Y_t= {Y}^*_t(X_t)\mid X_t) =  W(X_t).
$
This gives \eqref{eq-optimalpower-class}.
\end{proof}

\subsubsection{Proof of Corollary~\ref{optimalityclassifiid}}\label{proof:optimalityclassifiid}

To prove Corollary~\ref{optimalityclassifiid}, we use Theorem~\ref{th:rate:classif}. 
Consider any $\delta'>\delta$, $\bar q=\bar q(\delta')\in [0,B^*)$ and $T_0=T_0(\delta')\geq 1$ such that $\P(\Omega)\geq 1-\delta'$ for the event $\Omega=\cap_{t\geq T_0}\Omega_t$ with $\Omega_t=\set{\forall k\in [T_0,t], q_k\leq \bar q}$. 

From Theorem~\ref{th:rate:classif} \eqref{equ-convclassif} (and its proof) recall that we have for the rate $\rho_t=(\log t)^{1/2} t^{1/2-\beta} + t^{\beta-1} $, 
\begin{equation}\label{equintermomegat}
\P(|\IER_t(\mathcal{R}^{\OSCI}) - \alpha|+ \abs{\qt_t^{\OSCI}-q^{\bench}_{\alpha}} \leq C_0 \eta_{t})\mid \Omega_t)\geq 1- C_0/t,
\end{equation}
where we denoted $\eta_t=\rho_t+(\log t)^{\nu_0/2} t^{-\nu/2}$ for short. 
Let us now prove (i). 
\et{
Recall that by Lemma~\ref{lem:concentrespJt} \eqref{eqeintermlemmaJtconc}, we have for $t\geq T_0$ (for some $T_0\geq 1$), 
\begin{align}\label{eqintermproofforgotcond}
\abs{\probp{W_t(X_t)>q_t\mid \mathcal{F}_{t-1}}- \probp{W(X_t)>q_t  \mid \mathcal{F}_{t-1}}}&\leq C r_t + \varepsilon_t(\mathcal{F}_{t-1}),
\end{align}
for which we recall that  $r_t=\sqrt{C_1(\log t)^{\nu_0} t^{-\nu}}$ and $\probp{\varepsilon_t(\mathcal{F}_{t-1})> r_t}\leq \frac{3}{\pi^2t^3}$ by \eqref{equ:concvarepilonclassif} and \eqref{ass:esticlassif}. Note that the latter holds even outside the event $\Omega_t$. Now, since $\Omega_t$ is $\mathcal{F}_{t-1}$-measurable, we have 
 \begin{align}
    &\P(Y_t\neq \hat{Y}_t(X_t),\Omega_t\mid S^{\OSCI}_t(X_t)=1)\nonumber\\
    &= \e{ \ind{\Omega_t}\frac{\P(S^{\OSCI}_t(X_t)=1\mid \mathcal{F}_{t-1})}{\P(S^{\OSCI}_t(X_t)=1)}\P(Y_t\neq \hat{Y}_t(X_t) \mid S^{\OSCI}_t(X_t)=1, \mathcal{F}_{t-1})} \nonumber\\
    &= \e{ \IER_t(\mathcal{R}^{\OSCI})\ind{\Omega_t} \frac{\P(W_t(X_t)>q_t\mid \mathcal{F}_{t-1})}{\P(W_t(X_t)>q_t)}}\nonumber\\
    &\leq \e{ \IER_t(\mathcal{R}^{\OSCI})\ind{\Omega_t} \frac{\P(W_t(X_t)>q_t\mid \mathcal{F}_{t-1})}{\P(W_t(X_t)>q_t)} \ind{\abs{\qt_t^{\OSCI}-q^{\bench}_{\alpha}} \leq C_0 \eta_{t}} \ind{\varepsilon_t(\mathcal{F}_{t-1})\leq r_t}}\nonumber\\ 
    &\:\:\:\:+ \frac{\frac{3}{\pi^2t^3} + C_0/t}{\P(W_t(X_t)>q_t)}\label{equintermforcond}
    \end{align}
    Now observe that by \eqref{eqintermproofforgotcond} on the event $\varepsilon_t(\mathcal{F}_{t-1})\leq r_t$ and $\abs{\qt_t^{\OSCI}-q^{\bench}_{\alpha}} \leq C_0 \eta_{t}$, we have for $t$ large enough,
    \begin{align}\label{equboundcrucial}
    |\P(W_t(X_t)>q_t\mid \mathcal{F}_{t-1})-C_2|\leq C_0 \eta_{t},
    \end{align}
     for $C_{2}=\P( W(X)>q^{\bench}_{\alpha})\in (0,1)$, because $r_t$ is of order smaller than $\eta_{t}$, because the density of $W(X)$ is upper bounded and up to update $C_0$ accordingly. This in particular implies  that 
\begin{align}
    \P(W_t(X_t)>q_t)&\geq \E[\ind{\varepsilon_t(\mathcal{F}_{t-1})\leq r_t}\ind{\abs{\qt_t^{\OSCI}-q^{\bench}_{\alpha}} \leq C_0 \eta_{t}}\P(W_t(X_t)>q_t \:|\: \mathcal{F}_{t-1})]  \nonumber\\
    &\geq (C_2 - C_0\eta_t)(1 - 3/(\pi^2t^3) - C_0/t - \delta') \geq (C_2 - C_0\eta_t) (1-\delta'- C_0/t),\label{eqlowerboundprobareject}
\end{align}
     by adjusting $C_0$ and for $t$ large enough.
Relations \eqref{equboundcrucial} and \eqref{eqlowerboundprobareject} imply that on the events $\varepsilon_t(\mathcal{F}_{t-1})\leq r_t$ and $\abs{\qt_t^{\OSCI}-q^{\bench}_{\alpha}} \leq C_0 \eta_{t}$, we have
\begin{align}\label{eqneededbound}
    \frac{\P(W_t(X_t)>q_t\mid \mathcal{F}_{t-1})}{\P(W_t(X_t)>q_t)}\leq  \frac{1}{1-\delta'- C_0/t}\frac{C_2+C_0 \eta_{t}}{C_2- C_0\eta_t}\leq \frac{1+C_0\eta_t}{1-\delta'- C_0/t},
\end{align}
up to change $C_0$.
Now, coming back to \eqref{equintermforcond} and using \eqref{eqlowerboundprobareject} and \eqref{eqneededbound} above, we have
     \begin{align}
    &\P(Y_t\neq \hat{Y}_t(X_t)\mid S^{\OSCI}_t(X_t)=1)\nonumber\\
    & = \P(Y_t\neq \hat{Y}_t(X_t),\Omega_t\mid S^{\OSCI}_t(X_t)=1) +  \P(Y_t\neq \hat{Y}_t(X_t),(\Omega_t)^c\mid S^{\OSCI}_t(X_t)=1)\nonumber\\
&\leq \frac{1+C_0\eta_t}{1-\delta'- C_0/t}\E( \IER_t(\mathcal{R}^{\OSCI})\ind{\Omega_t})+ \frac{2}{C_2(1-\delta')}(3/(\pi^2t^3) + C_0/t)+  \frac{2}{C_2(1-\delta')}\P((\Omega_t)^c)\label{equplane},
    \end{align} 
    for a $t$ large enough. 
   Since by \eqref{equintermomegat} we have
   $$
   \E( \IER_t(\mathcal{R}^{\OSCI})\ind{\Omega_t}) \leq \alpha + C_0 \eta_{t} + C_0/t,
   $$
   and $\P(\Omega^c)\leq \delta'$ this finally implies (i) up to adjust again the constants.
 }
 
Let us now prove (ii). \et{We have 
\begin{align*}
&\E\big[ \ind{Y_t= \hat{Y}_t(X_t)}S^{\OSCI}_t(X_t) \ind{\Omega_t}\big] \\
&=\e{\ind{\Omega_t} \probp{S^{\OSCI}_t(X_t)=1 \mid \mathcal{F}_{t-1}}\probp{Y_t= \hat{Y}_t(X_t) \mid S^{\OSCI}_t(X_t)=1, \mathcal{F}_{t-1}}}\\
&=\e{\ind{\Omega_t} \probp{S^{\OSCI}_t(X_t)=1 \mid \mathcal{F}_{t-1}}(1-\IER_t(\mathcal{R}^{\OSCI}))}\\
&\geq \e{\ind{\Omega_t\cap \mathcal{G}_t} \probp{S^{\OSCI}_t(X_t)=1 \mid \mathcal{F}_{t-1}}(1-\IER_t(\mathcal{R}^{\OSCI}))}\\
&\geq (1 - \alpha - C_0 \eta_{t})\probp{\Omega_t, \mathcal{G}_t,S^{\OSCI}_t(X_t)=1},
  \end{align*}
  by considering the `good' event $\mathcal{G}_t := \{1-\IER_t(\mathcal{R}^{\OSCI})\geq 1 - \alpha - C_0 \eta_{t}\}$. 
Now, by using \eqref{equintermomegat} we have $\P(\Omega_t\cap (\mathcal{G}_t)^c)\leq \P(\Omega_t)C_0/t$, so that the latter is at least
\begin{align*}
& \P(S^{\OSCI}_t(X_t)=1, \Omega_t)(1-\alpha  -C_0\eta_t)- C_0 \P(\Omega_t)/t\\
&\geq \P(S^{\OSCI}_t(X_t)=1, \Omega_t) \probp{Y_t= {Y}^*(X_t)\mid S^*(X_t)=1} -C_0\eta_t \\
&= \P(S^{\OSCI}_t(X_t)=1, \Omega_t) \probp{Y_t= {Y}^*(X_t)\mid S^*(X_t)=1,\Omega_t} -C_0\eta_t,
  \end{align*}
  by \eqref{equqalphainfo} and because $\Omega_t$ is independent of $(X_t,Y_t)$.}

We now study two cases. First, if $\P(S^*(X_t)=1, \Omega_t)-\P(S^{\OSCI}_t(X_t)=1, \Omega_t)\geq 0$, then 
\begin{align*}
&0\geq\\
&\brac{\P(S^{\OSCI}_t(X_t)=1, \Omega_t) -\probp{S^*(X_t)=1, \Omega_t}}\brac{\probp{S^*(X_t)=1, \Omega_t}-\probp{Y_t= {Y}^*(X_t),S^*(X_t)=1,\Omega_t}}\\
&=\probp{Y_t= {Y}^*(X_t),S^*(X_t)=1,\Omega_t}\probp{S^*(X_t)=1, \Omega_t}-\probp{S^*(X_t)=1, \Omega_t}^2\\
&+\probp{S^{\OSCI}_t(X_t)=1, \Omega_t}\probp{S^*(X_t)=1, \Omega_t}\\
&-\probp{S^{\OSCI}_t(X_t)=1, \Omega_t}\probp{Y_t= {Y}^*(X_t),S^*(X_t)=1,\Omega_t}.\\
\end{align*}
Dividing by $\probp{S^*(X_t)=1, \Omega_t}$, and putting the last term in the right hand side to the left hand side we get,
\begin{align*}
    &\probp{S^{\OSCI}_t(X_t)=1, \Omega_t}\probp{Y_t= {Y}^*(X_t)\mid S^*(X_t)=1,\Omega_t}\geq\\
    &\probp{Y_t= {Y}^*(X_t),S^*(X_t)=1,\Omega_t}
    -\et{(\probp{S^*(X_t)=1, \Omega_t} -\probp{S^{\OSCI}_t(X_t)=1, \Omega_t})}.
\end{align*}
Second, if $\P(S^*(X_t)=1, \Omega_t)-\P(S^{\OSCI}_t(X_t)=1, \Omega_t)< 0$, then
\begin{align*}
    &\P(S^{\OSCI}_t(X_t)=1, \Omega_t) \probp{Y_t= {Y}^*(X_t)\mid S^*(X_t)=1,\Omega_t}\\
    &>\P(S^*(X_t)=1, \Omega_t) \probp{Y_t= {Y}^*(X_t)\mid S^*(X_t)=1,\Omega_t}\\
    &=\probp{Y_t= {Y}^*(X_t),S^*(X_t)=1,\Omega_t}.
\end{align*}
Those two inequalities finally lead to,
\begin{align*}
\E\big[ \ind{Y_t= \hat{Y}_t(X_t)}S^{\OSCI}_t(X_t) \ind{ \Omega_t}\big] 
&\geq  \e{\ind{Y_t= {Y}^*(X_t)} S^*(X_t) \ind{ \Omega_t}}-C_0 \eta_t\\
&-\paren{\P(S^*(X_t)=1, \Omega_t)-\P(S^{\OSCI}_t(X_t)=1, \Omega_t)}\vee 0.
  \end{align*}

\et{Now, similarly to the case (i), we have $\P(S^*(X_t)=1, \Omega_t)-\P(S^{\OSCI}_t(X_t)=1, \Omega_t)\leq C_0 \eta_t + C_0 r_t$ and thus }
\begin{align*}
 \e{\ind{Y_t= {Y}^*(X_t)} S^*(X_t)\ind{ \Omega_t}}
&\leq \E\big[ \ind{Y_t= \hat{Y}_t(X_t)}S^{\OSCI}_t(X_t) \ind{ \Omega_t}\big] +  C_0\eta_t + C_0r_t,
  \end{align*}
  which provides by independence between $\Omega_t$ and $(X_t,Y_t)$, 
\et{
\begin{align*}
 \e{\ind{Y_t= {Y}^*(X_t)} S^*(X_t)}
&\leq (\E\big[ \ind{Y_t= \hat{Y}_t(X_t)}S^{\OSCI}_t(X_t)\big] +  C_0\eta_t + C_0r_t)/\P(\Omega_t)\\
&\leq (1-\delta')^{-1}(\E\big[ \ind{Y_t= \hat{Y}_t(X_t)}S^{\OSCI}_t(X_t)\big] +  C_0\eta_t + C_0r_t),
  \end{align*}
  by using $\P(\Omega_t^c)\leq \delta'$. 
}
  This gives (ii).

  Finally, we prove (iii). 
By Proposition~\ref{prop:optimalclassif}, for any $(S_t(X_t),\tilde{Y}_t(X_t))_{t\geq 1}$ such that for all $t\geq T_0$, $\P(Y_t\neq \tilde{Y}_t(X_t)\mid S_t(X_t)=1)\leq \alpha+C\delta'+\upsilon_t$, we have
\begin{align*}
        \E\big[ \ind{Y_t= \tilde{Y}_t(X_t)}S_t(X_t)\big]
        &\leq \e{W(X_t)\ind{ W(X_t)>q_{\alpha}^\bench}} +  q^\bench_\alpha (C\delta'+\upsilon_t)/(1-\alpha-q^\bench_\alpha).
  \end{align*}
\et{This finally gives 
$$
 \E\big[ \ind{Y_t= \tilde{Y}_t(X_t)}S_t(X_t)\big]\leq  (1-\delta')^{-1}(\E\big[ \ind{Y_t= \hat{Y}_t(X_t)}S^{\OSCI}_t(X_t)\big] +  C_0\eta_t + C_0r_t) + C (\delta'+\upsilon_t),
$$
where $C$ does not depend on $\delta'$.}
This provides \eqref{eq-OSCIpower-class} because $r_t=O(\eta_t)$.

\subsection{Proofs for conformal testing in \S~\ref{sec:selectiveNDiid}}

\subsubsection{Proof of Theorem~\ref{th:rate:ND}}\label{proof:th:rate:ND}

The proof is analogue to the proof of Theorem~\ref{th:rate:classif}.  For $r_t:=\sqrt{C_1 (\log t)^{\nu_0} t^{-\nu}}$, we let
\begin{equation}
    \label{equassumptiontocheckND}
    \varepsilon_t(\mathcal{F}_{t-1}):=\probp{\abs{\widehat{\lfdr_t}(X_t)-\lfdr(X_t)}> r_t \mid  \mathcal{F}_{t-1}}.
\end{equation}
By \eqref{ass:estiND} and Markov's inequality,
\begin{align}
\probp{\varepsilon_t(\mathcal{F}_{t-1})> r_t}\leq \probp{\e{\abs{\widehat{\lfdr_t}(X_t)-\lfdr(X_t)}\mid \mathcal{F}_{t-1}} > r^2_t} \leq \frac{3}{\pi^2t^3}\label{equ:concvarepilonND}.
\end{align}

\paragraph*{Bounding $D_t$}

Let us show that for all $t\geq T_0$ (with an updated constant $T_0$),
\begin{equation}
    \label{boundDtND}
    \probp{D_t > C_0 r_t, \Omega_t} \leq  \frac{3}{\pi^2t^3}.
\end{equation}
For this, we have for all $q\in [0,\bar q]$,  
\begin{align*}
    \Pi_t(q)
&= \probp{Y_t=0 \mid 1-\widehat{\lfdr_t}(X_t)>q, \mathcal{F}_{t-1}}= \e{ \lfdr(X_t) \mid 1-\widehat{\lfdr_t}(X_t)>q, \mathcal{F}_{t-1}}
\end{align*}
By definition, we have 
\begin{align*}
    1-\Pi_t(q) 
    &=\e{ 1-\lfdr(X_t) \mid 1-\widehat{\lfdr_t}(X_t)>q, \mathcal{F}_{t-1}}.
    \end{align*}
    Hence, we can proceed exactly as in the proof of Theorem~\ref{th:rate:classif} with $W(X_t):=1-\lfdr(X_t)$ and $W_t(X_t):=1-\widehat{\lfdr_t}(X_t)$ to obtain for $t\geq T_0$,  on the event $\Omega_t$
\begin{equation*}
  |\Pi_t(q_t)-\Pi^\bench(q_t)| \ind{q_t\in [0,B]} \ind{\varepsilon_t(\mathcal{F}_{t-1})\leq r_t} \leq C_0 r_t,
\end{equation*}
which gives \eqref{boundDtND} with \eqref{equ:concvarepilonND}.
 
\paragraph*{Bounding $J(t)$}
\et{
Again, we proceed exactly as in the proof of Theorem~\ref{th:rate:classif} to prove that, for $t\geq T_0$ (with an updated constant $T_0$) 
\begin{equation*}
   \P\big( J(t) <\underline{J}(t), \Omega_t\big)\leq 2/t,
\end{equation*}
for 
\begin{align*}
    \underline{J}(t)&:=\lfloor t\:\probp{1-\lfdr(X_t)>\bar q  }/5\rfloor  .
\end{align*}
}
    
\paragraph*{Conclusion}

We apply  Theorem~\ref{rem:coroonanEvent} (case (i)), to obtain \eqref{equ-convND}.

\subsubsection{Optimality of the oracle procedure}

\begin{proposition}\label{prop:optimalND}
Consider the iid online testing model \eqref{model:ND} satisfying  Assumption~\ref{ass:support} with the test statistics $W(x):=1-{\lfdr}(x)$ \eqref{equ:lfdr} and support values $0\leq \kappa<B^*\leq 1$, and an error level $\alpha\in (0,1)$, such that $1-B^*<\alpha<\P(Y=0)$. Then any online testing procedure $(S_t(X_t),\{1\})_{t\geq 1}$ (for which $S_t(\cdot)$ is $\mathcal{F}_{t-1}$-measurable) is dominated by the oracle procedure $(S^*(X_t),\{1\})_{t\geq 1}$ given by \eqref{equ:NDoracleselect} in the following sense:  if for all $t\geq T_0$, $\P(Y_t=0\mid S_t(X_t)=1)\leq \alpha+z$, for $z\geq 0$, then for all $t\geq T_0$
    \begin{equation}\label{eq-optimalpower-ND}
        \E\big[ \ind{Y_t= 1}S_t(X_t)\big]\leq \E\big[ \ind{Y_t=1}S^*(X_t)\big]+z q^\bench_\alpha/(1-\alpha-q^\bench_\alpha).
    \end{equation}
\end{proposition}
For $z=0$, this result is (essentially) proved in \cite{SC2007}. Hence, Proposition~\ref{prop:optimalND} is an online extension, accommodating the case where $z>0$.

\begin{proof}
The proof is analogue to the proof of Proposition~\ref{prop:optimalclassif}, and we provide it here for completeness.
First note that by \eqref{equqalphainfo}, we have $\P(Y_t=0\mid S^*(X_t)=1)=\e{\lfdr(X_t)\mid S^*(X_t)=1} = \alpha$ and $q^\bench_\alpha<1-\alpha$.
Also, by considering $f$ in \eqref{equfSC}, we have
\begin{align*}
S^*(X_t)&=\ind{\lfdr(X_t)<1-q^\bench_\alpha} \\
&=\ind{f(\lfdr(X_t))<c^*}=\ind{\lfdr(X_t)-\alpha - c^* (1-\lfdr(X_t)))<0},    
\end{align*}
for $c^*:=f(1-q^\bench_\alpha)>0$. 
Now consider any online testing procedure $(S_t(X_t),\{1\})_{t\geq 1}$, for which $S_t(\cdot)$ is $\mathcal{F}_{t-1}$-measurable.
Hence, almost surely, we have
$$
(\lfdr(X_t)-\alpha - c^* (1-\lfdr(X_t)))(S^*(X_t)-S_t(X_t))\leq 0.
$$
Rearranging this inequality, we obtain
\begin{align*}
&c^* ((1-\lfdr(X_t))) S^*(X_t)- (1-\lfdr(X_t))S_t(X_t))\\&\geq (\lfdr(X_t)-\alpha)S^*(X_t) - (\lfdr(X_t)-\alpha)S_t(X_t).
\end{align*}
By integrating the latter in $X_t$ and $\mathcal{F}_{t-1}$, and by assuming $\P(Y_t=0\mid S_t(X_t)=1)\leq \alpha+z$, we obtain
\begin{align*}
&c^* (\E((1-\lfdr(X_t)) S^*(X_t))- \E((1-\lfdr(X_t)) S_t(X_t))\\
&\geq \E((\lfdr(X_t)-\alpha)S^*(X_t)) - \E((\lfdr(X_t)-\alpha)S_t(X_t)),
\end{align*}
that is,
\begin{align*}
c^* (\E(\ind{Y_t= 1} S^*(X_t) )- \E(\ind{Y_t= 1} S_t(X_t)))\geq 0  - z \P(S_t(X_t)=1)\geq -z,
\end{align*}
because
$
\P(Y_t=0\mid X_t,\mathcal{F}_{t-1}) = \lfdr(X_t)$.
This gives \eqref{eq-optimalpower-ND}.
\end{proof}

\subsubsection{Proof of Corollary~\ref{optimalityNDiid}}\label{proof:optimalityNDiid}

The proof is analogue to the proof of Corollary~\ref{optimalityclassifiid}, 
by using Theorem~\ref{th:rate:ND} and Proposition~\ref{prop:optimalND}.

\section{Auxiliary results}\label{sec:auxiliary}

\subsection{Specific lemmas}\label{sec:auxiliaryspe}

\begin{proposition}\label{prop:unavoidable}
    Consider any procedure $\cR=((S_t)_{t\geq 1}, (\cC_t)_{t\geq 1})$ \et{which is non-trivial in the sense that it produces non-trivial prediction sets on the selection, that is, such that $\cC_t\neq \cY$ whenever $S_t=1$.} Then it has the following intrinsic adversarial limitations (with the coverage error \eqref{coverror}):
    \begin{itemize}
        \item[(i)] If an FCP bound of the type \eqref{FCPbound_decre} is true, then the procedure should stop selection after some time.
        More formally, assume that for some $\alpha\in (0,1)$ and any sequence $(X_{t},Y_{t})_{t\geq 1}$, we have 
        $\FP_t(\mathcal{R})\leq \alpha + \Psi(J(t))$ for all $t\geq 1$, with $J(t)=\sum_{i=1}^{t-1} S_{i} + 1$ and $\lim_{j\to\infty}\Psi(j)=0$. Then there exists a sequence $(X_{t},Y_{t})_{t\geq 1}$ such that $\sup_{t\geq 1}J(t)<+\infty$;
        \item[(ii)] 
        No FCP bound of the type \eqref{FCPbound_decre} that always converges to some $\alpha\in (0,1)$ can hold. More formally, assume that for some $\alpha\in (0,1)$ and any sequence $(X_{t},Y_{t})_{t\geq 1}$, 
      we have  $\FP_t(\mathcal{R})\leq \alpha + \varepsilon(t)$ for all $t\geq 1$, for some sequence $(\varepsilon(t))_{t\geq 1}$. Then there exists a sequence $(X_{t},Y_{t})_{t\geq 1}$ such that either $S_t=0$ for all $t\geq 1$, or $(\varepsilon(t))_{t\geq 1}$ does not converge to $0$.
    \end{itemize}
\end{proposition}

\begin{proof}
    By assumption, since the procedure is non-trivial at each time, we can recursively build a sequence $(X_{t},Y_{t})_{t\geq 1}$   such that $Y_t\notin \cC_t (X_t,(X_i,Y_i)_{1\leq i\leq t-1})$ for each time point $t\geq 1$ with selection and thus $\FCP_t(\cR)=1$ for all time points $t$ larger than the first selection time. Hence, (i) and (ii) follow.
    \end{proof}

\begin{lemma}\label{lem:deltapositive}
    Consider  the \method\  threshold sequence $(q_t^{\OSCI})_{t\geq 1}$ \eqref{equACIquantile} and let $t_0\geq 1$ be the smallest (deterministic) integer for which $q_1+(1-\alpha)(\gamma_1+\dots+\gamma_{t_0})\geq B$. Then we have $q_T^\OSCI=q_{t_0+1}^\OSCI$ for all $T\geq t_0+1$ on the event where one always wrongly selects at time $\leq t_0$, that is,
\begin{equation}
    \label{equ:unluckyevent}
    \cap_{1\leq t\leq t_0} \set{W_t(X_t)>q'_t, Y_t\notin \cC_t(X_t,q'_t)},
\end{equation}
for $q'_1:=q_1$ and $q'_t:=q_1+(1-\alpha)(\gamma_1+\dots+\gamma_{t-1})$, $2\leq t\leq t_0$.
Moreover, the probability of \eqref{equ:unluckyevent} is positive in the iid model, if $W_t(\cdot)$ and $\cC_t(\cdot,\cdot)$ are  deterministic functions, if for all $1\leq t\leq t_0$, $\P(W_t(X)>q)>0$ for all $q<B$  and $\P(Y_t\notin \cC_t(X_t,q)\mid X_t)$ is positive almost surely.
\end{lemma}

\begin{proof}
    The first statement is obvious from the definition of \method. The second one is true because by the iid assumption, and since each $q'_t$ is deterministic
    \begin{align*}
        &\probp{\cap_{1\leq t\leq t_0} \set{W_t(X_t)>q'_t, Y_t\notin \cC_t(X_t,q'_t)}}\\
        &=\prod_{t=1}^{t_0} \probp{W_t(X)>q'_t, Y\notin \cC_t(X,q'_t)}\\
        &= \prod_{t=1}^{t_0} \e{\probp{Y\notin \cC_t(X,q'_t) \mid X}\ind{W_t(X)>q'_t}}.
    \end{align*}
    Now, if there exists $t\in [1,t_0]$ with $\e{\probp{Y\notin \cC_t(X,q'_t) \mid X}\ind{W_t(X)>q'_t}}=0$, we have almost surely that either $W_t(X)\leq q'_t$ or $\probp{Y\notin \cC_t(X,q'_t) \mid X}=0$. Since $q'_t<B$ for $t\leq t_0$, this is impossible with the assumption of the lemma. This shows that each term is positive and thus the probability of \eqref{equ:unluckyevent} is also positive.
\end{proof}

\begin{lemma}[Concentration of $J(t)$]\label{lem:JtconcentrMart}
Consider the threshold $(q_t^{\OSCI})_{t\geq 1}$ \eqref{equACIquantile}, used with any selection rule $(S_t(\cdot))_{t\geq 1}$. Then, without any model assumption, we have for all $t\geq 1$, $x>0$,
\begin{align}
    \probp{\et{J(t+1)-1}
    - \sum_{k=1}^t\probp{S_k(X_k,q^{\OSCI}_k)=1\mid \mathcal{F}_{k-1}}\geq x}\leq e^{-2x^2/t}\label{equJtconcentrMart};\\
    \probp{\et{J(t+1)-1}
    - \sum_{k=1}^t\probp{S_k(X_k,q^{\OSCI}_k)=1\mid \mathcal{F}_{k-1}}\leq- x}\leq e^{-2x^2/t}\label{equJtconcentrMartotherside},
\end{align}
where $J(t+1)$ is defined by \eqref{def:Jt}.
\end{lemma}
\begin{proof}
By definition, we have
    \begin{align*}
    \et{J(t+1)-1}&= \sum_{k=1}^t \ind{S_k(X_k,q^{\OSCI}_k)=1}\\
    &= \sum_{k=1}^t \paren{\ind{S_k(X_k,q^{\OSCI}_k)=1} - \probp{S_k(X_k,q^{\OSCI}_k)=1\mid \mathcal{F}_{k-1}}} \\
    &\:\:\:+ \sum_{k=1}^t\probp{S_k(X_k,q^{\OSCI}_k)=1\mid \mathcal{F}_{k-1}}\\
    &=M_t + \sum_{k=1}^t\probp{S_k(X_k,q^{\OSCI}_k)=1\mid \mathcal{F}_{k-1}}.
\end{align*}
The process $(M_t)_{t\geq 1}$ is a square integrable martingale adapted to the filtration $(\mathcal{F}_t)_{t\geq 1}$ with $\langle M \rangle_t=\sum_{k=1}^t \e{\Delta M_k^2\mid \mathcal{F}_{t-1}}\leq t/4$, $-\probp{S_k(X_k)=1\mid \mathcal{F}_{k-1}}=:A_k\leq \Delta M_t \leq B_k:=1-\probp{S_k(X_k)=1\mid \mathcal{F}_{k-1}}$ and $\mathcal{D}_t=\sum_{k=1}^t (B_k-A_k)^2\leq t$. Hence, Lemma~\ref{lem:concentrationMart} provides (using that $ 2\langle M \rangle_n+\mathcal{D}_n\leq y=3t/2$ almost surely) that for all $t\geq 1$, $x>0$, $$\probp{M_t\geq x}\leq e^{-2x^2/t}.$$
This leads to  \eqref{equJtconcentrMart}. The inequality \eqref{equJtconcentrMartotherside} is obtained similarly.
\end{proof}

\begin{lemma}\label{lem:concentrespJt}
Consider any sequence threshold $(q_t)_{t\geq 1}$ such that each $q_t$ is $\mathcal{F}_{t-1}$ measurable, and an informative selection rule of the form $S_t(x,q)=\ind{W_t(x)>q}$ for some $\mathcal{F}_{t-1}$-measurable function $W_t:x\in \mathcal{X}\mapsto W_t(x)\in \R$. In the iid model, consider any (deterministic) function $W:x\in \mathcal{X}\mapsto W(x)\in \R$ and assume that the random variable $W(X)$ has a density with respect to the Lebesgue measure bounded uniformly by $C>0$.
Then, we have for all $\delta>0$,
\begin{align}
    &\abs{\probp{W_k(X_k)>q_k\mid \mathcal{F}_{k-1}}- \probp{W(X_k)>q_k  \mid \mathcal{F}_{k-1}}}\nonumber\\&\leq C\delta + \probp{|W_k(X_k)-W(X_k)|>\delta\mid \mathcal{F}_{k-1}} \label{eqeintermlemmaJtconc}.
\end{align}
\end{lemma}

\begin{proof}
Let $\delta>0$. On the one hand, 
\begin{align*}
    &\probp{W_k(X_k)>q_k\mid \mathcal{F}_{k-1}}\\&\leq \probp{W(X_k)>q_k -|W_k(X_k)-W(X_k)| \mid \mathcal{F}_{k-1}} \\
    &\leq \probp{W(X_k)>q_k -\delta \mid \mathcal{F}_{k-1}}+ \probp{|W_k(X_k)-W(X_k)|>\delta\mid \mathcal{F}_{k-1}}\\
    &\leq \probp{W(X_k)>q_k \mid \mathcal{F}_{k-1}} + C\delta+ \probp{|W_k(X_k)-W(X_k)|>\delta\mid \mathcal{F}_{k-1}},
\end{align*}
because, by assumption, $\probp{q_k -\delta \leq W(X_k)\leq q_k  \mid \mathcal{F}_{k-1}} \leq C \delta$.
On the other hand,
\begin{align*}
    &\probp{W_k(X_k)>q_k\mid \mathcal{F}_{k-1}}\\&\geq \probp{W(X_k)>q_k +|W_k(X_k)-W(X_k)| \mid \mathcal{F}_{k-1}} \\
    &\geq \probp{W(X_k)>q_k +\delta \mid \mathcal{F}_{k-1}}- \probp{|W_k(X_k)-W(X_k)|>\delta\mid \mathcal{F}_{k-1}}\\
    &\geq \probp{W(X_k)>q_k  \mid \mathcal{F}_{k-1}}- C\delta- \probp{|W_k(X_k)-W(X_k)|>\delta\mid \mathcal{F}_{k-1}} .
\end{align*}
This establishes \eqref{eqeintermlemmaJtconc}. 
\end{proof}

\begin{lemma}\label{lem:forDT}
Consider the residual-based online adaptive score $V_t(x,y)$  \eqref{regressionVadaptive} with $B=1$ and a model where the distribution of $Y_t$ given $X_t=x, \mathcal{F}_{t-1}$ is given by $\mu(x)+\sigma(x) \xi$, where $\xi$ has a symmetric continuous density $f$ with respect to the Lebesgue measure which is positive on its support with $M_1:=\sup_{x\in \R} f(x)<\infty$ and   $M_2:=\sup_{x\in \R} \{|x|f(x)\}<\infty$. Then we have 
\begin{align}\label{equ:controlcondx}
&\sup_{q\in [0,B]}\abs{\probp{V_t(X_t,Y_t) >q \:|\: X_t=x, \mathcal{F}_{t-1}} - \Pi^0(q)}\\
&\:\:\:\:\:\:\leq M_1  (1/{\sigma}(x)) |\hat{\mu}_t(x)-{\mu}(x)| + 2(1\vee M_2) |\hat{\sigma}_t(x)  / \sigma(x)-1|.   \nonumber
\end{align}
\end{lemma}
\begin{proof}
To establish \eqref{equ:controlcondx}, we use that the model assumption to write
\begin{align*}
& \probp{V_t(X_t,Y_t) >q \:|\: X_t=x, \mathcal{F}_{t-1}} \\&=\Pc{ |\mu(x) - \hat{\mu}_t(x)+ \sigma(x) \xi|> \hat{\sigma}_t(x) \bar{\Phi}^{-1}(1/2-q/2)  }{\hat{\mu}_t(x),\hat{\sigma}_t(x)}\\
&=\bar{F}\Bigg(\frac{\hat{\sigma}_t(x)}{\sigma(x)} \bar{\Phi}^{-1}(1/2-q/2)+\frac{|\hat{\mu}_t(x)-{\mu}(x)|}{\sigma(x)}\Bigg)\nonumber\\
&\:\:\:+\bar{F}\Bigg(\frac{\hat{\sigma}_t(x)}{\sigma(x)} \bar{\Phi}^{-1}(1/2-q/2)-\frac{|\hat{\mu}_t(x)-{\mu}(x)|}{\sigma(x)}\Bigg),
\end{align*}
    because the noise as a symmetric distribution.
      We now use Lemma~\ref{lemma:cdfcontrol} to deduce 
  \begin{align*}
  &\bar{F}\Bigg(\frac{\hat{\sigma}_t(x)}{\sigma(x)} \bar{\Phi}^{-1}(1/2-q/2)+\frac{|\hat{\mu}_t(x)-{\mu}(x)|}{\sigma(x)}\Bigg)-\bar{F}\Bigg( \bar{\Phi}^{-1}(1/2-q/2)\Bigg) \\&\leq M_1  {\sigma}^{-1}(x) |\hat{\mu}_t(x)-{\mu}(x)| + M_2 |\hat{\sigma}_t(x)  / \sigma(x)-1|,  
  \end{align*}
and thus \eqref{equ:controlcondx}.
\end{proof}

\subsection{General lemmas}

\begin{lemma}\label{lemsuminte}
  For $\beta\in (1/2,1)$, we have as $n$ tends to infinity
  \begin{align}
 \sum_{j=n-\Delta_n}^{n}j^{-\beta} &\sim \Delta_n n^{-\beta} \:;\label{equiv1}\\
\sum_{j=n-\Delta_n}^{n}j^{-2\beta} &\sim  \Delta_n n^{-2\beta}\label{equiv2}\:,
 \end{align}
 provided that $n^{\beta}\ll \Delta_n\ll n$ (and $\Delta_n$ is an integer).
 Note that $\Delta_n n^{-\beta}$ tends to infinity in \eqref{equiv1} and $\Delta_n n^{-2\beta}$ tends to $0$ in \eqref{equiv2}.
 In addition, for any $c\in (0,1)$, we have for all $\beta>1/2$,
 \begin{align}
 \sum_{j=\lfloor nc \rfloor }^{n}j^{-\beta} &\sim (1-\beta)^{-1}(1-c^{1-\beta}) n^{1-\beta} \:;\label{equiv3}\\
\sum_{j=\lfloor nc \rfloor}^{n}j^{-2\beta} &\sim (2\beta-1)^{-1}(c^{1-2\beta}-1) n^{1-2\beta} \label{equiv4}\:,
 \end{align}
\end{lemma}
\begin{proof}
{ 
Let $\gamma>0$ with $\gamma\neq 1$. Let $1<m\leq M$ two integers. We have classically for all $j\in\range{m,M}$, $\int_{j}^{j+1}u^{-\gamma}\dd u\leq j^{-\gamma}\leq \int_{j-1}^ju^{-\gamma}\dd u$ since $\gamma>0$. By summing this inequality we get:
\begin{equation}\label{eq:SumPower}
\frac{1}{1-\gamma}\paren{(M+1)^{1-\gamma}-m^{1-\gamma}}\leq \sum_{j=m}^{M}j^{-\gamma}\leq \frac{1}{1-\gamma}\paren{M^{1-\gamma}-(m-1)^{1-\gamma}}.
\end{equation}
Let us first prove \eqref{equiv3} and \eqref{equiv4}. We use~\eqref{eq:SumPower} with $M=n$ and $m=\lfloor{nc\rfloor}$. We have for the left-hand-side and right-hand side of \eqref{eq:SumPower}:
\begin{align*}
  \frac{1}{1-\gamma}\paren{(n+1)^{1-\gamma}-\lfloor{nc\rfloor}^{1-\gamma}} 
  \sim\frac{n^{1-\gamma}}{1-\gamma}\paren{1-c^{1-\gamma}};\\  
\frac{1}{1-\gamma}\paren{n^{1-\gamma}-\lfloor{nc-1\rfloor}^{1-\gamma}} 
\sim\frac{n^{1-\gamma}}{1-\gamma}\paren{1-c^{1-\gamma}},
\end{align*}
which, by using~\eqref{eq:SumPower} gives \eqref{equiv3} and \eqref{equiv4} since for $\beta\in(1/2,1)$, $2\beta>1$.
}

{
Let us prove \eqref{equiv1} and \eqref{equiv2}. We now take $M=n$ and $m=n-\Delta_n$. We study the right hand side of \eqref{eq:SumPower}. We have $(1-\gamma)^{-1}\paren{n^{1-\gamma}-(n-\Delta_n-1)^{1-\gamma}}={n^{1-\gamma}}(1-\gamma)^{-1}{\paren{1-(1-\Delta_nn^{-1}-n^{-1})^{1-\gamma}}}$, and since $\Delta_nn^{-1}\rightarrow 0$, one can use the Taylor expansion of $x\mapsto(1+x)^{1-\gamma}$ near $0$ to obtain the following expansion:
\begin{align*}
    (1-\gamma)^{-1}\paren{n^{1-\gamma}-(n-\Delta_n-1)^{1-\gamma}}&=\frac{n^{1-\gamma}}{1-\gamma}\paren{(1-\gamma)(\Delta_nn^{-1}+n^{-1})+o(\Delta_nn^{-1}+n^{-1}))}\\
    &=\Delta_nn^{-\gamma}+n^{-\gamma}+o(\Delta_nn^{-\gamma}+n^{-\gamma}),
\end{align*}
and since $n^{-\gamma}=o(\Delta_nn^{-\gamma})$ (because $\Delta_n\rightarrow +\infty$), we have that the right-hand side of \eqref{eq:SumPower} is equivalent, when $n$ tends to infinity, to $\Delta_n n^{-\gamma}$. With the same computation, we have that the left hand side is also equivalent to $\Delta_nn^{-\gamma}$, which concludes the proof.
}
\end{proof}

\begin{lemma}\label{lem:decreasing}
Let $U$ be any random variable with values in $[A_1,A_2]$, for $0\leq A_1<A_2\leq 1$ \et{and $\P(U=1)=0$}. Consider the function $\Psi:u\in [0,1]\mapsto \E(U  \:|\: U<u)$, with the convention $\Psi(u)=A_1$ when $\P(U<u)=0$. Then the following holds: 
\begin{itemize}
    \item[(i)]  $\Psi$ is nondecreasing on $[0,1]$ with $\Psi(0)=A_1$ and $\Psi(1)=\E(U)$;
    \item[(ii)] if $U$ has a continuous density $g$ on $[A_1,A_2]$, which is positive on $(A_1,A_2)$, then $\Psi$ is also increasing on $[A_1,A_2]$ and continuous on $[0,1]$. In addition, $\Psi$ is differentiable on $(A_1,A_2)$ with $\Psi'(u)=\frac{g(u)(u -\Psi(u))}{\P(U<u)}$, $u\in (A_1,A_2)$, and $\sup_{u\in [u_0,v_0]}\Psi'(u)<+\infty$, $\inf_{u\in [u_0,v_0]}\Psi'(u)>0$ for any $A_1<u_0<v_0<A_2$.
\end{itemize}
\end{lemma}

\begin{proof}
Let us first check point (i). 
    Take $A_1\leq u<u'\leq A_2$ and show that $\Psi(u)\leq \Psi(u')$. If $\P(U<u)=0$, then $\Psi(u)=A_1$ and the inequality is trivial. Now suppose $\P(U<u)>0$ and let $\gamma=\P(U< u\:|\: U< u')= \frac{\P(U< u)}{\P(U< u')}\in [0,1]$. We have
    \begin{align*}
       \Psi(u')&=\E(U \:|\: U<u') \\
       & = \E( U \ind{U<u} \:|\: U<u') + \E( U \ind{u\leq U<u'} \:|\: U<u')\\
       &= \gamma \E( U  \:|\: U<u) + (1-\gamma)\E( U  \:|\: u\leq U<u')\\
       &\geq \gamma \Psi(u) + (1-\gamma) u \geq \Psi(u),
    \end{align*}
    because $\Psi(u)\leq u$.     This gives the first statement.    
    
    Let us now check (ii). 
    Since in that case $\Psi(u)=\frac{\e{U\ind{U<u}}}{\probp{U<u}}$ for $u\in (A_1,1]$, it is clear that $\Psi$ is continuous on $(A_1,1]$. In addition, $\Psi(u)=A_1$ for $u<A_1$ and $A_1\leq \Psi(A_1+\epsilon)\leq A_1+\epsilon$ for $\epsilon>0$ small enough. Hence, $\Psi$ is also continuous in $A_1$ and thus on $[0,1]$.

Consider $A_1<u<A_2$. We get that $\Psi(u)=\int_{A_1}^uxg(x)\dd x/\P(U\leq u)$ with $g$ the density of $U$. By continuity of $g$ on $[A_1,A_2]$ we get that $\Psi$ is differentiable on $(A_1,A_2)$ with derivative:
\begin{align*}
\Psi'(u)&=\frac{ug(u)\P(U< u)-g(u)\int_{A_1}^u xg(x)d x}{\P(U< u)^2}\\
&=g(u)\frac{u-\Psi(u)}{\P(U< u)}.
\end{align*}
Take $A_1<u_0\leq v_0<A_2$ and let $u\in(u_0,v_0)$.
We first prove the lower bound. We have that $u-\Psi(u)=\E((u-U)\ind{U< u})/\P(U<u)\geq 0$. If this inequality is an equality, since $(u-U)\ind{U< u}\geq 0$ a.s. we get that $u=U$ on the event $\set{U< u}$. However this is impossible since $g$ is continuous on $(A_1,A_2)$ and $u\in(A_1,A_2)$. Hence $u-\Psi(u)>0$ and so $\min_{u\in[u_0,v_0]}(u-\Psi(u))>0$, which implies that for all $u\in[u_0,v_0]$,
\begin{equation*}
    \Psi'(u)\geq \min_{x\in[u_0,v_0]}g(x)\min_{x\in[u_0,v_0]}(x-\Psi(x))>0,
\end{equation*}
with the last inequality holding since $g$ is  positive on $(A_1,A_2)$.  This leads to $\min_{u\in[u_0,v_0]}\Psi'(u)>0$.
We now turn to prove the upper bound. Let $u\in[u_0,v_0]$,
\begin{align*}
    \Psi'(u)&\leq \frac{g(u)u}{\P(U\leq u)}\\
    &\leq \frac{\sup_{x\in[u_0,v_0]}g(x)\times 1}{\P(U\leq u_0)}.
\end{align*}
$g$ being continuous on $(A_1,A_2)$ implies that $g$ is bounded on any compact of $(A_1,A_2)$, and since $u_0>A_1$ we get that $\P(U\leq u_0)>0$, which concludes the proof.
\end{proof}

\begin{lemma}\label{lemma:cdfcontrol}
Let $f$ be a bounded density (with respect to the Lebesgue measure) such that $\sup_{x\in\R}\abs{x f(x) }<\infty$. 
Then the corresponding upper-tail function $\bar F(x)=\int_{x}^{\infty} f(t)dt$ is such that for all $\sigma_1,\sigma_2>0$, for all $x\in\R$ and $h\in \R$,
\begin{equation}\label{equ:cdfcontrol}
\abs{\bar F\paren{\frac{\sigma_1}{\sigma_2}x+\frac{h}{\sigma_2}}-\bar F(x)}\leq\frac{\norm{f}_\infty}{\sigma_2}\abs{h}+2(1\vee\sup_{x\in\R}\abs{x f(x) })\cdot \abs{\sigma_1/\sigma_2-1}. 
\end{equation}
\end{lemma}

\begin{proof}
Let $x\in\R$ and $h\in\R$. We have:
\begin{align*}
    \abs{\bar F\paren{\frac{\sigma_1}{\sigma_2}x+\frac{h}{\sigma_2}}-\bar F(x)}&\leq\abs{\bar F\paren{\frac{\sigma_1}{\sigma_2}x+\frac{h}{\sigma_2}}-\bar F\paren{\frac{\sigma_1}{\sigma_2}x}}+\abs{\bar F\paren{\frac{\sigma_1}{\sigma_2}x}-\bar F\paren{x}}\\
    &= \abs{\bar F\paren{\frac{\sigma_1}{\sigma_2}\brac{x+\frac{h}{\sigma_1}}}-\bar F\paren{\frac{\sigma_1}{\sigma_2}x}}+\abs{\bar F\paren{\frac{\sigma_1}{\sigma_2}x}-\bar F\paren{x}}.
\end{align*}
By the mean value inequality, we have
\begin{equation*}
     \abs{\bar F\paren{\frac{\sigma_1}{\sigma_2}\brac{x+\frac{h}{\sigma_1}}}-\bar F\paren{\frac{\sigma_1}{\sigma_2}x}}\leq \frac{\norm{f}_\infty}{\sigma_2} \abs{h}.
\end{equation*}
 Let us bound the second term. Denote $\gamma=\sigma_1/\sigma_2> 0$ for short. We first assume that $\gamma\geq 1/2$. By the mean value inequality with the function $u\in[0,\infty)\mapsto \bar F(ux)$, we have:
 \begin{align*}
&\abs{\bar F\paren{\gamma x}-\bar F\paren{x}}\\&\leq \sup_{u\in[\gamma\wedge 1,\gamma\vee 1]}\abs{x}f(ux)\abs{\gamma-1}\leq\paren{\sup_{x\in\R}\abs{x}f(x)}
\frac{\abs{\gamma-1}}{\gamma\wedge 1}\leq 2 \paren{\sup_{x\in\R}\abs{x}f(x)} \abs{\gamma-1}.
 \end{align*}
In addition, for $\gamma\in (0,1/2)$, we have $\abs{\bar F\paren{\gamma x}-\bar F\paren{x}}\leq 1\leq 2\abs{\gamma-1}$. This gives, for all $\gamma> 0$,
 \begin{equation*}
     \abs{\bar F\paren{\gamma x}-\bar F\paren{x}}\leq 2\paren{1\vee\sup_{x\in\R}\abs{x}f(x)}\abs{\gamma-1},
 \end{equation*}
and finishes the proof. 
\end{proof}

\begin{lemma}[\citealp{bercu2015concentration}]\label{lem:concentrationMart}
    Let $(M_n)_{n\geq 1}$ be a square integrable martingale adapted to a filtration $(\mathcal{F}_n)_{n\geq 1}$ with $M_0=0$ and consider $\langle M \rangle_n=\sum_{k=1}^n \e{\Delta M_k^2\mid \mathcal{F}_{k-1}}$ where $\Delta M_n=M_n-M_{n-1}$.
    Let $n\geq 1$ and assume that for $k\in [n]$, we have $A_k\leq \Delta M_k\leq B_k$ almost surely, where $A_k$ and $B_k$ are bounded $\mathcal{F}_{k-1}$-measurable random variables. 
    Let $\mathcal{D}_n=\sum_{k=1}^n (B_k-A_k)^2$.
    Then for $x,y>0$,
    $$
\probp{M_n\geq x, 2\langle M \rangle_n+\mathcal{D}_n\leq y}\leq e^{-3x^2/y}.
$$
\end{lemma}

\begin{lemma}[\citealp{dedecker2015subgaussian}]\label{lem:concentrationJtAR}
Consider $(Y_t)_{t\in\mbz}$ the autoregressive model \eqref{regressionAR:model}, where the variables $(Z_t)_{t\in\mbz}$ are iid $\mathcal{N}(0,1)$, $(Y_t)_{t\in\mbz}$ is stationary with $\E(Y_t)=0$ for all $t\in\mbz$, and the polynomial $P(z)=1-\sum_{k=1}^d \varphi_k z^k$ has no root of modulus smaller or equal to $1$. Let $t\geq d+1$ and $X_k=(Y_{k-1},\dots,Y_{k-d})$ for all $k\leq t$.
Consider a function $K(x_1,\dots,x_{t-1})$ which is separately bounded in the sense that for all $k\in\range{t-1}$, $x_1,\dots,x_{t-1}$ and $x'_k$ in $\R^d$,
\begin{equation*}
    \abs{K(x_1,\dots,x_{k-1},x_k,x_{k+1},\dots,x_{t-1})-K(x_1,\dots,x_{k-1},x'_k,x_{k+1},\dots,x_{t-1})}\leq L_k,
\end{equation*}
for some sequence $L_k>0$, $k\in [t-1]$. Then there exists a constant $M_0>0$ (independent of $t$) such that for all $x>0$, we have
\begin{equation}
   \label{equ:concentrationJtAR}
   \P( \abs{K(X_1,\dots,X_{t-1})-\e{K(X_1,\dots,X_{t-1})}}>x)\leq 2 e^{-M_0^{-1} x^2/\sum_{k\in [t-1]}L_k^2}.
\end{equation}
\end{lemma}

\begin{proof}
  By Theorem 1' in \cite{Mokkadem1988mixing},  the process $(X_s)_{s\in\mathbb{Z}}$ is geometrically completely regular and a geometrically ergodic Markov chain. It is also irreducible and aperiodic, so we can apply Theorem 2 in \cite{dedecker2015subgaussian}.
\end{proof}

\section{Algorithms}\label{appendix-algorithms-for-applications}

\begin{algorithm}[h!]
\SetKwInOut{Input}{Input}
\Input{Targeted level $\alpha \in (0,1)$,  decreasing step size $(\gamma_j)_{j\geq 1}$, initial threshold $q_1$, initial training data $\mathcal{D}_0$, 
prediction algorithm $\mathcal{A}$ taking as input labeled observations and returning an estimator of $\mu(x)=\E (Y_t\mid X_t=x)$ and $\sigma(x)=\V^{1/2} (Y_t\mid X_t=x)$, selection rule $S_t(\cdot,\cdot)$.}
$J\gets 0$\;
$\mathcal{D}\gets \mathcal{D}_0$\;
\For{$t\geq 1$}{
$\hat{\mu}_t(\cdot),\hat{\sigma}_t(\cdot) \gets\mathcal{A}\paren{\cD}$ (estimation at time $t$)\;
Observe $X_t$\;
Compute $S_t\gets S_t(X_t,q_t)$\;
\If{
$S_t=1$
}{
$J \gets J+1$\;
\If{$q_t\in (0, 1)$}{
$c \gets \Phi^{-1}(1/2 + q_t/2)  $\;
$\mathcal{C}_t\gets 
[\hat{\mu}_t(X_t)\pm \hat{\sigma}_t(X_t) c] $ (prediction region)}
\If{$q_t\leq 0$}{
$\mathcal{C}_t\gets 
\{\hat{\mu}_t(X_t)\} $ (default prediction region when $q_t\leq 0$)}
\Else{
$\mathcal{C}_t\gets \R$ (default prediction region when $q_t\geq B$)
}
Observe $Y_t$\;
\If{$q_t\geq 0$ {\bf{and}} $Y_t\in \mathcal{C}_t $ }{
$q_{t+1}\gets q_t-\gamma_{J}\alpha$}
\Else{$q_{t+1}\gets q_t+\gamma_{J}(1-\alpha)$}
}
\Else{
$S_t\gets 0$ (do not select $X_t$)\;
$q_{t+1}\gets q_t$\;
Observe $Y_t$ (if adaptive $\mathcal{A}$);
}
$\mathcal{D}\gets\mathcal{D}\cup\set{\paren{X_t,Y_t}}$ (if adaptive $\mathcal{A}$)\;
}
\SetKwInOut{Output}{Output}
\Output{At each time $t$, selection indicator $S_t$ and prediction set $\mathcal{C}_t$ when $S_t=1$. }
\caption{Online selective regression} \label{alg:PredictXoriented}
\end{algorithm}

\begin{algorithm}[h!]
\SetKwInOut{Input}{Input}
\Input{Targeted level $\alpha \in (0,1)$, decreasing step size $(\gamma_j)_{j\geq 1}$, initial threshold $q_1$, initial training data $\cD_0$, classification algorithm $\mathcal{A}$ taking as input labeled observations and returning estimated probabilities of being in each class.}
$J\gets 0$\;
$\mathcal{D}\gets \mathcal{D}_0$\;
\For{$t\geq 1$}{
$\hat{\pi}_t\gets\mathcal{A}\paren{\cD}$, with $\hat{\pi}_t(y|x)$ being an estimation of ${\pi}_t(y|x)=\P(Y_t=y\:|\: X_t=x)$\;
Observe $X_t$\;
$\mathcal{C}_t=[K]$ (default prediction set)\;
$\{\wh{y}_t\}\gets\operatorname{argmax}_{y\in\range{K}}\set{\hat{\pi}_t(y|X_t)}$ (point prediction)\;
\If{$q_t<1$}{
$\mathcal{C}_t=\{\wh{y}_t\}$ (prediction set)\;
}
\If{$\hat{\pi}_t(\wh{y}_t|X_t)>q_t$}{
$S_t\gets 1$ (select $X_t$)\;
$J\gets J+1$\;
Observe $Y_t$\;
\If{$q_t \geq 0$ {\bf{and}} $Y_t=\wh{y}_t$}{
$q_{t+1}\gets q_t-\gamma_{J}\alpha$\;
}
\Else{
$q_{t+1}\gets q_t+\gamma_{J}(1-\alpha)$\;
}
}
\Else{
$S_t\gets 0$ (do not select $X_t$)\;
$q_{t+1}\gets q_t$\;
Observe $Y_t$ (if adaptive $\mathcal{A}$)\;
}
$\mathcal{D}\gets\mathcal{D}\cup\set{\paren{X_t,Y_t}}$ (if adaptive $\mathcal{A}$)\;
}
\SetKwInOut{Output}{Output}
\Output{At each time $t$, selection indicator $S_t$ and point prediction $\mathcal{C}_t=\{\wh{y}_t\}$ for the selected times. }
\caption{Online selective classification} \label{alg:BinaryClassification}
\end{algorithm}

\begin{algorithm}[h!]
\SetKwInOut{Input}{Input}
\Input{Targeted level $\alpha \in (0,1)$; step size sequence $(\gamma_j)_{j\geq 1}$; initial threshold $q_1$; initial labeled training data $\mathcal D_0$; an algorithm that takes as input some training labeled data and the covariates of a new example, and returns the estimated $\lfdr$, $\widehat{\lfdr}$, for the new example.}
$J\gets 0$\;
$\mathcal D \gets \mathcal D_0$\;
\For{$t\geq 1$}{
Observe $X_t$\;
$\widehat{\lfdr}_t \gets \widehat{\lfdr}(X_t; \mathcal D)$ (estimated lfdr for the example at time point $t$)\;
\If{$\widehat{\lfdr}_t< 1-q_t$}{
$S_t\gets 1$ (reject ``$Y_t=0$'')\;
$J\gets J+1$\;
Observe $Y_t$\;
\If{$q_t\geq 0$ {\bf{and}} $Y_t=1$ }{
$q_{t+1}\gets q_t-\gamma_{J}\alpha$\;
}
\Else{
$q_{t+1}\gets q_t+\gamma_{J}(1-\alpha)$\;
}
}
\Else{
$S_t\gets 0$ (do not reject ``$Y_t=0$'')\;
$q_{t+1}\gets q_t$\;
Observe $Y_t$ (if adaptive $\mathcal{A}$)\;
}
$\mathcal D \gets \mathcal D\cup \{(X_t,Y_t) \}$ (if adaptive $\mathcal{A}$)\;
}
\SetKwInOut{Output}{Output}
\Output{At each time $t$, $S_t$ the rejection indicator.}
\caption{Online conformal testing
} 
\label{alg:onlineBH}
\end{algorithm}

\end{document}